\documentclass{elsarticle}

\usepackage[margin=1in]{geometry}

\usepackage{bbm}
\usepackage{graphicx}
\usepackage{pstool}
\usepackage{wrapfig}
\usepackage{caption}
\usepackage{subcaption}
\usepackage[section]{placeins} 

\usepackage{fancyhdr} 
\usepackage{lastpage} 
\usepackage{extramarks} 

\usepackage{xcolor} 
\usepackage{enumerate} 
\usepackage{paralist} 
\usepackage{amsmath, amsthm, amssymb, mathtools}
\usepackage{mathabx, pifont, stmaryrd} 
\usepackage[explicit]{titlesec} 
\usepackage{etoolbox} 
\usepackage{bibentry} 
\makeatletter\let\saved@bibitem\@bibitem\makeatother 
\usepackage[colorlinks, bookmarksopen, bookmarksnumbered,
            citecolor=red,urlcolor=red]{hyperref} 
\makeatletter\let\@bibitem\saved@bibitem\makeatother 

\usepackage{algorithm}
\usepackage{algorithmic}

\newtheorem{remark}{Remark}

\theoremstyle{definition}

\usepackage{bm} 
\usepackage{amsfonts} 

\DeclareMathOperator*{\argmax}{arg\,max}

\newcommand{\norm}[1]{\ensuremath{\left\| #1 \right\|}}

\newcommand{\suchthat}{\mathrel{}\middle|\mathrel{}}

\newcommand{\argoptunc}[2]{\underset{#1}{\arg\min} ~~ #2}

\newcommand{\pder}[2]{\ensuremath{\frac{\partial #1}{\partial #2}}}

\newcommand{\Ccal}{\ensuremath{\mathcal{C}}}
\newcommand{\Dcal}{\ensuremath{\mathcal{D}}}
\newcommand{\Ecal}{\ensuremath{\mathcal{E}}}

\newcommand{\Ical}{\ensuremath{\mathcal{I}}}

\newcommand{\Ocal}{\ensuremath{\mathcal{O}}}
\newcommand{\Pcal}{\ensuremath{\mathcal{P}}}
\newcommand{\Qcal}{\ensuremath{\mathcal{Q}}}
\newcommand{\Rcal}{\ensuremath{\mathcal{R}}}

\newcommand{\Ucal}{\ensuremath{\mathcal{U}}}
\newcommand{\Vcal}{\ensuremath{\mathcal{V}}}
\newcommand{\Wcal}{\ensuremath{\mathcal{W}}}

\newcommand{\Rbb}{\ensuremath{\mathbb{R} }}

\newcommand\Ubm{{\ensuremath{\bm{U}}}}

\newcommand\rbm{{\ensuremath{\bm{r}}}}

\newcommand\alphabold{{\ensuremath{\boldsymbol{\alpha}}}}

\newcommand\betabold{{\ensuremath{\boldsymbol{\beta}}}}

\newcommand\rhobold{{\ensuremath{\boldsymbol{\rho}}}}

\newcommand\psibold{{\ensuremath{\boldsymbol{\psi}}}}

\usepackage{tikz}
\usepackage{pgfplots}
\usepackage{pgfplotstable, filecontents, booktabs}
\pgfplotsset{compat=1.9}

\usetikzlibrary{pgfplots.groupplots}
\usepgfplotslibrary{fillbetween}
\usetikzlibrary{calc,fit,matrix,arrows,automata,positioning,shapes}
\usetikzlibrary{arrows.meta}

\pgfplotsset{select coords between index/.style 2 args={
    x filter/.code={
        \ifnum\coordindex<#1\fi
        \ifnum\coordindex>#2\fi
    }
}}

\tikzset{
 invisible/.style={opacity=0},
 visible on/.style={alt={#1{}{invisible}}},
 alt/.code args={<#1>#2#3}{%
   \alt<#1>{\pgfkeysalso{#2}}{\pgfkeysalso{#3}}
 },
}

\newcommand{\colorbarMatlabParula}[5]{
\begin{tikzpicture}
\begin{axis}[
   hide axis, scale only axis,
   height=0pt, width=0pt,
   colormap={parula}{rgb255=(62,38,168) rgb255=(62,39,172) rgb255=(63,40,175) rgb255=(63,41,178) rgb255=(64,42,180) rgb255=(64,43,183) rgb255=(65,44,186) rgb255=(65,45,189) rgb255=(66,46,191) rgb255=(66,47,194) rgb255=(67,48,197) rgb255=(67,49,200) rgb255=(67,50,202) rgb255=(68,51,205) rgb255=(68,52,208) rgb255=(69,53,210) rgb255=(69,55,213) rgb255=(69,56,215) rgb255=(70,57,217) rgb255=(70,58,220) rgb255=(70,59,222) rgb255=(70,61,224) rgb255=(71,62,225) rgb255=(71,63,227) rgb255=(71,65,229) rgb255=(71,66,230) rgb255=(71,68,232) rgb255=(71,69,233) rgb255=(71,70,235) rgb255=(72,72,236) rgb255=(72,73,237) rgb255=(72,75,238) rgb255=(72,76,240) rgb255=(72,78,241) rgb255=(72,79,242) rgb255=(72,80,243) rgb255=(72,82,244) rgb255=(72,83,245) rgb255=(72,84,246) rgb255=(71,86,247) rgb255=(71,87,247) rgb255=(71,89,248) rgb255=(71,90,249) rgb255=(71,91,250) rgb255=(71,93,250) rgb255=(70,94,251) rgb255=(70,96,251) rgb255=(70,97,252) rgb255=(69,98,252) rgb255=(69,100,253) rgb255=(68,101,253) rgb255=(67,103,253) rgb255=(67,104,254) rgb255=(66,106,254) rgb255=(65,107,254) rgb255=(64,109,254) rgb255=(63,110,255) rgb255=(62,112,255) rgb255=(60,113,255) rgb255=(59,115,255) rgb255=(57,116,255) rgb255=(56,118,254) rgb255=(54,119,254) rgb255=(53,121,253) rgb255=(51,122,253) rgb255=(50,124,252) rgb255=(49,125,252) rgb255=(48,127,251) rgb255=(47,128,250) rgb255=(47,130,250) rgb255=(46,131,249) rgb255=(46,132,248) rgb255=(46,134,248) rgb255=(46,135,247) rgb255=(45,136,246) rgb255=(45,138,245) rgb255=(45,139,244) rgb255=(45,140,243) rgb255=(45,142,242) rgb255=(44,143,241) rgb255=(44,144,240) rgb255=(43,145,239) rgb255=(42,147,238) rgb255=(41,148,237) rgb255=(40,149,236) rgb255=(39,151,235) rgb255=(39,152,234) rgb255=(38,153,233) rgb255=(38,154,232) rgb255=(37,155,232) rgb255=(37,156,231) rgb255=(36,158,230) rgb255=(36,159,229) rgb255=(35,160,229) rgb255=(35,161,228) rgb255=(34,162,228) rgb255=(33,163,227) rgb255=(32,165,227) rgb255=(31,166,226) rgb255=(30,167,225) rgb255=(29,168,225) rgb255=(29,169,224) rgb255=(28,170,223) rgb255=(27,171,222) rgb255=(26,172,221) rgb255=(25,173,220) rgb255=(23,174,218) rgb255=(22,175,217) rgb255=(20,176,216) rgb255=(18,177,214) rgb255=(16,178,213) rgb255=(14,179,212) rgb255=(11,179,210) rgb255=(8,180,209) rgb255=(6,181,207) rgb255=(4,182,206) rgb255=(2,183,204) rgb255=(1,183,202) rgb255=(0,184,201) rgb255=(0,185,199) rgb255=(0,186,198) rgb255=(1,186,196) rgb255=(2,187,194) rgb255=(4,187,193) rgb255=(6,188,191) rgb255=(9,189,189) rgb255=(13,189,188) rgb255=(16,190,186) rgb255=(20,190,184) rgb255=(23,191,182) rgb255=(26,192,181) rgb255=(29,192,179) rgb255=(32,193,177) rgb255=(35,193,175) rgb255=(37,194,174) rgb255=(39,194,172) rgb255=(41,195,170) rgb255=(43,195,168) rgb255=(44,196,166) rgb255=(46,196,165) rgb255=(47,197,163) rgb255=(49,197,161) rgb255=(50,198,159) rgb255=(51,199,157) rgb255=(53,199,155) rgb255=(54,200,153) rgb255=(56,200,150) rgb255=(57,201,148) rgb255=(59,201,146) rgb255=(61,202,144) rgb255=(64,202,141) rgb255=(66,202,139) rgb255=(69,203,137) rgb255=(72,203,134) rgb255=(75,203,132) rgb255=(78,204,129) rgb255=(81,204,127) rgb255=(84,204,124) rgb255=(87,204,122) rgb255=(90,204,119) rgb255=(94,205,116) rgb255=(97,205,114) rgb255=(100,205,111) rgb255=(103,205,108) rgb255=(107,205,105) rgb255=(110,205,102) rgb255=(114,205,100) rgb255=(118,204,97) rgb255=(121,204,94) rgb255=(125,204,91) rgb255=(129,204,89) rgb255=(132,204,86) rgb255=(136,203,83) rgb255=(139,203,81) rgb255=(143,203,78) rgb255=(147,202,75) rgb255=(150,202,72) rgb255=(154,201,70) rgb255=(157,201,67) rgb255=(161,200,64) rgb255=(164,200,62) rgb255=(167,199,59) rgb255=(171,199,57) rgb255=(174,198,55) rgb255=(178,198,53) rgb255=(181,197,51) rgb255=(184,196,49) rgb255=(187,196,47) rgb255=(190,195,45) rgb255=(194,195,44) rgb255=(197,194,42) rgb255=(200,193,41) rgb255=(203,193,40) rgb255=(206,192,39) rgb255=(208,191,39) rgb255=(211,191,39) rgb255=(214,190,39) rgb255=(217,190,40) rgb255=(219,189,40) rgb255=(222,188,41) rgb255=(225,188,42) rgb255=(227,188,43) rgb255=(230,187,45) rgb255=(232,187,46) rgb255=(234,186,48) rgb255=(236,186,50) rgb255=(239,186,53) rgb255=(241,186,55) rgb255=(243,186,57) rgb255=(245,186,59) rgb255=(247,186,61) rgb255=(249,186,62) rgb255=(251,187,62) rgb255=(252,188,62) rgb255=(254,189,61) rgb255=(254,190,60) rgb255=(254,192,59) rgb255=(254,193,58) rgb255=(254,194,57) rgb255=(254,196,56) rgb255=(254,197,55) rgb255=(254,199,53) rgb255=(254,200,52) rgb255=(254,202,51) rgb255=(253,203,50) rgb255=(253,205,49) rgb255=(253,206,49) rgb255=(252,208,48) rgb255=(251,210,47) rgb255=(251,211,46) rgb255=(250,213,46) rgb255=(249,214,45) rgb255=(249,216,44) rgb255=(248,217,43) rgb255=(247,219,42) rgb255=(247,221,42) rgb255=(246,222,41) rgb255=(246,224,40) rgb255=(245,225,40) rgb255=(245,227,39) rgb255=(245,229,38) rgb255=(245,230,38) rgb255=(245,232,37) rgb255=(245,233,36) rgb255=(245,235,35) rgb255=(245,236,34) rgb255=(245,238,33) rgb255=(246,239,32) rgb255=(246,241,31) rgb255=(246,242,30) rgb255=(247,244,28) rgb255=(247,245,27) rgb255=(248,247,26) rgb255=(248,248,24) rgb255=(249,249,22) rgb255=(249,251,21) },
   colorbar horizontal,
   point meta min=#1, point meta max=#5,
   colorbar style={width=10cm, xtick={#1,#2,#3,#4,#5}}
]
\addplot [draw=none] coordinates {(0,0)};
\end{axis}
\end{tikzpicture}
}

\usepackage{setspace}
\usepackage{enumitem}
\usepackage{ulem}

\newbool{fastcompile}
\setbool{fastcompile}{false}

\begin{document}
\title{An adaptive model reduction method leveraging \\ locally supported basis functions}

\author[rvt1]{Han Gao\fnref{fn1}}
\ead{hgao1@seas.harvard.edu}

\author[rvt2]{Matthew J. Zahr\fnref{fn2}\corref{cor1}}
\ead{mzahr@nd.edu}

\address[rvt1]{School of Engineering and Applied Science, Harvard University, Cambridge, MA  02138, United States}

\address[rvt2]{Department of Aerospace and Mechanical Engineering, University
	of Notre Dame, Notre Dame, IN 46556, United States}

\cortext[cor1]{Corresponding author}

\fntext[fn1]{Postdoctoral Scholar, School of Engineering and Applied Science, Harvard University}
\fntext[fn2]{Assistant Professor, Department of Aerospace and Mechanical
             Engineering, University of Notre Dame}


\begin{abstract}
We propose a new method, the continuous Galerkin method with globally and locally supported basis functions (CG-GL), to address the parametric robustness issues of reduced-order models (ROMs) by incorporating solution-based adaptivity with locally supported finite element basis functions. The CG-GL method combines the accuracy of locally supported basis functions with the efficiency of globally supported data-driven basis functions. Efficient output-based dual-weighted residual error estimates are derived and implemented for the CG-GL method and used to drive efficient online trial space adaptation. An empirical quadrature procedure is introduced for rapid evaluation of nonlinear terms that does not require retraining throughout the adaptation process. Two numerical experiments demonstrate the potential of the CG-GL method to produce accurate approximations with limited training and its tunable tradeoff between accuracy and computational cost.
\end{abstract}
    
\maketitle

\section{Introduction}
\label{sec:intro}

Projection-based reduced-order models (ROMs), which project the dynamics of a large-scale system onto a low-dimensional subspace derived from experimental or numerical data, have been considered a game-changing technology needed to enable rapid, high-fidelity, many-query analyses of mission-critical challenge problems \cite{amsallem2008interpolation}. While reduced-order models have demonstrated exceptional reduction in computational cost, up to six orders of magnitude relative to the full fidelity simulation \cite{washabaugh2016faster}, they lack parametric robustness and often require substantial training. That is, at least for practical problems that commonly arise in engineering applications, they have shown limited prediction potential outside the parameter set at which they were trained. This lack of parametric robustness, common to even state-of-the-art methods, has prevented ROMs from delivering on their game-changing potential for relevant many-query analysis scenarios and remains an open problem in the community.

A majority of model reduction techniques are built on an offline-online decomposition of cost to realize computational efficiency in a parametric setting \cite{benner2015survey}. In the offline phase, the low-dimensional subspace is constructed by meticulously collecting basis functions, e.g., using a greedy method to sample the parametrized solution over the parameter space \cite{buffa_priori_2012,veroy2003posteriori}. In the online phase, the nonlinear dynamics are integrated in the low-dimensional space to provide rapid predictions of the original system at unseen parameter configurations. Despite this decomposition, most model reduction methods suffer from a costly bottleneck in the online phase associated with projection of nonlinear terms onto the low-dimensional subspace, which involves high-dimensional operations and must be done at every iteration or time step. A slew of hyperreduction techniques have been introduced to eliminate this bottleneck, including the Discrete Empirical Interpolation Method \cite{chaturantabut2010nonlinear}, the Energy Conserving Sampling and Weighting method \cite{farhat2014dimensional}, and the empirical quadrature procedure (EQP) method \cite{yano2019discontinuous}.

This offline-online decomposition is, in part, responsible for the lack of parametric robustness of ROMs because they can only predict linear combinations of states encountered during the offline, training phase. However, because of their online efficiency, many approaches have been developed to improve parametric robustness in this setting.  One such class of methods aim to localize the approximation in either state space \cite{amsallem2012nonlinear}, parameter space \cite{haasdonk_training_2011}, or time \cite{dihlmann_model_2011}. These methods can lead to smaller, more accurate ROMs; however, they do not significantly improve parametric robustness because they are still restricted to linear combinations of states encountered in the offline phase. Recently, nonlinear manifolds \cite{nair_transported_2019,lee2020model,taddei_registration_2020,barnett2022quadratic,mirhoseini2023model} have been used in place of the traditional affine model reduction approximation to improve the parametric robustness of ROMs for convection-dominated problem. While these approaches are significantly more complicated than traditional ROMs, they have proven effective at increasing the predictive potential of ROMs in an offline-online setting. 


Adaptive approaches that break the traditional offline-online framework have emerged to improve the parametric robustness of ROM. The salient feature of such approaches is the subspace of ROM constructed in the offline phase is modified during the online phase, thus allowing the ROM to adapt based on online information. One class of methods sample the solution of the large-scale system when an error indicator exceeds a threshold, which has been used successfully to accelerate optimization problems \cite{arian_trust-region_2000,agarwal_trust-region_2013,zahr2015progressive,wen2023globally}. Because they require additional expensive large-scale system solves during the online phase, they have not been adopted in more general settings.
Adaptive ROMs \cite{kerfriden2011bridging,peherstorfer_model_2020,huang2023predictive,zucatti2023adaptive} modify the affine approximation space with the underlying high-dimensional discretization at selected locations in the spatial domain during the online phase. A similar approach \cite{carlberg_adaptive_2015} adapts the ROM by splitting the reduced basis vectors into multiple vectors with disjoint support, a form of $h$-refinement in the context of model reduction. These methods exhibit excellent parametric prediction but have limited speedup potential. Finally, domain decomposition approaches \cite{lucia_reduced_2003,riffaud2021dgdd} use a standard high-fidelity discretization (e.g., finite volume or finite element) in regions of the domain not amenable to low-dimensional representation and the ROM elsewhere. These approaches have greater parametric robustness than traditional ROM but require offline knowledge of the appropriate domain decomposition, which is not usually available for parametric or unsteady problems. The method proposed in this work aims to overcome this key limitation of domain decomposition approaches using efficient error estimation and adaptivity, as well as an approach to hyperreduction that is amenable to adaptation without retraining.


We propose a new framework for model reduction that will remedy the parametric robustness issues of ROMs and provide a framework with a controllable tradeoff between accuracy and computational cost. Inspired by the accuracy and robustness of numerical discretizations that utilize locally supported basis functions, e.g., finite element method, and the low cost of reduced-order models that use a globally supported, data-driven basis, the proposed continuous Galerkin method with globally and locally supported basis functions (CG-GL) aims to combine the advantages of each. In the spirit of model reduction, basis functions with global support will be used to represent the solution up to some localized error of parametrized features. Only a small number of global modes will be required to approximate the global features of the solution since they will be derived from data from training simulations \cite{sirovich1987turbulence}. In the spirit of finite element discretizations, high-order basis functions with local support will compensate for information missing from the global modes. We also develop an efficient empirical quadrature procedure to evaluate integrals in CG-GL variational form over the portion of the domain where only global basis functions are employed. A key element of the proposed approach is an output-based dual-weighted residual (DWR) error estimation and adaptation framework \cite{becker2001optimal,fidkowski2011review} for the CG-GL method that provides a controllable tradeoff between accuracy and computational cost. We will show that this new framework allows for efficient online basis adaptation, a capability missing from current model reduction approaches. The main contributions of this work are summarized as:
\begin{itemize}
\item a variational framework for a finite element method that simultaneously employs both globally supported and locally supported basis functions, 
\item an empirical quadrature procedure to efficiently integrate terms in the variational formulation over the region of the domain with only globally supported basis functions, and
\item an error estimation and adaptation framework for the CG-GL method.
\end{itemize}

The remainder of this paper is organized as follows. Section~\ref{sec:mthd} introduces the governing conservation law and the proposed CG-GL method. Section~\ref{sec:hpr} introduces an approach to efficiency integrate terms in the CG-GL formulation over the region where only global basis functions are employed. Section~\ref{sec:errest} specializes the well-known DWR error estimation framework \cite{fidkowski2011review} to the CG-GL setting and proposes a strategy to adapt the CG-GL trial space based on the error estimates. Section~\ref{sec:rslt} demonstrates the merits of the proposed CG-GL framework for parametrized, linear and nonlinear partial differential equations (PDEs). Section~\ref{sec:conclu} concludes the paper.

%

\section{Governing equations and  global-local discretization}
\label{sec:mthd}

\subsection{Parametrized system of conservation laws}
\label{sec:mthd:claw}
We are interested in a $\mu$-parameterized system of $m$ conservation laws in $d$-dimensions of the form
\begin{equation}
	\nabla\cdot f(u,\nabla u; \mu) = s(u,\nabla u; \mu)\quad\text{in}\quad\Omega,
	\label{eqn:claw}
\end{equation}
where $\mu\in \Dcal$ is the paramter, $\Dcal$ is the parameter domain, $u(x;\mu)\in\Rbb^m$ is the solution at a point $x\in\Omega\subset\Rbb^d$, $\nabla\;:=[\partial_{x_1},...,\partial_{x_d}]$ is the gradient on the domain $\Omega$ such that $\nabla z\;=[\partial_{x_1}z,...,\partial_{x_d}z]$, $f:\Rbb^m\times\Rbb^{m\times d}\times \Dcal\rightarrow\Rbb^{m\times d}$ is the flux function, and $s:\Rbb^m\times\Rbb^{m\times d}\times \Dcal\rightarrow\Rbb^m$ is the source term. The boundary of the domain is $\partial\Omega$ with outward unit normal $n:\partial\Omega\rightarrow\Rbb^d$. The formulation of the conservation law in \eqref{eqn:claw} is sufficiently general to encapsulate steady second-order PDEs in a $d$-dimensional spatial domain such as Poisson's equation and the Navier-Stokes (NS) equations.

\subsection{Global-local discretization of the conservation law}
\label{sec:mthd:cggl}

\begin{figure}[htp]
	\centering
	\begin{tikzpicture}
\begin{axis}[
axis equal image,
axis line style={gray},
axis x line*=bottom,
axis y line*=left,
width=0.7\textwidth,
xtick=\empty,
ytick=\empty,
grid=major,
ymax=3.3,
xmax=10.3,
ylabel=$x_2$,
xmin=-5.3,
xlabel=$x_1$,
ymin=-3.3]
\addplot [opacity=0.6, fill=black!30!white, opacity=0.6, forget plot]
coordinates {
( 1.00000000e+01, -3.00000000e+00)
( 1.00000000e+01,  3.00000000e+00)
(-5.00000000e+00,  3.00000000e+00)
(-5.00000000e+00, -3.00000000e+00)
( 1.00000000e+01, -3.00000000e+00)};

\addplot [opacity=0.6, fill=red, forget plot]
coordinates {
( 2.50000000e+00, -1.00000000e+00)
( 2.50000000e+00,  1.00000000e+00)
(-1.25000000e+00,  1.00000000e+00)
(-1.25000000e+00, -1.00000000e+00)
( 2.50000000e+00, -1.00000000e+00)};

\node[above]    at    (axis cs:-4, -2) {$\bar\Omega_l$};
\node[above]    at    (axis cs:0, -1) {$\Omega_l$};
\end{axis}
\end{tikzpicture}
	\caption{Decomposition of the conservation law domain $\Omega$ into $\bar\Omega_l$,  the region where the solution is approximated solely with basis functions with global support, and $\Omega_l$, the region where the solution is approximated as a linear combination of globally and locally supported basis functions.  }
	\label{fig:domain_decomp}
\end{figure}
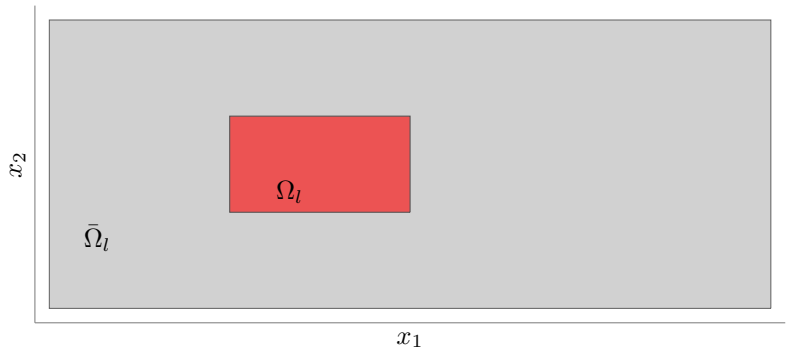

The construction of the proposed CG-GL method begins as a standard continuous Galerkin discretization on a special non-overlapping decomposition of the domain: $\Omega=\Omega_l\cup\bar\Omega_l$ where $\bar\Omega_l = \Omega \setminus \Omega_l$ (Figure~\ref{fig:domain_decomp}). The solution of \eqref{eqn:claw}
will be approximated solely using globally supported basis functions in $\bar\Omega_l$ and a mix of global and local basis functions in $\Omega_l$.  The boundary of the local domain $\partial\Omega_l$, can be decomposed into three parts: $\partial\Omega_l = \partial\Omega_{l,\mathrm{int}} \cup \partial\Omega_{l,\mathrm{nat}} \cup \partial\Omega_{l,\mathrm{ess}}$,  where $\partial\Omega_{l,\mathrm{int}} = \partial\Omega_l\cap\partial\bar\Omega_l$ is the boundary of $\Omega_l$ that intersects the boundary of $\bar\Omega_l$, $\partial\Omega_{l,\mathrm{ess}}$ is the boundary of $\Omega_l$ that intersects the boundary of $\Omega$ that has essential boundary conditions ($\partial\Omega_\mathrm{ess}$) , and $\partial\Omega_{l,\mathrm{nat}} = \partial\Omega_l \setminus (\partial\Omega_{l,\mathrm{int}}\cup\partial\Omega_{l,\mathrm{ess}})$ represents the boundary that has natural boundary conditions.

Let $\Ecal_{l,h}\subset 2^\Omega$ denote a discretization of $\Omega_l$ into $N_l^\mathrm{e}$ non-overlapping, curved elements. To establish the finite-dimensional CG formulation, we first introduce the CG approximation (trial) space of continuous piecewise polynomials associated with the mesh $\Ecal_{l,h}$,
\begin{equation}
	\bar\Vcal_{l,h,p}=\Bigl\{ \phi\in[H^1(\Omega_l)]^m\;\Big|\; \phi|_K\in[\Qcal_p(K)]^m,\forall K\in\Ecal_{l,h}, \left.\phi\right|_{\partial\Omega_{l,\mathrm{int}} \cup \partial\Omega_{l,\mathrm{ess}}} = 0  \Bigr\},
\end{equation}
where $\Qcal_p(K)$ is the space of polynomial functions of degree at most $p\geq1$ on the element $K$. Furthermore, in the context of CG formulation, we define the space of continuous piecewise polynomials of degree $p$ associated with the mesh $\Ecal_{l,h}$ via zero extension of $\bar\Vcal_{l,h,p}$ to the domain $\Omega$ as $\Vcal_{l,h,p}$, 
\begin{equation}
	\Vcal_{l,h,p} \coloneqq \left\{v \in [H^1({\Omega})]^m \Big|\; \left.v\right|_{\Omega_l} \in \bar\Vcal_{l,h,p}, \left.v\right|_{\bar\Omega_l} = 0 \right\},
\end{equation}
and introduce $\{\phi_1,\dots,\phi_{N_l}\}$ as a nodal basis of the space $\Vcal_{l,h,p}$.

Next, we introduce another space of continuous functions over the domain $\Omega$,
\begin{equation}
	\Wcal_k := \Bigg\{ w \in[H^1(\Omega)]^m \Bigg|\;  w\in\mathrm{span}(\psi_1,\dots,\psi_k),\quad\psi_i|_{\partial\Omega_\mathrm{ess}}=0,\quad i=1,\dots,k  \Bigg\},
	\label{eqn:global_space}
\end{equation}
where $\{\psi_1,\dots,\psi_k\}$ is a globally supported basis of $\Wcal_k$. Finally, we introduce a function $\bar{\psi}\in[H^1(\Omega)]^m$  that satisfies the essential boundary conditions. The complete trial space is taken as
\begin{equation}
	\Ucal_{l,h,p,k} \coloneqq \Vcal_{l,h,p} + \Wcal_k + \bar\psi
\end{equation}
and the corresponding space for homogeneous essential boundary conditions is
\begin{equation}
	\bar\Ucal_{l,h,p,k} \coloneqq \Vcal_{l,h,p} + \Wcal_k.
\end{equation}
An element of the trial space $v\in\Ucal_{l,h,p,k}$ takes the form
\begin{equation}
	v(x) = \bar\psi + \sum_{I=1}^{N_l}\phi_I(x)\alpha_I + \sum_{J=1}^k\psi_J(x)\beta_J,
	\label{eqn:trail-form}
\end{equation}
where $\alphabold = (\alpha_1,\dots,\alpha_{N_l}) \in \Rbb^{N_l}$ are the coefficients corresponding to the local basis functions $\{\phi_I\}_{I=1}^{N_l}$ and $\betabold = (\beta_1,\dots,\beta_k) \in\Rbb^k$ are the coefficients corresponding to the global basis functions $\{\psi_J\}_{J=1}^k$ . The local basis functions $\phi_i$ are constructed from Lagrange polynomials using standard finite element techniques and the global basis functions $\psi_i$ are constructed using the method of snapshots, i.e., an approximation to the PDE solution $u$ is obtained at several parameter configures by solving a standard finite element discretization on a sufficiently refined mesh and, if necessary, compressed using e.g., proper orthogonal decomposition \cite{sirovich1987turbulence}.

The variational formulation of the CG-GL method is: given $\mu\in\Dcal$, find $u_{l,h,p,k}\in\Ucal_{l,h,p,k}$ such that
\begin{equation}
 r(u_{l,h,p,k}, w; \mu) = 0,
\end{equation}
for all $w \in \bar\Ucal_{l,h,p,k}$, where $r:\Ucal_{l,h,p,k}\times\bar\Ucal_{l,h,p,k}\times \Dcal\rightarrow \Rbb$ is the CG-GL variational residual
\begin{equation}
 r(v, w; \mu) \coloneqq \int_{\partial\Omega} w\cdot f(v,\nabla v;\mu)n \, dS - \int_{\Omega}\nabla w:f(v,\nabla v;\mu) \, dV-\int_{\Omega} w\cdot s(v,\nabla v;\mu) \, dV = 0.
\end{equation}
The residual can be split into contributions from the local and global domains as
\begin{equation}
	\begin{split}
		r(v,w;\mu)\coloneq r_{\Omega_l}(v,w;\mu)+r_{\bar\Omega_l}(v,w;\mu),
	\end{split}
	\label{eqn:glres}
\end{equation}
where $r_\Rcal:\Ucal_{l,h,p,k}\times\bar\Ucal_{l,h,p,k}\times \Dcal\rightarrow \Rbb$ is the contribution of the region $\Rcal\subset\Omega$
to the residual
\begin{equation} \label{eqn:globlocres}
	r_\Rcal(v,w;\mu) \coloneq \int_{\partial\Omega\cap\partial \Rcal}w\cdot f(v,\nabla v;\mu)n \, dS-\int_\Rcal\nabla w:f(v,\nabla v;\mu) \, dV-\int_\Rcal w\cdot s(v,\nabla v;\mu) \, dV.
\end{equation}
\begin{remark}
 In the limit of $\Omega_l = \Omega$ and $k = 0$, the CG-GL formulation is a standard continuous Galerkin discretization of the governing equations. In the limit of $\Omega_l = \emptyset$ and $N_l=0$, the CG-GL formulation is a Galerkin discretization employing global basis functions, e.g., a Ritz method or reduced-order model.
\end{remark}

The contribution of the local domain to the residual can be evaluated as a summation of elemental contributions
\begin{equation}
 r_{\Omega_l}(v,w;\mu) = \sum_{K \in \Ecal_{l,h}} r_K(v,w;\mu).
\end{equation}
As written, the same approach cannot be used to evaluate the global contribution to the residual, $r_{\bar\Omega_l}$, because, by construction, there is no mesh of this domain. Therefore, we introduce a \textit{quadrature mesh} of the global region that will be solely used to evaluate integrals over $\bar\Omega_l$. We use $\bar\Ecal_l$ to denote this mesh of $\bar{N}_l^\mathrm{e} = |\bar\Ecal_l|$ non-overlapping elements such that $\bigcup_{K\in\bar\Ecal_l} = \bar\Omega_l$. Then the contribution of the global domain to the residual can be evaluated as a summation
of elemental contributions
\begin{equation}
 r_{\bar\Omega_l}(v,w;\mu) = \sum_{K \in \bar\Ecal_l} r_K(v,w;\mu).
\end{equation}
Because $\bar\Ecal_l$ is solely used for quadrature, it does not need to be conforming and can be coarser than grids used to approximate the solution if high-order quadrature rules are used.

To obtain an algebraic form of the problem, we introduce an operator $V_{l,h,p} : \Rbb^{N_l}\rightarrow\Vcal_{l,h,p}$ which associates generalized coordinate $\alphabold\in\Rbb^{N_l}$ to a function $V_{l,h,p}(\alphabold)=\sum_{I=1}^{N_l}\alpha_I \phi_I\in\Vcal_{l,h,p}$. We next introduce an operator $W_k:\Rbb^k\rightarrow\Wcal_k$ which associates generalized coordinate $\betabold\in\Rbb^k$ to a function $W_k(\betabold)=\sum_{J=1}^k\beta_J \psi_J\in\Wcal_k$. Finally,  we define the operator $U_{l,h,p,k}:\Rbb^{N_l+k}\rightarrow\Ucal_{l,h,p,k}$ which associates a generalized coordinate $ \Ubm=(\alphabold,\betabold)\in\Rbb^{N_l+k}$ ($\alphabold\in\Rbb^{N_l},\betabold\in\Rbb^{k}$)  to a function 
\begin{equation}
	U_{l,h,p,k}(\Ubm)\coloneqq \bar\psi+ V_{l,h,p}(\alphabold)+W_k(\betabold)\in\Ucal_{l,h,p,k}
	\label{eqn:discretized2continuous}
\end{equation}
 to maintain a connection between the algebraic and functional representation of the CG-GL solution.

The weak formulation reduces to a system of nonlinear algebraic equations with residual
$\rbm_{l,h,p,k} : \Rbb^{N_l+k} \times \Dcal \rightarrow \Rbb^{N_l+k}$
by introducing a basis $\{\tilde{\phi}_i\}_{i=1}^{N_l+k}$ for the test space $\bar\Ucal_{l,h,p,k}$
\begin{equation}
 \left(\rbm_{l,h,p,k}(\Ubm;\mu)\right)_i \coloneqq r(U_{l,h,p,k}(\Ubm),\tilde{\phi}_i;\mu) =  r_{\Omega_l}(U_{l,h,p,k}(\Ubm),\tilde{\phi}_i;\mu) + r_{\bar\Omega_l}(U_{l,h,p,k}(\Ubm),\tilde{\phi}_i;\mu),
\end{equation}
for $i = 1,\dots,N_l+k$. We take $\tilde\phi_i = \phi_i$ for $i=1,\dots,N_l$ and $\tilde\phi_i = \psi_{i-N_l}$ for $i = N_l+1,\dots,N_l+k$, which allows us to write the residual as
\begin{equation}
 \left(\rbm_{l,h,p,k}(\Ubm;\mu)\right)_i =
 \begin{cases}
  r_{\Omega_l}(U_{l,h,p,k}(\Ubm),\phi_i;\mu) + r_{\bar\Omega_l}(U_{l,h,p,k}(\Ubm),\phi_i;\mu) & i = 1,\dots,N_l \\
  r_{\Omega_l}(U_{l,h,p,k}(\Ubm),\psi_j;\mu) + r_{\bar\Omega_l}(U_{l,h,p,k}(\Ubm),\psi_j;\mu) & i = N_l+1,\dots,N_l+k,
 \end{cases}
\end{equation}
where $j = i-N_l$.
First, we observe that $r_{\bar\Omega_l}(U_{l,h,p,k}(\Ubm),\phi_i;\mu) = 0$ because the local basis functions $\phi_i$ are zero in $\bar\Omega_l$, by construction. Next, we introduce the algebraic form of the remaining three terms as
\begin{equation} \label{eqn:cggl_res}
 \begin{aligned}
  (\rbm_{l,h,p,k}^{\mathrm{\mathrm{ll}}}(\Ubm;\mu))_i &\coloneqq r_{\Omega_l}(U_{l,h,p,k}(\Ubm),\phi_i;\mu) \\
  (\rbm_{l,h,p,k}^{\mathrm{\mathrm{gl}}}(\Ubm;\mu))_i &\coloneqq r_{\Omega_l}(U_{l,h,p,k}(\Ubm),\psi_i;\mu) \\
  (\rbm_k^{\mathrm{\mathrm{gg}}}(\betabold;\mu))_j &\coloneqq r_{\bar\Omega_l}(\bar\psi+W_k(\betabold),\psi_i;\mu),
 \end{aligned}
\end{equation}
for $i=1,\dots,N_l$ and $j=1,\dots,k$, where
$\rbm_{l,h,p,k}^{\mathrm{ll}} : \Rbb^{N_l+k}\times\Dcal \rightarrow \Rbb^{N_l}$ is the algebraic form of the residual over the local domain corresponding to equations from the local test functions,
$\rbm_{l,h,p,k}^{\mathrm{gl}} : \Rbb^{N_l+k}\times\Dcal \rightarrow \Rbb^{N_l}$ is the algebraic form of the residual over the local domain corresponding to equations from the global test functions, and
$\rbm_k^{\mathrm{gg}} : \Rbb^k\times\Dcal \rightarrow \Rbb^k$ is the algebraic form of the residual over the global domain.
With these definitions, the complete residual can be written as
\begin{equation} \label{eqn:cggl_res2}
	\rbm_{l,h,p,k}(\Ubm;\mu)= (\rbm_{l,h,p,k}^{\mathrm{ll}}(\Ubm;\mu), \rbm_{l,h,p,k}^{\mathrm{gl}}(\Ubm;\mu)+\rbm_k^{\mathrm{gg}}(\betabold;\mu))
\end{equation} 
and the primal problem reads: given $\mu\in\Dcal$, find $\Ubm^\star\in\Rbb^{N_l+k}$ such that
\begin{equation} \label{eqn:cggl_eqn}
	\rbm_{l,h,p,k}(\Ubm^\star;\mu)=0
\end{equation}
and the continuous solution is obtained as $u_{l,h,p,k}(\,\cdot\,;\mu)=U_{l,h,p,k}(\Ubm^\star(\mu))$.

We close this section by expanding the algebraic residuals in terms of elemental contributions
\begin{equation} \label{eqn:cggl_res_sum}
 \begin{aligned}
  \rbm_{l,h,p,k}^{\mathrm{ll}}(\Ubm;\mu) &= \sum_{K\in\Ecal_{l,h}} \hat\rbm^{\mathrm{ll}}_{l,h,p,k}(\Ubm;\mu,K) \\
  \rbm_{l,h,p,k}^{\mathrm{gl}}(\Ubm;\mu) &= \sum_{K\in\Ecal_{l,h}} \hat\rbm_{l,h,p,k}^{\mathrm{gl}}(\Ubm;\mu,K)  \\
  \rbm_k^{\mathrm{gg}}(\betabold;\mu) &= \sum_{K\in\bar\Ecal_l} \hat\rbm_k^{\mathrm{gg}}(\betabold;\mu,K),
 \end{aligned}
\end{equation}
where $\hat\rbm_{l,h,p,k}^{\mathrm{ll}} : \Rbb^{N_l+k}\times\Dcal\times 2^{\Omega_l} \rightarrow \Rbb^{N_l}$ and
$\hat\rbm_{l,h,p,k}^{\mathrm{gl}} : \Rbb^{N_l+k}\times\Dcal\times 2^{\Omega_l} \rightarrow \Rbb^k$ are subdomain
contributions to the corresponding residual, and
$\hat\rbm_k^{\mathrm{gg}} : \Rbb^k\times\Dcal\times 2^{\bar\Omega_l} \rightarrow \Rbb^k$ is a subdomain contribution to $\rbm_k^{\mathrm{gg}}$. These elemental contributions are defined as
\begin{equation} \label{eqn:cggl_res_el}
\begin{aligned}
(\hat\rbm_{l,h,p,k}^{\mathrm{ll}}(\Ubm;\mu,K))_i &\coloneqq r_K(U_{l,h,p,k}(\Ubm),\phi_i;\mu) \\
(\hat\rbm_{l,h,p,k}^{\mathrm{gl}}(\Ubm;\mu,K))_i &\coloneqq r_K(U_{l,h,p,k}(\Ubm),\psi_i;\mu) \\
 (\hat\rbm_k^{\mathrm{gg}}(\betabold;\mu,K))_i &\coloneqq r_K(\bar\psi+W_k(\betabold),\psi_i;\mu).
\end{aligned}
\end{equation}

\begin{remark}
The evaluation of $\rbm_k^{\mathrm{gg}}$ requires $\Ocal(\bar{N}_l^\mathrm{e})$ operations, which could
spoil the efficiency gained by the use of global basis functions if $\bar{N}_l^\mathrm{e} \gg k$. We propose a hyperreduction
approach in Section~\ref{sec:hpr} to circumvent the potential bottleneck associated with evaluating $\rbm_k^{\mathrm{gg}}$.
\end{remark}
%

\subsection{Quantity of interest}
\label{sec:mthd:qoi}
In the PDE setting, we are usually interested in quantities integrated over the volume or boundary. We write a general quantity of interest in terms of a volume and boundary term as
\begin{equation}
	q(v,\mu) \coloneqq \int_\Omega q_v(v,\nabla v; \mu) \, dV  + \int_{\partial\Omega} q_b(v,\nabla v; \mu, n) \, dS
\end{equation}
where $q_v : \Rbb^m \times \Rbb^{m\times d} \times \Dcal \rightarrow \Rbb$ and $q_b : \Rbb^m \times \Rbb^{m\times d} \times \Dcal \times \Rbb^d \rightarrow \Rbb$ are the volumetric and boundary integrands, respectively. 
We define the algebraic representation, $J_{l,h,p,k}:\Rbb^{N_l+k}\times\Dcal\rightarrow\Rbb$, of the quantity of interest as
\begin{equation}
	J_{l,h,p,k}(\Ubm, \mu) \coloneqq q(U_{l,h,p,k}(\Ubm),\mu).
	\label{eqn:JH}
\end{equation}


\section{Hyperreduction}
\label{sec:hpr}
In this section, we propose an approach to accelerate integral evaluations over the global region $\bar\Omega_l$, which is a potential bottleneck in the evaluation of the CG-GL residual or quantities of interest as it scales with $\Ocal(\bar{N}_l^\mathrm{e})$. To eliminate this bottleneck, we adapt the EQP method \cite{yano2019discontinuous} to the CG-GL framework in such a way that avoids retraining the weights when the local-global partition is adapted, e.g., to enrich the approximation space based on error estimates. 

\subsection{Empirical quadrature procedure over subdomains}
\label{sec:hpr:eqp}
If EQP is directly applied to a mesh of $\bar\Omega_l$, the quadrature weights would need to be recomputed whenever the local-global partition is modified, which prevents the CG-GL method from being efficiently embedded in an adaptive framework. To avoid this issue, we partition the domain $\Omega$ into $P$ non-overlapping subdomains $\Pcal_1,...,\Pcal_{P}\subset\Omega$ such that $\Pcal_1\cup\Pcal_2\cup\dots\cup\Pcal_{P}= \Omega$. Let $\Ecal_h^{\Pcal_i}$ be a mesh of the patch $\Pcal_i$ with ordered elements $\Ecal_h^{\Pcal_i} = \{\Pcal_{i,e}\}_{e=1}^{N_{\Pcal_i}}$, where $N_{\Pcal_i} = |\Ecal_h^{\Pcal_i}|$ is the number of elements in the mesh of patch $\Pcal_i$ and $\Pcal_{i,e} \subset\Pcal_i$ is the $e$th element of the mesh. 
Then we define the global domain as the union of a subset of the patches $\bar\Omega_l = \bigcup_{i \in \Ical} \Pcal_i$, where $\Ical \subset \{1,\dots,P\}$ defines which patches make up the global domain and $\Omega_l = \Omega\setminus\bar\Omega_l$. The meshes of the local and global domains can then be defined in terms of the patch meshes as
\begin{equation}
 \Ecal_{l,h} = \bigcup_{i\in\{1,\dots,P\}\setminus\Ical} \Ecal_h^{\Pcal_i}, \qquad
 \bar\Ecal_l = \bigcup_{i\in\Ical} \Ecal_h^{\Pcal_i}.
\end{equation}
Based on these definitions, the global residual can be written as
\begin{equation}
 \rbm_k^{\mathrm{gg}}(\betabold; \mu) = \sum_{i \in \Ical} \hat\rbm_k^{\mathrm{gg}}(\betabold; \mu, \Pcal_i),
\end{equation}
where $\hat\rbm_k^{\mathrm{gg}}(\betabold; \mu, \Pcal_i)$ is the contribution of patch $\Pcal_i$ to the residual (\ref{eqn:cggl_res}),
which can in turn be expanded as in terms of elemental contributions as
\begin{equation} \label{eqn:cggl_res_hrom0}
  \hat\rbm_k^{\mathrm{gg}}(\betabold; \mu,\Pcal_i) = \sum_{e=1}^{N_{\Pcal_i}} \hat\rbm_k^{\mathrm{gg}}(\betabold; \mu,\Pcal_{i,e}).
\end{equation}

Finally, we apply the EQP approach of \cite{yano2019discontinuous} to each patch individually. That is, for a given patch $\Pcal_i$, we approximate the global residual contribution with the weighted residual $\tilde\rbm_k^{\mathrm{gg}} : \Rbb^k \times\Dcal\times 2^{\Omega}\times\Rbb_{\geq 0}^{N_{\Pcal_i}}\rightarrow\Rbb^k$ as
\begin{equation} \label{eqn:cggl_res_hrom}
 \hat\rbm_k^{\mathrm{gg}}(\betabold;\mu,\Pcal_i) \approx \tilde\rbm_k^{\mathrm{gg}}(\betabold;\mu,\Pcal_i,\rhobold_i) \coloneqq
 \sum_{e=1}^{N_{\Pcal_i}} (\rhobold_i)_e \hat\rbm_k^{\mathrm{gg}}(\betabold; \mu,\Pcal_{i,e}),
\end{equation}
where $\rhobold_i \in \Rbb_{\geq 0}^{N_{\Pcal_i}}$ is a vector of elemental weights. For any element $\Pcal_{i,e}\in\Ecal_h^{\Pcal_i}$ where $(\rhobold_i)_e = 0$, the term in the summation can be skipped. Thus, if the vector of weights is highly sparse, computational efficiency is achieved. The advantage of this approach is that adaptations to the local-global decomposition (Section~\ref{sec:errest}) can be accomplished by moving a patch $\Pcal_i$ from $\bar\Omega_l$ to $\Omega_l$. Because empirical quadrature rules are associated with \textit{patches}, the new empirical quadrature rule will come directly from the quadrature rules associated with its remaining patches $\Ical \setminus \{i\}$.

\subsection{Training}
\label{sec:hpr:train}
To promote sparsity of the weight vectors and ensure accuracy of the hyperreduced residuals for a patch $\Pcal_i$, we train the EQP weights by solving an $l_1$ minimization problem that includes constraints requiring the global residual in \eqref{eqn:cggl_res_hrom0} and hyperreduced global residual in \eqref{eqn:cggl_res_hrom} to be sufficiently close on some training set. Let $\Xi\subset\Dcal$ be a collection of EQP training parameters and define the EQP weights $\rhobold_{i}^\star$ as the solution of the following linear program
\begin{equation} \label{eqn:eqp_solve}
	\rhobold_{i}^\star   = \argoptunc{\rhobold\in\Ccal_i  }{\sum_{e=1}^{N_{\Pcal_i}} (\rhobold)_e     }
\end{equation}
where $\Ccal_i\coloneqq \Ccal^{  \mathrm{nn}  }_i \cap \Ccal_i^{\mathrm{dv}} \cap \Ccal_i^{\mathrm{rp}}$ with 
\begin{equation}
	\begin{split}
		&\Ccal^{  \mathrm{nn}  }_i\coloneqq \Bigg\{\rhobold\in\Rbb^{N_{\Pcal_i}}\Bigg|  (\rhobold)_e \geq 0,  \quad e = 1,\dots,N_{\Pcal_i}  \Bigg\}, \\
		&\Ccal^{\mathrm{dv}}_i \coloneqq\Bigg\{ \rhobold\in\Rbb^{N_{\Pcal_i}}\Bigg|\; \Bigg|  
		\sum_{e=1}^{N_{\Pcal_i}}\Bigg((\rhobold)_e|\Omega_{e}^{\Pcal_i}| - |\Omega_{e}^{\Pcal_i}|\Bigg)
		\Bigg|    <\delta_{\mathrm{dv}}  \Bigg\},   \\
		&\Ccal^{\mathrm{rp}}_{i}\coloneqq\Bigg\{ \rhobold\in\Rbb^{N_{\Pcal_i}}\Bigg|\;
		\big|\big|\tilde{\rbm}_k^{\mathrm{gg}}\big(\betabold^\star(\mu); \mu, \Pcal_i, \rhobold\big)  -  \hat\rbm_k^{\mathrm{gg}}\big(\betabold^\star(\mu); \mu, \Pcal_i\big)\big|\big|_\infty<\delta_{\mathrm{rp}},\forall \mu\in\Xi
		\Bigg\},
	\end{split}
\end{equation} 
$\delta_\mathrm{dv} > 0$ and $\delta_\mathrm{rp}>0$ are user-defined tolerances, and $\betabold^\star(\mu)$ is the solution defined in \eqref{eqn:cggl_eqn} with $\Omega = \bar\Omega_l$ and $N_l = 0$ (only global basis functions). In this work we take $\delta_\mathrm{dv} = \delta_\mathrm{rp} = 10^{-8}$ and choose $\Xi$ to be the parameter set used to construct the global basis functions.

\subsection{Hyperreduced CG-GL formulation}
\label{sec:hpr:form}
Finally, we define a weighted version of the CG-GL residual $\tilde\rbm_{l,h,p,k} : \Rbb^{N_l+k}\times\Dcal\rightarrow\Rbb^{N_l+k}$ by replacing the residual in (\ref{eqn:cggl_eqn}) with its weighted version
\begin{equation}		
	\tilde\rbm_{l,h,p,k}(\Ubm;\mu)= (\rbm_{l,h,p,k}^{\mathrm{ll}}(\Ubm;\mu), \rbm_{l,h,p,k}^{\mathrm{gl}}(\Ubm;\mu)+\sum_{i\in \Ical} \tilde\rbm_k^{\mathrm{gg}}(\betabold;\mu,\Pcal_i,\rhobold_i^\star)),
\end{equation} 
where $\rhobold_1^\star,\dots,\rhobold_{P}^\star$ are the trained, sparse weight vectors from (\ref{eqn:eqp_solve}). Then, the weighted CG-GL problem reads: given $\mu\in\Dcal$, find $\tilde\Ubm^\star\in\Rbb^{N_l+k}$ such that
\begin{equation} \label{eqn:cggl_eqn_eqp}
	\tilde\rbm_{l,h,p,k}(\tilde\Ubm^\star;\mu)=0
\end{equation}
and the continuous solution is obtained as $U_{l,h,p,k}(\tilde\Ubm^\star)$.

\section{Goal-oriented error estimation and solution adaptation}
\label{sec:errest}
In this section we extend standard error estimation and adaptation technologies to the CG-GL setting to yield a fully adaptive method. We detail output-based DWR \textit{a posteriori} error estimation (Section~\ref{sec:errest:err}) including choice of the fine space (Section~\ref{sec:errest:fine}), localization of the error estimates to subdomains (patches) (Section~\ref{sec:errest:local}), and a strategy to refine the CG-GL discretization from the error estimates (\ref{sec:errest:refine}). For brevity, we drop the dependence of all terms on $\mu$ in this section.

\subsection{Error estimation}
\label{sec:errest:err}
We begin with standard construct of DWR error estimates using enriched spaces in the context of the CG-GL discretization. To this end, given two function spaces $\Ucal_{l_1,h_1,p_1, k_1}$, $ \Ucal_{l_2,h_2,p_2,k_2}$ such that $\Ucal_{l_1,h_1,p_1,k_1}\subseteq \Ucal_{l_2,h_2,p_2,k_2}$, which requires
\begin{equation}
 h_2 \leq h_1, \quad
 p_2 \geq p_1, \quad
 \Wcal_{k_1} \subset \Wcal_{k_2}, \quad
 \Omega^{l_1} \subset \Omega^{l_2}.
\end{equation}
For notational brevity, we define $\mathtt{H} = (l_1,h_1,p_1, k_1)$ and $\mathtt{h} = (h_2,p_2, l_2,k_2)$, and write
$\Ucal_\mathtt{H}$ and $\Ucal_\mathtt{h}$ in place of $\Ucal_{l_1,h_1,p_1, k_1}$ and $ \Ucal_{h_2,p_2, l_2,k_2}$.
Then we define a prolongation operator, 	$I_\mathtt{H}^\mathtt{h}:\Rbb^{N_{\mathrm{lb}_1}+k_1}\rightarrow\Rbb^{N_{\mathrm{lb}_2}+k_2}$ such that
\begin{equation}
U_\mathtt{h}(I_\mathtt{H}^\mathtt{h}(\Ubm))=U_\mathtt{H}(\Ubm).
\end{equation}
for any $\Ubm\in\Rbb^{N_{\mathrm{lb}_1}+k_1}$. That is, the prolongation operator converts the generalized coordinates
$\Ubm$ of a function in the coarse space $\Ucal_\mathtt{H}$ to the generalized coordinates of the same function
in the fine space $\Ucal_\mathtt{h}$.

Next, we define the mean of residual Jacobian $\overline{\partial{\rbm}}_{l,h,p,k}:\Rbb^{N_{\mathrm{lb}}+k}\times\Rbb^{N_{\mathrm{lb}}+k}\rightarrow \Rbb^{(N_{\mathrm{lb}}+k)\times(N_{\mathrm{lb}}+k)}$,
\begin{equation}
	\overline{\partial\rbm}_{l,h,p,k}(\Ubm_1,\Ubm_2) = \int_0^1 \pder{\rbm_{l,h,p,k}}{\Ubm}\big(\Ubm_1+ \theta(\Ubm_2-\Ubm_1) \big)\,d\theta
	\label{eqn:mean_r}
\end{equation} 
and the mean of QoI Jacobian $\bar{J}_{l,h,p,k}:\Rbb^{N_l+k}\times\Rbb^{N_l+k}\rightarrow \Rbb^{1\times(N_l+k)}$,
\begin{equation}
	\overline{\partial J}_{l,h,p,k}(\Ubm_1,\Ubm_2) = \int_0^1 \pder{J_{l,h,p,k}}{\Ubm}\big(\Ubm_1+ \theta(\Ubm_2-\Ubm_1)\big)\,d\theta.
		\label{eqn:mean_j}
\end{equation} 
With these definitions, the DWR \textit{a posteriori} error takes the form
\begin{equation}
	\begin{split}
		\overline{\partial\rbm}_{\mathtt{h}}\bigl(\Ubm_\mathtt{h}^\star, I_\mathtt{H}^\mathtt{h}(\Ubm_\mathtt{H}^\star) \bigr)^T 
		 \psibold_\mathtt{h} =- \overline{\partial J}_\mathtt{h} \bigl( \Ubm_\mathtt{h}^\star, I_\mathtt{H}^\mathtt{h}(\Ubm_\mathtt{H}^\star) \bigr)^T, \\
		J_\mathtt{H}(\Ubm_\mathtt{H}^\star) - J_\mathtt{h}(\Ubm_\mathtt{h}^\star)  = -\psibold_\mathtt{h}^\mathrm{mv}\cdot\rbm_\mathtt{h}\bigr( I_\mathtt{H}^\mathtt{h}(\Ubm_\mathtt{H}^\star) \bigl),
	\end{split}
	\label{eqn:dwr_truth}
\end{equation}
where $\Ubm_\mathtt{H}^\star$ and $\Ubm_\mathtt{h}^\star$ are the coarse and fine space, respectively, CG-GL solutions (\ref{eqn:cggl_eqn}), and $\psibold_\mathtt{h}^\mathrm{mv}$ is the adjoint solution used to weight the entries of the residual to define the error. This is the \textit{mean value} method and provides an exact formula to compute the error $J_\mathtt{H}(\Ubm_\mathtt{H}^\star) - J_\mathtt{h}(\Ubm_\mathtt{h}^\star)$; however, it is inefficient owing to the dependence on the fine state $\Ubm_\mathtt{h}^\star$ and the repeated Jacobian evaluations required to evaluate the mean values. To improve the efficiency of error estimation, we use the following DWR error estimate
\begin{equation}
	\begin{split}
		\pder{\rbm_\mathtt{h}}{\Ubm}\bigl(I_\mathtt{H}^\mathtt{h}(\Ubm_\mathtt{H}^\star) \bigr)^T 
		 \psibold_{\mathtt{H},\mathtt{h}} =- \pder{J_\mathtt{h}}{\Ubm}\bigl(I_\mathtt{H}^\mathtt{h}(\Ubm_\mathtt{H}^\star) \bigr)^T, \\
		J_\mathtt{H}(\Ubm_\mathtt{H}^\star) - J_\mathtt{h}(\Ubm_\mathtt{h}^\star)  \approx E_{\mathtt{H},\mathtt{h}} \coloneqq -\psibold_{\mathtt{H},\mathtt{h}}\cdot\rbm_\mathtt{h}\bigr( I_\mathtt{H}^\mathtt{h}(\Ubm_\mathtt{H}^\star) \bigl),
	\end{split}
	\label{eqn:errorestimation}
\end{equation}
where the mean value integrals were approximated with a one-point quadrature rule at $\theta = 1$, $\psibold_{\mathtt{H},\mathtt{h}}$ is the fine space adjoint solution corresponding to the prolonged coarse state, and $\Ecal_{\mathtt{H},\mathtt{h}} \in \Rbb$ is the DWR estimate of the error $J_\mathtt{H}(\Ubm_\mathtt{H}^\star) - J_\mathtt{h}(\Ubm_\mathtt{h}^\star)$. Similar to the primal CG-GL solution, the adjoint can be partitioned into local and global contributions $\psibold_{\mathtt{H},\mathtt{h}} = (\psibold_{\mathtt{H},\mathtt{h}}^l,\psibold_{\mathtt{H},\mathtt{h}}^g)$, where $\psibold_{\mathtt{H},\mathtt{h}}^l \in \Rbb^{N_{\mathrm{lb}_2}}$ and $\psibold_{\mathtt{H},\mathtt{h}}^g\in\Rbb^{k_2}$,
which means the error estimate in (\ref{eqn:errorestimation}) can be split into three terms when using the decomposition of the residual in (\ref{eqn:cggl_res2})
\begin{equation} \label{eqn:errest2}
 \psibold_{\mathtt{H},\mathtt{h}}\cdot\rbm_\mathtt{h}\bigr( I_\mathtt{H}^\mathtt{h}(\Ubm_\mathtt{H}^\star) \bigl) =
 \psibold_{\mathtt{H},\mathtt{h}}^l\cdot\rbm_\mathtt{h}^{\mathrm{ll}}\bigr( I_\mathtt{H}^\mathtt{h}(\Ubm_\mathtt{H}^\star) \bigl) +
 \psibold_{\mathtt{H},\mathtt{h}}^g\cdot\rbm_\mathtt{h}^{\mathrm{gl}}\bigr( I_\mathtt{H}^\mathtt{h}(\Ubm_\mathtt{H}^\star) \bigl) +
 \psibold_{\mathtt{H},\mathtt{h}}^g\cdot\rbm_\mathtt{h}^{\mathrm{gg}}\bigr( I_\mathtt{H}^\mathtt{h}(\Ubm_\mathtt{H}^\star) \bigl).
\end{equation}

\subsection{Choice of fine space}
\label{sec:errest:fine}
The error estimation framework in the previous section only requires $\Ucal_\mathtt{H} \subset \Ucal_\mathtt{h}$; however, the quality of the error estimate depends on the choice of the fine space $\Ucal_\mathtt{h}$. To ensure we obtain a reasonable error estimate in the entire domain, we choose $l_2$ such that $\Omega_{l_2} = \Omega$, i.e., there is no region where only global basis functions are used. This will facilitate adaptation of the global-local decomposition. For simplicity, we choose our fine space in the local region $\Omega_{l_1}\subset\Omega$ via uniform $h$-refinement ($h_2 < h_1$) for fixed polynomial degree ($p_2 = p_1$) and global basis function ($\Wcal_{k_1}=\Wcal_{k_2}$). While other choices are possible, we show this approach is effective and straightforward.

\subsection{Error localization}
\label{sec:errest:local}
To identify the regions of the domain that require refinement, we localize the error estimate to elements of the fine mesh using (\ref{eqn:cggl_res_sum}).
Owing to our choice of the fine space that takes $\Omega_{l_2} = \Omega$, the expression in (\ref{eqn:errest2}) reduces to
\begin{equation}
 \psibold_{\mathtt{H},\mathtt{h}}\cdot\rbm_\mathtt{h}\bigr( I_\mathtt{H}^\mathtt{h}(\Ubm_\mathtt{H}^\star) \bigl) =
 \sum_{K\in\Ecal_{h_2}^{l_2}} \left(\psibold_{\mathtt{H},\mathtt{h}}^l\cdot\rbm_\mathtt{h}^{\mathrm{ll}}\bigr( I_\mathtt{H}^\mathtt{h}(\Ubm_\mathtt{H}^\star); K \bigl) + \psibold_{\mathtt{H},\mathtt{h}}^g\cdot\rbm_\mathtt{h}^{\mathrm{gl}}\bigr( I_\mathtt{H}^\mathtt{h}(\Ubm_\mathtt{H}^\star); K \bigl)\right).
\end{equation}
Then for each element $K\in\Ecal_{l_2,h_2}$, the local error estimate $\eta_K \in \Rbb$ is naturally
\begin{equation}
 \eta_K = \psibold_{\mathtt{H},\mathtt{h}}^l\cdot\rbm_\mathtt{h}^{\mathrm{ll}}\bigr( I_\mathtt{H}^\mathtt{h}(\Ubm_\mathtt{H}^\star); K \bigl) + \psibold_{\mathtt{H},\mathtt{h}}^g\cdot\rbm_\mathtt{h}^{\mathrm{gl}}\bigr( I_\mathtt{H}^\mathtt{h}(\Ubm_\mathtt{H}^\star); K \bigl).
\end{equation}
From these estimates on the fine mesh, we obtain estimates over a region of the domain $\Rcal\subset\Omega$ as
\begin{equation}
 \bar\eta_{\Rcal} = \sum_{\substack{K \in \Ecal_{l_2,h_2} \\ K \subset \Rcal}} \eta_K.
\end{equation}
Then the error estimate over an element $K$ of the coarse mesh is $\bar\eta_K$ and over a patch $\Pcal_i$ 
is $\bar\eta_{\Pcal_i}$.

\subsection{Adaptive refinement}
\label{sec:errest:refine}
With the localized error estimates available, we propose an approach to use these estimates to adapt the CG-GL trial space. There are four potential ways to improve resolution of the trial space: 1) reduce the size of elements in the local domain $h' < h_1$,
2) increase the polynomial degree in the local domain $p' > p_1$, 3) add global basis function $\Wcal_{k_1}\subset\Wcal_{k'}$,
or 4) reduce the portion of the domain where only global basis functions used $\Omega_{l_1} \subset \Omega_{l'}$. We will not consider the third option because new global basis functions usually are expensive to construct, e.g., from simulations on finer grids. Additionally, for simplicity, we will not consider refinements to the polynomial degree ($p' = p_1$), which leaves $h$-refinement of the local domain and modifications to the local-global decomposition as the available refinement mechanisms.

First, let us define $\Pcal^\star$ as the patch with the largest error estimate
\begin{equation}
 \Pcal^\star = \argmax_{\Pcal \in \{\Pcal_1,\dots,\Pcal_{P}\}} |\eta_\Pcal|.
\end{equation}
If $\Pcal^\star \subset \Omega\setminus \Omega^{l_1}$, the largest error lies in the global region and we increase resolution in this patch by enriching with local basis functions. In practice, this is achieved by moving this patch to the local domain, which leads to the following update to the local-global decomposition
\begin{equation}
 \Omega_{l'} =
 \begin{cases}
  \Omega_{l_1} & \Pcal^\star \subset \Omega_{l_1} \\
  \Omega_{l_1} \cup \Pcal^\star & \Pcal^\star \subset \Omega\setminus\Omega_{l_1}.
 \end{cases}
\end{equation}
On the other hand, if $\Pcal^\star$ lies in the local domain, we use $h$-adaption to refine elements of that patch, i.e.,
\begin{equation}
 \Ecal_{h'}^{\Pcal_i} =
 \begin{cases}
  \Ecal_{h_1}^{\Pcal_i} & \Pcal_i \neq \Pcal^\star \text{~or~} \Pcal_i = \Pcal^\star \subset \Omega\setminus \Omega_{l_1} \\
  \mathtt{refine}(\Ecal_{h_1}^{\Pcal_i}) & \Pcal_i = \Pcal^\star \subset \Omega_{l_1},
 \end{cases}
\end{equation}
where $\mathtt{refine}: 2^\Omega \rightarrow 2^\Omega$ applies a chosen $h$-refinement strategy to its argument. The result is a new trial space $\Ucal_{l',h',p_1,k_1}$ that is embedded in the original space, i.e., $\Ucal_{l',h',p_1,k_1} \subset \Ucal_{l_1,h_1,p_1,k_1}$. Once a single error estimation and refinement iteration is complete, it is repeated with $\mathtt{H}' = (l',h',p_1,k_1)$ in place of $\mathtt{H}$ until a maximum number of refinements are performed or the error estimate $E_{\mathtt{H},\mathtt{h}}$ is sufficiently small.

\begin{remark}
Because the empirical quadrature weights are trained patch-by-patch, they only need to be updated when $h$-refinement is performed on a patch. However, because the refinement strategy only performs $h$-refinement on patches in the local domain and empirical quadrature is only used in the global domain, the empirical quadrature weights never need to be retrained.
\end{remark}

\begin{remark}
The adaptive CG-GL framework can be initialized with any suitable global-local trial space. In this work, we initialize it with a traditional reduced-order model, i.e., $\bar\Omega_l = \Omega$, and a coarse mesh of each patch.
\end{remark}

\section{Numerical experiments}
\label{sec:rslt}
In this section, we investigate the performance of the proposed CG-GL method on two computational benchmark problems. For the first problem (Poisson's equation in a square domain), we study the convergence of the proposed CG-GL method as the reduced basis is refined, the performance of parameter interpolation and extrapolation with limited training, and the scalability of  EQP in the CG-GL setting. The second problem (incompressible flow over a backward-facing step) demonstrates the ability of the CG-GL method to provide accurate approximations with a limited sampling of the parameter space. 

For each problem, we will study the solution and QoI accuracy of the adaptive CG-GL method starting from a standard reduced-order model. The solution accuracy will be measured by the relative $L^2(\Omega)$ error with respect a ``truth'' solution $U_t(\,\cdot\,;\mu)$ computed on a highly refined mesh as
\begin{equation}
	e_\mathrm{sln}(\mu) \coloneq \sqrt{\frac{\int_\Omega||U_t(x;\mu) -   U_{\mathrm{H}_r}(x;\mu)||_2^2 \, dV}{\int_\Omega||U_t(x;\mu) ||_2^2 \, dV}},
\end{equation}
where $U_{\mathtt{H}_r}(\,\cdot\,;\mu)$ is the CG-GL solution after $r$ adaption iterations. The QoI accuracy will be measured by the relative error with respect to a ``truth'' QoI $q_t(\mu) \coloneqq J(U_t(,\cdot\,;\mu),\mu)$ as
\begin{equation}
	e_\mathrm{qoi}(\mu) \coloneq \frac{|q_{\mathrm{H}_r}(\mu)-q_t(\mu)|}{|q_t(\mu)|},
\end{equation} 
where $q_{\mathrm{H}_r}(\mu) = J(U_{\mathtt{H}_r}(\,\cdot\,;\mu))$ is the CG-GL QoI after $r$ adaptation iterations. The quality of the error estimate will be measured by the relative error with respect to the true QoI error as
\begin{equation}
	e_\mathrm{est}(\mu) \coloneqq \frac{\bigl||E_{\mathrm{H}_r,\mathrm{h}_r}|  -|q_{\mathrm{H}_r}(\mu)-q_t(\mu)|\bigr|}{|q_t(\mu)|}.
\end{equation} 


\subsection{Poisson's equation}
\label{sec:rslt:poi}
\begin{figure}
	\centering
	\begin{tikzpicture}
		\begin{groupplot}[
			group style={
				group size=4 by 1,
				horizontal sep=0.5cm
			},
			width=0.32\textwidth,
			axis equal image,
			xlabel={$x_1$},
			ylabel={$x_2$},
			xtick = {0.0, 4.0, 8.0},
			ytick = {0.0, 4.0, 8.0},
			xmin=0, xmax=8,
			ymin=0, ymax=8
			]
			\nextgroupplot[title={$\mu=(1 , 1,  4,     4)$}]
			\addplot graphics [xmin=0, xmax=8.0, ymin=0, ymax=8.0] {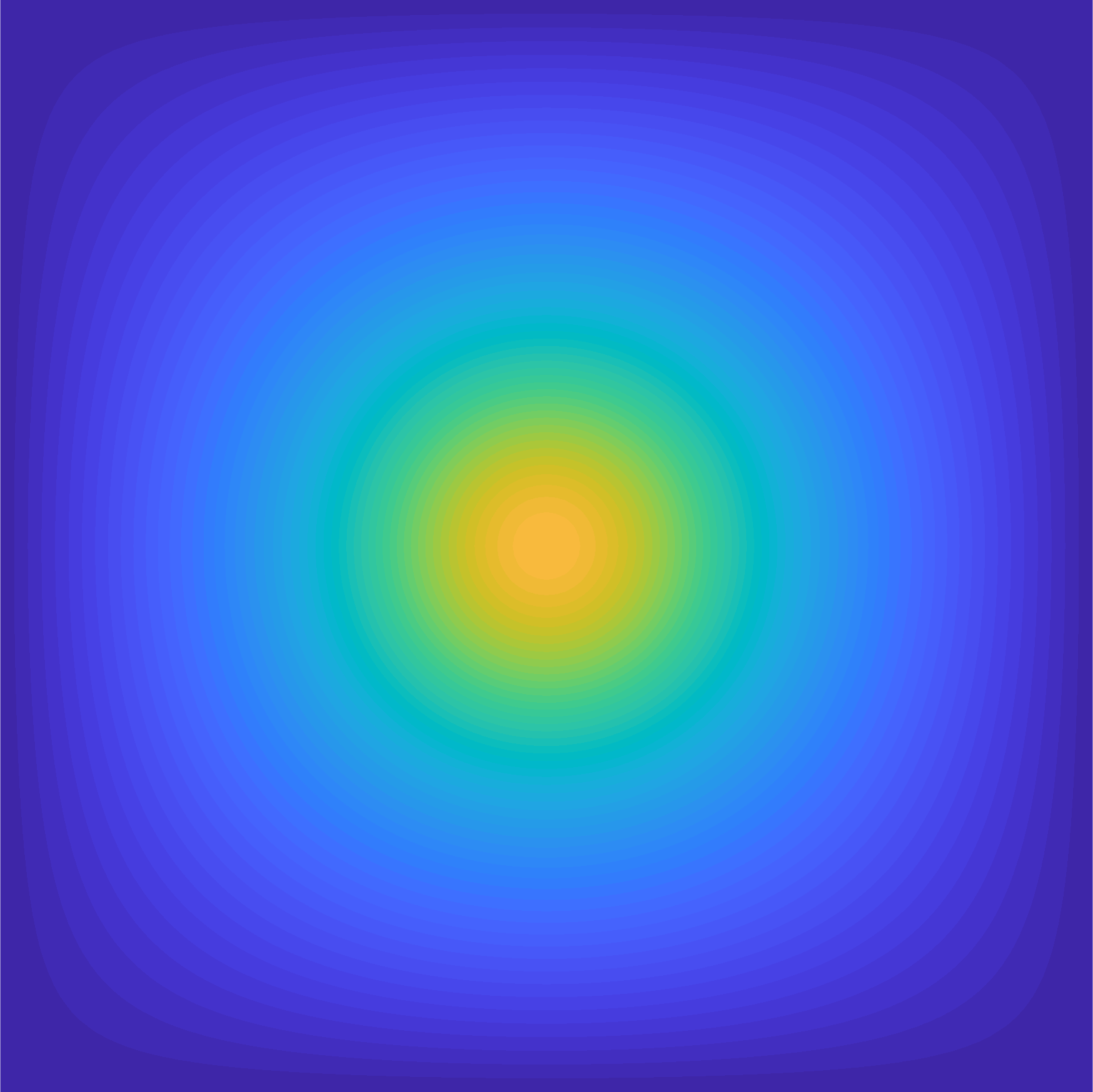};
			
			\nextgroupplot[title={$\mu=(3, 0.5, 3, 4, 4)$}, ylabel={}, ytick=\empty]
			\addplot graphics [xmin=0, xmax=8.0, ymin=0, ymax=8.0] {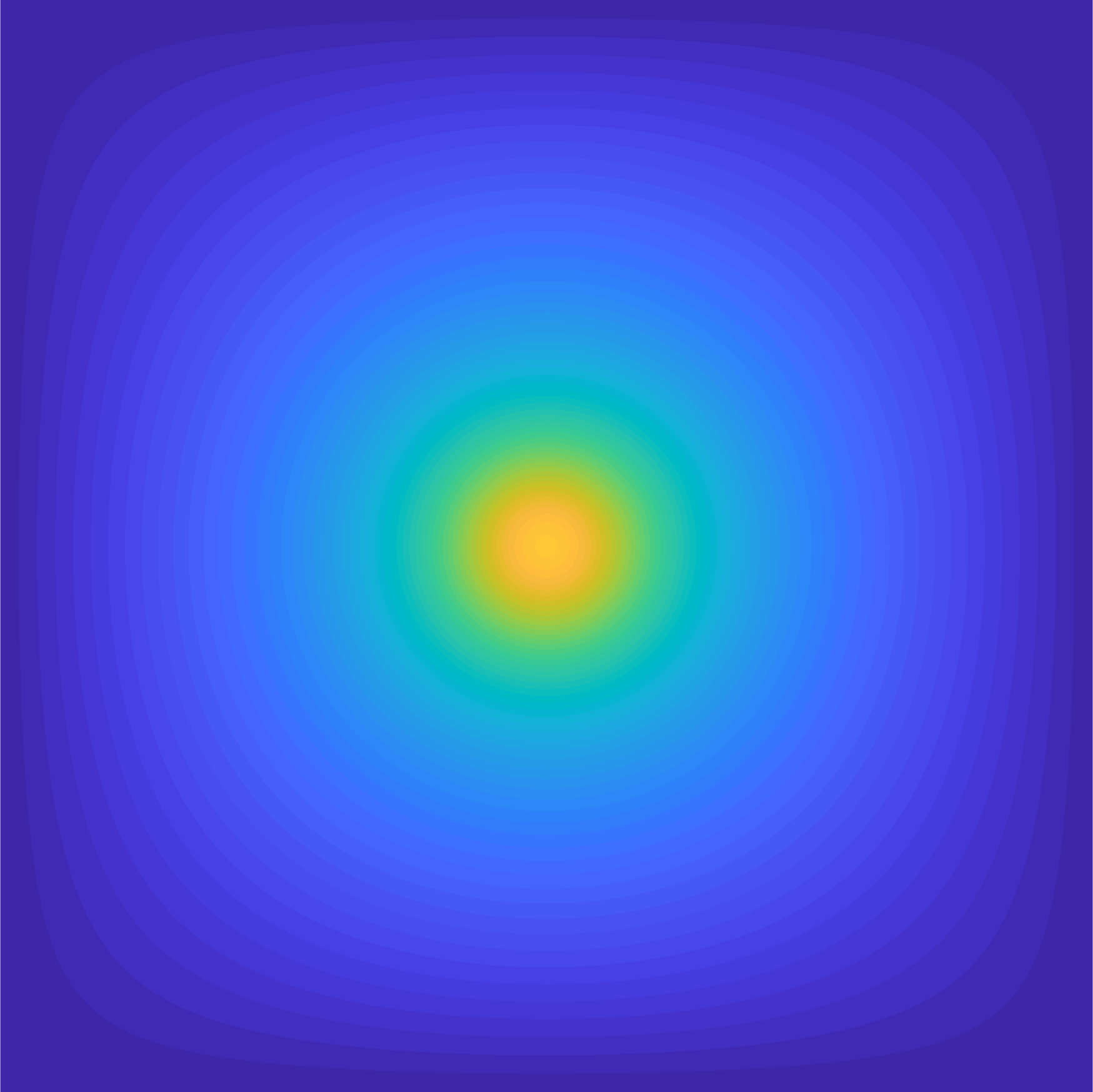};
			
			\nextgroupplot[title={$\mu=(2, 0.5, 4, 4)$}, ylabel={}, ytick=\empty]
			\addplot graphics [xmin=0, xmax=8.0, ymin=0, ymax=8.0] {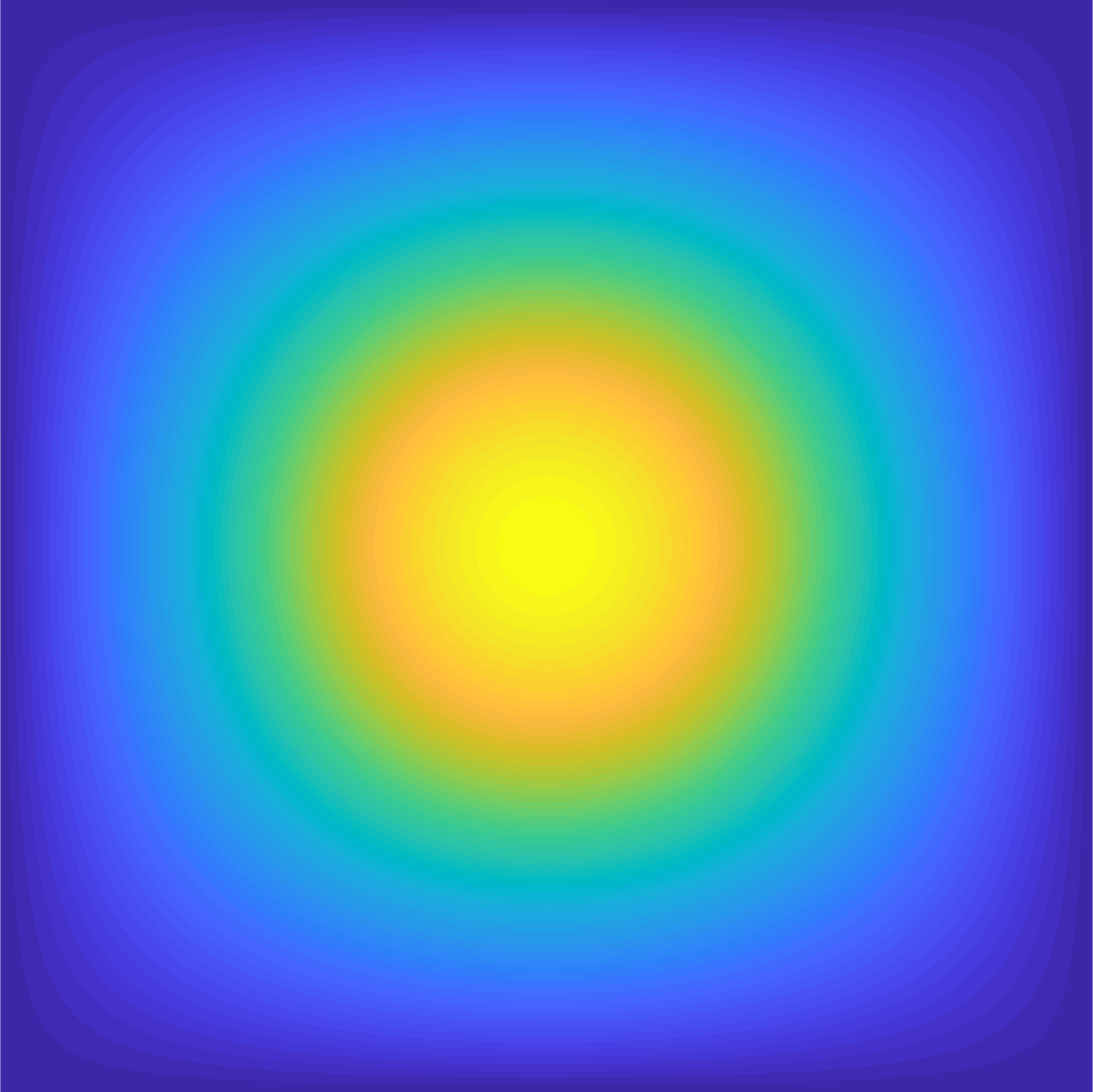};
			
			\nextgroupplot[title={$\mu=(4, 1, 0.5,0.5)$}, ylabel={}, ytick=\empty]
			\addplot graphics [xmin=0, xmax=8.0, ymin=0, ymax=8.0] {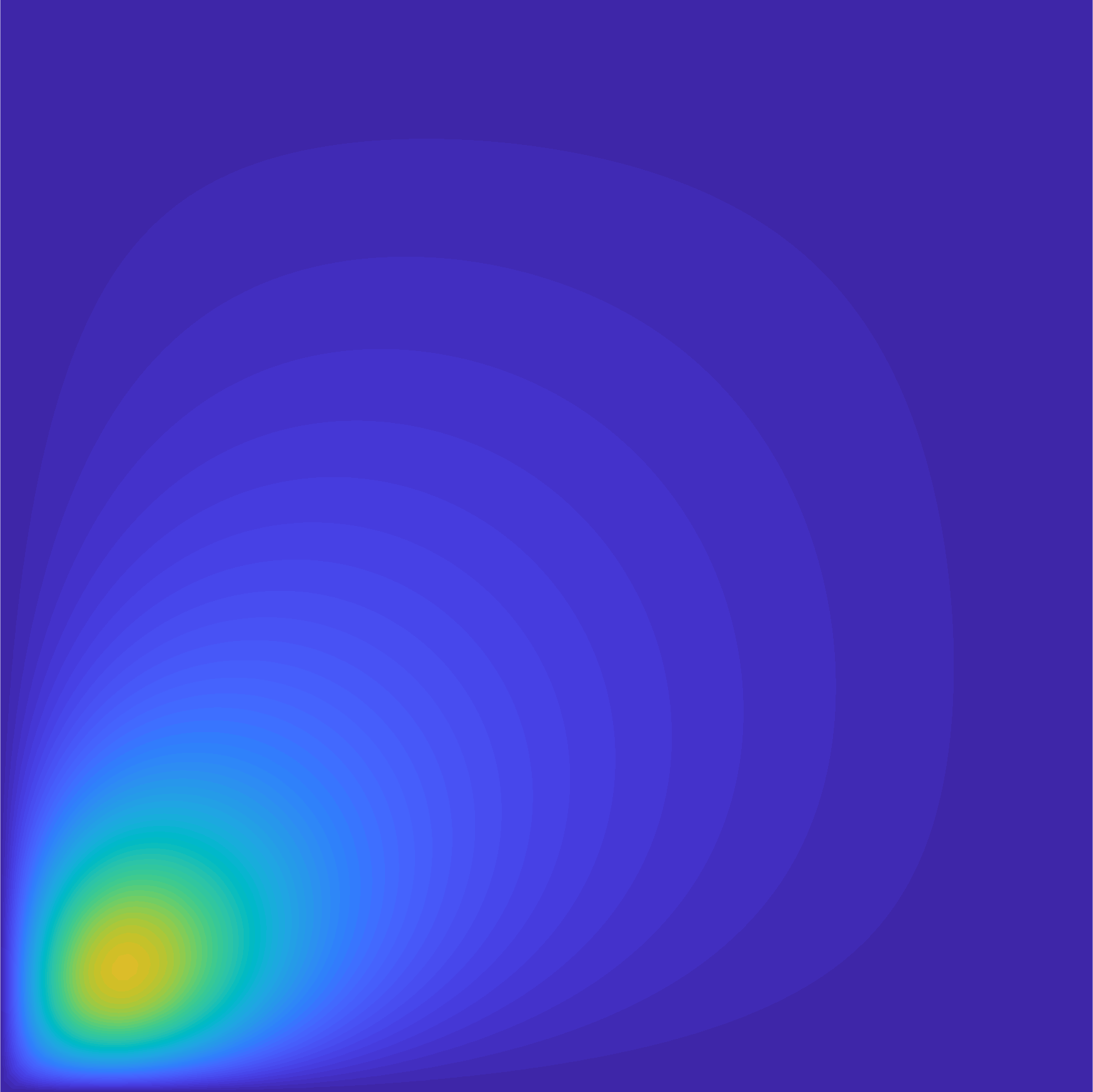};
		\end{groupplot}
	\end{tikzpicture}
	\colorbarMatlabParula{0}{0.2647}{0.5293}{0.7940}{1.0568}
	\caption{Solution of Poisson equation (\ref{eqn:poi}) for
		four different parameter configurations.}
	\label{fig:poi_hdm_sol_demo}
\end{figure}

We start by considering the parametrized linear Poisson equation in a two-dimensional domain $\Omega\coloneqq [0,8]^2$ 
\begin{equation} \label{eqn:poi}
 \Delta u = f\quad\text{in}\quad\Omega, \quad u=0\quad\text{on}\quad\partial\Omega,
\end{equation}
where $f: \Omega \times \Dcal \rightarrow \Rbb$ is the source term
\begin{equation}
f(x; (a,\sigma,c)) = a\exp\left(\frac{\norm{x-c}^2}{\sigma^2}\right),
\end{equation}
$u(x;\mu)\in\Rbb$ is the solution implicitly defined by (\ref{eqn:poi}), and $\mu\coloneqq (a, \sigma, c)\in\Dcal\coloneqq[1,10]\times[0.1,10]\times\Omega$ parameter vector. The reference solution is
obtained from a continuous Galerkin discretization on a structured grid of $5184$ quadratic
quadrilateral elements ($p=2$). Figure~\ref{fig:poi_hdm_sol_demo} shows the solution for four parameter configurations;
the first parameter ($a$) controls the magnitude of exponential source term, the second parameter ($\sigma$) controls the width of the source term peak, and the last two parameters $(c_1, c_2)$ controls the exponential kernel location. The QoI for this problem is the solution weighted by the exponential kernel
\begin{equation}
	q(v,\mu) = \int_\Omega v(x) f(x; \mu) \, dV.
\end{equation}

\subsubsection{Convergence under basis refinement, fixed training set}
\label{sec:rslt:poi:conv}
\begin{figure}
\centering
\begin{tikzpicture}
\begin{axis}[
width=0.75\textwidth,
xtick={0, 1, 2, 3, 4, 5, 6, 7, 8},
ytick={0, 1, 2, 3, 4, 5, 6, 7, 8},
xlabel=$c_1$,
ymax=8,
xmax=8,
ylabel=$c_2$,
xmin=0,
ymin=0]
\addplot [mark options={solid, thick}, mark=x, mark size=3, cyan, only marks]
coordinates {
( 2.00000000e-01,  2.37142857e+00)
( 2.00000000e-01,  5.62857143e+00)
( 1.28571429e+00,  2.00000000e-01)
( 1.28571429e+00,  3.45714286e+00)
( 1.28571429e+00,  6.71428571e+00)
( 2.37142857e+00,  1.28571429e+00)
( 3.45714286e+00,  2.00000000e-01)
( 3.45714286e+00,  7.80000000e+00)
( 4.54285714e+00,  6.71428571e+00)
( 5.62857143e+00,  1.28571429e+00)
( 6.71428571e+00,  2.00000000e-01)
( 6.71428571e+00,  3.45714286e+00)
( 6.71428571e+00,  6.71428571e+00)
( 7.80000000e+00,  1.28571429e+00)
( 7.80000000e+00,  4.54285714e+00)
( 7.80000000e+00,  7.80000000e+00)};\label{line:poi_para_ext}

\addplot [mark options={solid, thick}, mark=+, mark size=3, red, only marks]
coordinates {
( 2.37142857e+00,  2.37142857e+00)
( 2.37142857e+00,  3.45714286e+00)
( 2.37142857e+00,  4.54285714e+00)
( 2.37142857e+00,  5.62857143e+00)
( 3.45714286e+00,  2.37142857e+00)
( 3.45714286e+00,  3.45714286e+00)
( 3.45714286e+00,  4.54285714e+00)
( 3.45714286e+00,  5.62857143e+00)
( 4.54285714e+00,  2.37142857e+00)
( 4.54285714e+00,  3.45714286e+00)
( 4.54285714e+00,  4.54285714e+00)
( 4.54285714e+00,  5.62857143e+00)
( 5.62857143e+00,  2.37142857e+00)
( 5.62857143e+00,  3.45714286e+00)
( 5.62857143e+00,  4.54285714e+00)
( 5.62857143e+00,  5.62857143e+00)};\label{line:poi_para_int}

\addplot [mark options={solid, thick}, mark=*, mark size=3, blue, only marks]
coordinates {
( 2.00000000e+00,  2.00000000e+00)
( 2.00000000e+00,  4.00000000e+00)
( 2.00000000e+00,  6.00000000e+00)
( 4.00000000e+00,  2.00000000e+00)
( 4.00000000e+00,  4.00000000e+00)
( 4.00000000e+00,  6.00000000e+00)
( 6.00000000e+00,  2.00000000e+00)
( 6.00000000e+00,  4.00000000e+00)
( 6.00000000e+00,  6.00000000e+00)};\label{line:poi_para_train}

\end{axis}
\end{tikzpicture}
\caption{Training set $\Ccal_{\mathrm{train}}$ (\ref{line:poi_para_train}), interploated testing set $\Ccal_{\mathrm{interp}}$ (\ref{line:poi_para_int}), extrapolated testing set $\Ccal_{\mathrm{extrap}}$ (\ref{line:poi_para_ext}) used for the convergence study in Section~\ref{sec:rslt:poi:conv}.}
\label{fig:poi_train_para_0}
\end{figure}
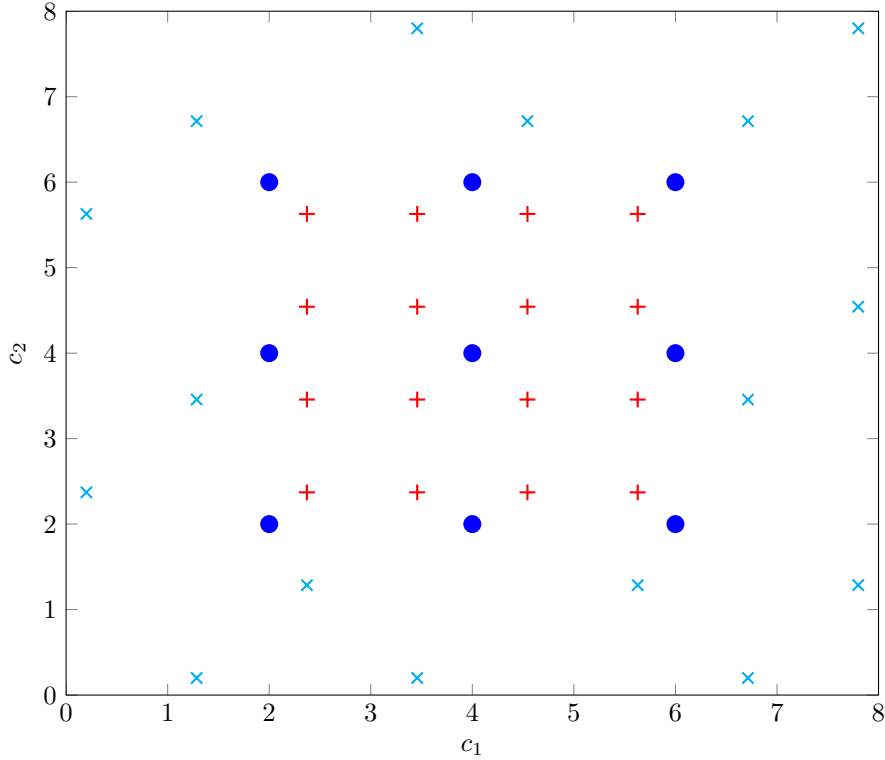

First, we demonstrate that the addition of pre-trained global basis functions leads to limited improvement of a traditional ROM in terms
of prediction accuracy, especially for extrapolated prediction. On the other hand, the CG-GL method provides a framework that can drive the error to small values with additional computational cost. For this study, we consider a reduced parameter space $\tilde{\Dcal}$ consisting of a two-dimensional subset of $\Dcal$ where only the third and fourth parameters $(c_1, c_2)$ vary
\begin{equation}
	\tilde{\Dcal}\coloneqq \left\{(1,1,c_1,c_2)  \suchthat (c_1,c_2)\in   [0,8]^2 \right\},
\end{equation}
which will cause the peak of the primary feature to move throughout the domain. Now, we define the collection of training parameters
$\Dcal_\mathrm{train}\subset\tilde\Dcal$ and testing parameters $\Dcal_\mathrm{test} \coloneqq \Dcal_\mathrm{interp}\cup\Dcal_\mathrm{extrap}$ as
\begin{equation}
 \begin{aligned}
  \Dcal_\mathrm{train} &\coloneqq \left\{(1,1,c_1,c_2) \suchthat (c_1,c_2) \in \Ccal_\mathrm{train}\right\}, \\
  \Dcal_\mathrm{interp} &\coloneqq \left\{(1,1,c_1,c_2) \suchthat (c_1,c_2) \in \Ccal_\mathrm{interp}\right\}, \\
  \Dcal_\mathrm{extrap} &\coloneqq \left\{(1,1,c_1,c_2) \suchthat (c_1,c_2) \in \Ccal_\mathrm{extrap}\right\},
 \end{aligned}
\end{equation}
where $\Ccal_\mathrm{train}$, $\Ccal_\mathrm{interp}$, $\Ccal_\mathrm{extrap}$ are shown in Figure~\ref{fig:poi_train_para_0}.

\begin{figure}
	\centering
	\begin{tikzpicture}
		\begin{groupplot}[
			group style={
				group size=3 by 3,
				horizontal sep=0.5cm,
				vertical sep=0.5cm
			},
			width=0.35\textwidth,
			axis equal image,
			xlabel={$x_1$},
			ylabel={$x_2$},
			xtick = {0.0, 4.0, 8.0},
			ytick = {0.0, 4.0, 8.0},
			xmin=0, xmax=8,
			ymin=0, ymax=8
			]
			\nextgroupplot[ ylabel={}, ytick=\empty, xlabel={}, xtick=\empty]
			\addplot graphics [xmin=0, xmax=8.0, ymin=0, ymax=8.0] {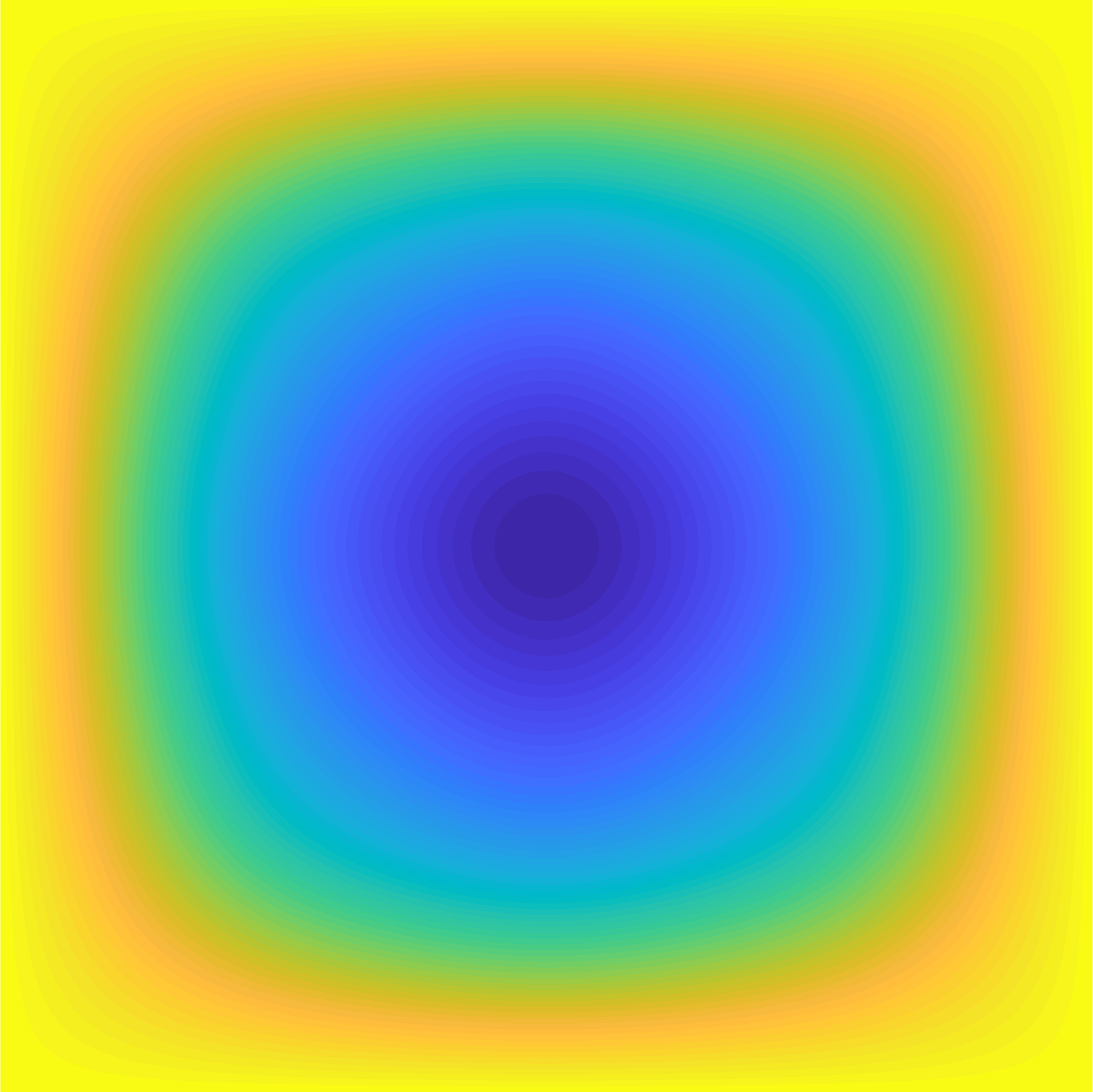};
			
			\nextgroupplot[ ylabel={}, ytick=\empty, xlabel={}, xtick=\empty]
			\addplot graphics [xmin=0, xmax=8.0, ymin=0, ymax=8.0] {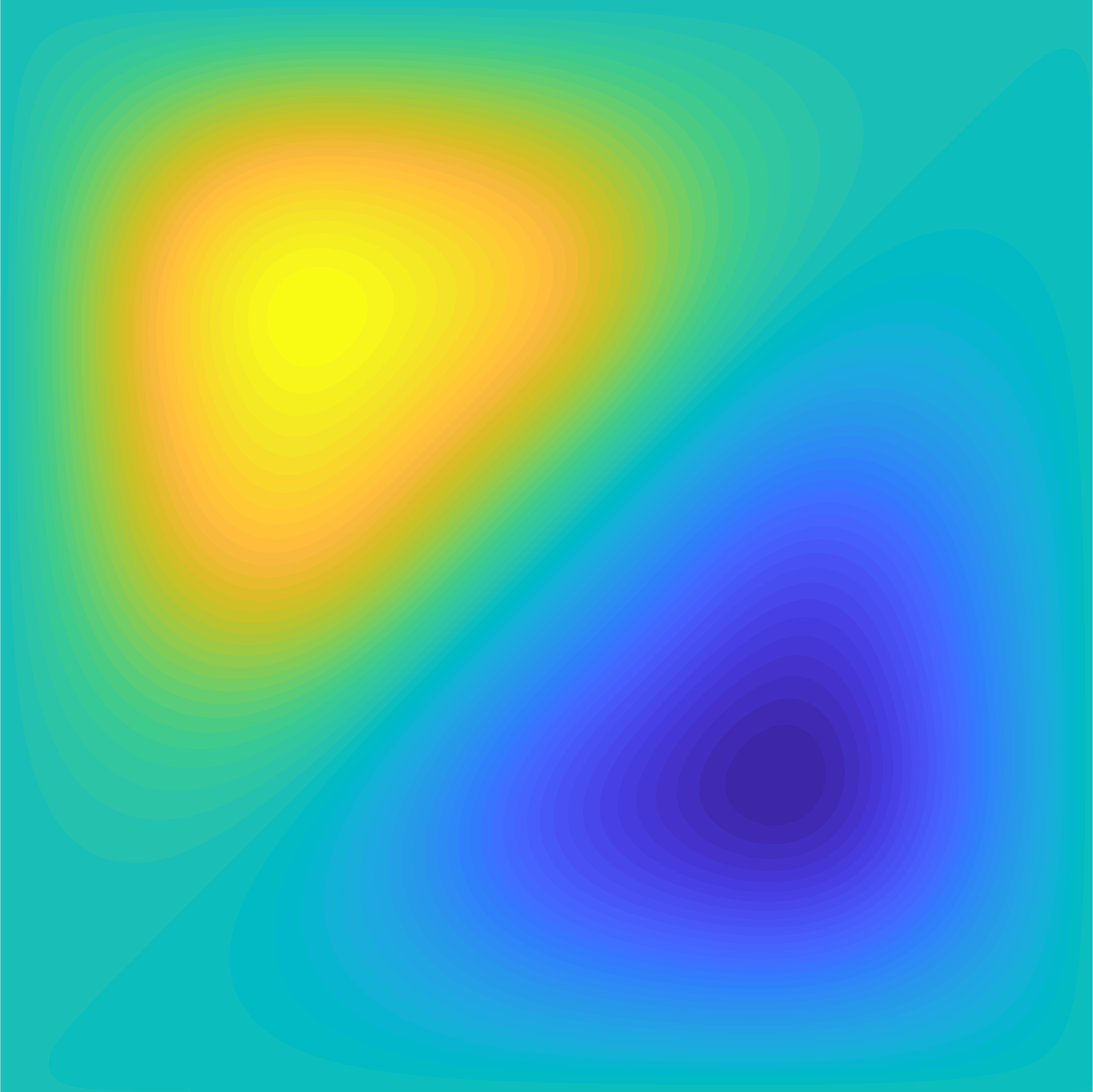};
			
			\nextgroupplot[ ylabel={}, ytick=\empty, xlabel={}, xtick=\empty]
			\addplot graphics [xmin=0, xmax=8.0, ymin=0, ymax=8.0] {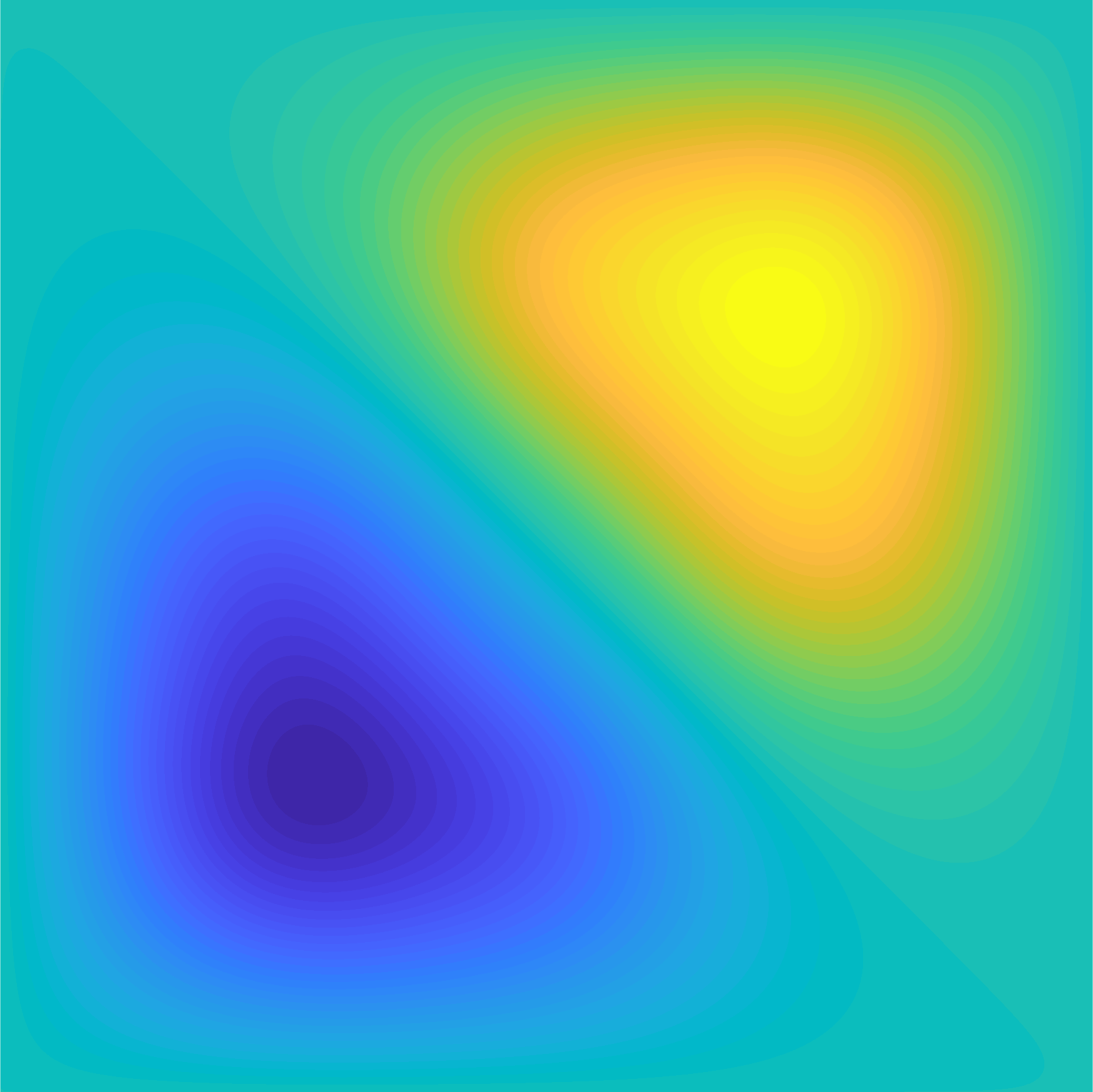};
			
			\nextgroupplot[ ylabel={}, ytick=\empty, xlabel={}, xtick=\empty]
			\addplot graphics [xmin=0, xmax=8.0, ymin=0, ymax=8.0] {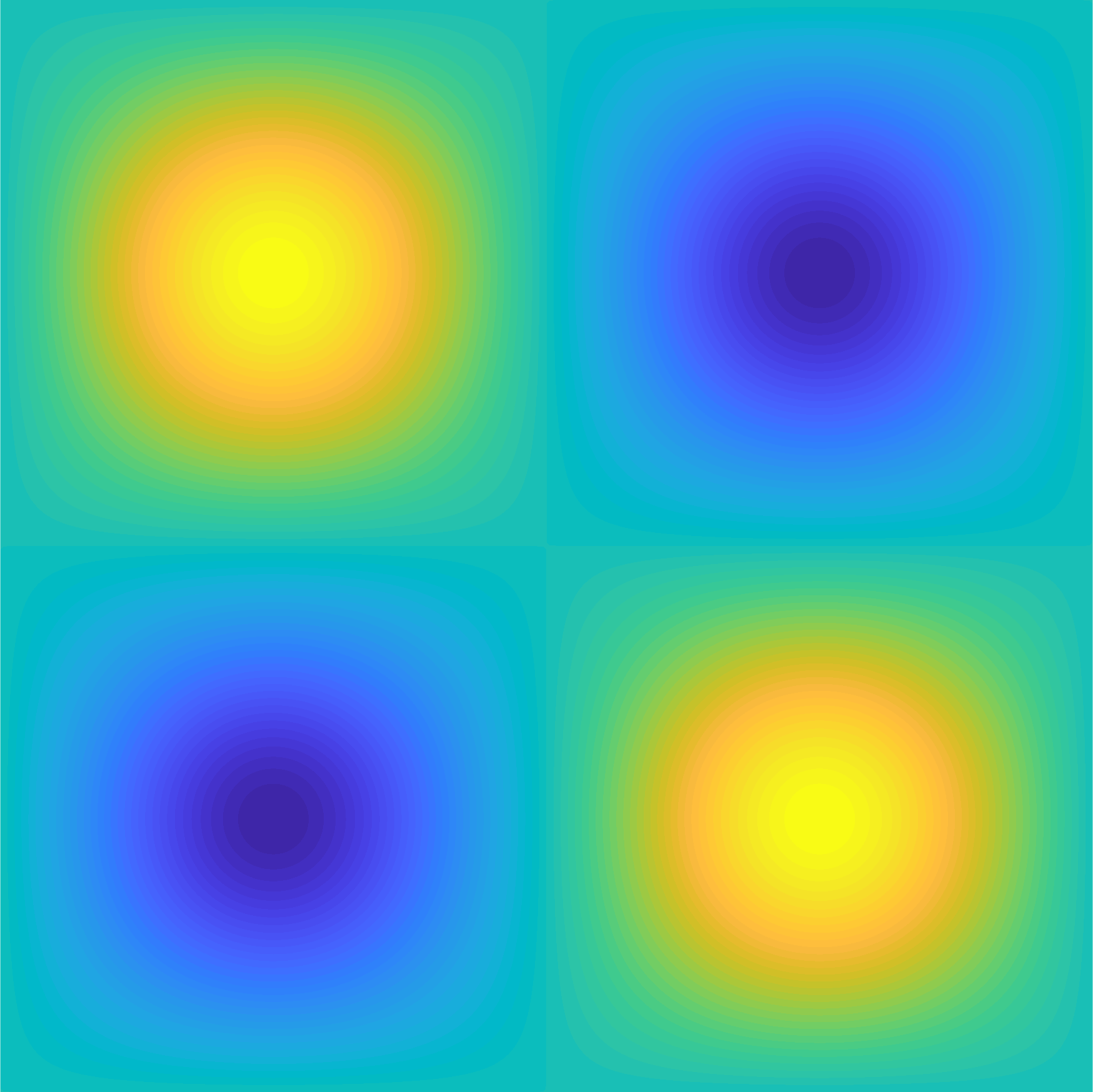};
			
			\nextgroupplot[ ylabel={}, ytick=\empty, xlabel={}, xtick=\empty]
			\addplot graphics [xmin=0, xmax=8.0, ymin=0, ymax=8.0] {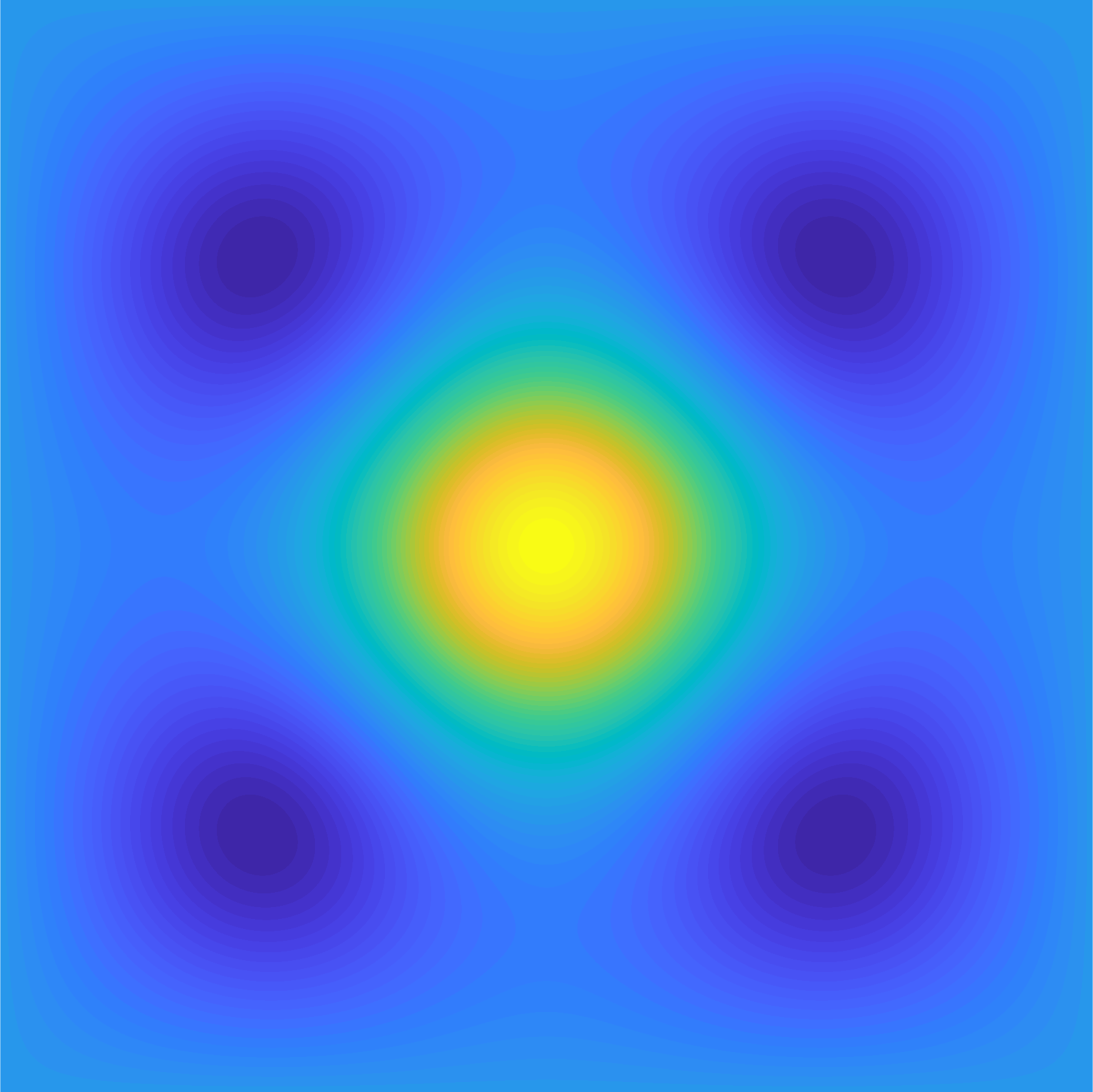};
			
			\nextgroupplot[ ylabel={}, ytick=\empty, xlabel={}, xtick=\empty]
			\addplot graphics [xmin=0, xmax=8.0, ymin=0, ymax=8.0] {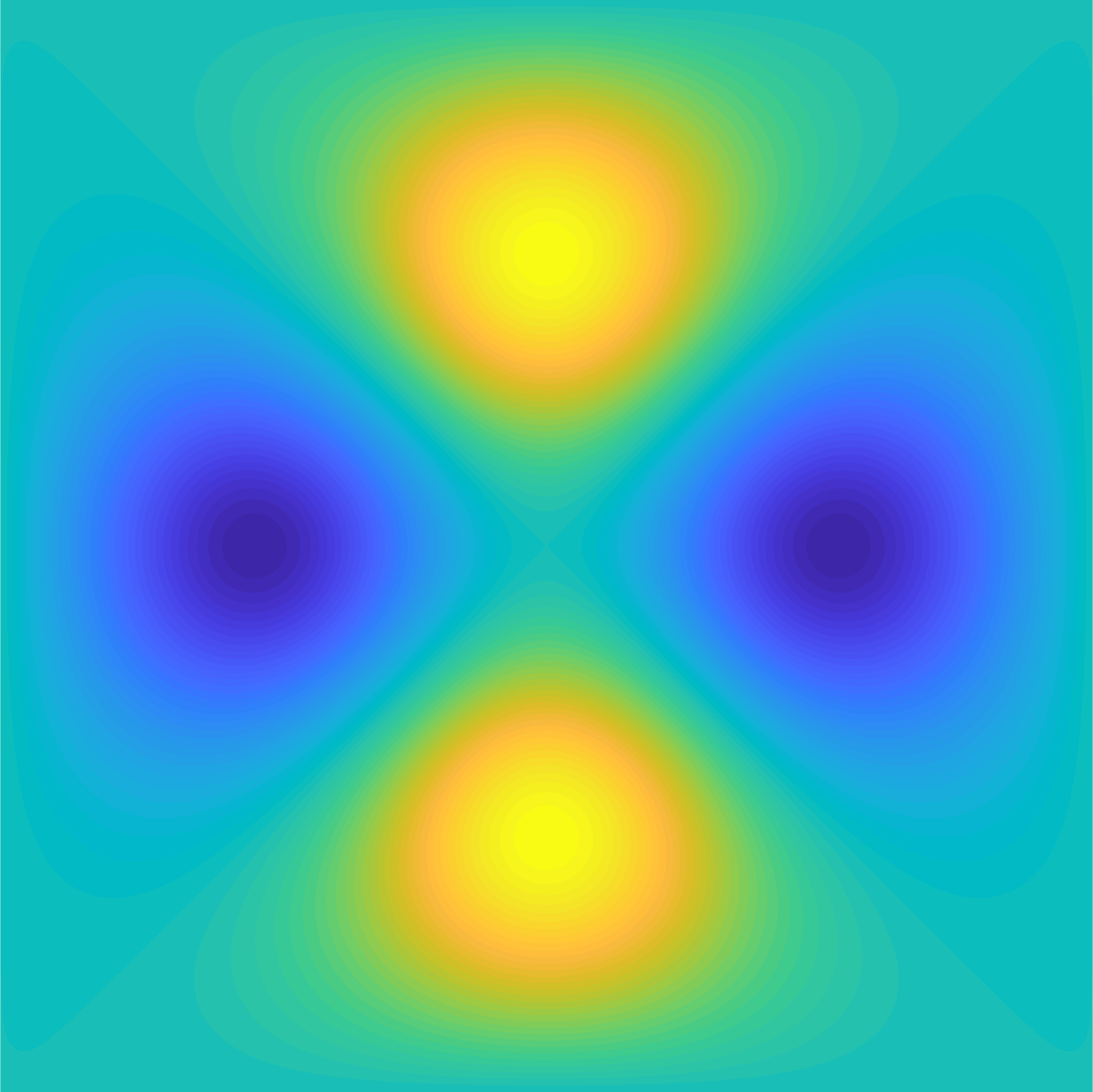};
			
			\nextgroupplot[ ylabel={}, ytick=\empty, xlabel={}, xtick=\empty]
			\addplot graphics [xmin=0, xmax=8.0, ymin=0, ymax=8.0] {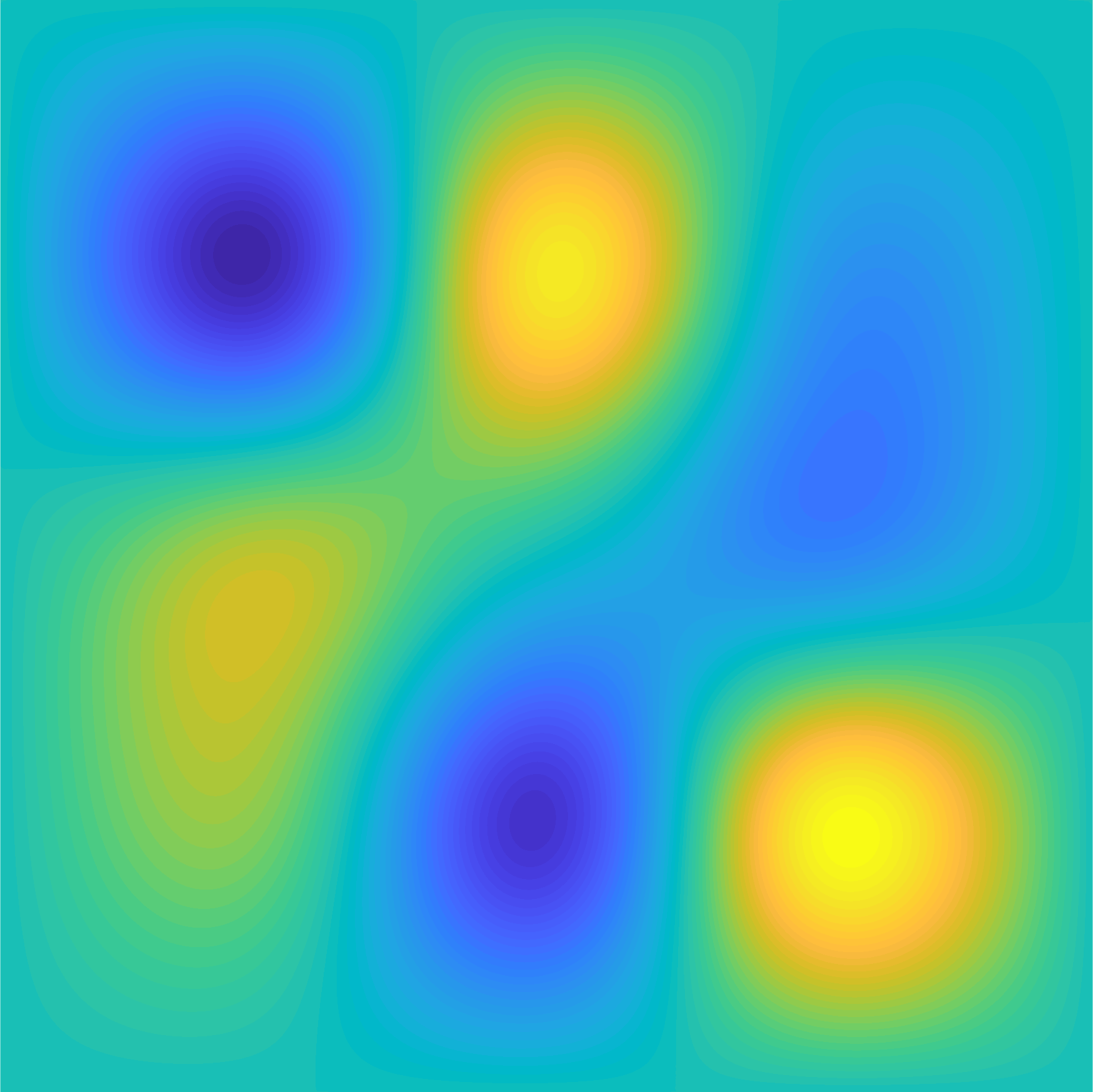};
			
			\nextgroupplot[ ylabel={}, ytick=\empty, xlabel={}, xtick=\empty]
			\addplot graphics [xmin=0, xmax=8.0, ymin=0, ymax=8.0] {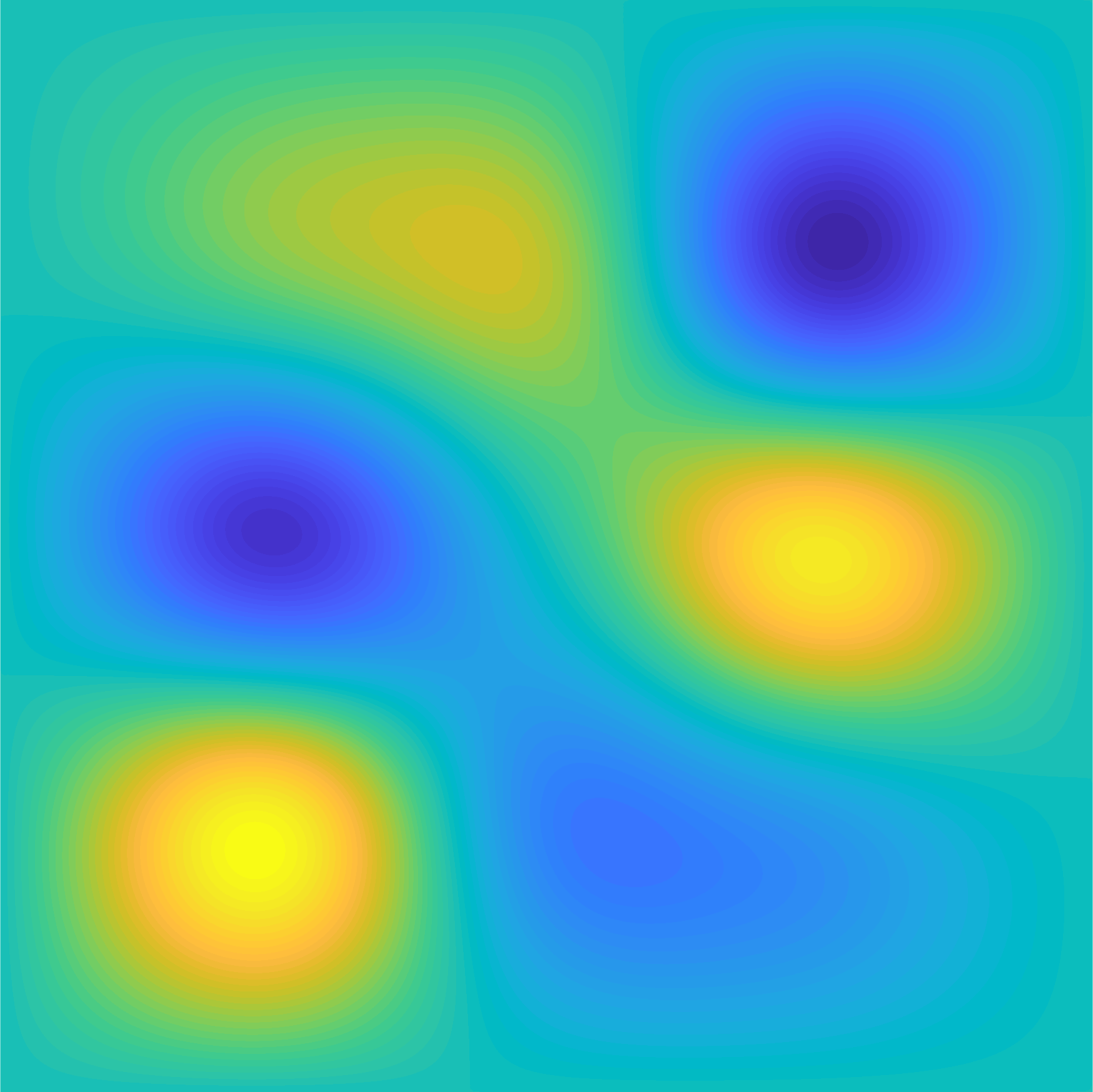};
			
			\nextgroupplot[ ylabel={}, ytick=\empty, xlabel={}, xtick=\empty]
			\addplot graphics [xmin=0, xmax=8.0, ymin=0, ymax=8.0] {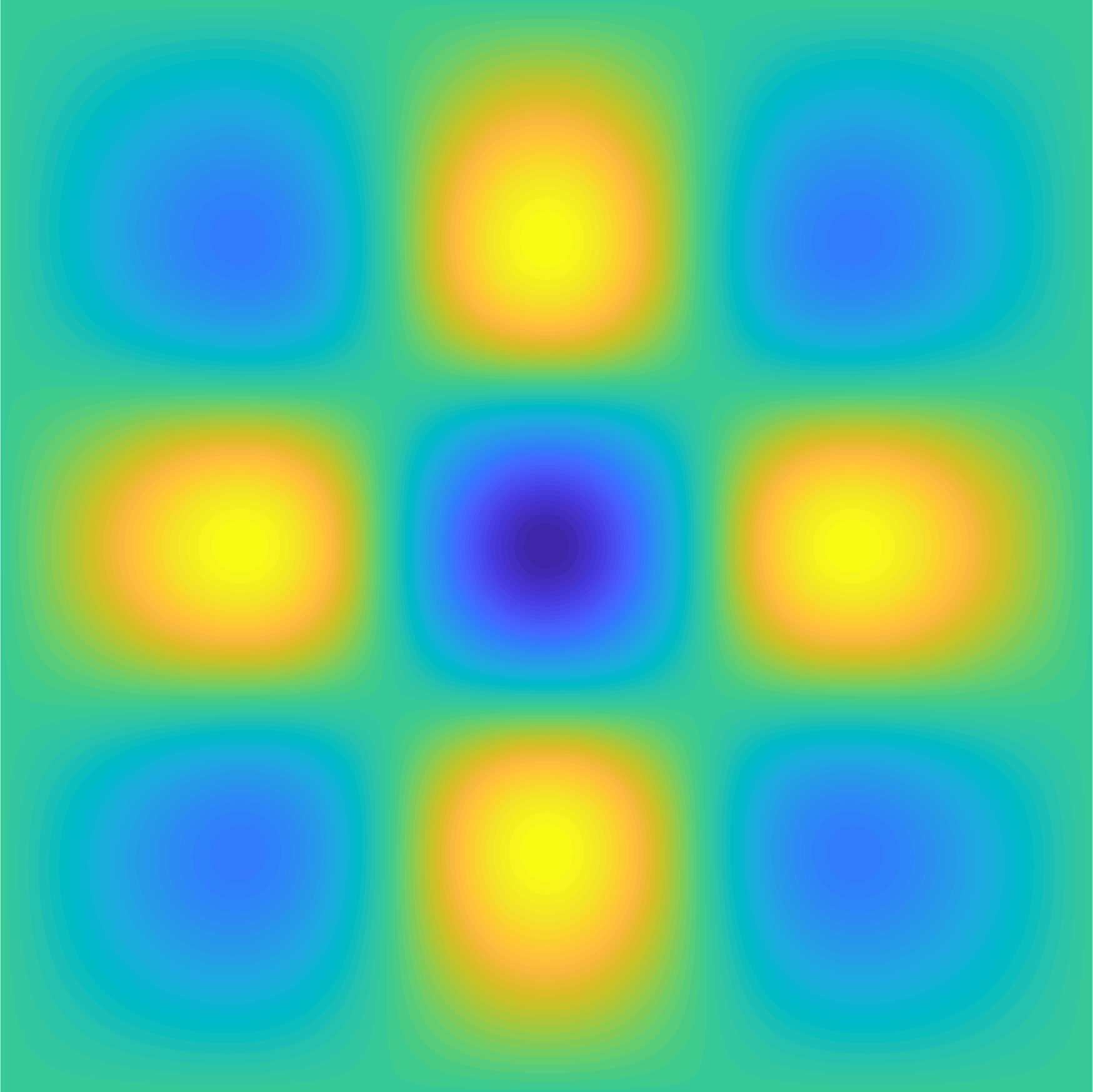};			
		\end{groupplot}
	\end{tikzpicture}
	\caption{All nine POD modes (\textit{left-to-right, top-to-bottom}) corresponding to the snapshots generated from the training set $\Dcal_{\mathrm{train}}$ in Section~\ref{sec:rslt:poi:conv}. Each figure uses a different colorbar, scaled to its range, to highlight features in the corresponding basis function.}
	\label{fig:poi_pod_0}
\end{figure}

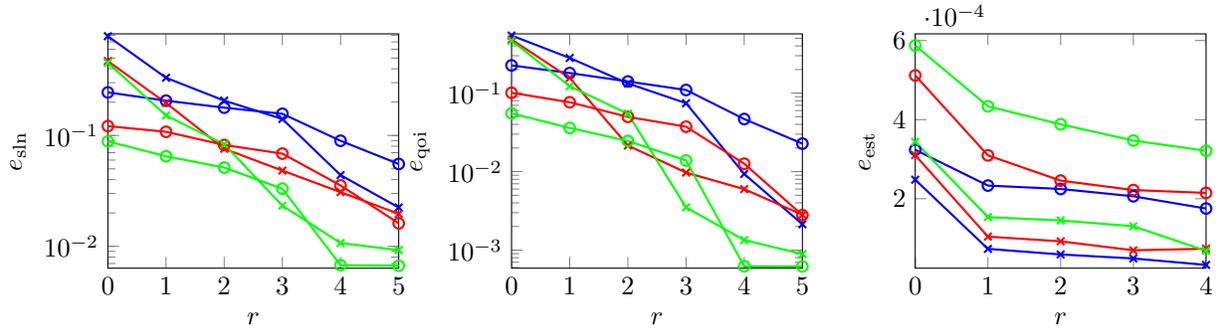
\begin{figure}
	\centering
	\begin{tikzpicture}
\begin{groupplot} [
group style={group size = 3 by 1, horizontal sep = 1.5cm}]
\nextgroupplot[width=0.33\textwidth, xtick={0,1,2,3,4,5}, xlabel=$r$, ymax=0.826569333117877, xmax=5, ylabel=$e_\mathrm{sln}$, xmin=0, ymin=0.006346261418677826,ymode=log]
\addplot [color=blue, mark=o, thick]
coordinates {
( 0.00000000e+00,  2.44875783e-01)
( 1.00000000e+00,  2.06218549e-01)
( 2.00000000e+00,  1.78360730e-01)
( 3.00000000e+00,  1.56998157e-01)
( 4.00000000e+00,  8.95952150e-02)
( 5.00000000e+00,  5.53275331e-02)};\label{poi_nrb3_int}

\addplot [color=blue, mark=x, thick]
coordinates {
( 0.00000000e+00,  7.87208889e-01)
( 1.00000000e+00,  3.31728488e-01)
( 2.00000000e+00,  2.06465065e-01)
( 3.00000000e+00,  1.40913031e-01)
( 4.00000000e+00,  4.38865219e-02)
( 5.00000000e+00,  2.24228992e-02)};\label{poi_nrb3_ext}

\addplot [color=red, mark=o, thick]
coordinates {
( 0.00000000e+00,  1.21532826e-01)
( 1.00000000e+00,  1.07997955e-01)
( 2.00000000e+00,  8.20916030e-02)
( 3.00000000e+00,  6.86059576e-02)
( 4.00000000e+00,  3.54140414e-02)
( 5.00000000e+00,  1.61168799e-02)};\label{poi_nrb6_int}

\addplot [color=red, mark=x, thick]
coordinates {
( 0.00000000e+00,  4.68900686e-01)
( 1.00000000e+00,  1.96579610e-01)
( 2.00000000e+00,  7.59425526e-02)
( 3.00000000e+00,  4.82253519e-02)
( 4.00000000e+00,  3.08103475e-02)
( 5.00000000e+00,  1.96298381e-02)};\label{poi_nrb6_ext}

\addplot [color=green, mark=o, thick]
coordinates {
( 0.00000000e+00,  8.85152381e-02)
( 1.00000000e+00,  6.49963373e-02)
( 2.00000000e+00,  5.13423317e-02)
( 3.00000000e+00,  3.30848996e-02)
( 4.00000000e+00,  6.75624019e-03)
( 5.00000000e+00,  6.68027518e-03)};\label{poi_nrb9_int}

\addplot [color=green, mark=x, thick]
coordinates {
( 0.00000000e+00,  4.52765688e-01)
( 1.00000000e+00,  1.51506795e-01)
( 2.00000000e+00,  8.32556948e-02)
( 3.00000000e+00,  2.32624050e-02)
( 4.00000000e+00,  1.06994928e-02)
( 5.00000000e+00,  9.23896095e-03)};\label{poi_nrb9_ext}

\nextgroupplot[width=0.33\textwidth, xtick={0,1,2,3,4,5}, xlabel=$r$, ymax=0.5693885767506327, xmax=5, ylabel=$e_\mathrm{qoi}$, xmin=0, ymin=0.0005884894905236263,ymode=log]
\addplot [color=blue, mark=o, thick, forget plot]
coordinates {
( 0.00000000e+00,  2.25703207e-01)
( 1.00000000e+00,  1.80376181e-01)
( 2.00000000e+00,  1.40832503e-01)
( 3.00000000e+00,  1.09730914e-01)
( 4.00000000e+00,  4.67606567e-02)
( 5.00000000e+00,  2.27312308e-02)};

\addplot [color=blue, mark=x, thick, forget plot]
coordinates {
( 0.00000000e+00,  5.42274835e-01)
( 1.00000000e+00,  2.80904902e-01)
( 2.00000000e+00,  1.32939755e-01)
( 3.00000000e+00,  7.44553435e-02)
( 4.00000000e+00,  9.34101627e-03)
( 5.00000000e+00,  2.13504469e-03)};

\addplot [color=red, mark=o, thick, forget plot]
coordinates {
( 0.00000000e+00,  1.01011586e-01)
( 1.00000000e+00,  7.65999076e-02)
( 2.00000000e+00,  4.98143361e-02)
( 3.00000000e+00,  3.73862374e-02)
( 4.00000000e+00,  1.26073250e-02)
( 5.00000000e+00,  2.79336002e-03)};

\addplot [color=red, mark=x, thick, forget plot]
coordinates {
( 0.00000000e+00,  4.81511874e-01)
( 1.00000000e+00,  1.56632651e-01)
( 2.00000000e+00,  2.15104917e-02)
( 3.00000000e+00,  9.72600725e-03)
( 4.00000000e+00,  5.97932342e-03)
( 5.00000000e+00,  2.81836014e-03)};

\addplot [color=green, mark=o, thick, forget plot]
coordinates {
( 0.00000000e+00,  5.52630871e-02)
( 1.00000000e+00,  3.59689103e-02)
( 2.00000000e+00,  2.45928835e-02)
( 3.00000000e+00,  1.38390682e-02)
( 4.00000000e+00,  6.22048442e-04)
( 5.00000000e+00,  6.19462622e-04)};

\addplot [color=green, mark=x, thick, forget plot]
coordinates {
( 0.00000000e+00,  4.65177417e-01)
( 1.00000000e+00,  1.22628651e-01)
( 2.00000000e+00,  5.45490393e-02)
( 3.00000000e+00,  3.50026688e-03)
( 4.00000000e+00,  1.34236559e-03)
( 5.00000000e+00,  8.90451932e-04)};

\nextgroupplot[width=0.33\textwidth, xtick={0,1,2,3,4}, xlabel=$r$, ymax=0.0006169588585348394, xmax=4, ylabel=$e_\mathrm{est}$, xmin=0, ymin=2.434191186343053e-05]
\addplot [color=blue, mark=o, thick, forget plot]
coordinates {
( 0.00000000e+00,  3.23652654e-04)
( 1.00000000e+00,  2.33261361e-04)
( 2.00000000e+00,  2.24684362e-04)
( 3.00000000e+00,  2.06174535e-04)
( 4.00000000e+00,  1.75192294e-04)};

\addplot [color=blue, mark=x, thick, forget plot]
coordinates {
( 0.00000000e+00,  2.48112998e-04)
( 1.00000000e+00,  7.31694081e-05)
( 2.00000000e+00,  5.90399078e-05)
( 3.00000000e+00,  4.88744141e-05)
( 4.00000000e+00,  3.26445072e-05)};

\addplot [color=red, mark=o, thick, forget plot]
coordinates {
( 0.00000000e+00,  5.11862168e-04)
( 1.00000000e+00,  3.09615495e-04)
( 2.00000000e+00,  2.45729209e-04)
( 3.00000000e+00,  2.21970049e-04)
( 4.00000000e+00,  2.14680117e-04)};

\addplot [color=red, mark=x, thick, forget plot]
coordinates {
( 0.00000000e+00,  3.10105676e-04)
( 1.00000000e+00,  1.04141315e-04)
( 2.00000000e+00,  9.22607284e-05)
( 3.00000000e+00,  6.95281114e-05)
( 4.00000000e+00,  7.35785213e-05)};

\addplot [color=green, mark=o, thick, forget plot]
coordinates {
( 0.00000000e+00,  5.87579865e-04)
( 1.00000000e+00,  4.33669103e-04)
( 2.00000000e+00,  3.88186894e-04)
( 3.00000000e+00,  3.47050040e-04)
( 4.00000000e+00,  3.21310591e-04)};

\addplot [color=green, mark=x, thick, forget plot]
coordinates {
( 0.00000000e+00,  3.43857346e-04)
( 1.00000000e+00,  1.53490000e-04)
( 2.00000000e+00,  1.45047059e-04)
( 3.00000000e+00,  1.30368712e-04)
( 4.00000000e+00,  6.82708145e-05)};

\end{groupplot}\end{tikzpicture}
	\caption{The mean error of the CG-GL method applied to the Poisson problem over the test set $\Dcal_\mathrm{interp}$ 
	with basis size $k=3$ (\ref{poi_nrb3_int}), $k=6$ (\ref{poi_nrb6_int}), $k=9$ (\ref{poi_nrb9_int}), and the test set
	$\Dcal_\mathrm{extrap}$ with basis size  $k=3$ (\ref{poi_nrb3_ext}) ,  $k=6$ (\ref{poi_nrb6_ext}),  $k=9$
	(\ref{poi_nrb9_ext}). Because the CG-GL method is initialized with $\Omega_l = \emptyset$ and $N_l=0$,
	$r= 0$ corresponds to a traditional ROM.}
	\label{fig:poi_rf_basis}
\end{figure}

To begin, we construct a traditional ROM by computing the reference solution for each $\mu\in\Dcal_\mathrm{train}$ and compressing them into a globally supported basis of dimension $k \in \{3, 6, 9\}$ using POD. The first POD mode captures the centered peak of the solution, and subsequent modes provide adjustments to its position and shape (Figure~\ref{fig:poi_pod_0}). We use these global modes to initialize the CG-GL method ($\Omega_l = \emptyset$, $N_l = 0$) with 16 patches of 36 quadratic quadrilateral elements. We use both the ROM and adaptive CG-GL method to predict the PDE solution and QoI at all testing points and directly compare the errors (Figure~\ref{fig:poi_rf_basis}).  For predictions at interpolated parameters, increasing the global basis reduces the solution and QoI error using the ROM alone. As expected, the ROM predictions have higher errors at extrapolated parameters than at interpolated parameters. The improvement observed by increasing the number of global functions is marginal for the ROM alone: as $k$ increases from $6$ to $9$, the solution error decreases from $46.8\%$ to $45.2\%$. On the other hand, for a fixed global basis dimension $k$, the CG-GL method enriches the POD basis with piecewise polynomial basis function, which steadily reduces the solution and QoI error for each parameter in the test set. The error decreases most rapidly with adaptation for the $k=9$ case suggesting the global basis functions are providing benefit over purely local adaption. In most cases, the error at extrapolated points drop close to or below the errors at interpolated points suggesting the adaptive CG-GL method is not sensitive to location of the testing points relative to the training. Finally, Figure~\ref{fig:poi_rf_basis} shows the error estimates produced are highly accurate for this problem at all stages of refinement.

\subsubsection{Performance with limited training}
\label{sec:rslt:poi:limited}
\begin{figure}
	\centering
	\begin{tikzpicture}
		\begin{groupplot}[
			group style={
				group size=3 by 1,
				horizontal sep=0.5cm
			},
			width=0.35\textwidth,
			axis equal image,
			xlabel={$x_1$},
			ylabel={$x_2$},
			xtick = {0.0, 4.0, 8.0},
			ytick = {0.0, 4.0, 8.0},
			xmin=0, xmax=8,
			ymin=0, ymax=8
			]
			\nextgroupplot[ ylabel={}, ytick=\empty, xlabel={}, xtick=\empty]
			\addplot graphics [xmin=0, xmax=8.0, ymin=0, ymax=8.0] {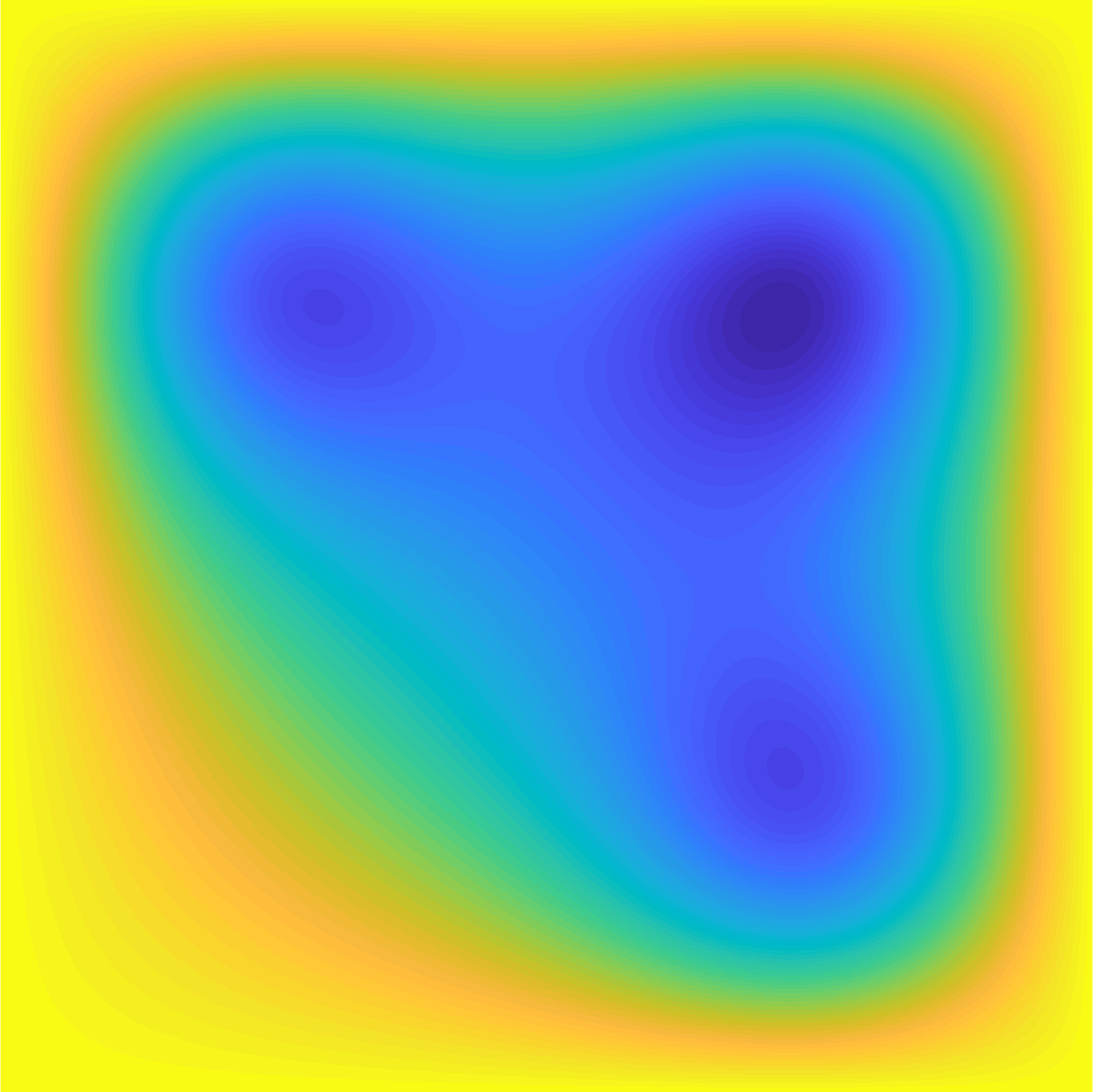};
			
			\nextgroupplot[ ylabel={}, ytick=\empty, xlabel={}, xtick=\empty]
			\addplot graphics [xmin=0, xmax=8.0, ymin=0, ymax=8.0] {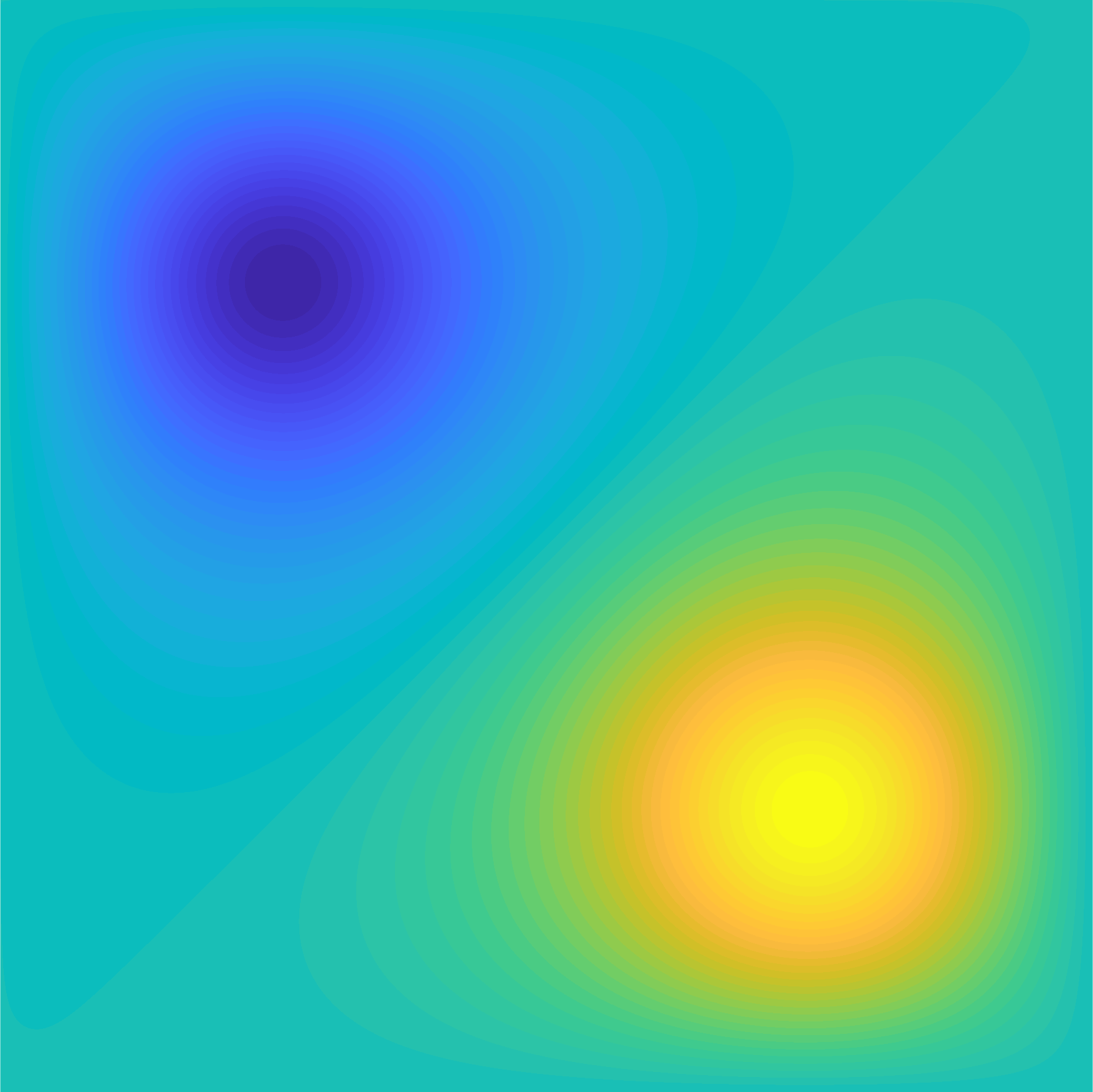};
			
			\nextgroupplot[ ylabel={}, ytick=\empty, xlabel={}, xtick=\empty]
			\addplot graphics [xmin=0, xmax=8.0, ymin=0, ymax=8.0] {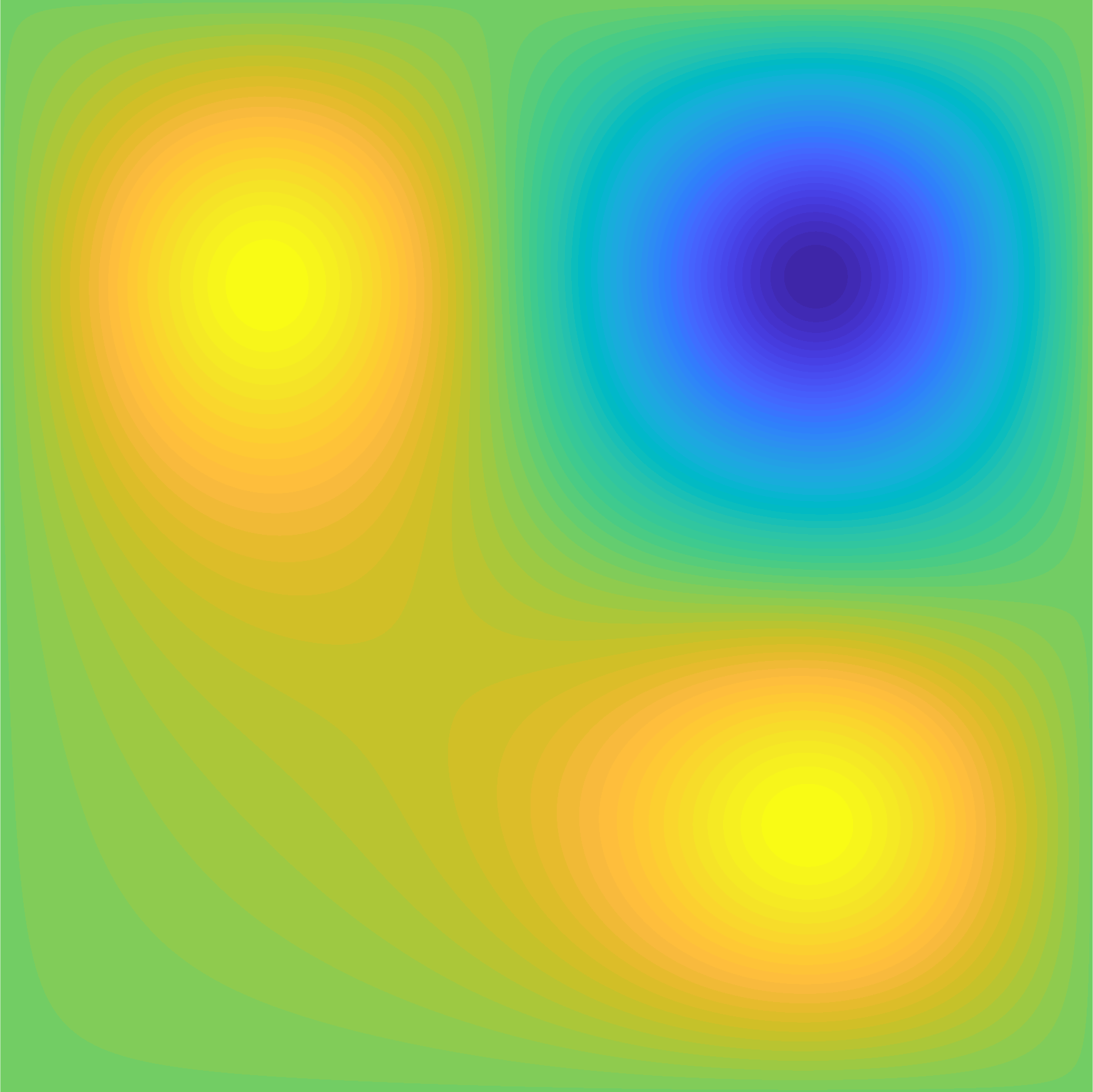};

		\end{groupplot}
	\end{tikzpicture}
	\caption{All three POD modes (\textit{left-to-right}) corresponding to the snapshots generated from the training set $\Dcal_{\mathrm{train}}$ in Section~\ref{sec:rslt:poi:limited}. Each figure uses a different colorbar, scaled to its range to highlight features in the corresponding basis function.}
	\label{fig:poi_pod_1}
\end{figure}

\begin{figure}
	\centering
	\begin{tikzpicture}
		\begin{groupplot}[
			group style={
				group size=4 by 3,
				horizontal sep=0.5cm,
				vertical sep=0.5cm
			},
			width=0.35\textwidth,
			axis equal image,
			xlabel={$x_1$},
			ylabel={$x_2$},
			xtick = {0.0, 4.0, 8.0},
			ytick = {0.0, 4.0, 8.0},
			xmin=0, xmax=8,
			ymin=0, ymax=8
			]
			\nextgroupplot[ ylabel={}, ytick=\empty, xlabel={}, xtick=\empty]
			\addplot graphics [xmin=0, xmax=8.0, ymin=0, ymax=8.0] {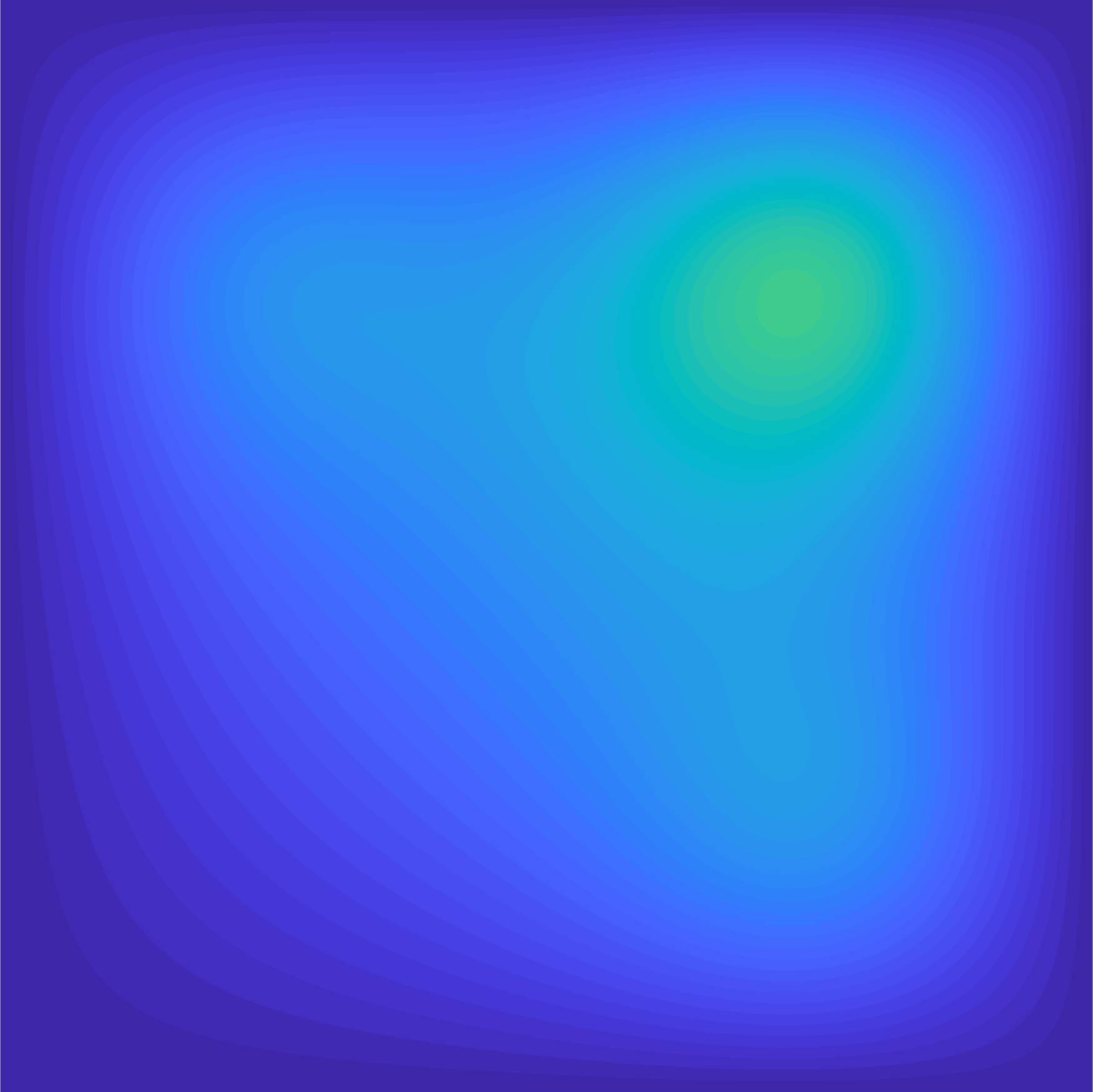};
			
			\nextgroupplot[ ylabel={}, ytick=\empty, xlabel={}, xtick=\empty]
			\addplot graphics [xmin=0, xmax=8.0, ymin=0, ymax=8.0] {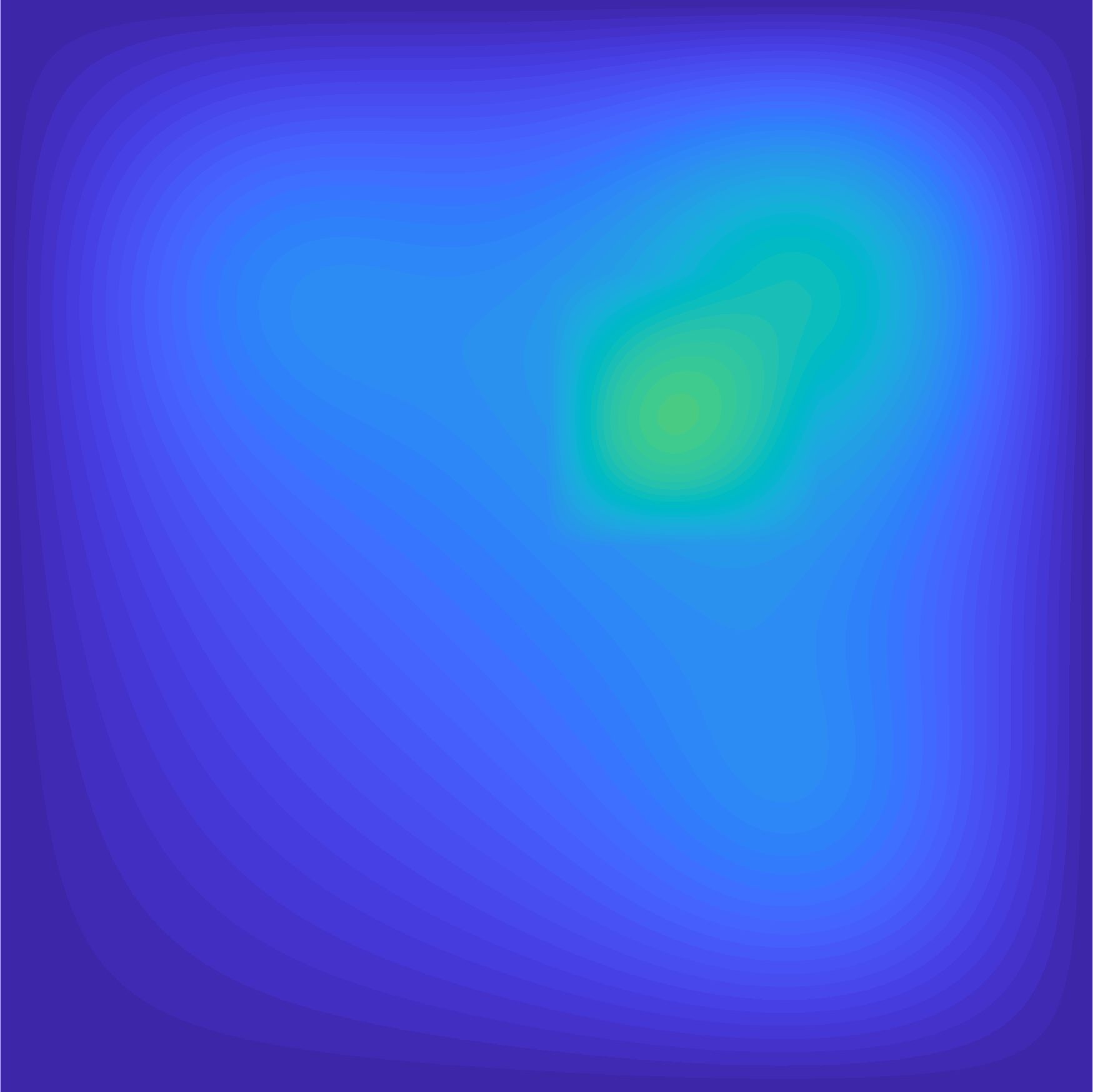};
			
			\nextgroupplot[ ylabel={}, ytick=\empty, xlabel={}, xtick=\empty]
			\addplot graphics [xmin=0, xmax=8.0, ymin=0, ymax=8.0] {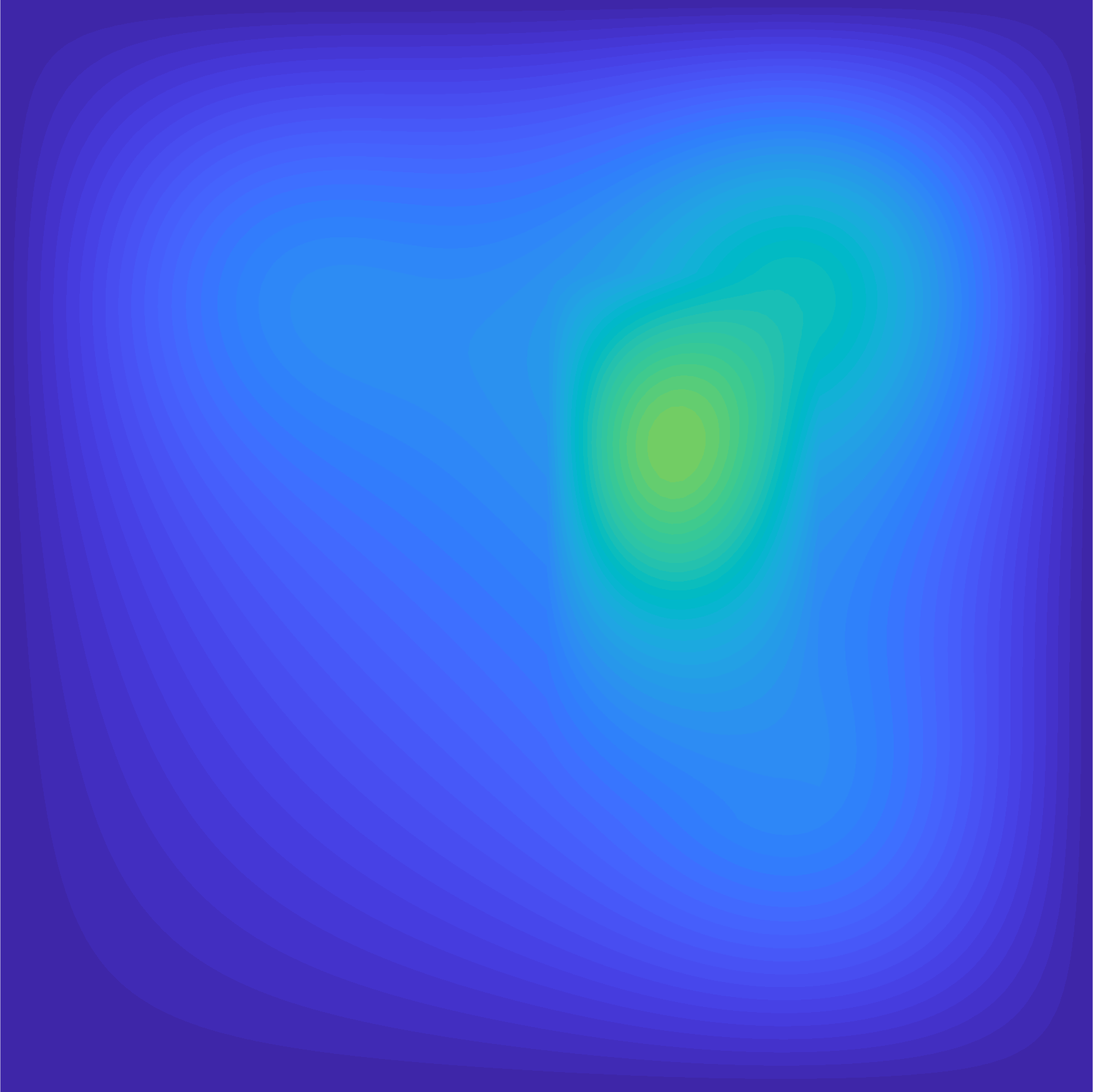};
			
			\nextgroupplot[ ylabel={}, ytick=\empty, xlabel={}, xtick=\empty]
			\addplot graphics [xmin=0, xmax=8.0, ymin=0, ymax=8.0] {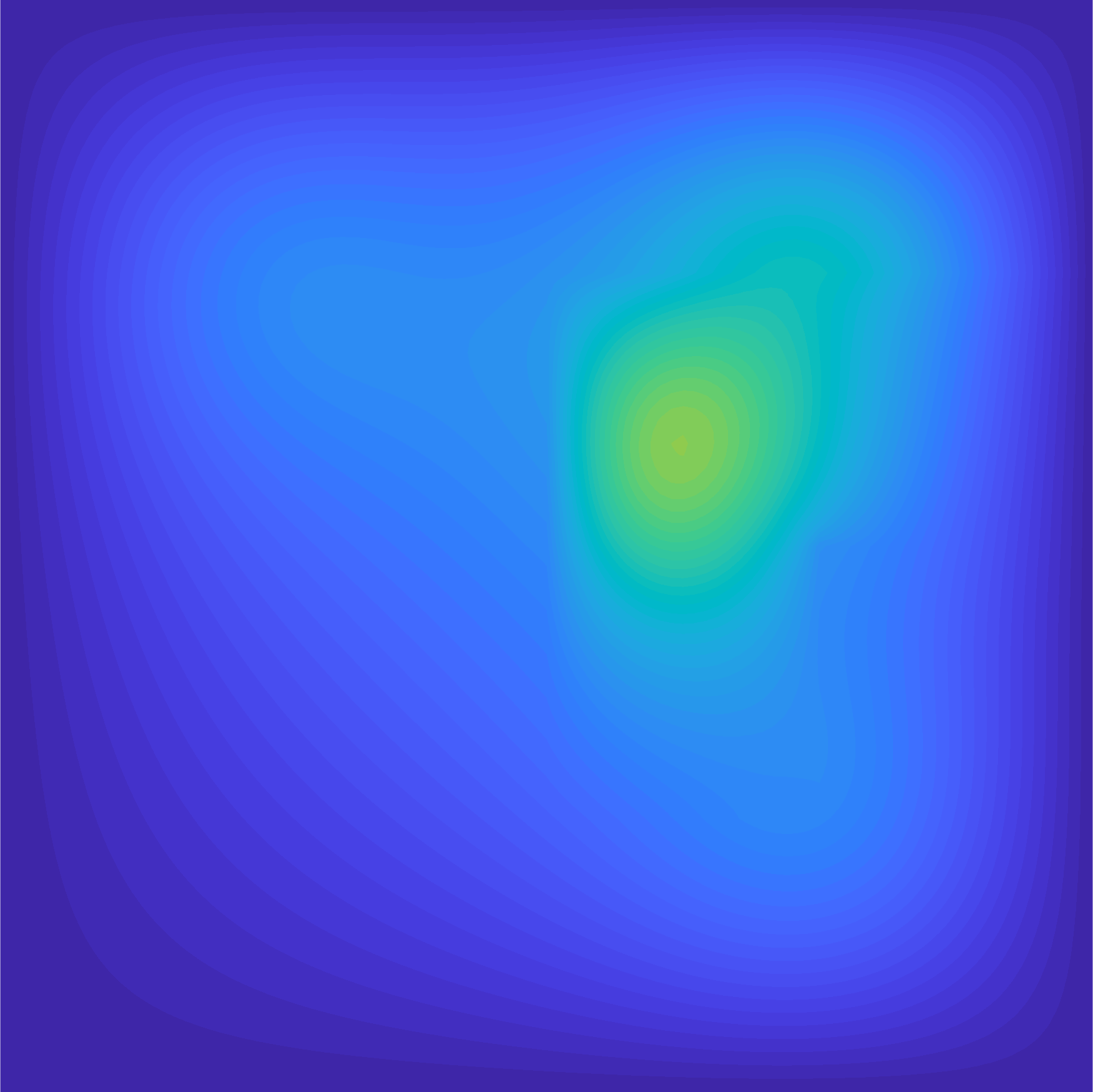};
			
			\nextgroupplot[ ylabel={}, ytick=\empty, xlabel={}, xtick=\empty]
			\addplot graphics [xmin=0, xmax=8.0, ymin=0, ymax=8.0] {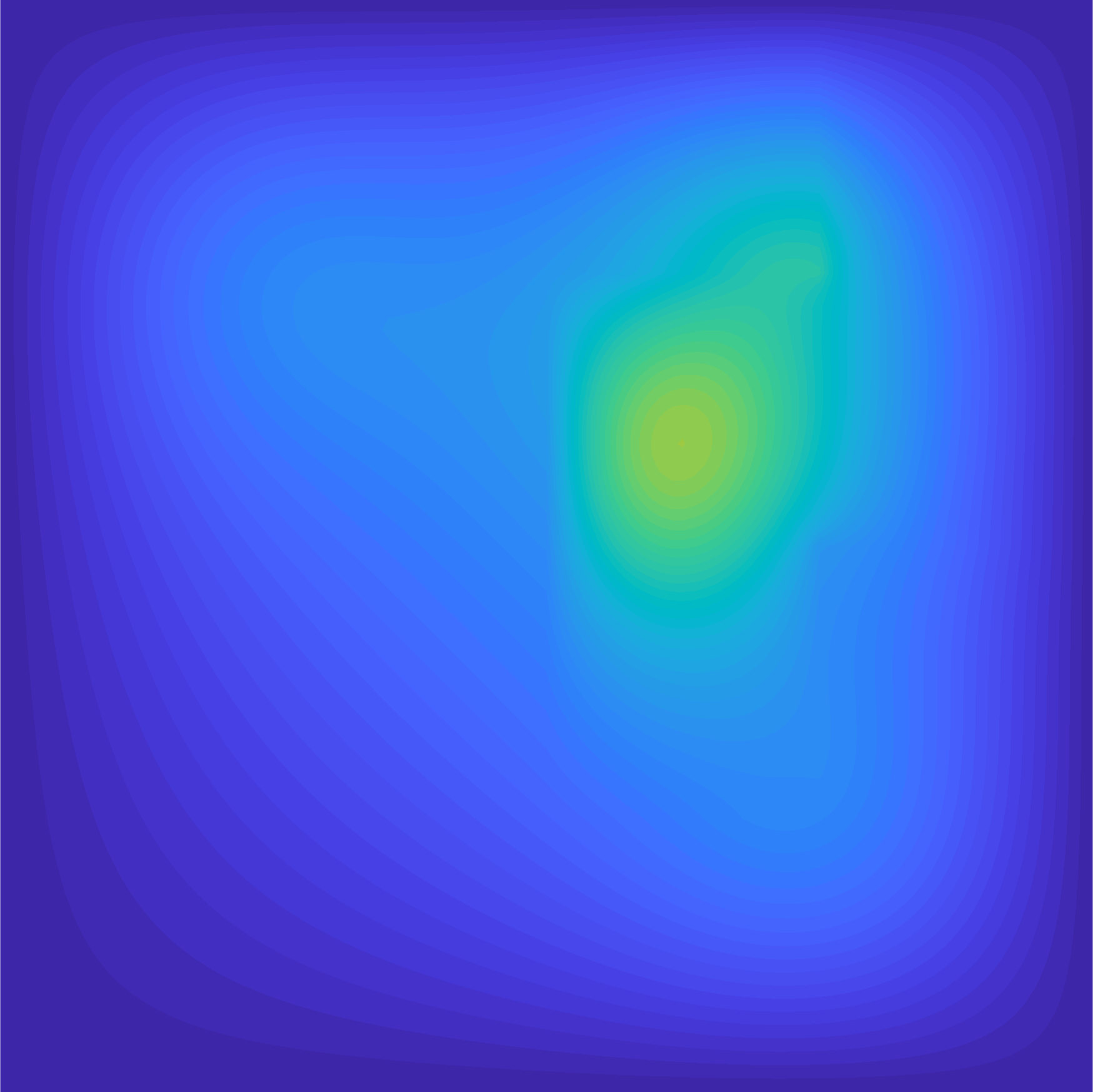};
			
			\nextgroupplot[ ylabel={}, ytick=\empty, xlabel={}, xtick=\empty]
			\addplot graphics [xmin=0, xmax=8.0, ymin=0, ymax=8.0] {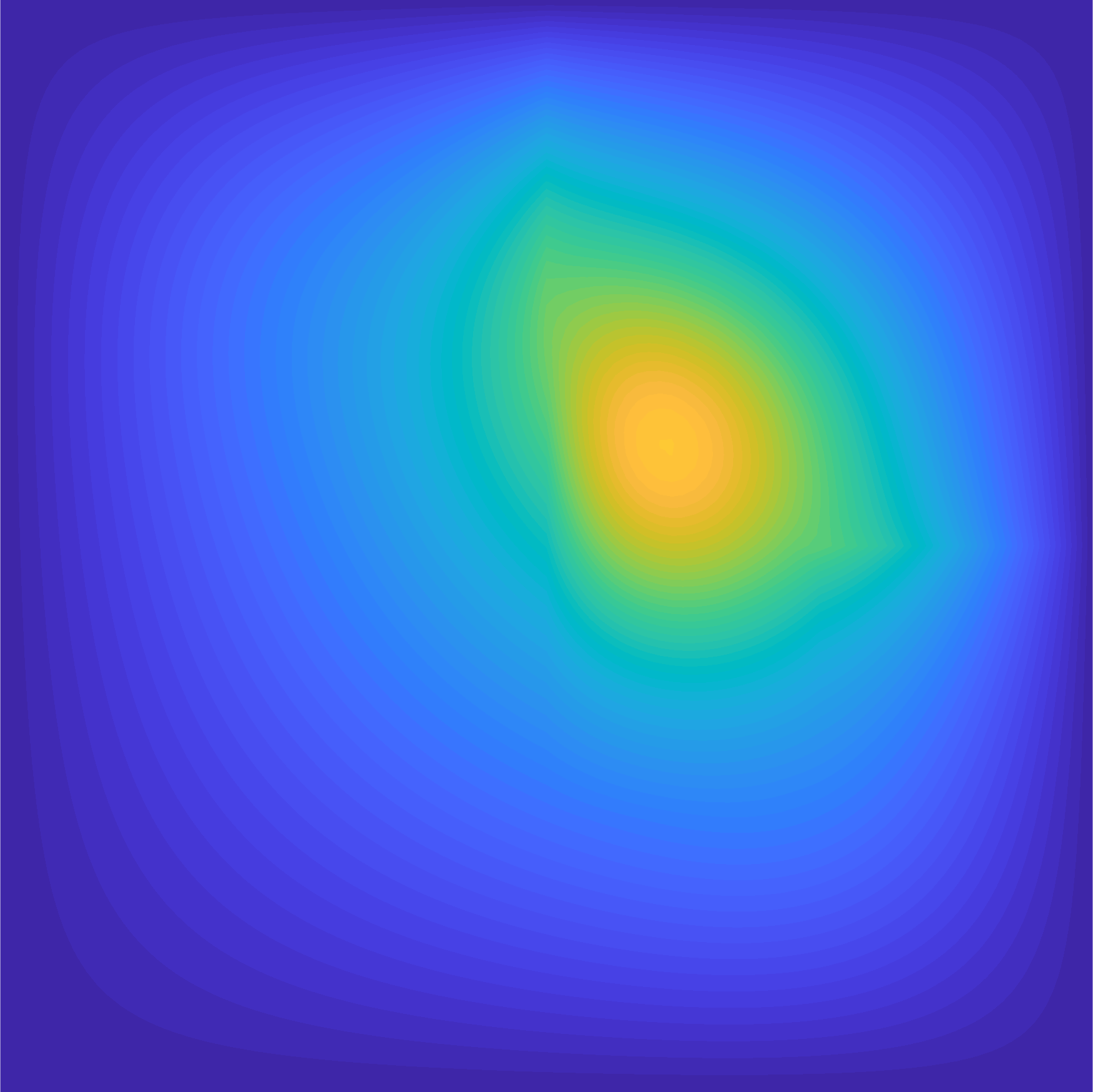};
			
			\nextgroupplot[ ylabel={}, ytick=\empty, xlabel={}, xtick=\empty]
			\addplot graphics [xmin=0, xmax=8.0, ymin=0, ymax=8.0] {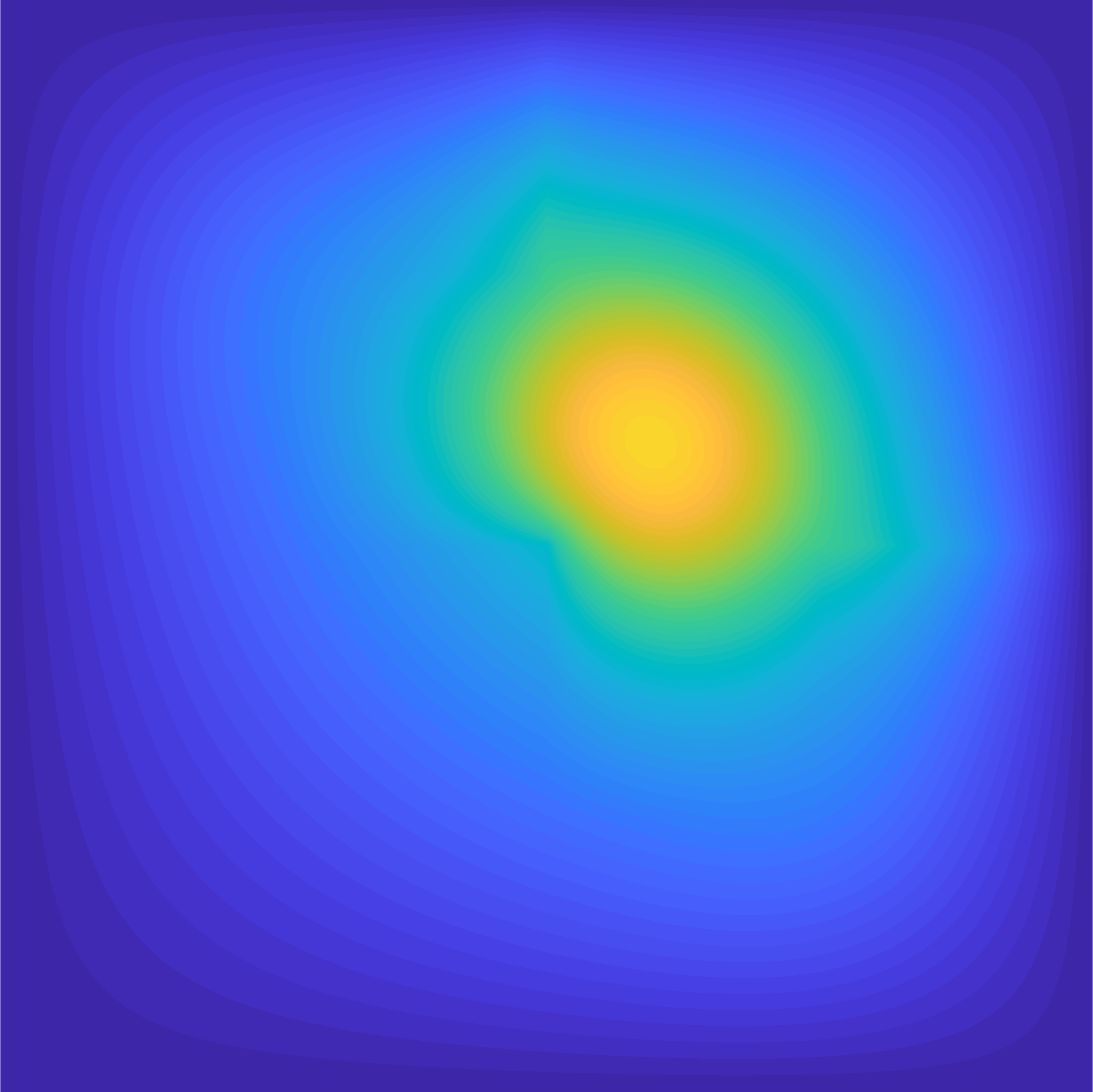};
			
			\nextgroupplot[ ylabel={}, ytick=\empty, xlabel={}, xtick=\empty]
			\addplot graphics [xmin=0, xmax=8.0, ymin=0, ymax=8.0] {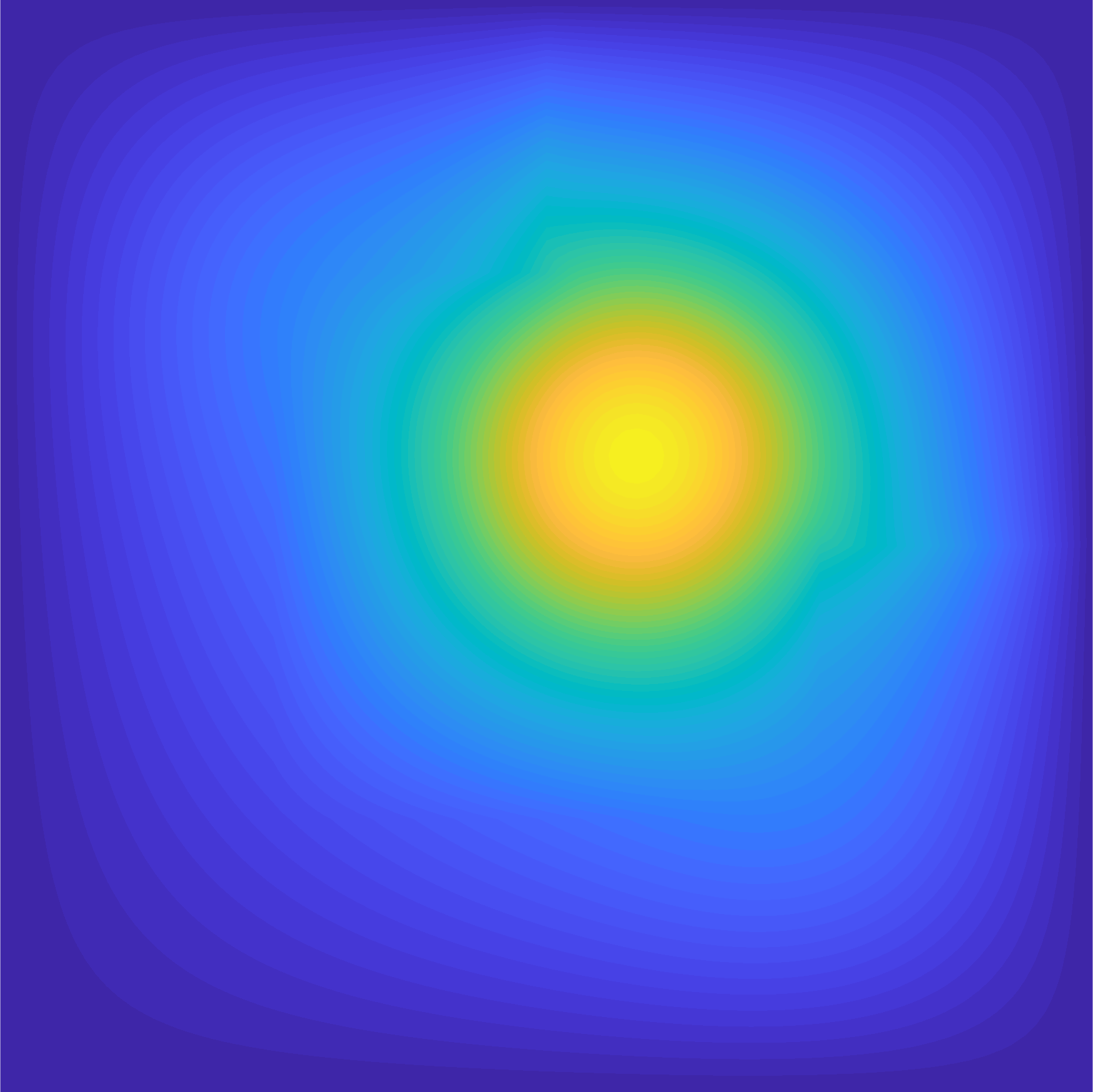};
			
			\nextgroupplot[ ylabel={}, ytick=\empty, xlabel={}, xtick=\empty]
			\addplot graphics [xmin=0, xmax=8.0, ymin=0, ymax=8.0] {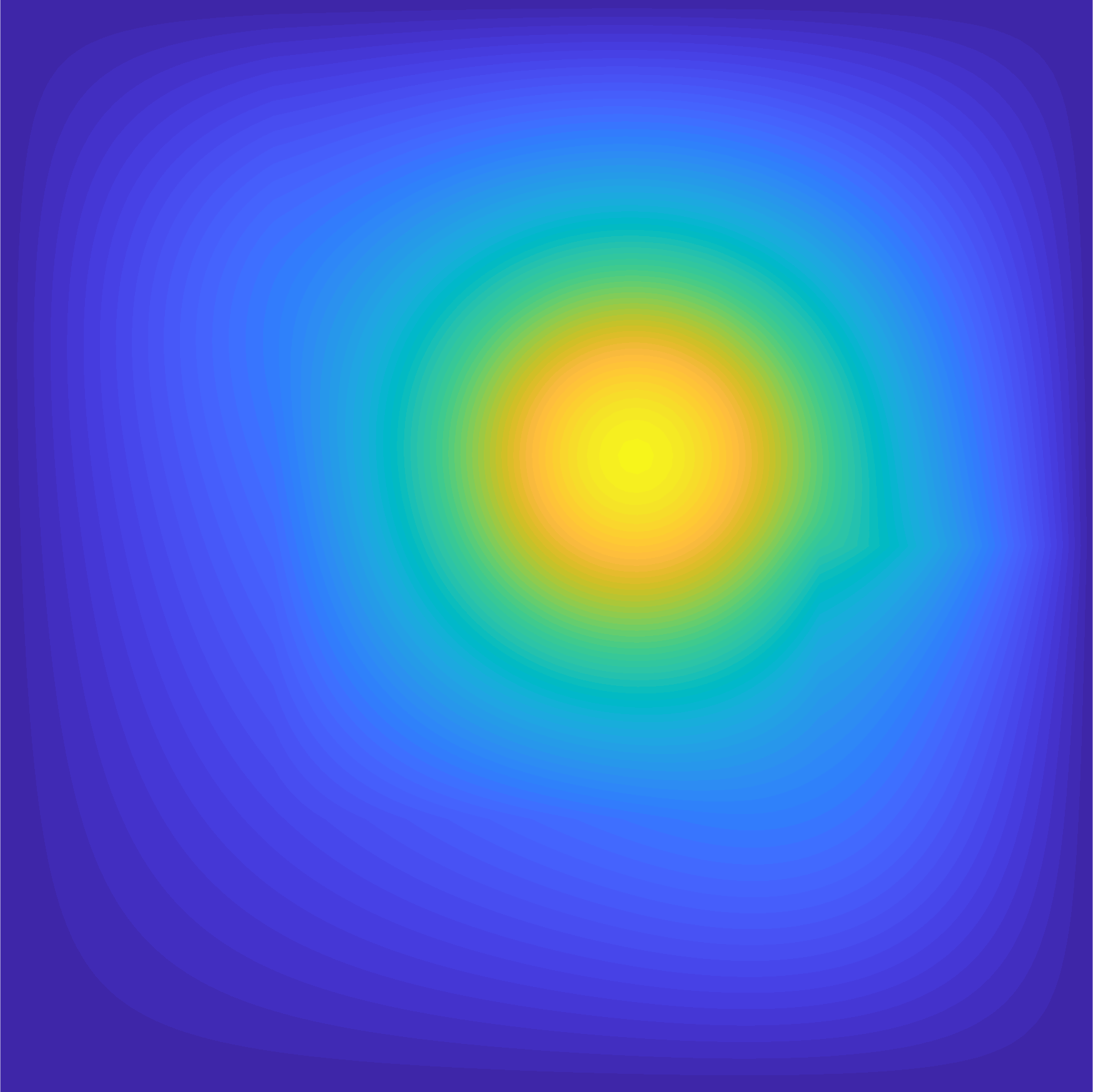};	
			
			\nextgroupplot[ ylabel={}, ytick=\empty, xlabel={}, xtick=\empty]
			\addplot graphics [xmin=0, xmax=8.0, ymin=0, ymax=8.0] {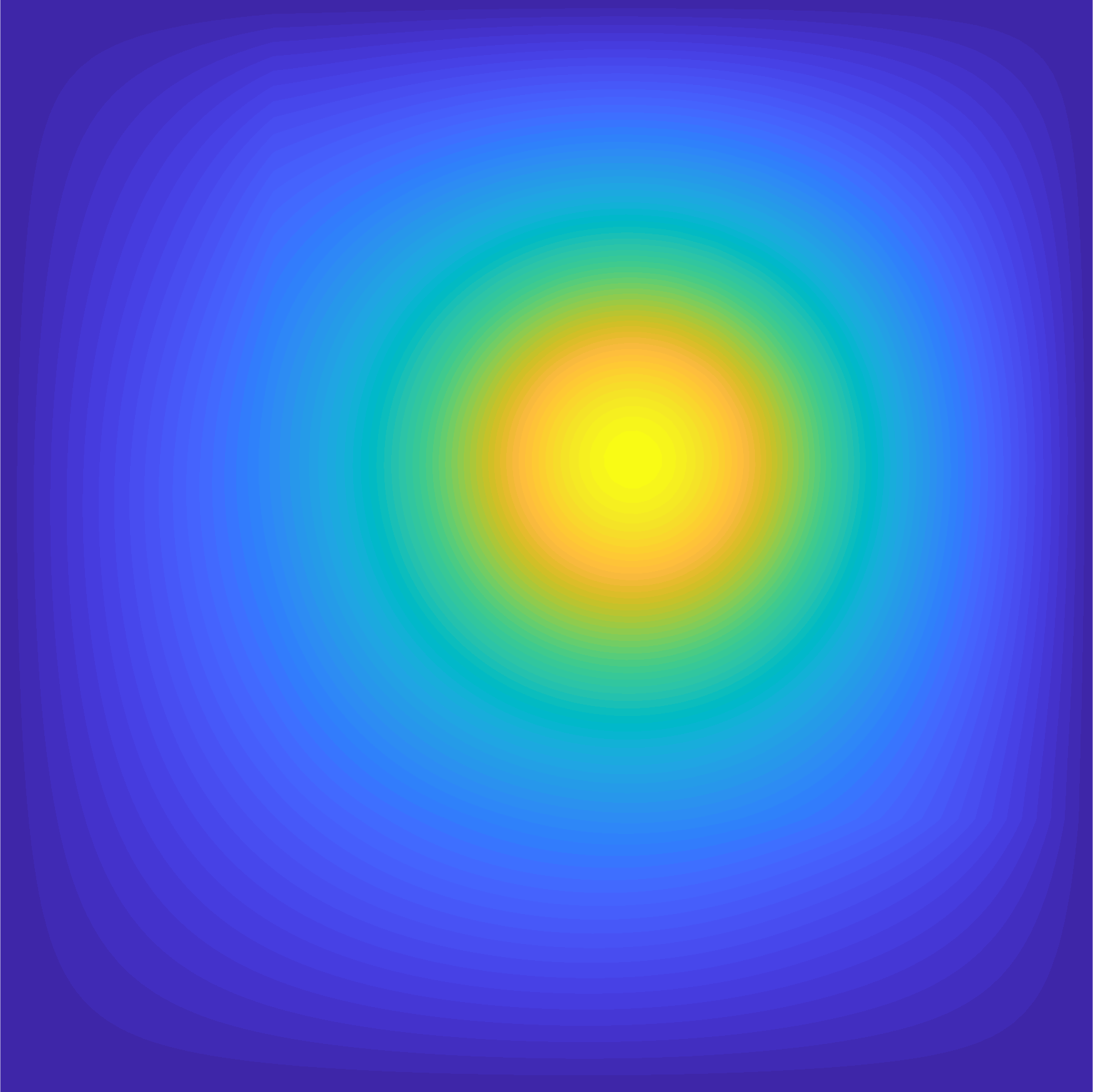};
			
			\nextgroupplot[ ylabel={}, ytick=\empty, xlabel={}, xtick=\empty]
			\addplot graphics [xmin=0, xmax=8.0, ymin=0, ymax=8.0] {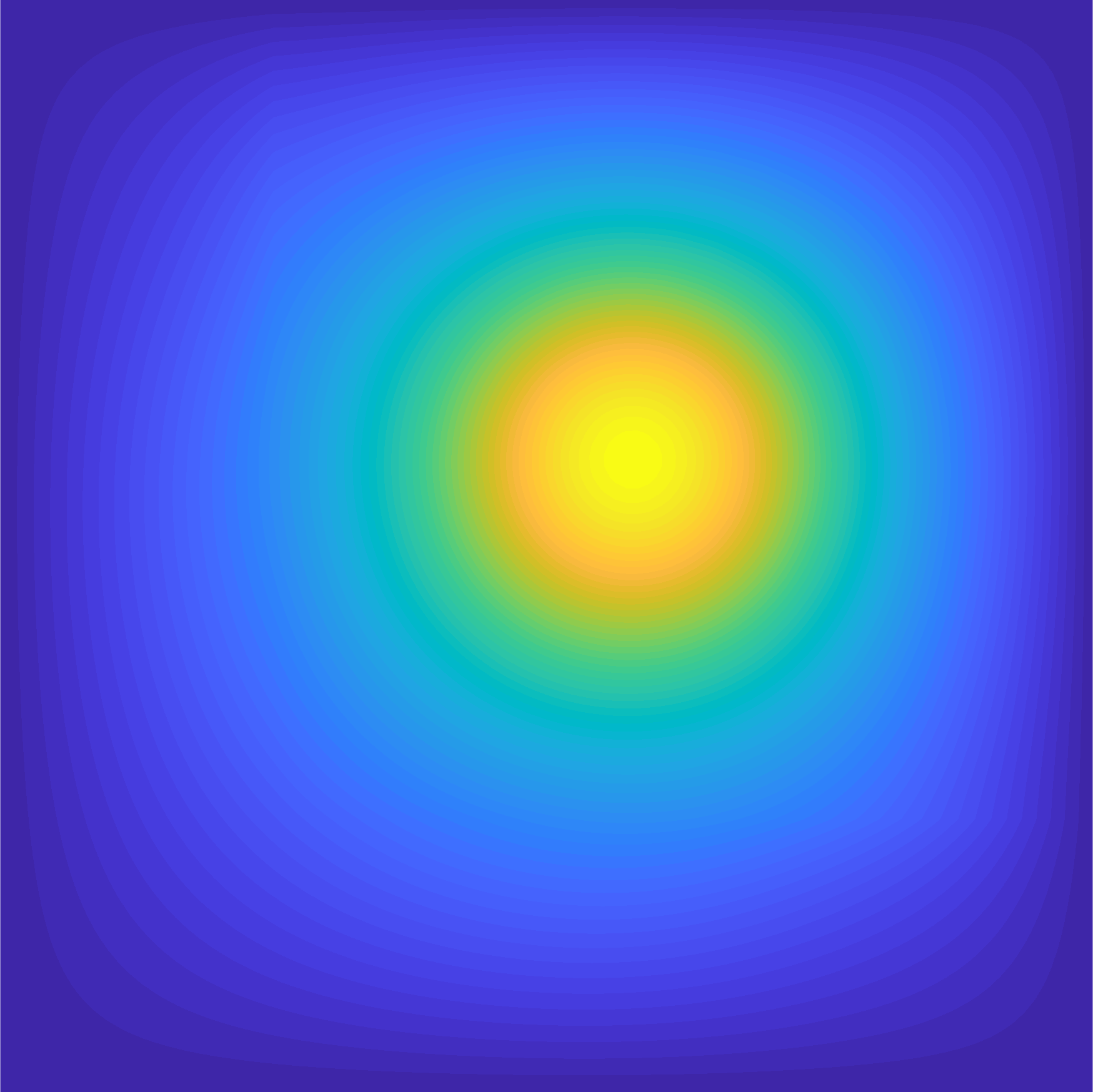};
						
			\nextgroupplot[ ylabel={}, ytick=\empty, xlabel={}, xtick=\empty]
			\addplot graphics [xmin=0, xmax=8.0, ymin=0, ymax=8.0] {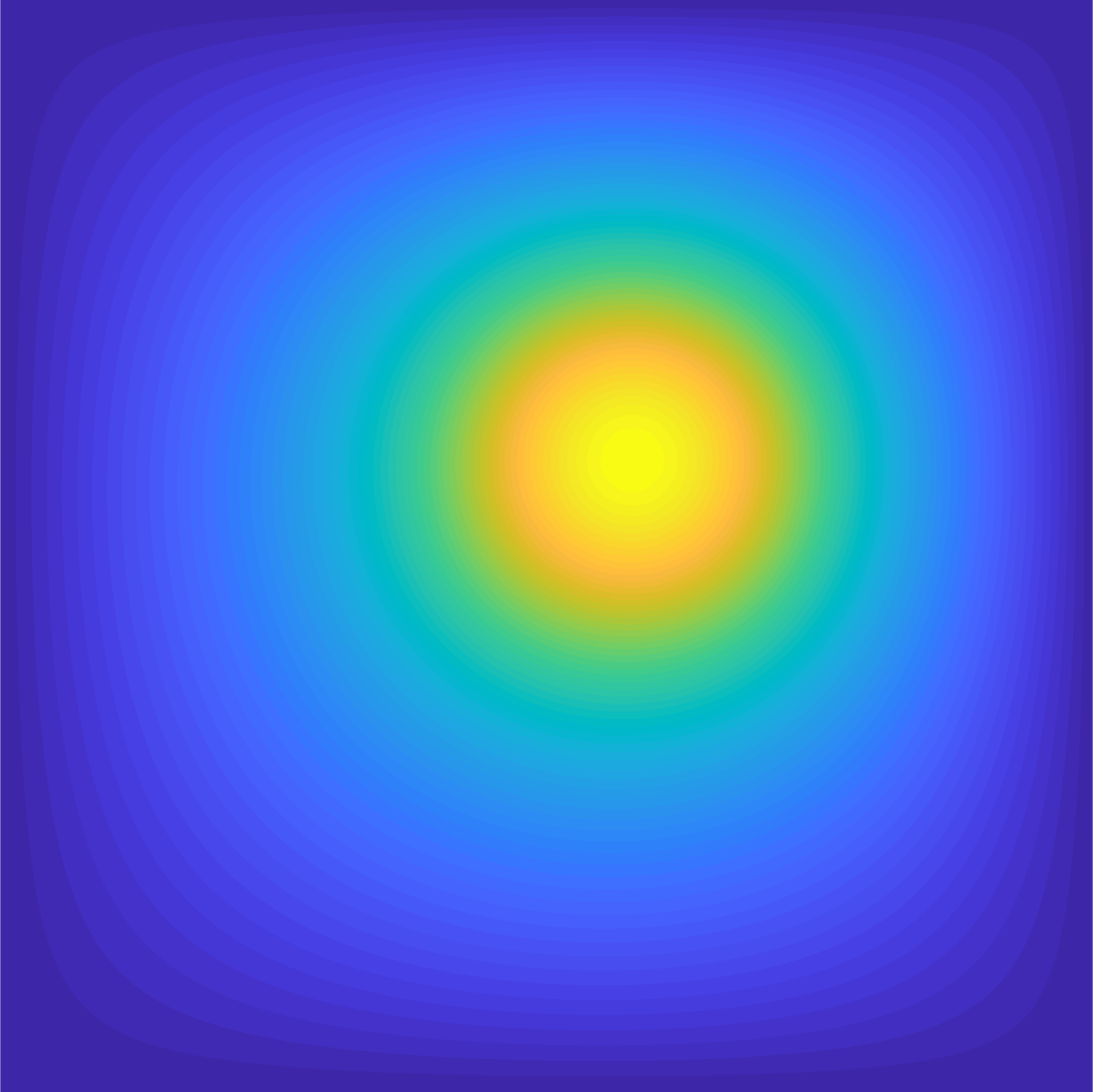};	
		\end{groupplot}
	\end{tikzpicture}
	\colorbarMatlabParula{0}{0.2130}{0.4260}{0.6391}{0.8521}
	\caption{Solution of the Poisson problem at $\Dcal_\mathrm{interp}$ (Section~\ref{sec:rslt:poi:limited}) using the CG-GL method at each adaptation iteration (\textit{left-to-right}, \textit{top-to-bottom}), as well as the reference solution (\textit{bottom right}). Because the CG-GL method is initialized with $\Omega_l = \emptyset$ and $N_l=0$, the \textit{top left} figure is a traditional ROM solution.}
	\label{fig:poi_interp_example_sol_sequential}
\end{figure}

\begin{figure}
	\centering
	\begin{tikzpicture}
		\begin{groupplot}[
			group style={
				group size=2 by 1,
				horizontal sep=0.1 cm,
			},
			width=0.32\textwidth,
			axis equal image,
			xlabel={$x_1$},
			ylabel={$x_2$},
			xtick = {0.0, 4.0, 8.0},
			ytick = {0.0, 4.0, 8.0},
			xmin=0, xmax=8,
			ymin=0, ymax=8
			]
			\nextgroupplot[ ylabel={}, ytick=\empty, xlabel={}, xtick=\empty]
			\addplot graphics [xmin=0, xmax=8.0, ymin=0, ymax=8.0] {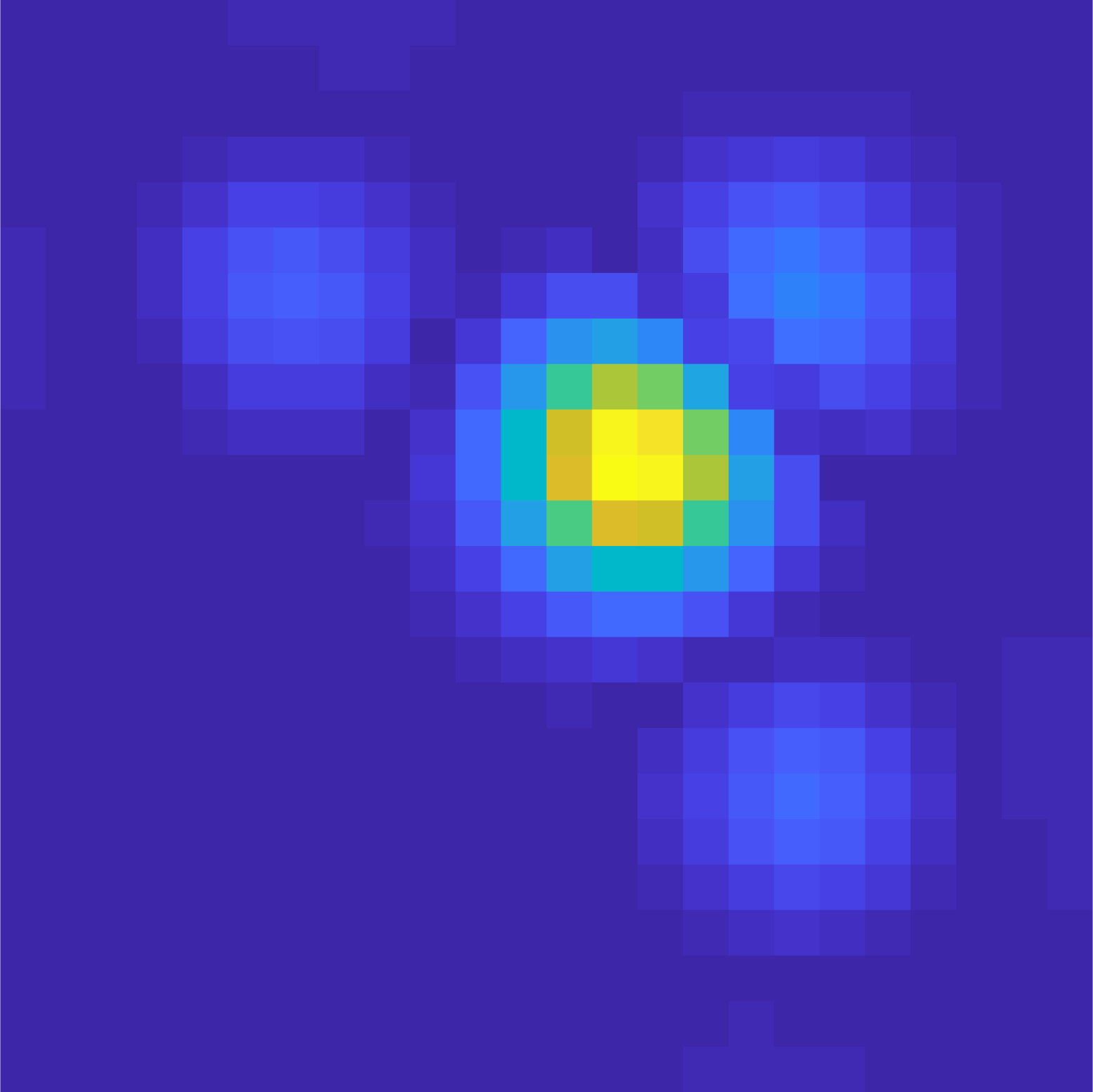};
			
			\nextgroupplot[ ylabel={}, ytick=\empty, xlabel={}, xtick=\empty]
			\addplot graphics [xmin=0, xmax=8.0, ymin=0, ymax=8.0] {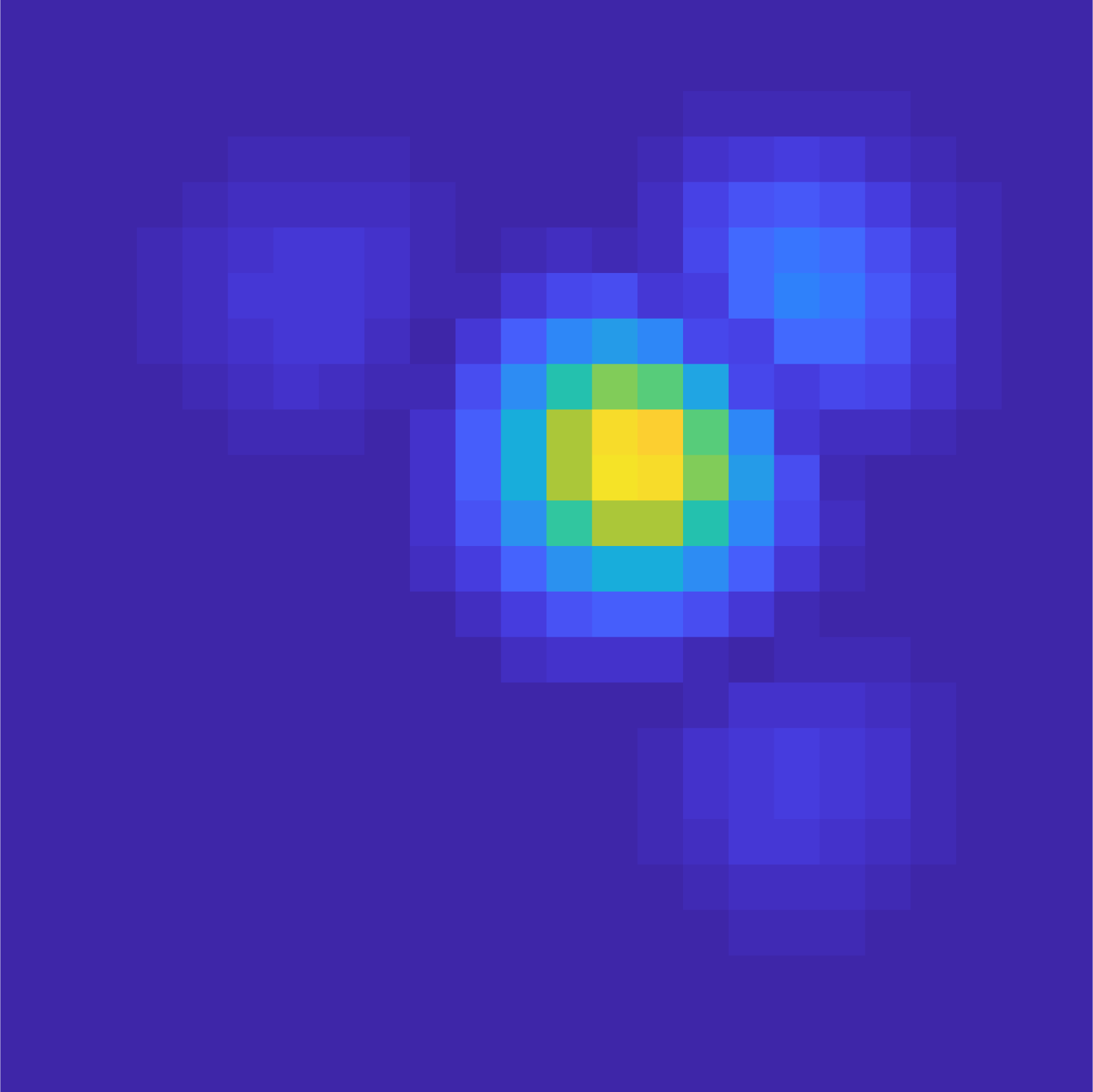};
		\end{groupplot}
	\end{tikzpicture}
\hspace{4pt}
	\begin{tikzpicture}
	\begin{groupplot}[
		group style={
			group size=2 by 1,
			horizontal sep=0.1 cm,
		},
		width=0.32\textwidth,
		axis equal image,
		xlabel={$x_1$},
		ylabel={$x_2$},
		xtick = {0.0, 4.0, 8.0},
		ytick = {0.0, 4.0, 8.0},
		xmin=0, xmax=8,
		ymin=0, ymax=8
		]
		\nextgroupplot[ ylabel={}, ytick=\empty, xlabel={}, xtick=\empty]
		\addplot graphics [xmin=0, xmax=8.0, ymin=0, ymax=8.0] {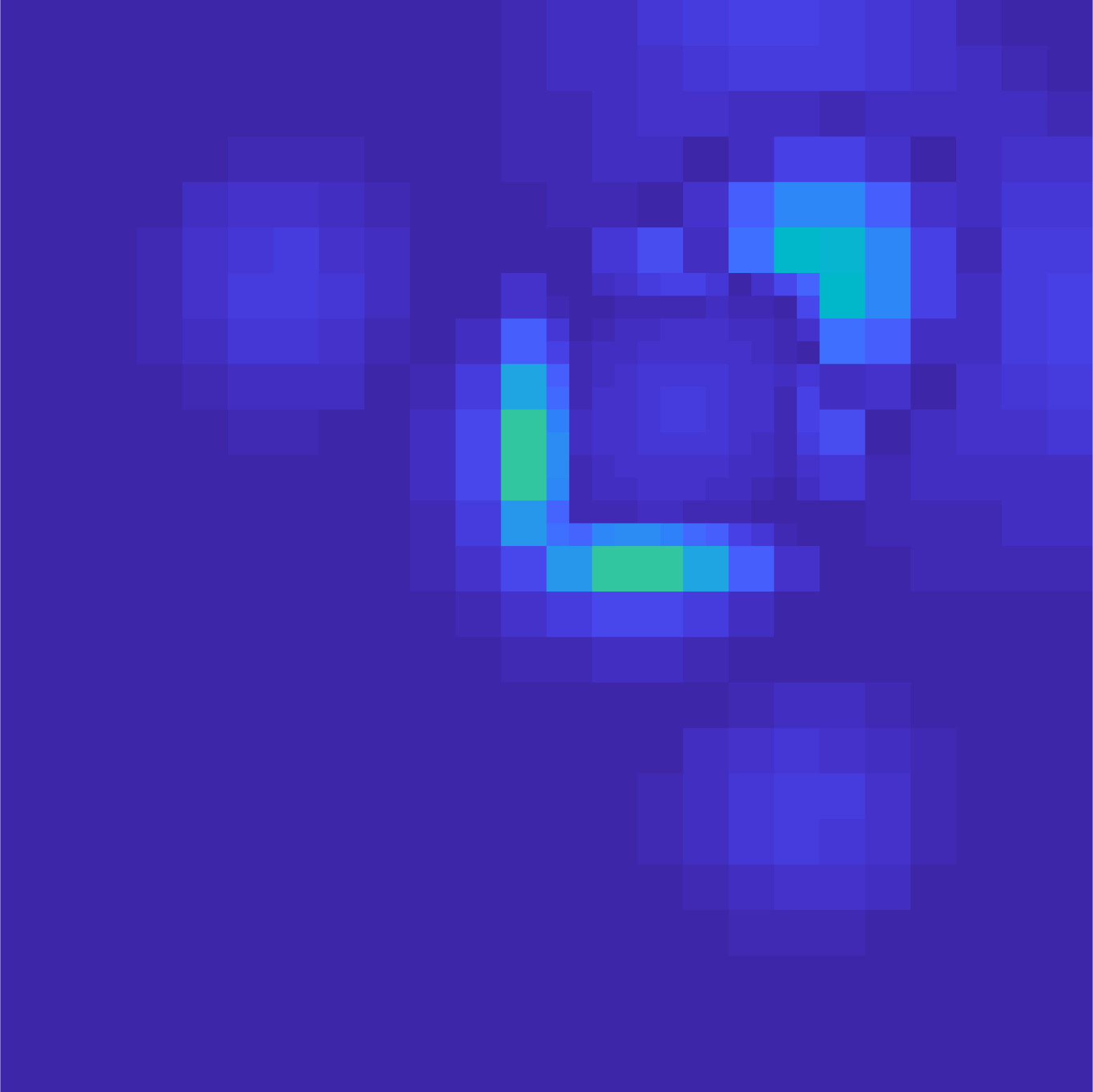};
		
		\nextgroupplot[ ylabel={}, ytick=\empty, xlabel={}, xtick=\empty]
		\addplot graphics [xmin=0, xmax=8.0, ymin=0, ymax=8.0] {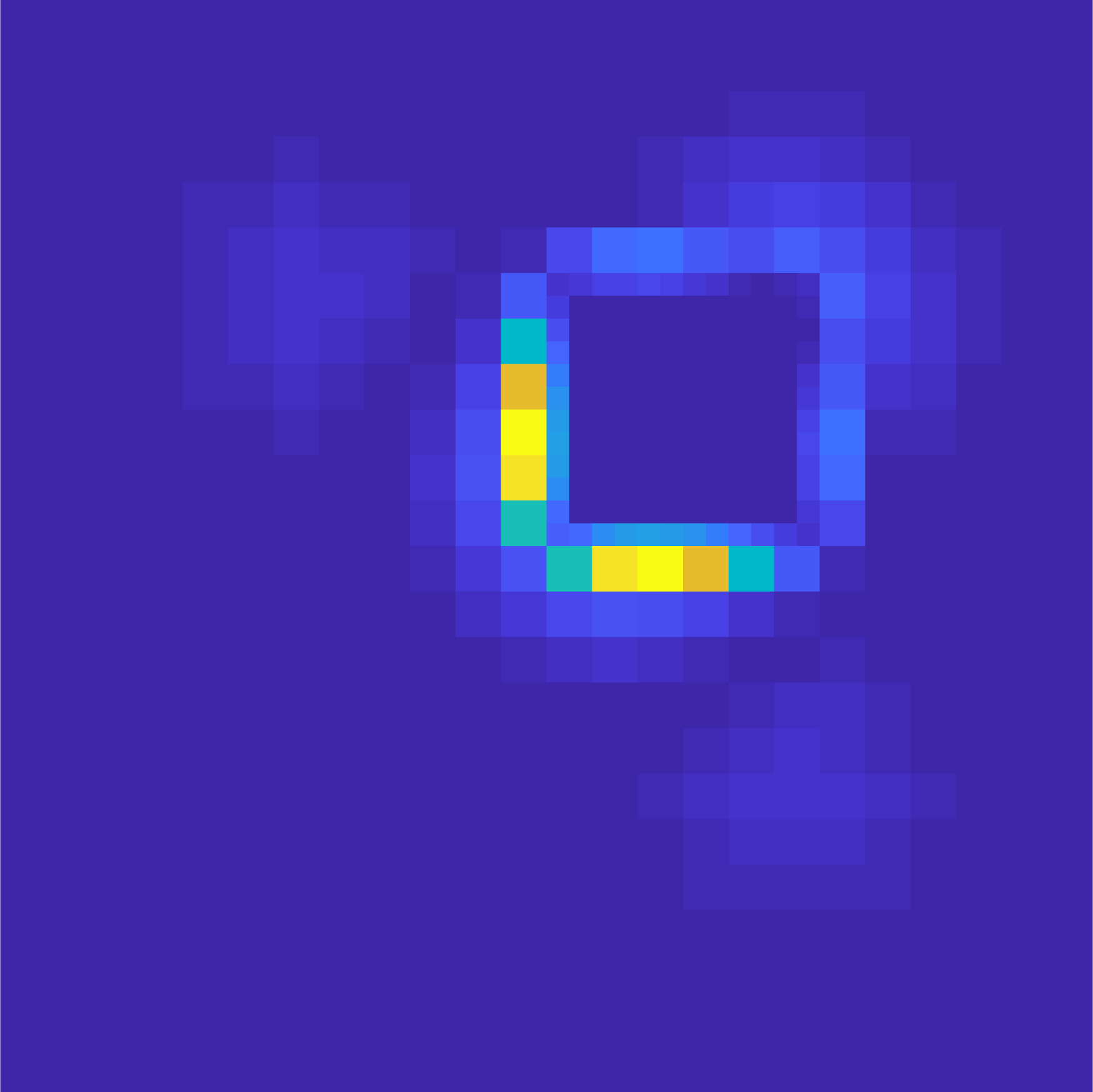};
	\end{groupplot}
\end{tikzpicture}\\
\vspace{8pt}
\begin{tikzpicture}
	\begin{groupplot}[
		group style={
			group size=2 by 1,
			horizontal sep=0.1 cm,
		},
		width=0.32\textwidth,
		axis equal image,
		xlabel={$x_1$},
		ylabel={$x_2$},
		xtick = {0.0, 4.0, 8.0},
		ytick = {0.0, 4.0, 8.0},
		xmin=0, xmax=8,
		ymin=0, ymax=8
		]
		\nextgroupplot[ ylabel={}, ytick=\empty, xlabel={}, xtick=\empty]
		\addplot graphics [xmin=0, xmax=8.0, ymin=0, ymax=8.0] {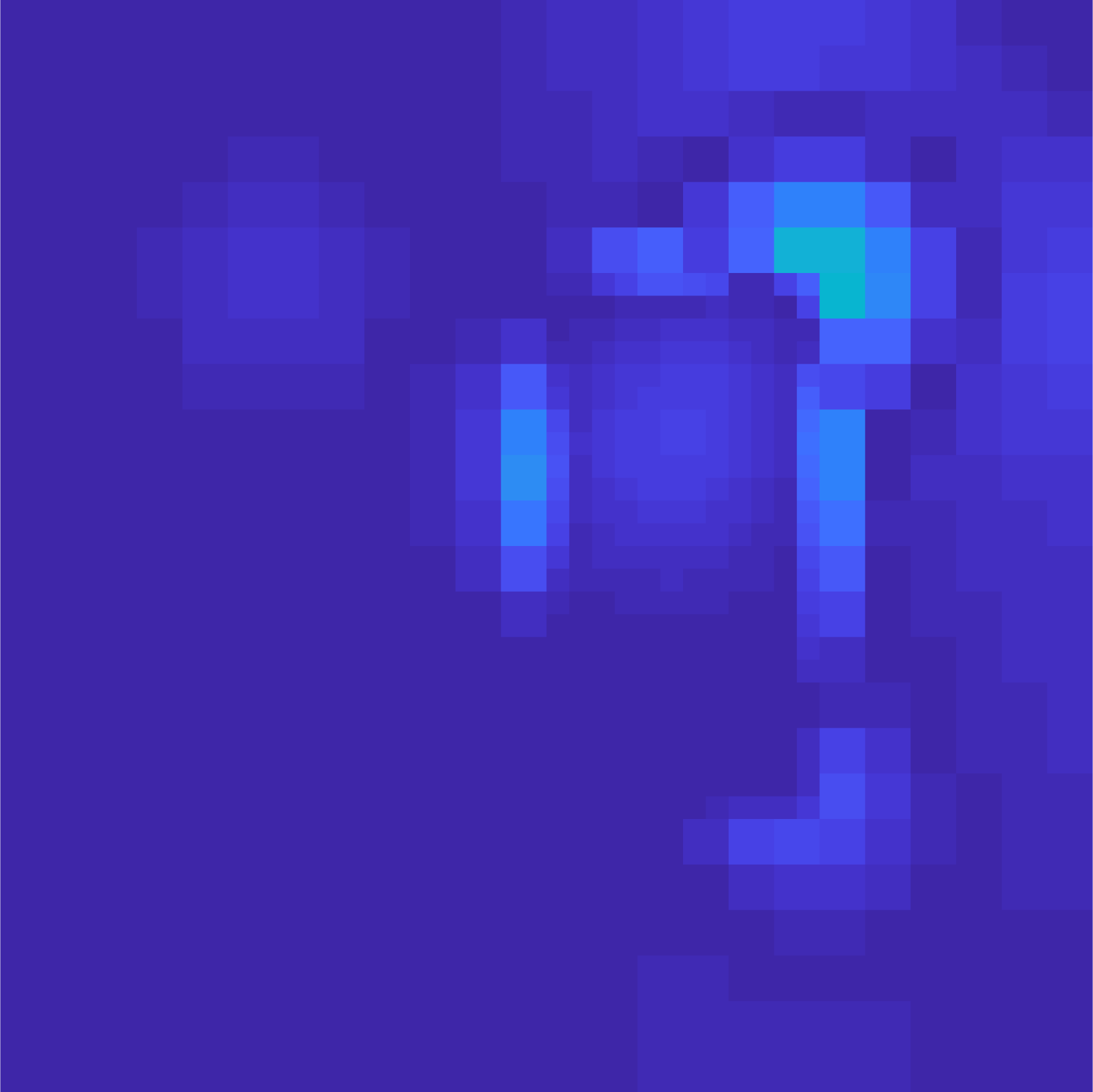};
		
		\nextgroupplot[ ylabel={}, ytick=\empty, xlabel={}, xtick=\empty]
		\addplot graphics [xmin=0, xmax=8.0, ymin=0, ymax=8.0] {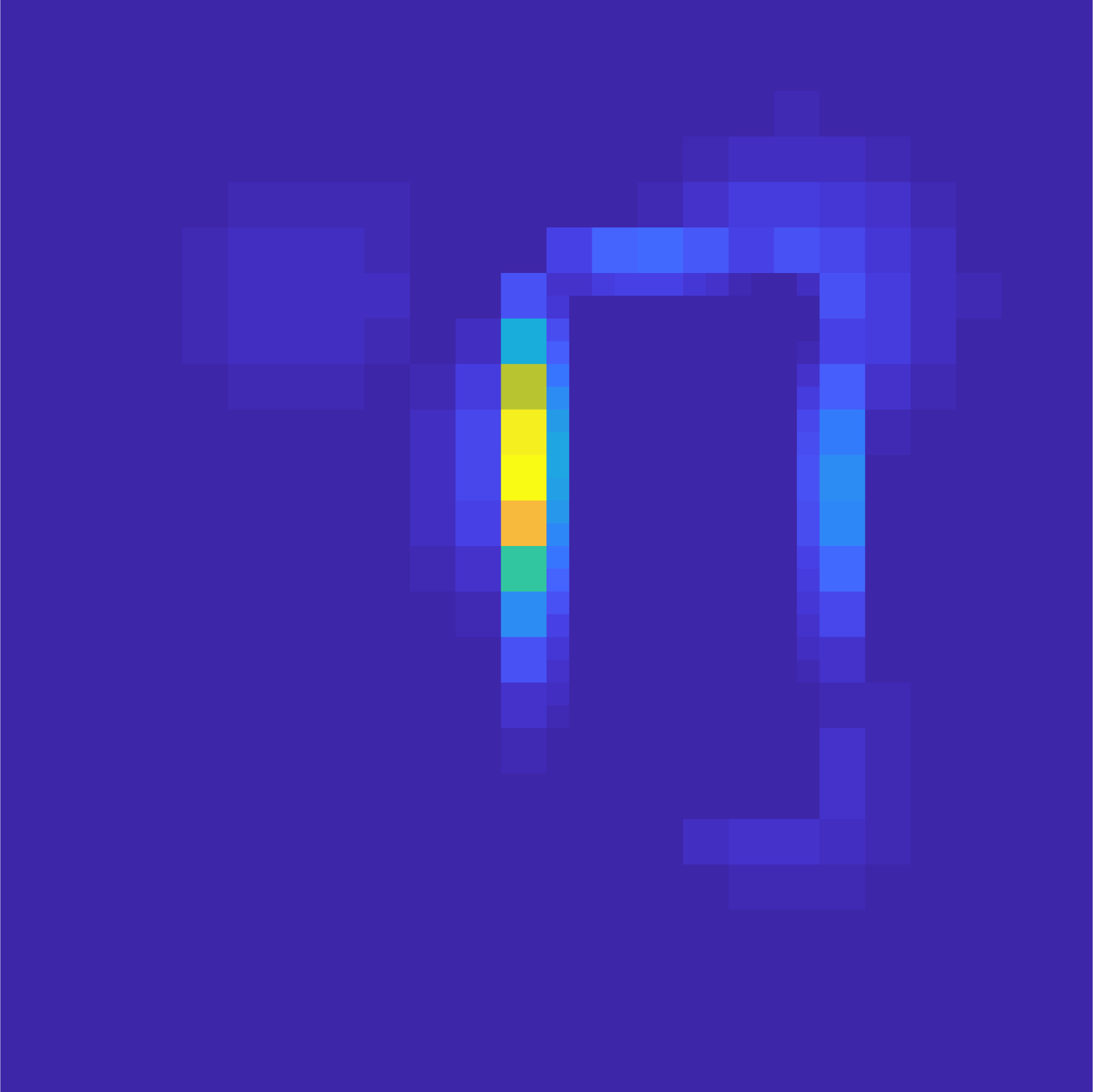};
	\end{groupplot}
\end{tikzpicture}
\hspace{4pt}
\begin{tikzpicture}
	\begin{groupplot}[
		group style={
			group size=2 by 1,
			horizontal sep=0.1 cm,
		},
		width=0.32\textwidth,
		axis equal image,
		xlabel={$x_1$},
		ylabel={$x_2$},
		xtick = {0.0, 4.0, 8.0},
		ytick = {0.0, 4.0, 8.0},
		xmin=0, xmax=8,
		ymin=0, ymax=8
		]
		\nextgroupplot[ ylabel={}, ytick=\empty, xlabel={}, xtick=\empty]
		\addplot graphics [xmin=0, xmax=8.0, ymin=0, ymax=8.0] {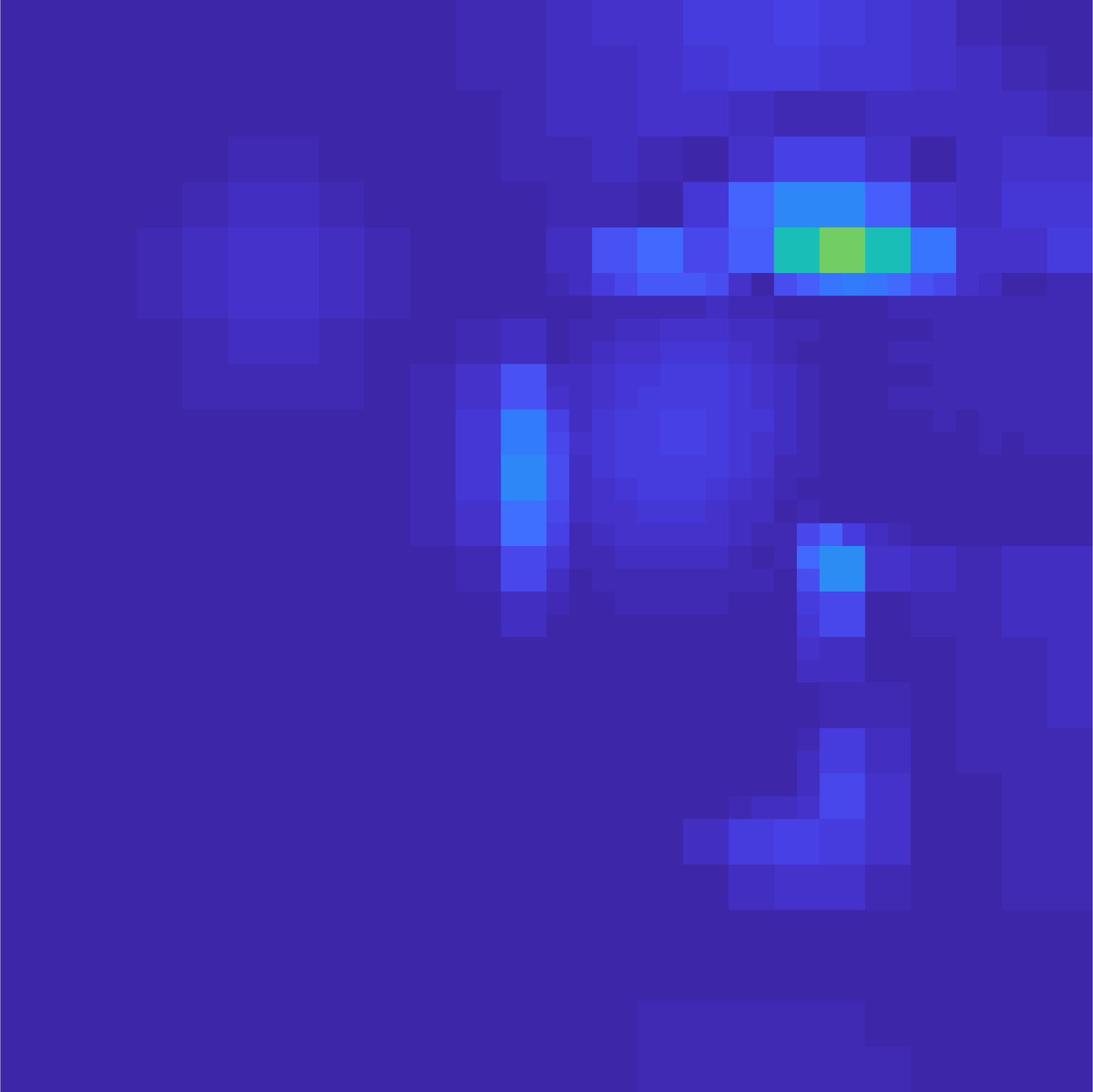};
		
		\nextgroupplot[ ylabel={}, ytick=\empty, xlabel={}, xtick=\empty]
		\addplot graphics [xmin=0, xmax=8.0, ymin=0, ymax=8.0] {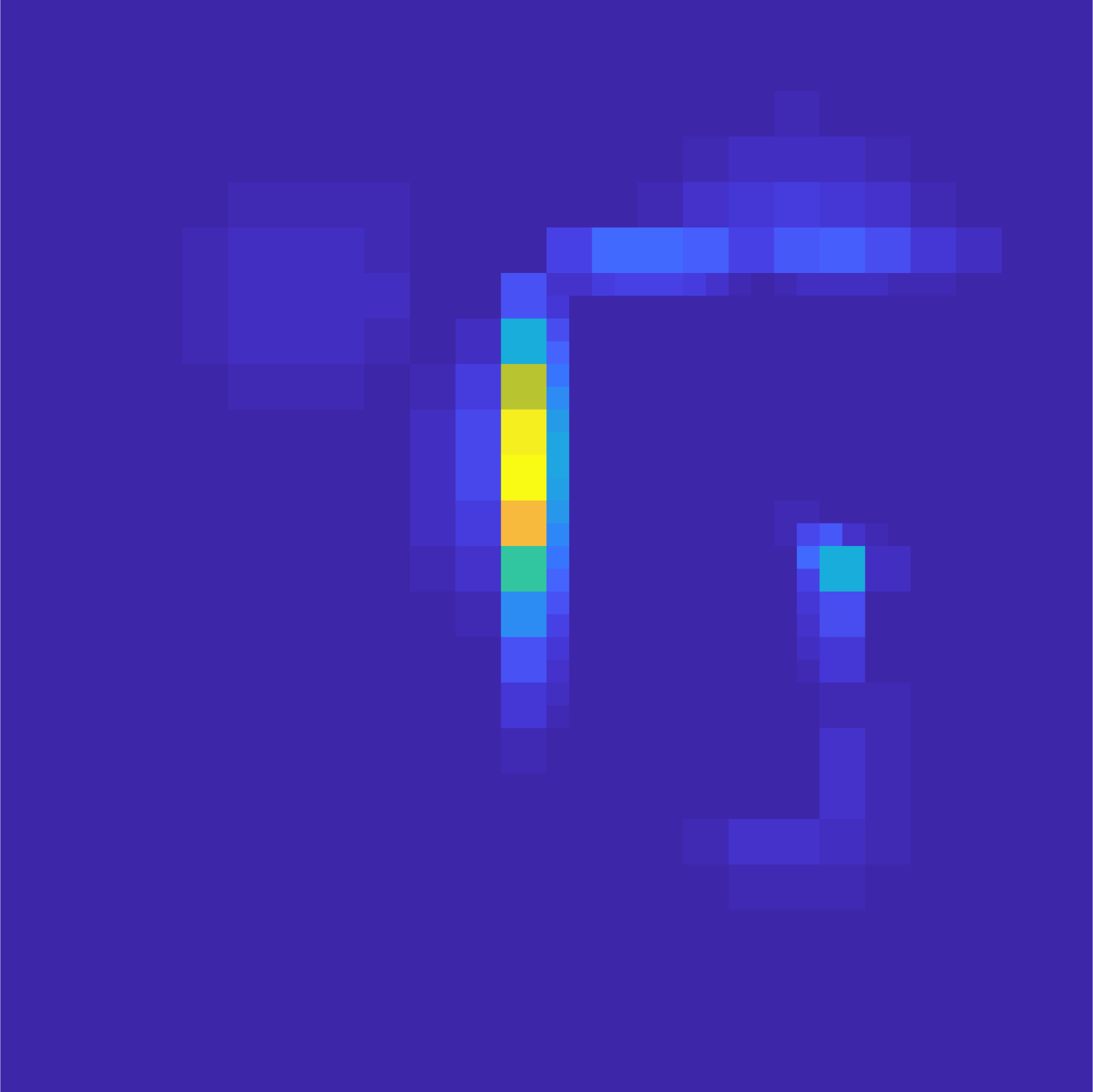};
	\end{groupplot}
\end{tikzpicture}

\vspace{8pt}
\begin{tikzpicture}
	\begin{groupplot}[
		group style={
			group size=2 by 1,
			horizontal sep=0.1 cm,
		},
		width=0.32\textwidth,
		axis equal image,
		xlabel={$x_1$},
		ylabel={$x_2$},
		xtick = {0.0, 4.0, 8.0},
		ytick = {0.0, 4.0, 8.0},
		xmin=0, xmax=8,
		ymin=0, ymax=8
		]
		\nextgroupplot[ ylabel={}, ytick=\empty, xlabel={}, xtick=\empty]
		\addplot graphics [xmin=0, xmax=8.0, ymin=0, ymax=8.0] {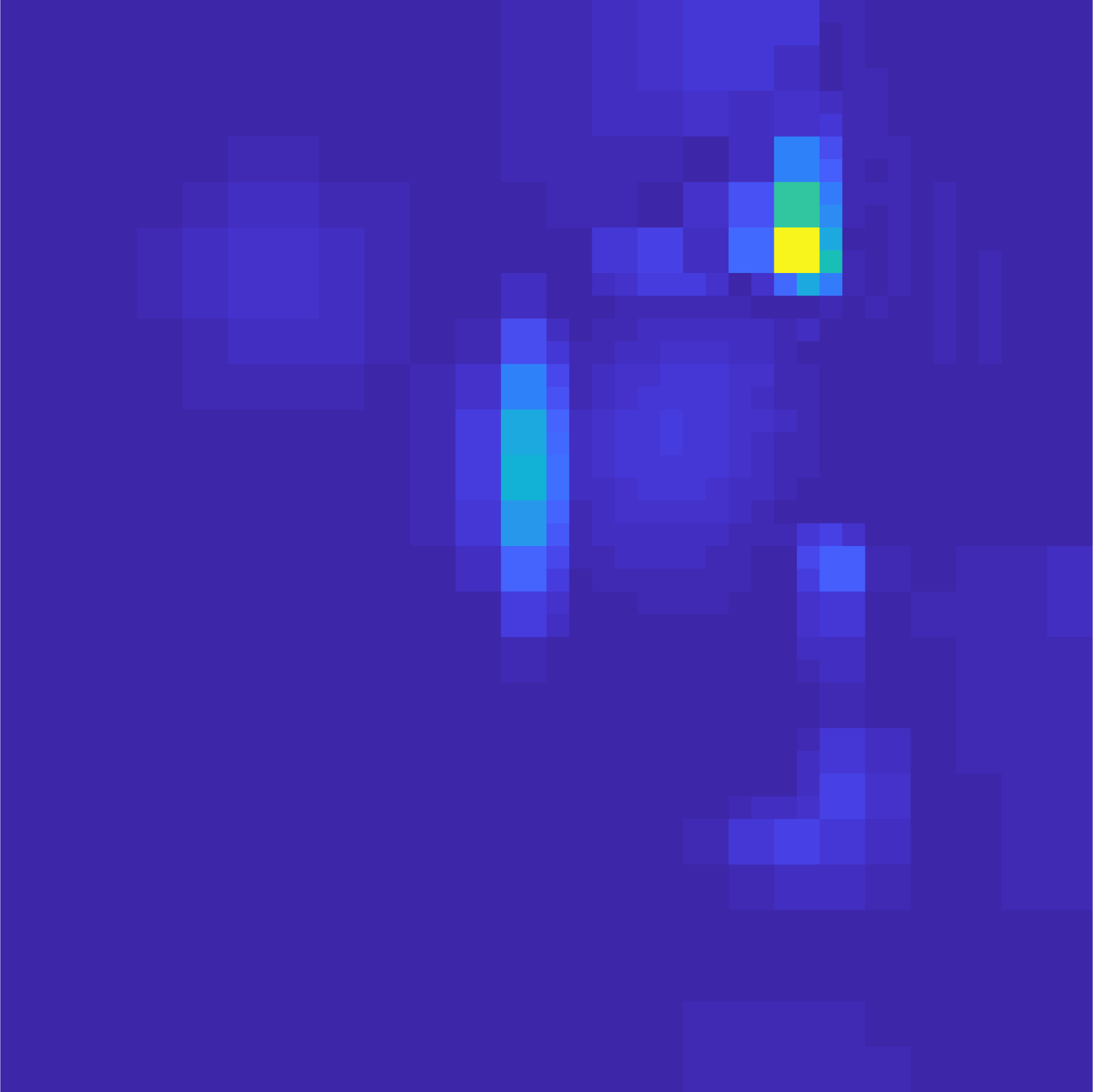};
		
		\nextgroupplot[ ylabel={}, ytick=\empty, xlabel={}, xtick=\empty]
		\addplot graphics [xmin=0, xmax=8.0, ymin=0, ymax=8.0] {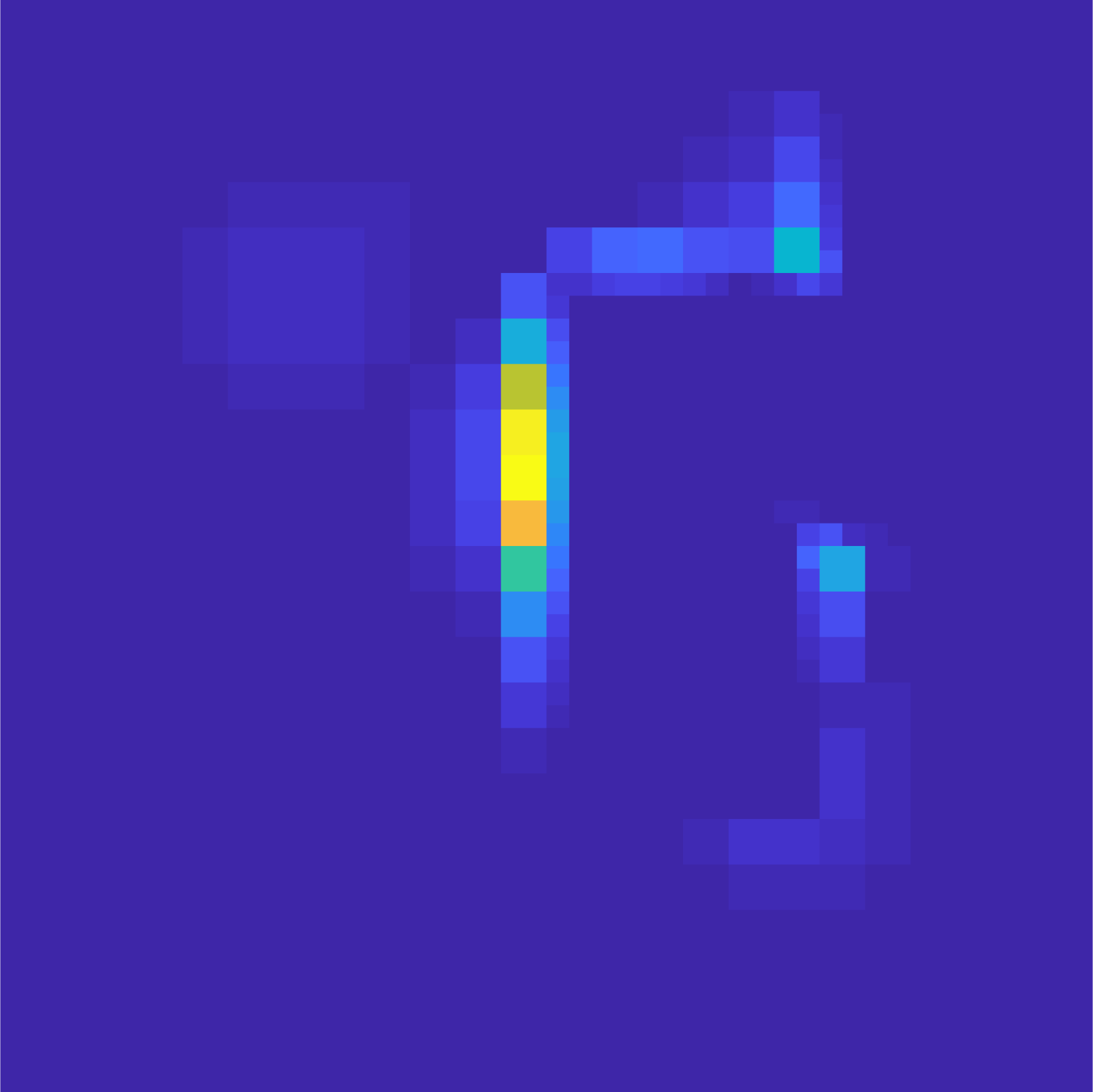};
	\end{groupplot}
\end{tikzpicture}
\hspace{4pt}
\begin{tikzpicture}
	\begin{groupplot}[
		group style={
			group size=2 by 1,
			horizontal sep=0.1 cm,
		},
		width=0.32\textwidth,
		axis equal image,
		xlabel={$x_1$},
		ylabel={$x_2$},
		xtick = {0.0, 4.0, 8.0},
		ytick = {0.0, 4.0, 8.0},
		xmin=0, xmax=8,
		ymin=0, ymax=8
		]
		\nextgroupplot[ ylabel={}, ytick=\empty, xlabel={}, xtick=\empty]
		\addplot graphics [xmin=0, xmax=8.0, ymin=0, ymax=8.0] {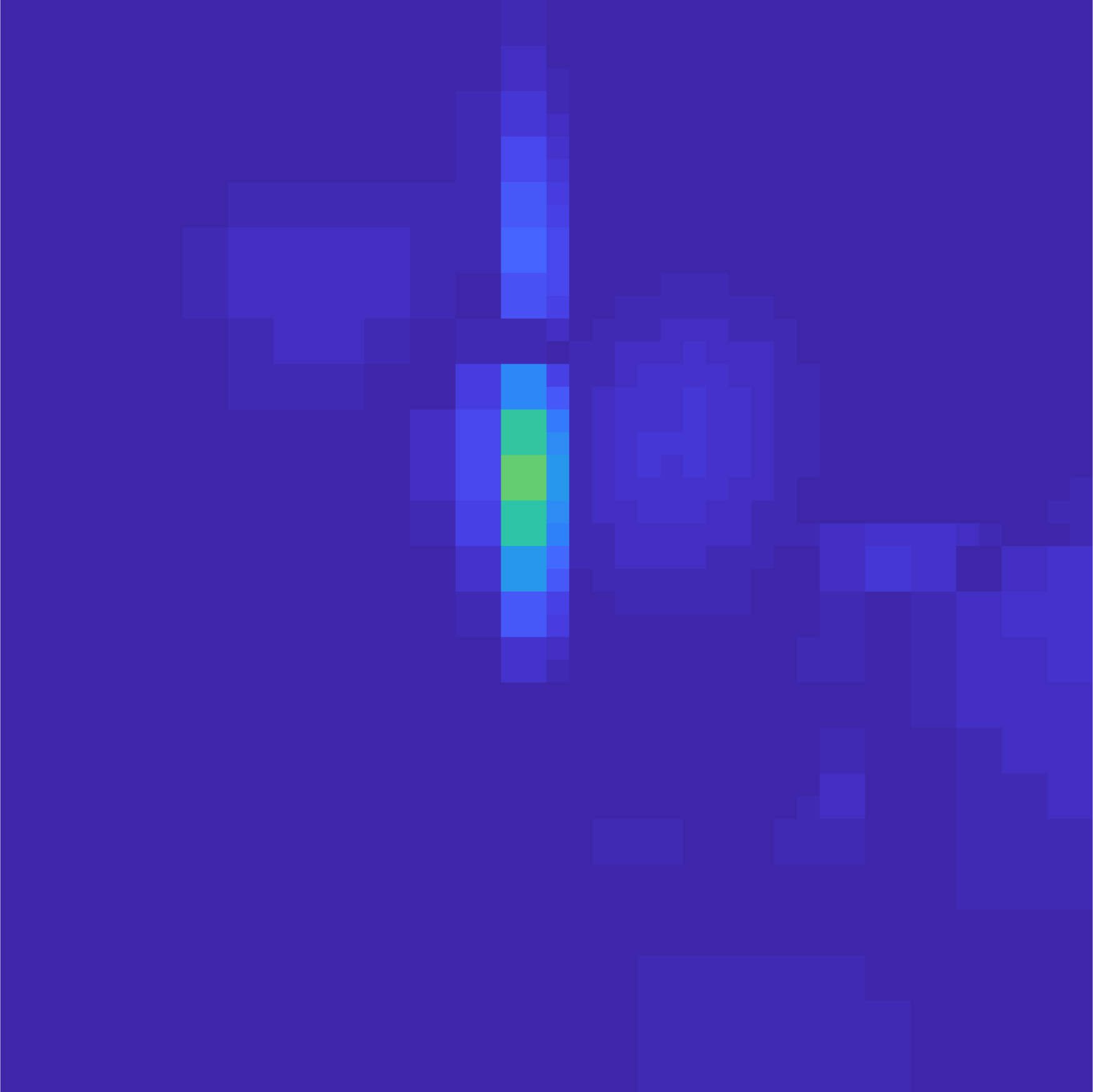};
		
		\nextgroupplot[ ylabel={}, ytick=\empty, xlabel={}, xtick=\empty]
		\addplot graphics [xmin=0, xmax=8.0, ymin=0, ymax=8.0] {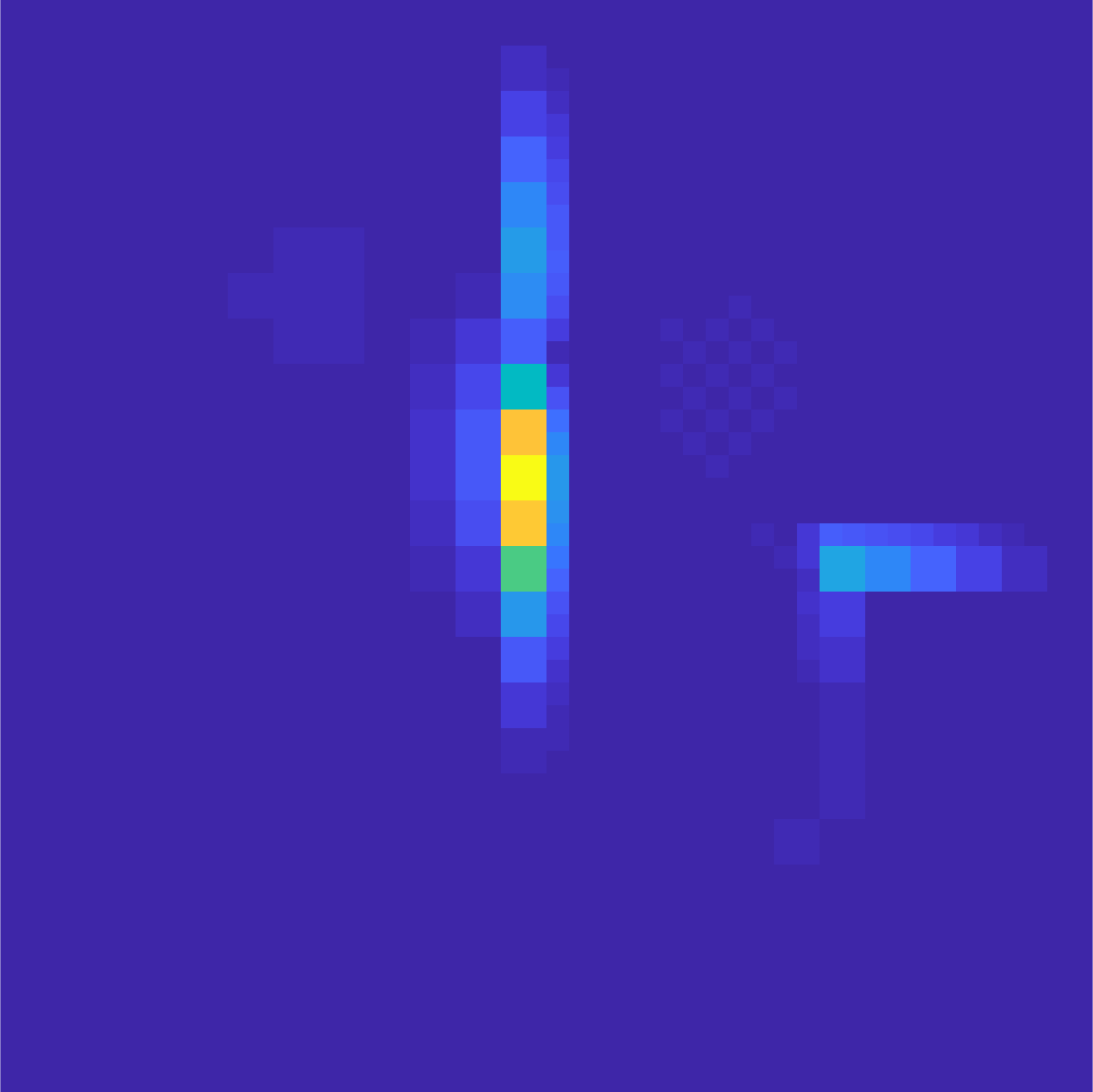};
	\end{groupplot}
\end{tikzpicture}

\vspace{8pt}
\begin{tikzpicture}
	\begin{groupplot}[
		group style={
			group size=2 by 1,
			horizontal sep=0.1 cm,
		},
		width=0.32\textwidth,
		axis equal image,
		xlabel={$x_1$},
		ylabel={$x_2$},
		xtick = {0.0, 4.0, 8.0},
		ytick = {0.0, 4.0, 8.0},
		xmin=0, xmax=8,
		ymin=0, ymax=8
		]
		\nextgroupplot[ ylabel={}, ytick=\empty, xlabel={}, xtick=\empty]
		\addplot graphics [xmin=0, xmax=8.0, ymin=0, ymax=8.0] {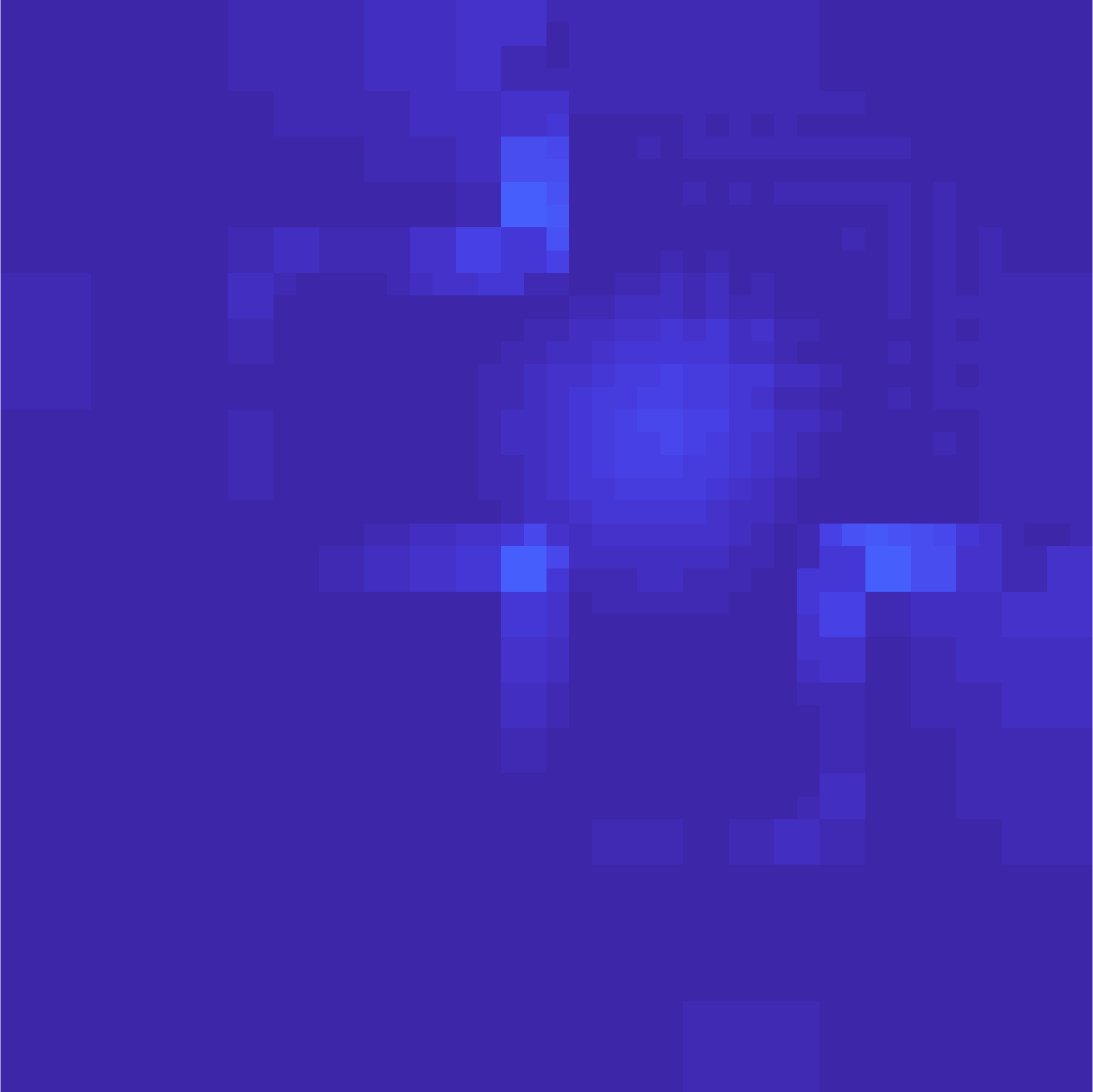};
		
		\nextgroupplot[ ylabel={}, ytick=\empty, xlabel={}, xtick=\empty]
		\addplot graphics [xmin=0, xmax=8.0, ymin=0, ymax=8.0] {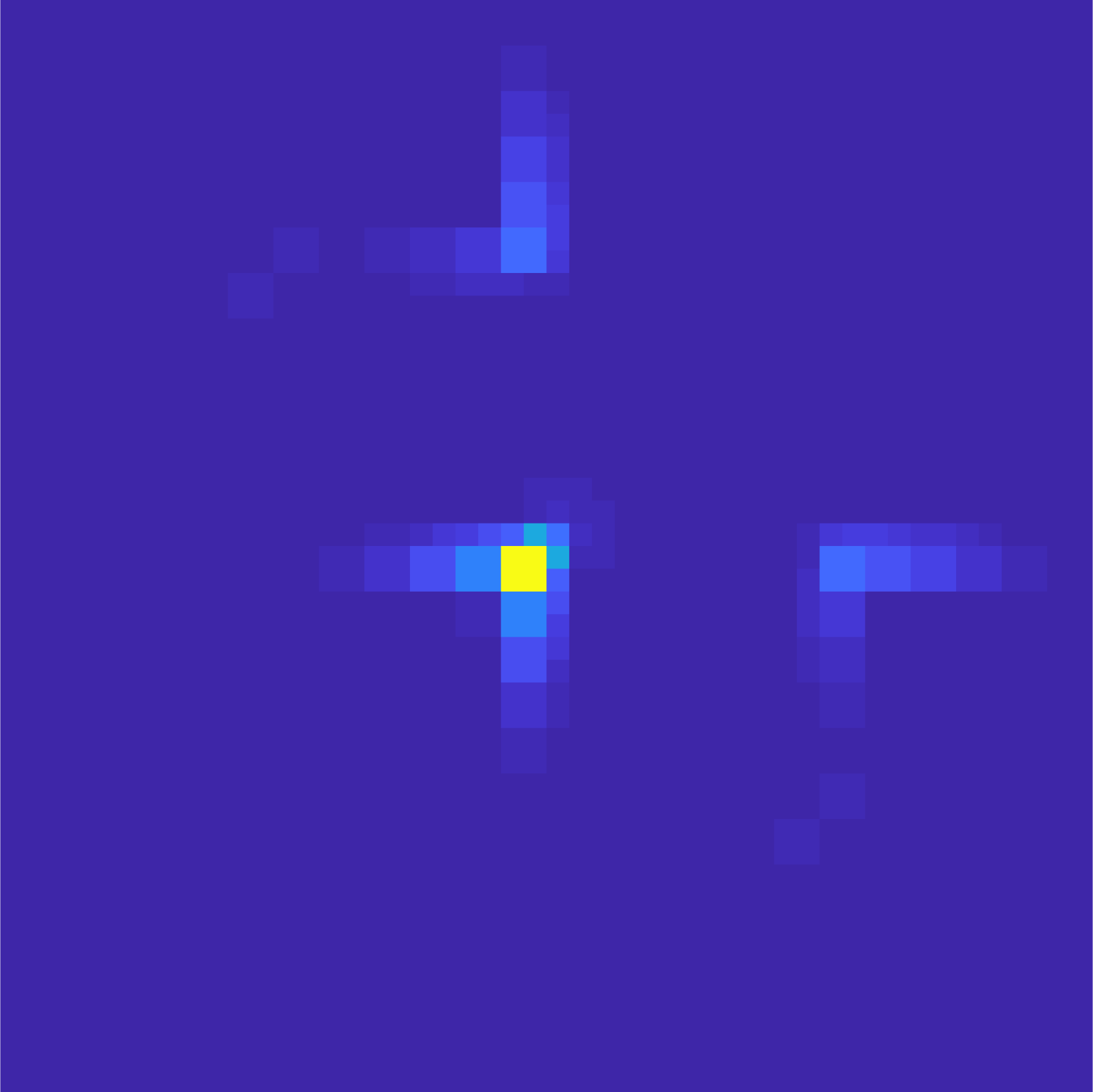};
	\end{groupplot}
\end{tikzpicture}
\hspace{4pt}
\begin{tikzpicture}
	\begin{groupplot}[
		group style={
			group size=2 by 1,
			horizontal sep=0.1 cm,
		},
		width=0.32\textwidth,
		axis equal image,
		xlabel={$x_1$},
		ylabel={$x_2$},
		xtick = {0.0, 4.0, 8.0},
		ytick = {0.0, 4.0, 8.0},
		xmin=0, xmax=8,
		ymin=0, ymax=8
		]
		\nextgroupplot[ ylabel={}, ytick=\empty, xlabel={}, xtick=\empty]
		\addplot graphics [xmin=0, xmax=8.0, ymin=0, ymax=8.0] {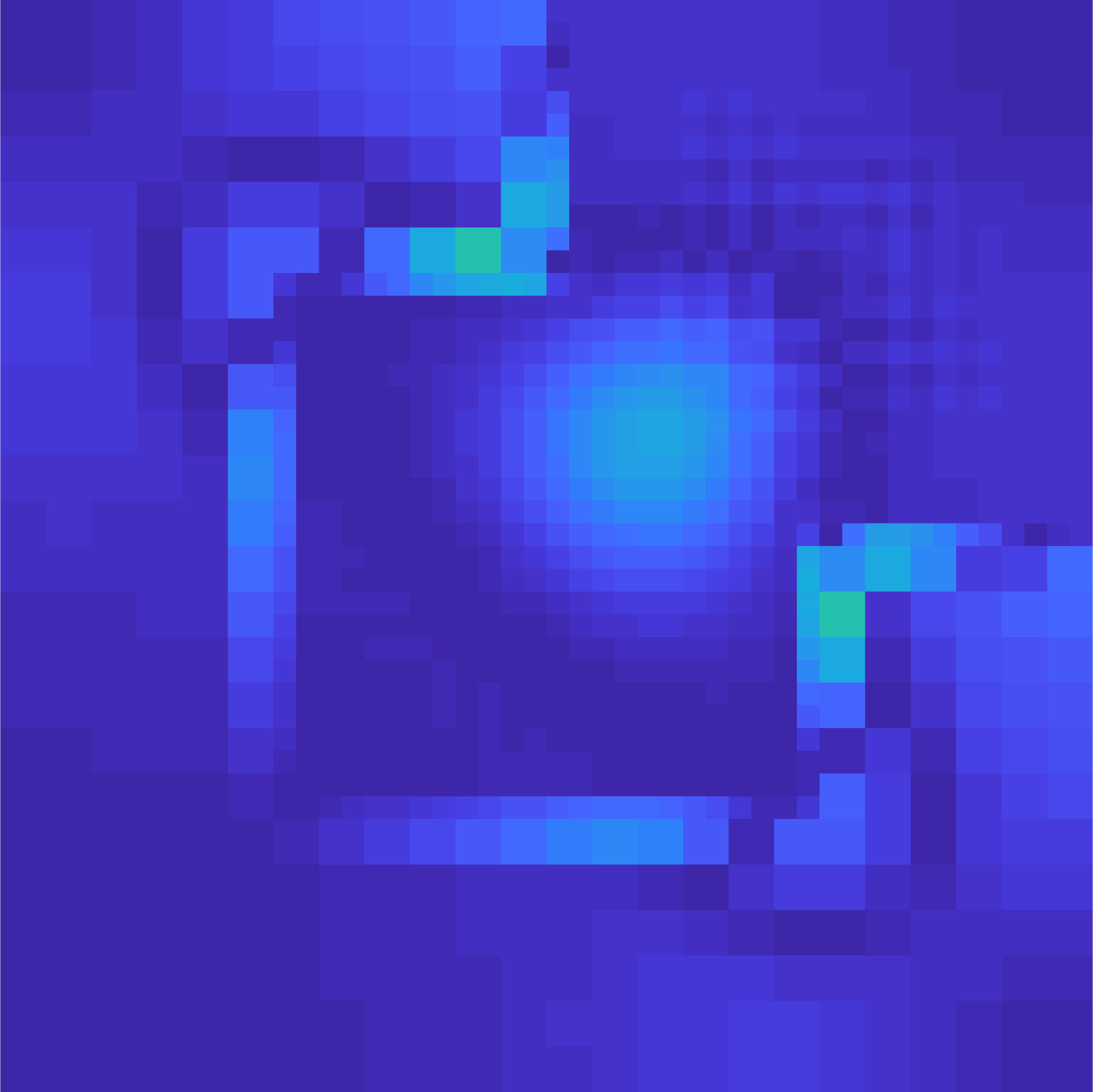};
		
		\nextgroupplot[ ylabel={}, ytick=\empty, xlabel={}, xtick=\empty]
		\addplot graphics [xmin=0, xmax=8.0, ymin=0, ymax=8.0] {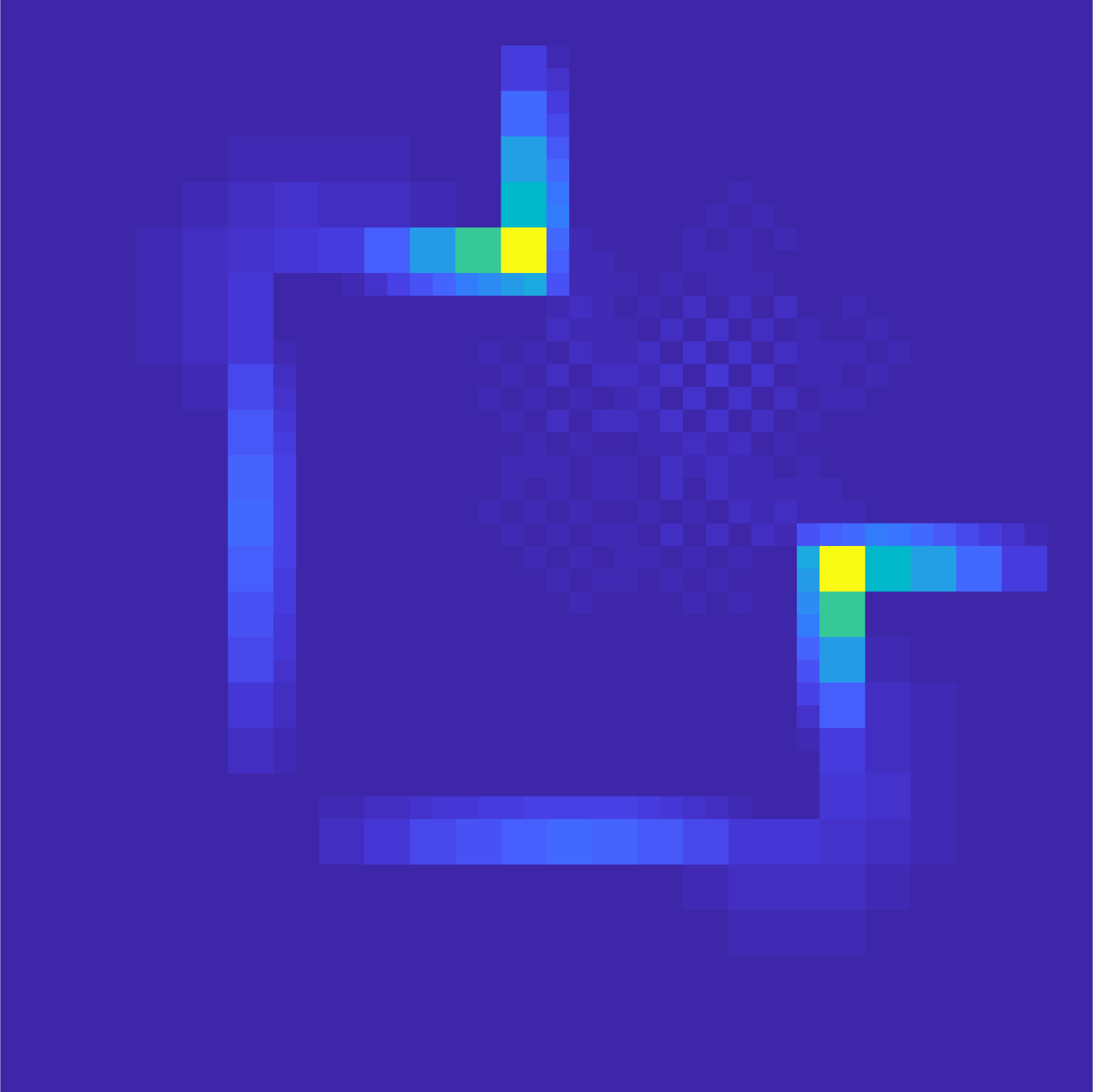};
	\end{groupplot}
\end{tikzpicture}

	\caption{The localized DWR error estimate (\textit{left}) and true error (\textit{right}) at the first eight CG-GL adaptation
	iterations (\textit{left-to-right, top-to-bottom}) for the Poisson problem at $\Dcal_\mathrm{interp}$ (Section~\ref{sec:rslt:poi:limited}).
	Each pair of figures corresponds to a different CG-GL iteration and use its own colorbar, scaled to its range to highlight similarity
	between the proposed error estimate and actual error.}
	\label{fig:error_contour_0}
\end{figure}

\begin{figure}
	\centering
	\begin{tikzpicture}
		\begin{groupplot}[
			group style={
				group size=4 by 3,
				horizontal sep=1cm,
			},
			width=0.32\textwidth,
			axis equal image,
			xlabel={$x_1$},
			ylabel={$x_2$},
			xtick = {0.0, 4.0, 8.0},
			ytick = {0.0, 4.0, 8.0},
			xmin=0, xmax=8,
			ymin=0, ymax=8
			]
			\nextgroupplot[ ylabel={}, ytick=\empty, xlabel={}, xtick=\empty]
			\addplot graphics [xmin=0, xmax=8.0, ymin=0, ymax=8.0] {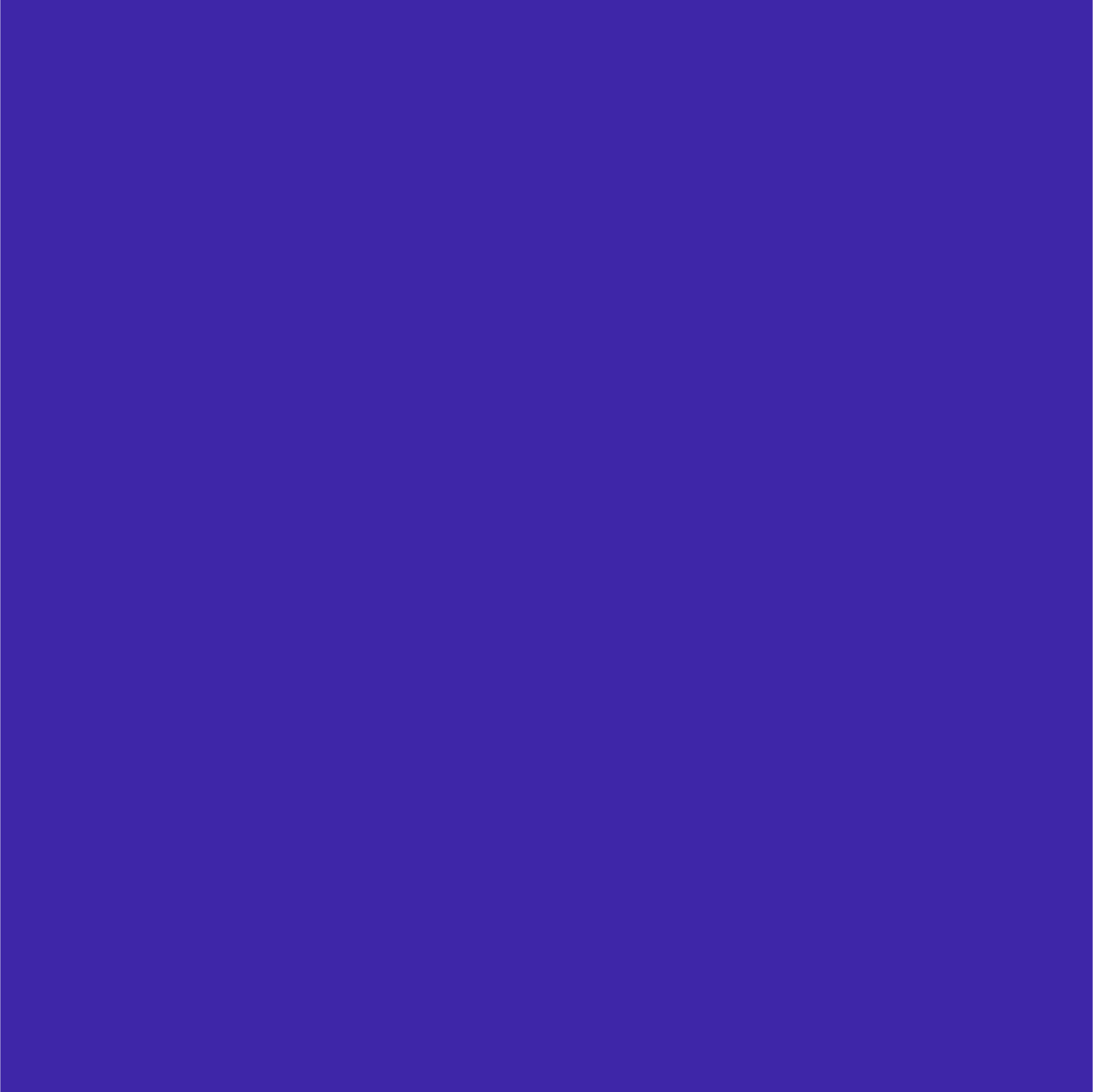};
			
			\nextgroupplot[ ylabel={}, ytick=\empty, xlabel={}, xtick=\empty]
			\addplot graphics [xmin=0, xmax=8.0, ymin=0, ymax=8.0] {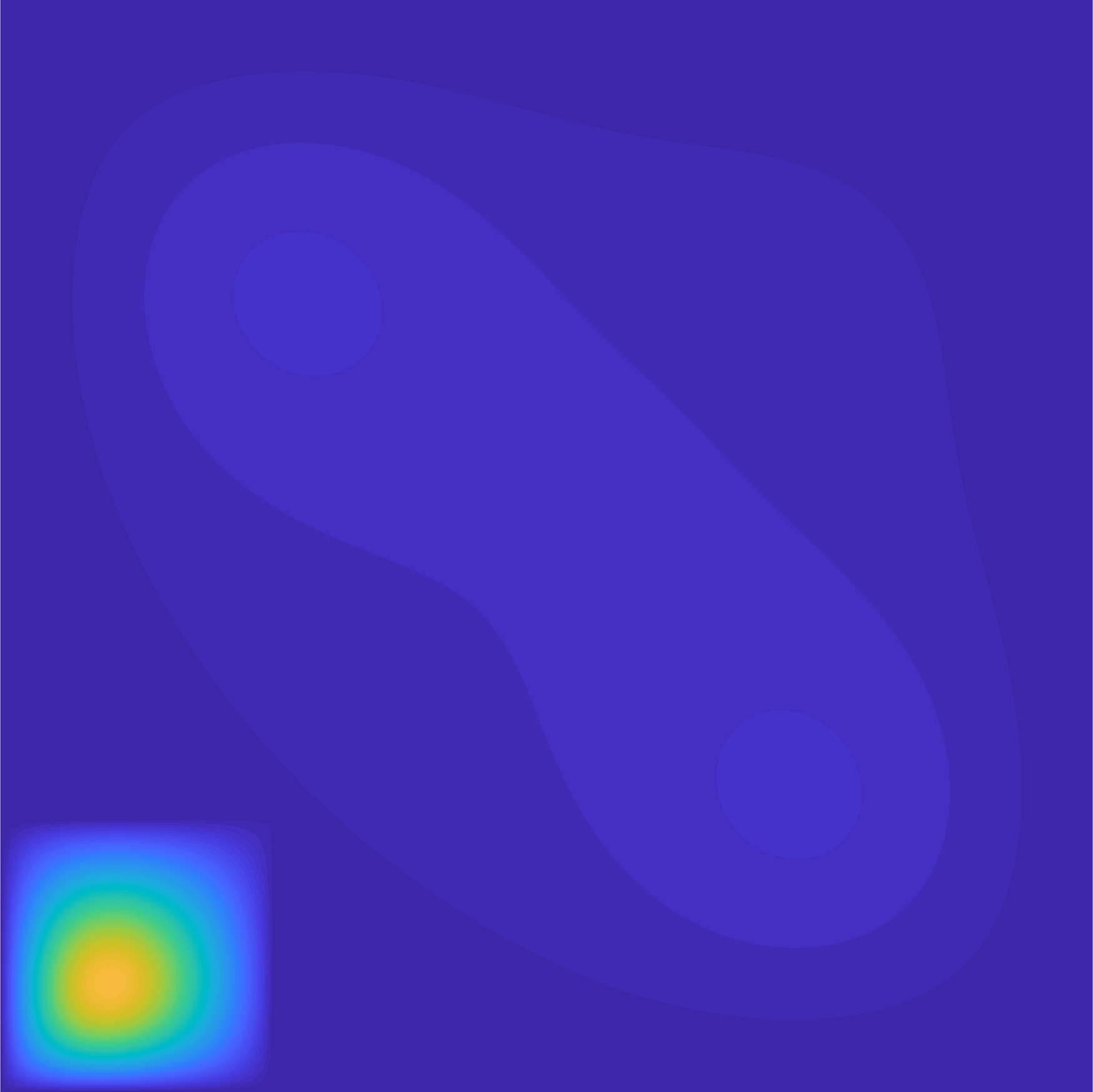};
			
			\nextgroupplot[ ylabel={}, ytick=\empty, xlabel={}, xtick=\empty]
			\addplot graphics [xmin=0, xmax=8.0, ymin=0, ymax=8.0] {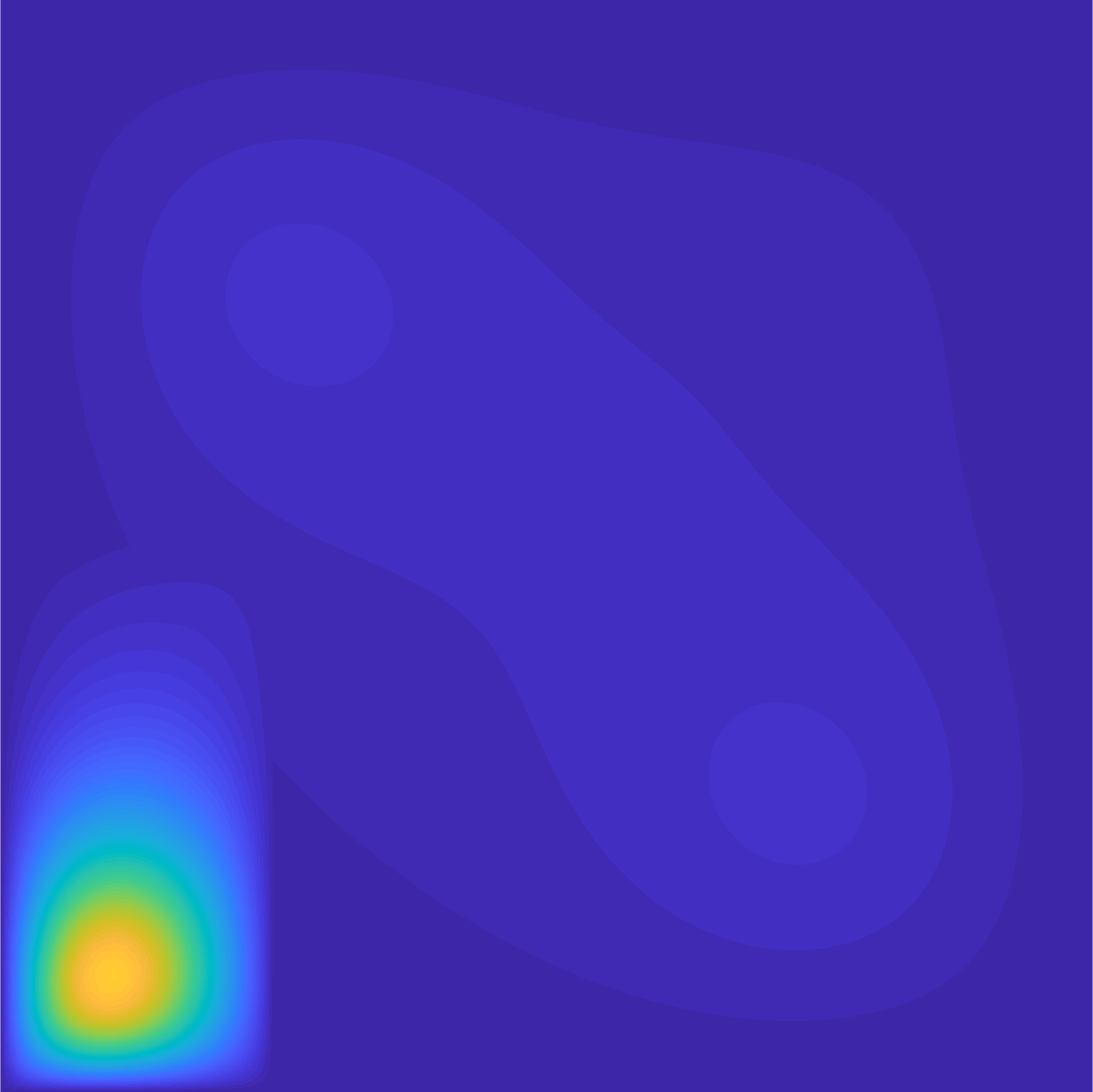};
			
			\nextgroupplot[ ylabel={}, ytick=\empty, xlabel={}, xtick=\empty]
			\addplot graphics [xmin=0, xmax=8.0, ymin=0, ymax=8.0] {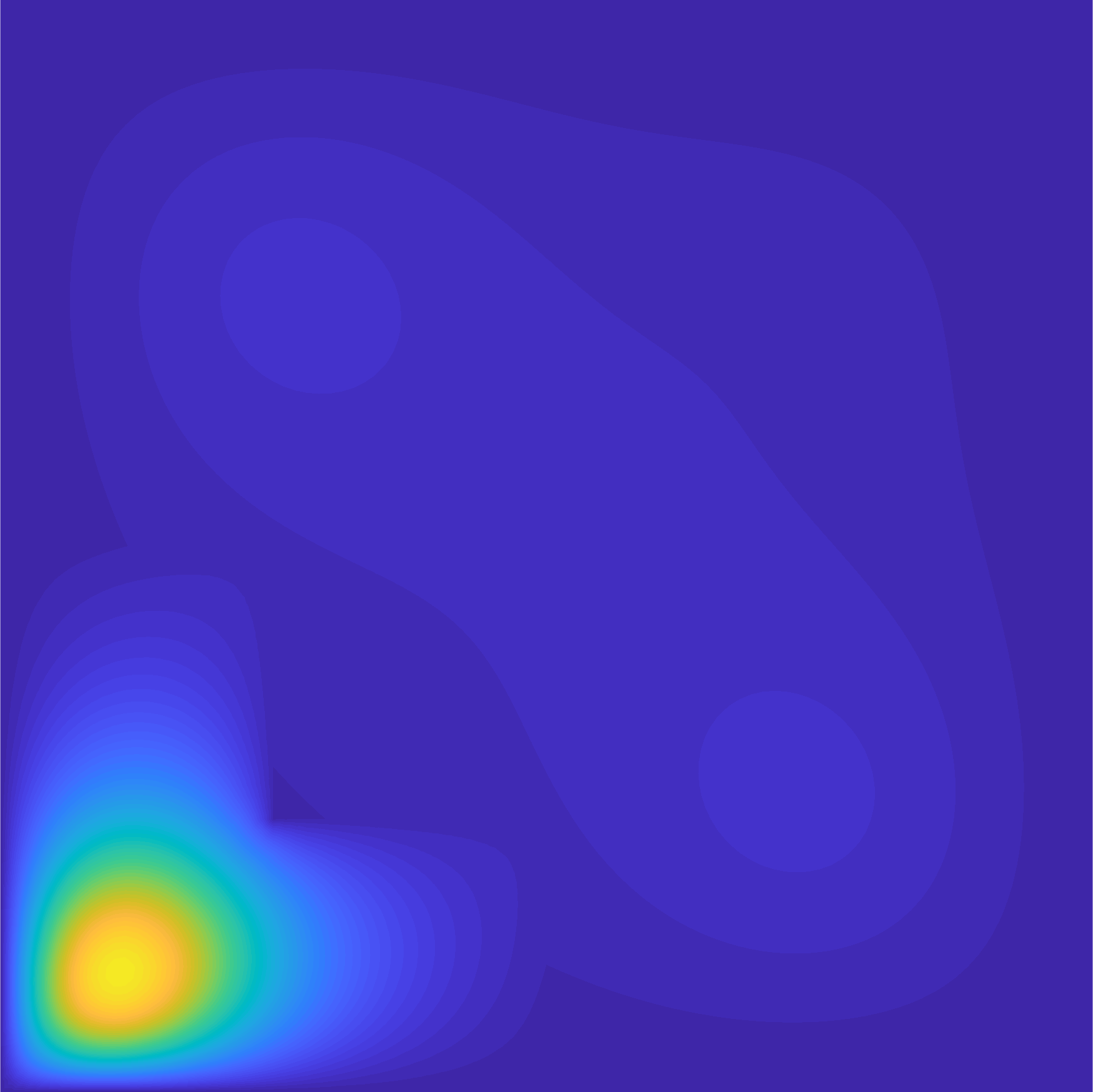};
			
			\nextgroupplot[ ylabel={}, ytick=\empty, xlabel={}, xtick=\empty]
			\addplot graphics [xmin=0, xmax=8.0, ymin=0, ymax=8.0] {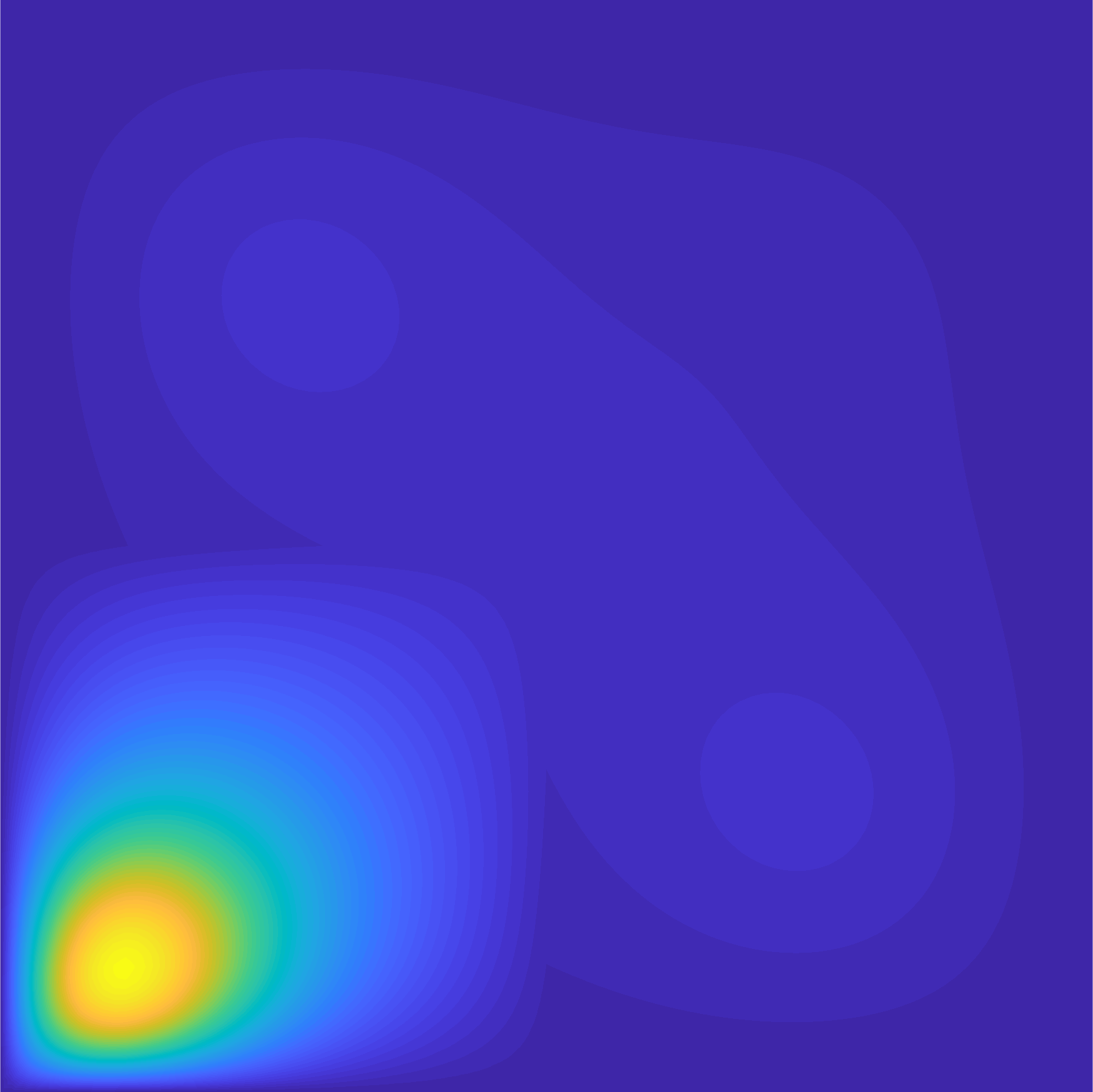};
			
			\nextgroupplot[ ylabel={}, ytick=\empty, xlabel={}, xtick=\empty]
			\addplot graphics [xmin=0, xmax=8.0, ymin=0, ymax=8.0] {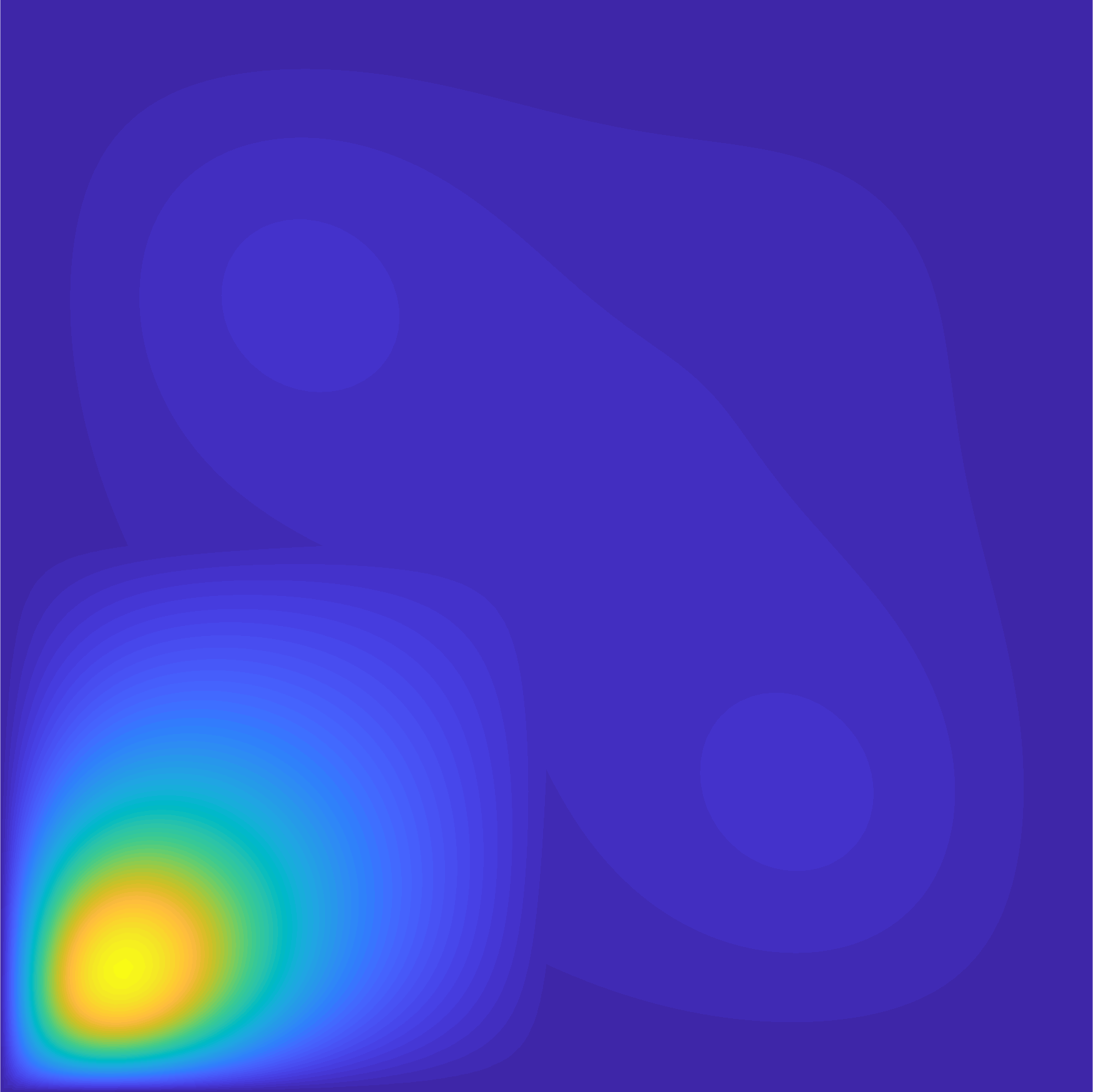};
			
			\nextgroupplot[ ylabel={}, ytick=\empty, xlabel={}, xtick=\empty]
			\addplot graphics [xmin=0, xmax=8.0, ymin=0, ymax=8.0] {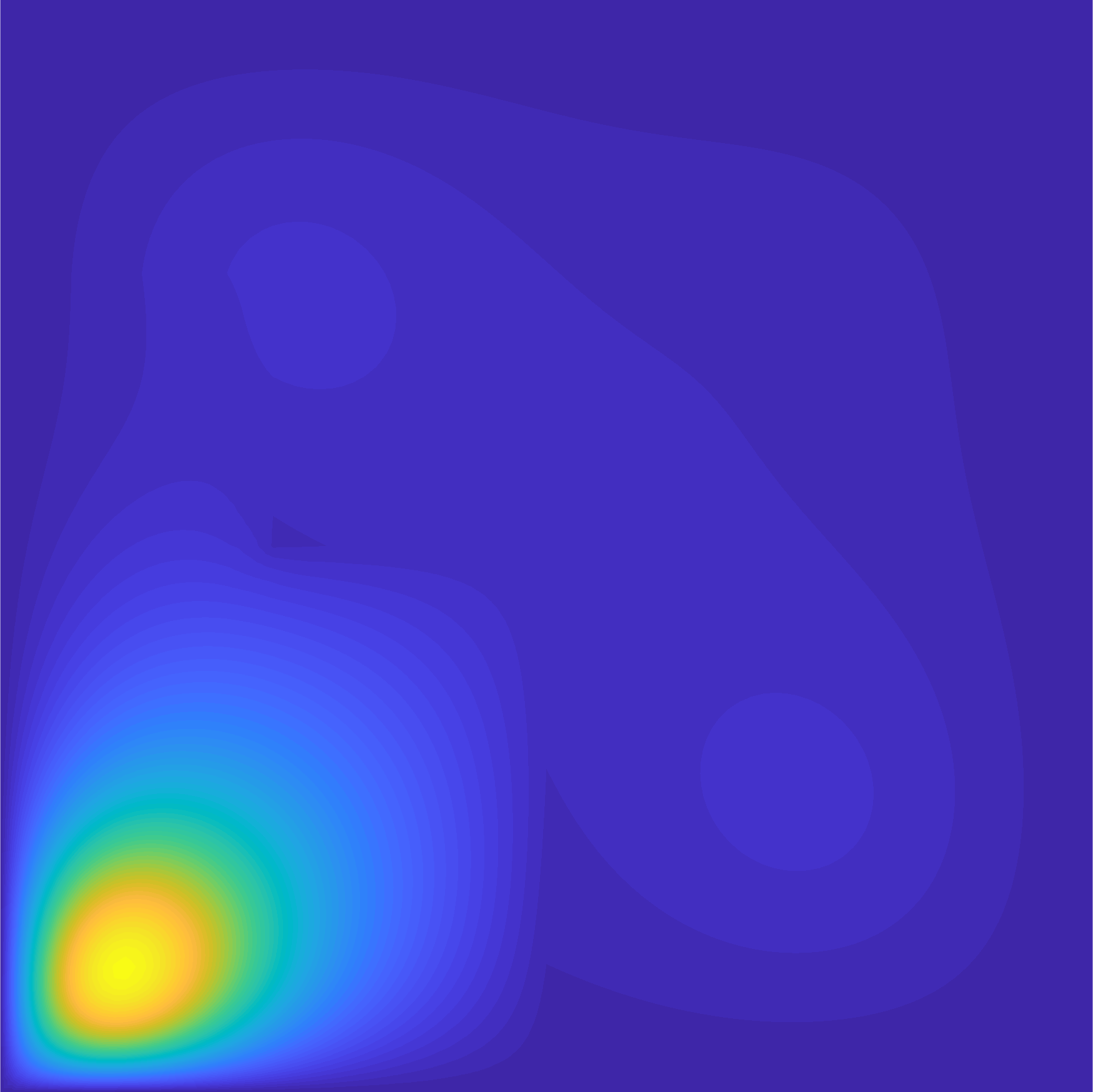};
			
			\nextgroupplot[ ylabel={}, ytick=\empty, xlabel={}, xtick=\empty]
			\addplot graphics [xmin=0, xmax=8.0, ymin=0, ymax=8.0] {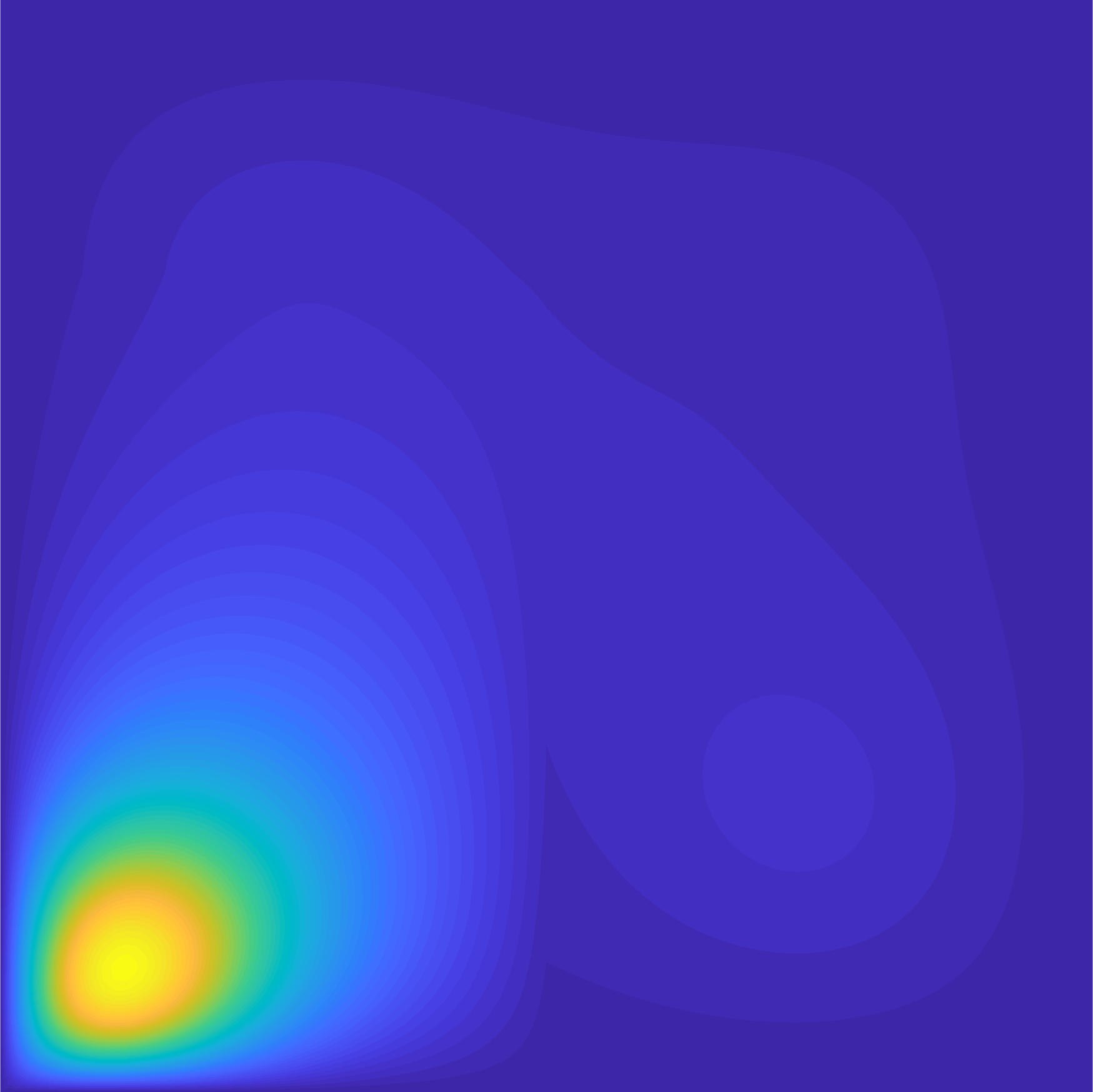};
			
			\nextgroupplot[ ylabel={}, ytick=\empty, xlabel={}, xtick=\empty]
			\addplot graphics [xmin=0, xmax=8.0, ymin=0, ymax=8.0] {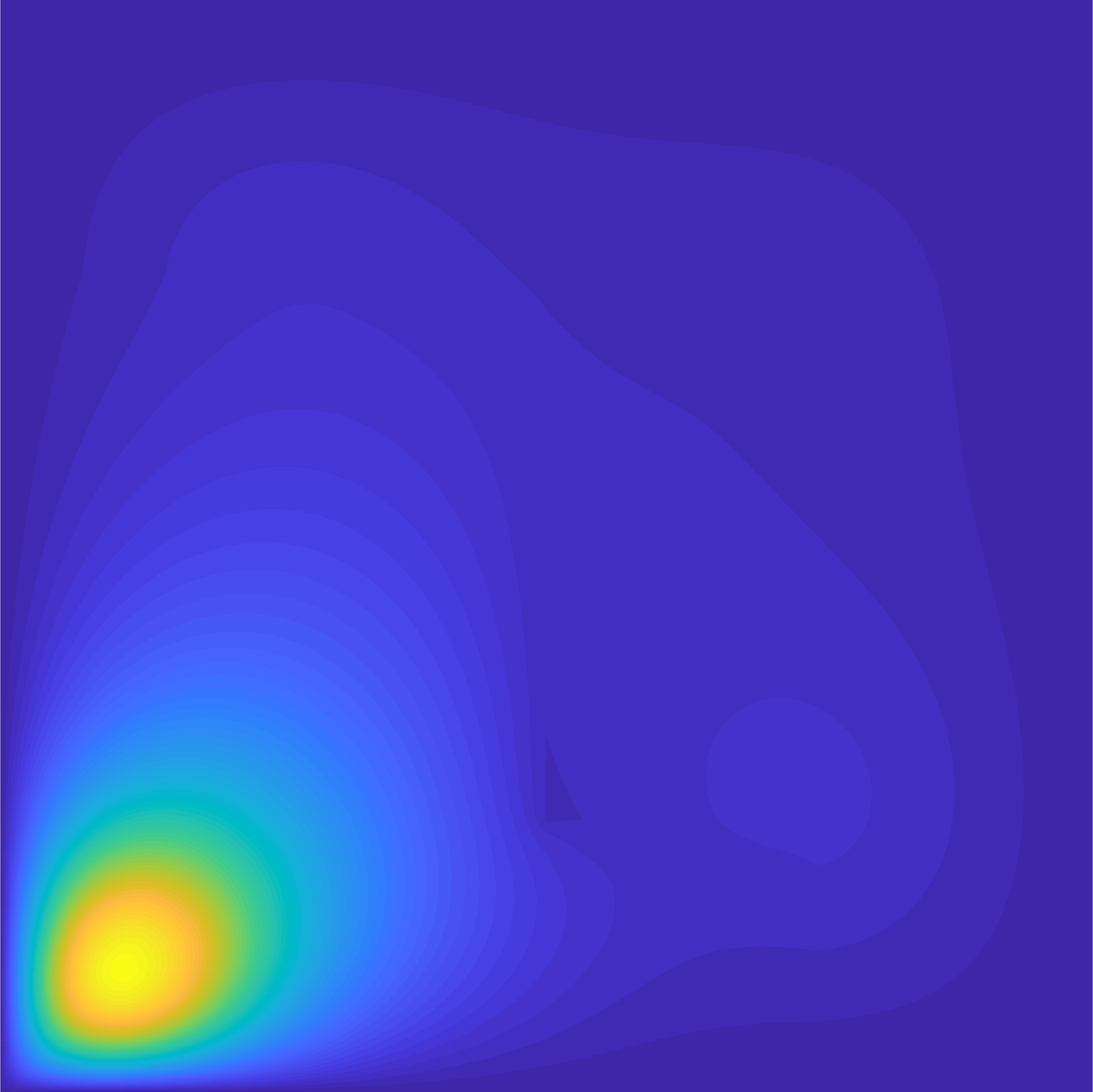};	
			
			\nextgroupplot[ ylabel={}, ytick=\empty, xlabel={}, xtick=\empty]
			\addplot graphics [xmin=0, xmax=8.0, ymin=0, ymax=8.0] {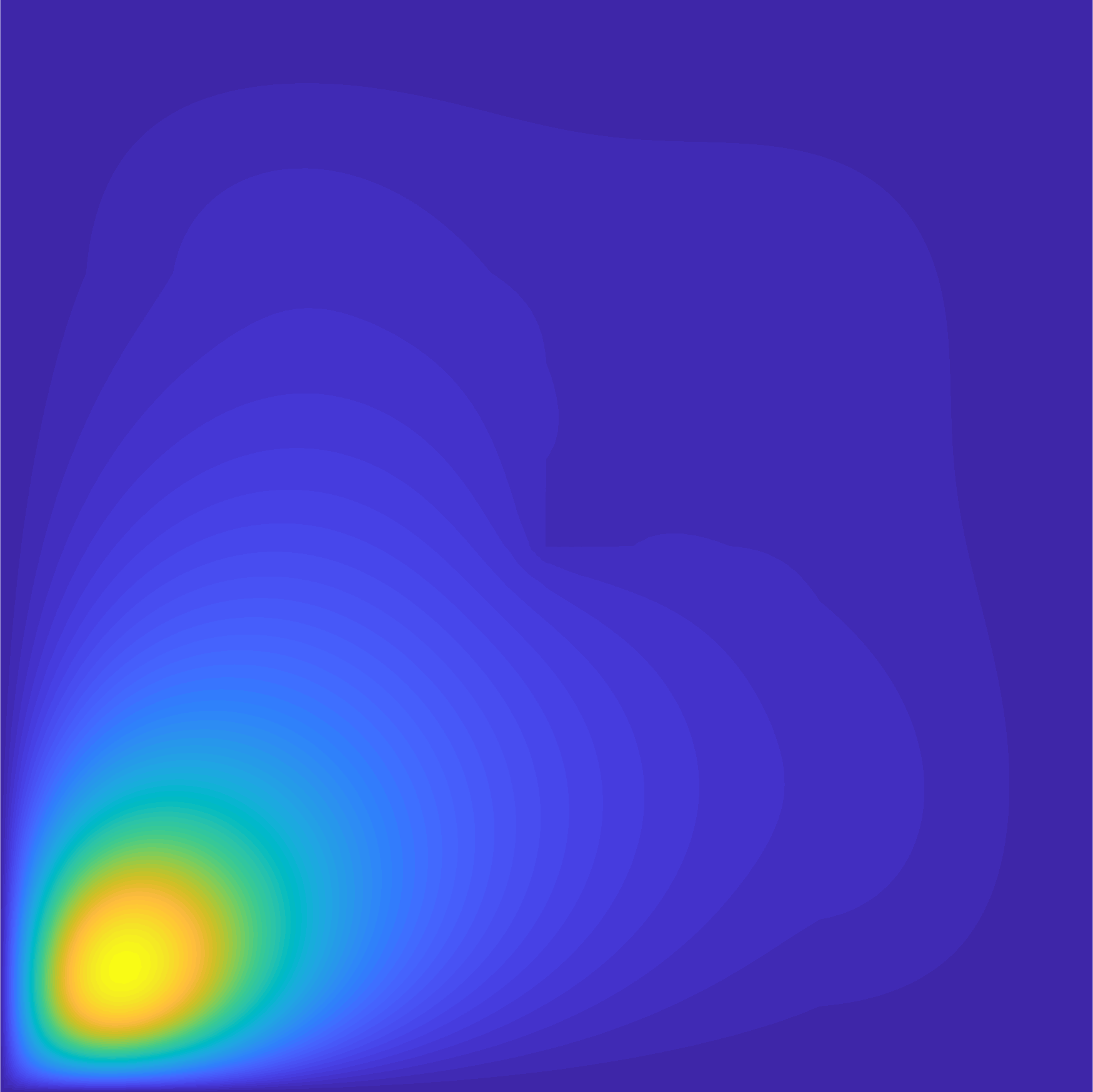};
			
			\nextgroupplot[ ylabel={}, ytick=\empty, xlabel={}, xtick=\empty]
			\addplot graphics [xmin=0, xmax=8.0, ymin=0, ymax=8.0] {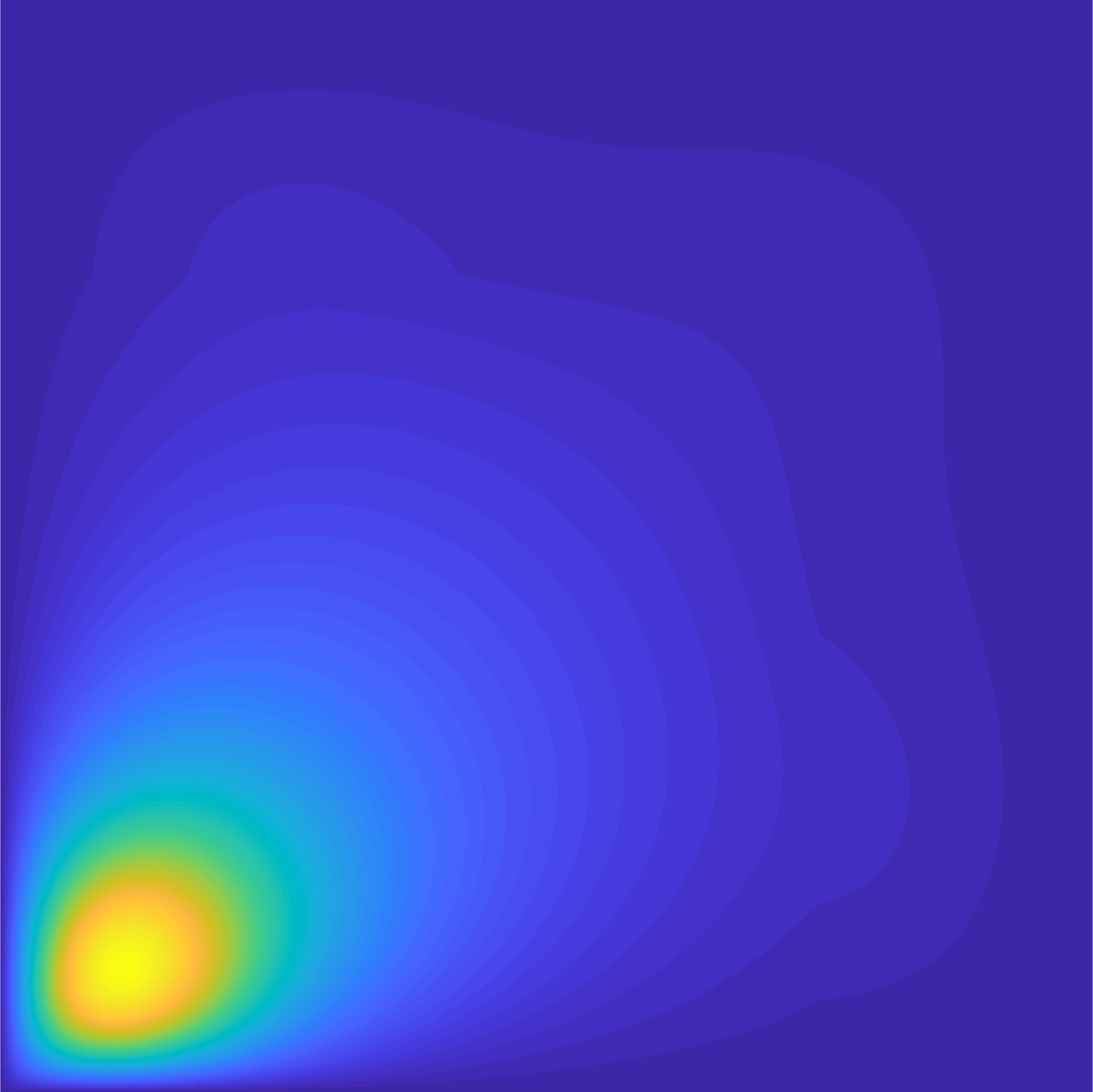};
			
			\nextgroupplot[ ylabel={}, ytick=\empty, xlabel={}, xtick=\empty]
			\addplot graphics [xmin=0, xmax=8.0, ymin=0, ymax=8.0] {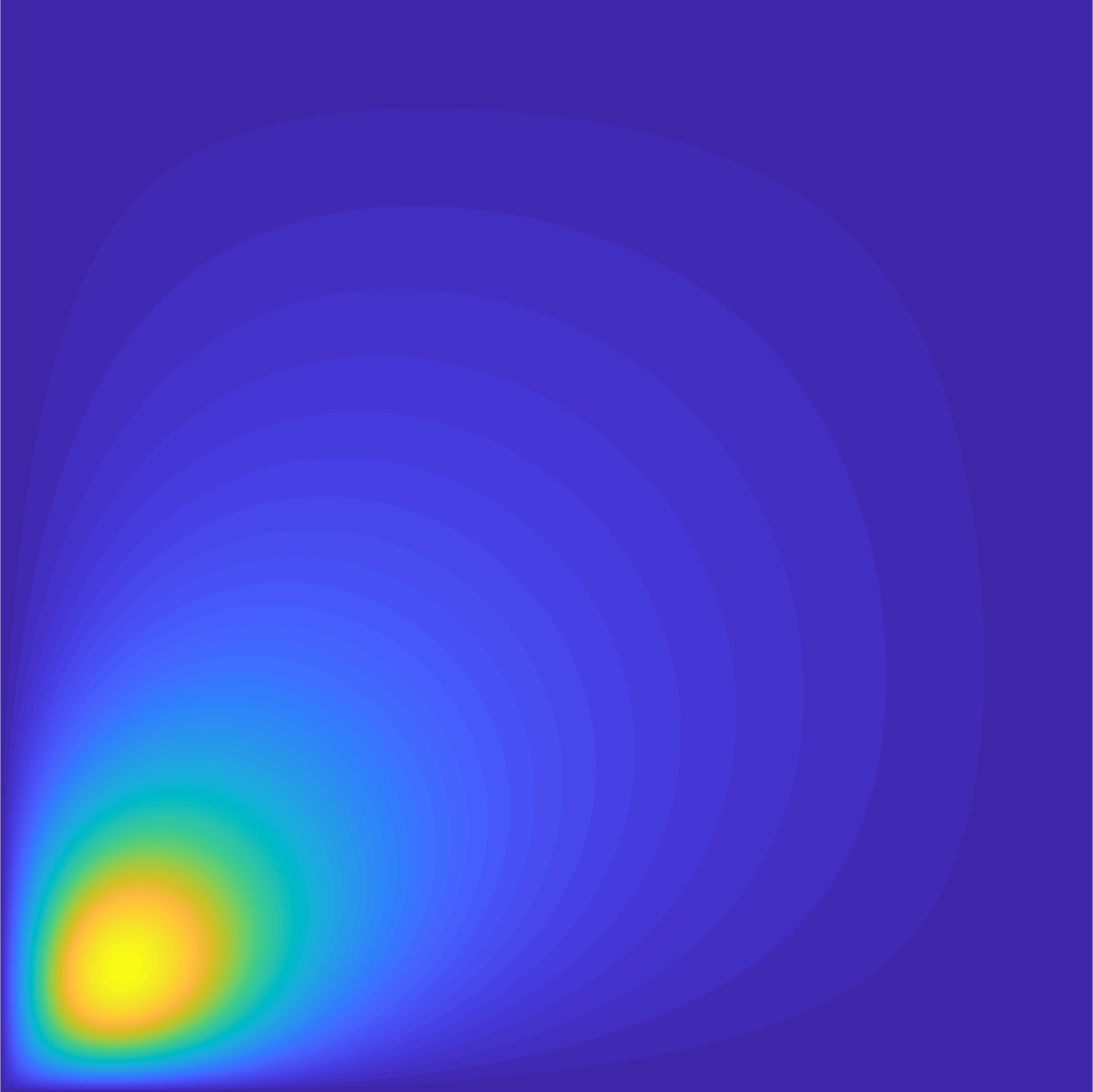};	
		\end{groupplot}
	\end{tikzpicture}
\colorbarMatlabParula{0}{0.0511}{0.1022}{0.1533}{0.2044}
	\caption{Solution of the Poisson problem at $\Dcal_\mathrm{extrap}$ (Section~\ref{sec:rslt:poi:limited}) using the CG-GL method at each adaptation iteration (\textit{left-to-right}, \textit{top-to-bottom}), as well as the reference solution (\textit{bottom right}). Because the CG-GL method is initialized with $\Omega_l = \emptyset$ and $N_l=0$, the \textit{top left} figure is a traditional ROM solution.}
	\label{fig:poi_extrap_example_sol_sequential}
\end{figure}

\begin{figure}
	\centering
	\begin{tikzpicture}
		\begin{groupplot}[
			group style={
				group size=2 by 1,
				horizontal sep=0.1 cm,
			},
			width=0.32\textwidth,
			axis equal image,
			xlabel={$x_1$},
			ylabel={$x_2$},
			xtick = {0.0, 4.0, 8.0},
			ytick = {0.0, 4.0, 8.0},
			xmin=0, xmax=8,
			ymin=0, ymax=8
			]
			\nextgroupplot[ ylabel={}, ytick=\empty, xlabel={}, xtick=\empty]
			\addplot graphics [xmin=0, xmax=8.0, ymin=0, ymax=8.0] {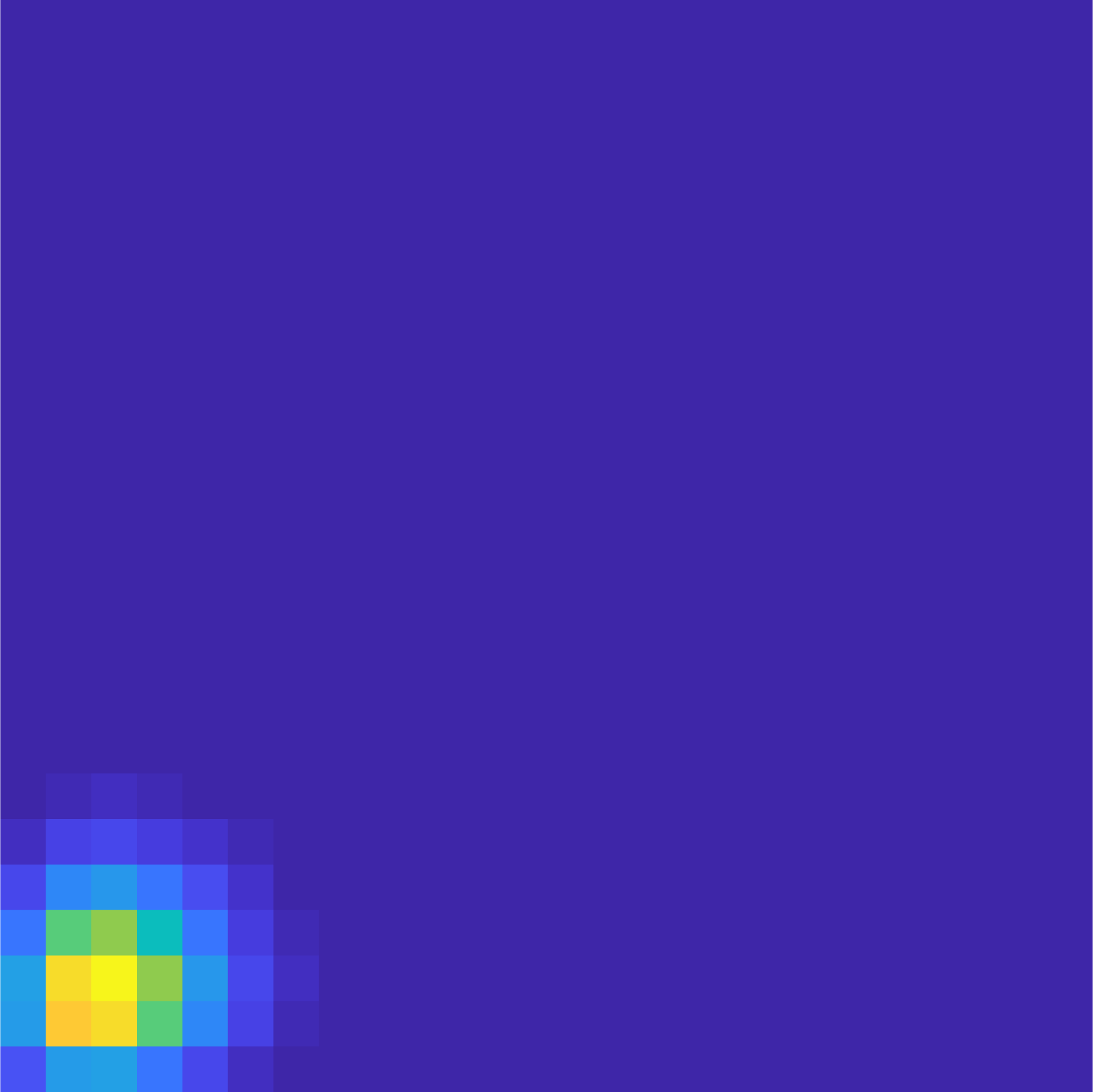};
			
			\nextgroupplot[ ylabel={}, ytick=\empty, xlabel={}, xtick=\empty]
			\addplot graphics [xmin=0, xmax=8.0, ymin=0, ymax=8.0] {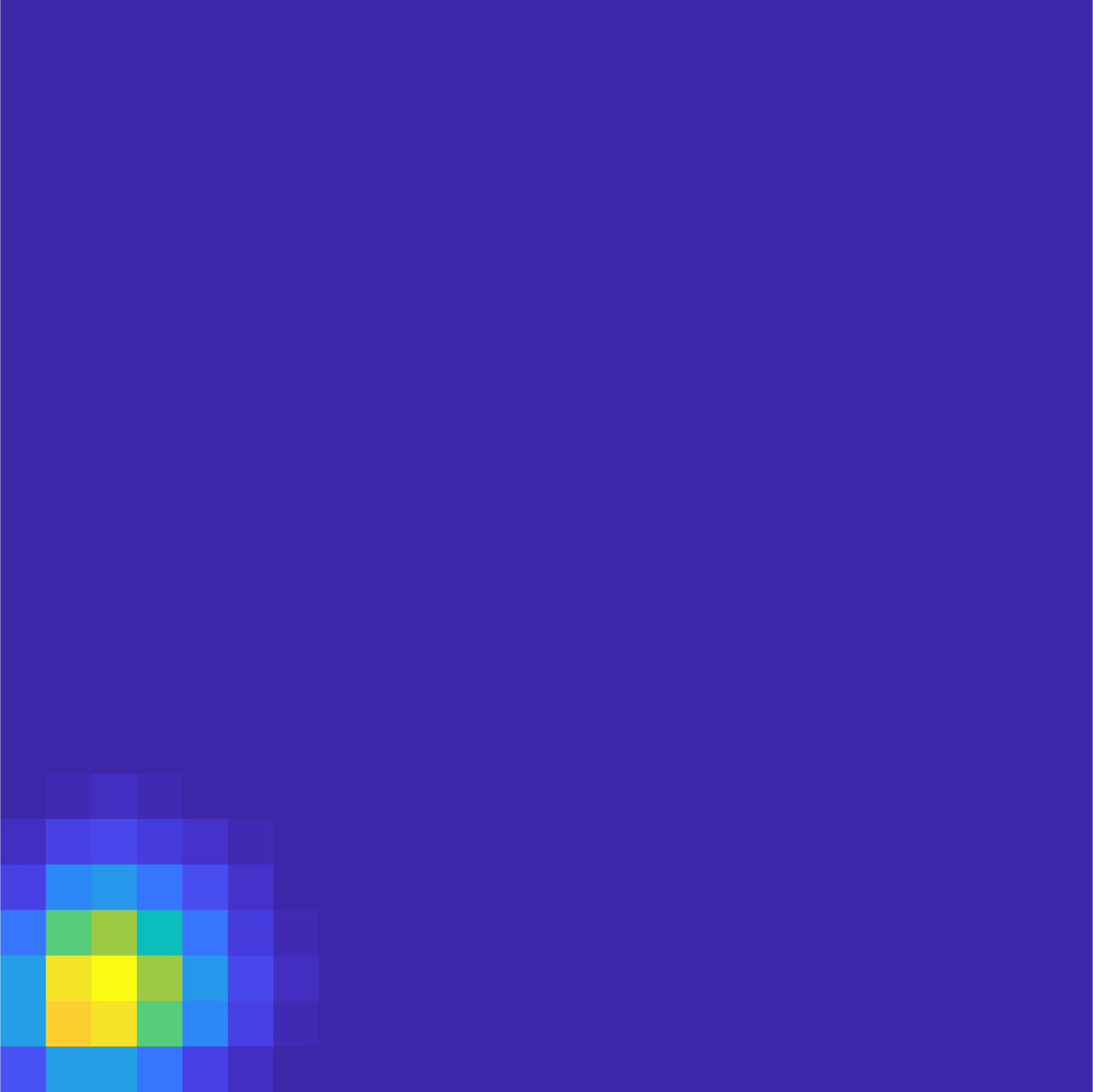};
		\end{groupplot}
	\end{tikzpicture}
	\hspace{4pt}
	\begin{tikzpicture}
		\begin{groupplot}[
			group style={
				group size=2 by 1,
				horizontal sep=0.1 cm,
			},
			width=0.32\textwidth,
			axis equal image,
			xlabel={$x_1$},
			ylabel={$x_2$},
			xtick = {0.0, 4.0, 8.0},
			ytick = {0.0, 4.0, 8.0},
			xmin=0, xmax=8,
			ymin=0, ymax=8
			]
			\nextgroupplot[ ylabel={}, ytick=\empty, xlabel={}, xtick=\empty]
			\addplot graphics [xmin=0, xmax=8.0, ymin=0, ymax=8.0] {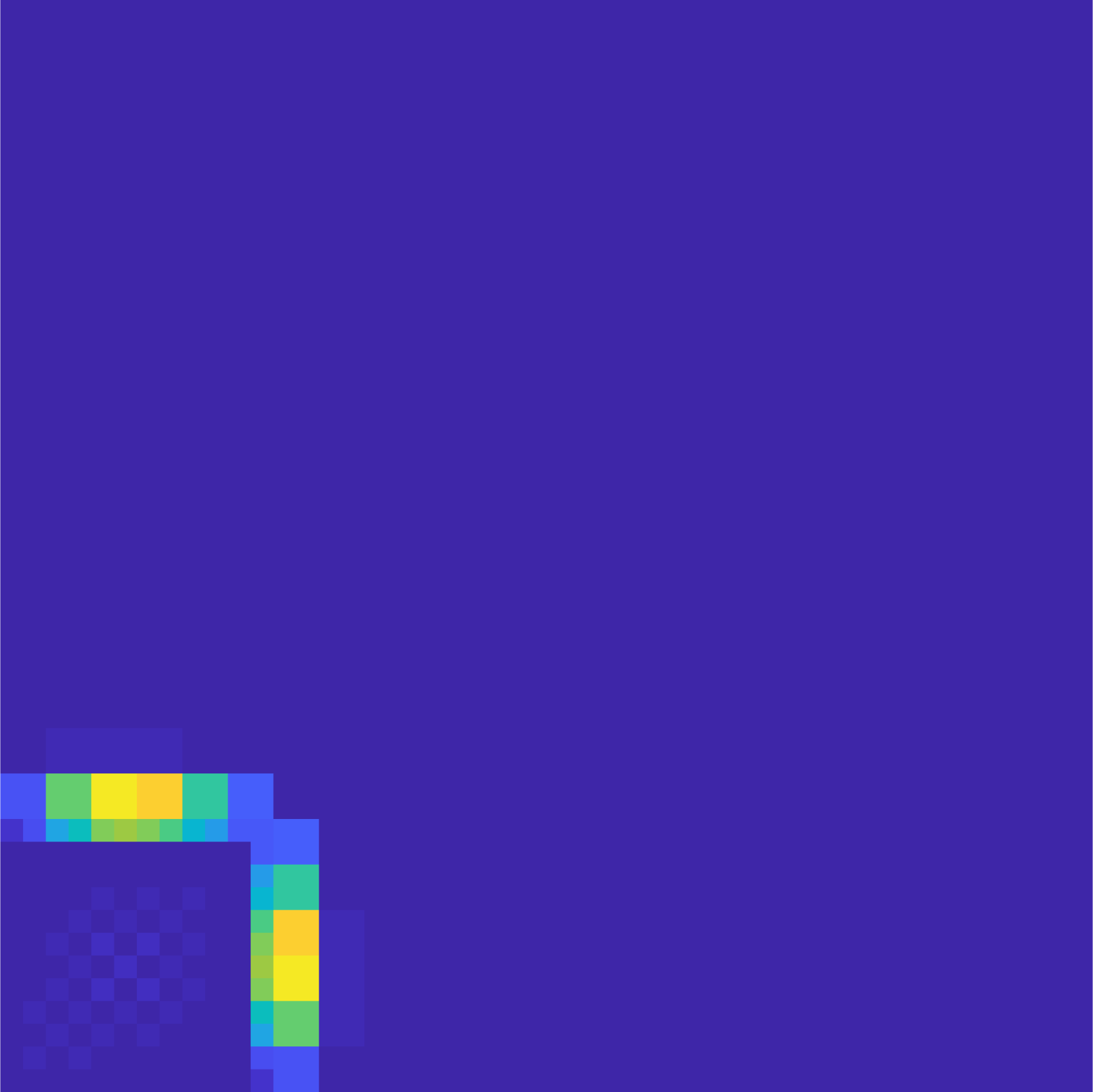};
			
			\nextgroupplot[ ylabel={}, ytick=\empty, xlabel={}, xtick=\empty]
			\addplot graphics [xmin=0, xmax=8.0, ymin=0, ymax=8.0] {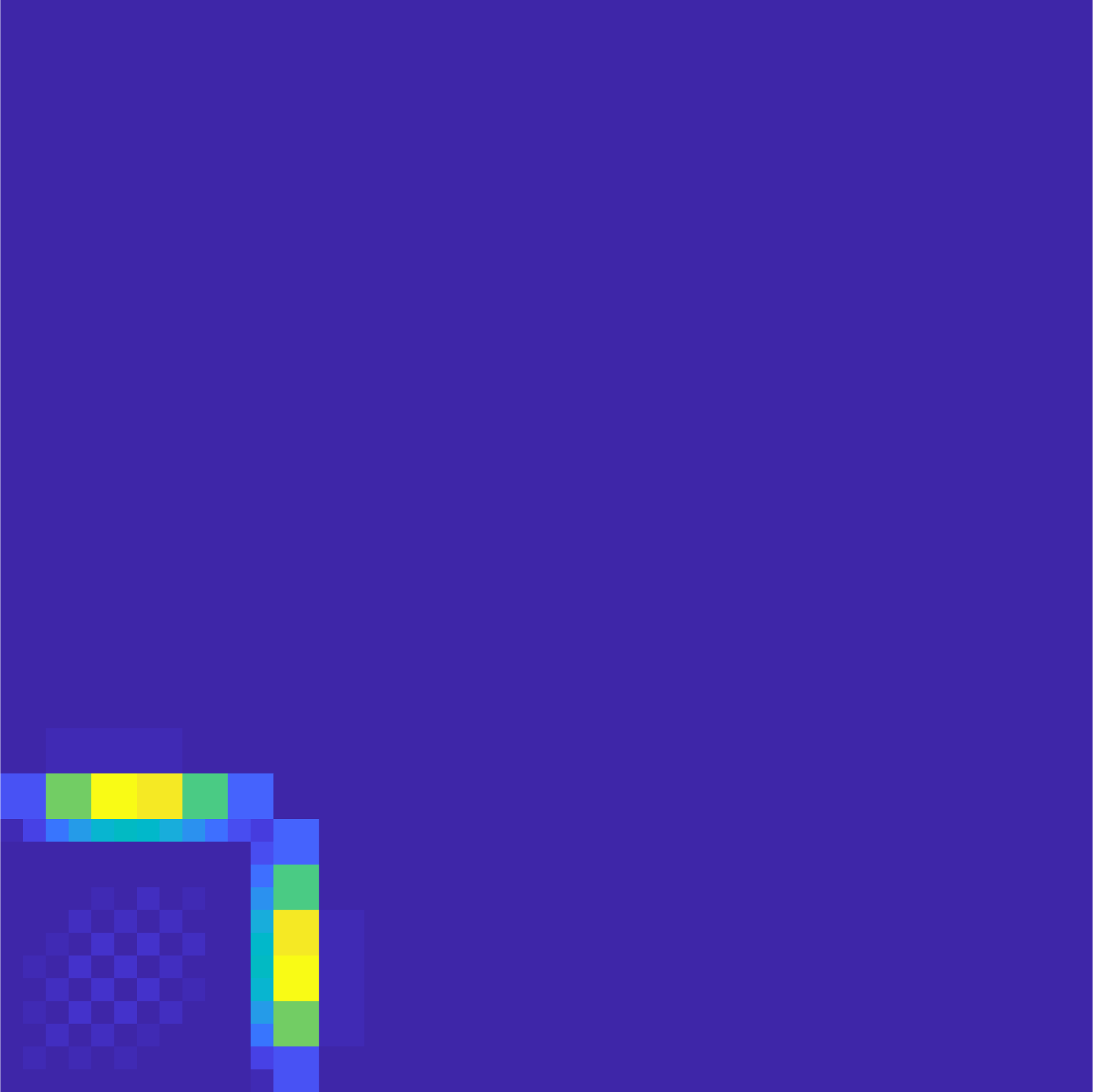};
		\end{groupplot}
	\end{tikzpicture}\\
	\vspace{8pt}
	\begin{tikzpicture}
		\begin{groupplot}[
			group style={
				group size=2 by 1,
				horizontal sep=0.1 cm,
			},
			width=0.32\textwidth,
			axis equal image,
			xlabel={$x_1$},
			ylabel={$x_2$},
			xtick = {0.0, 4.0, 8.0},
			ytick = {0.0, 4.0, 8.0},
			xmin=0, xmax=8,
			ymin=0, ymax=8
			]
			\nextgroupplot[ ylabel={}, ytick=\empty, xlabel={}, xtick=\empty]
			\addplot graphics [xmin=0, xmax=8.0, ymin=0, ymax=8.0] {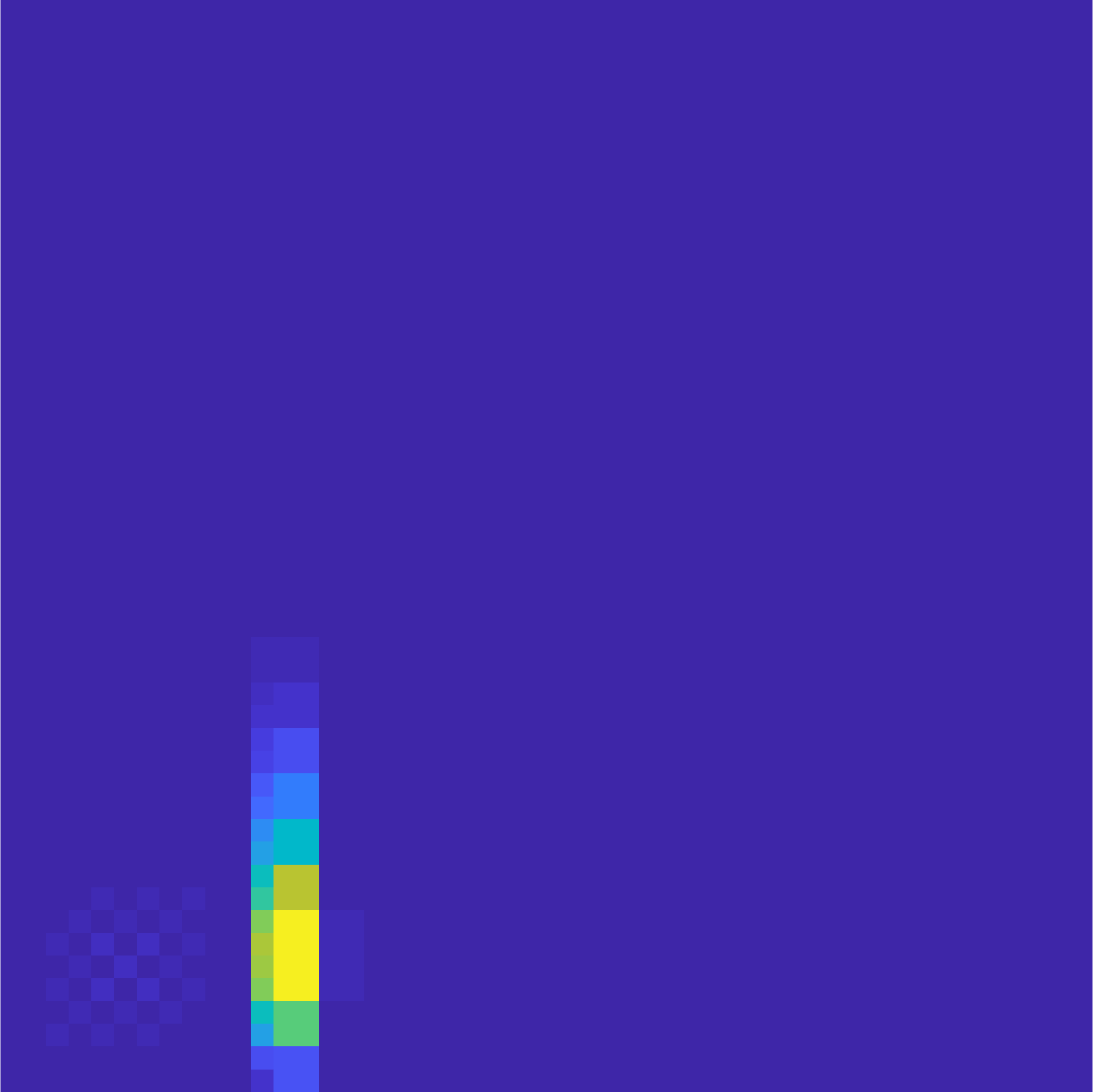};
			
			\nextgroupplot[ ylabel={}, ytick=\empty, xlabel={}, xtick=\empty]
			\addplot graphics [xmin=0, xmax=8.0, ymin=0, ymax=8.0] {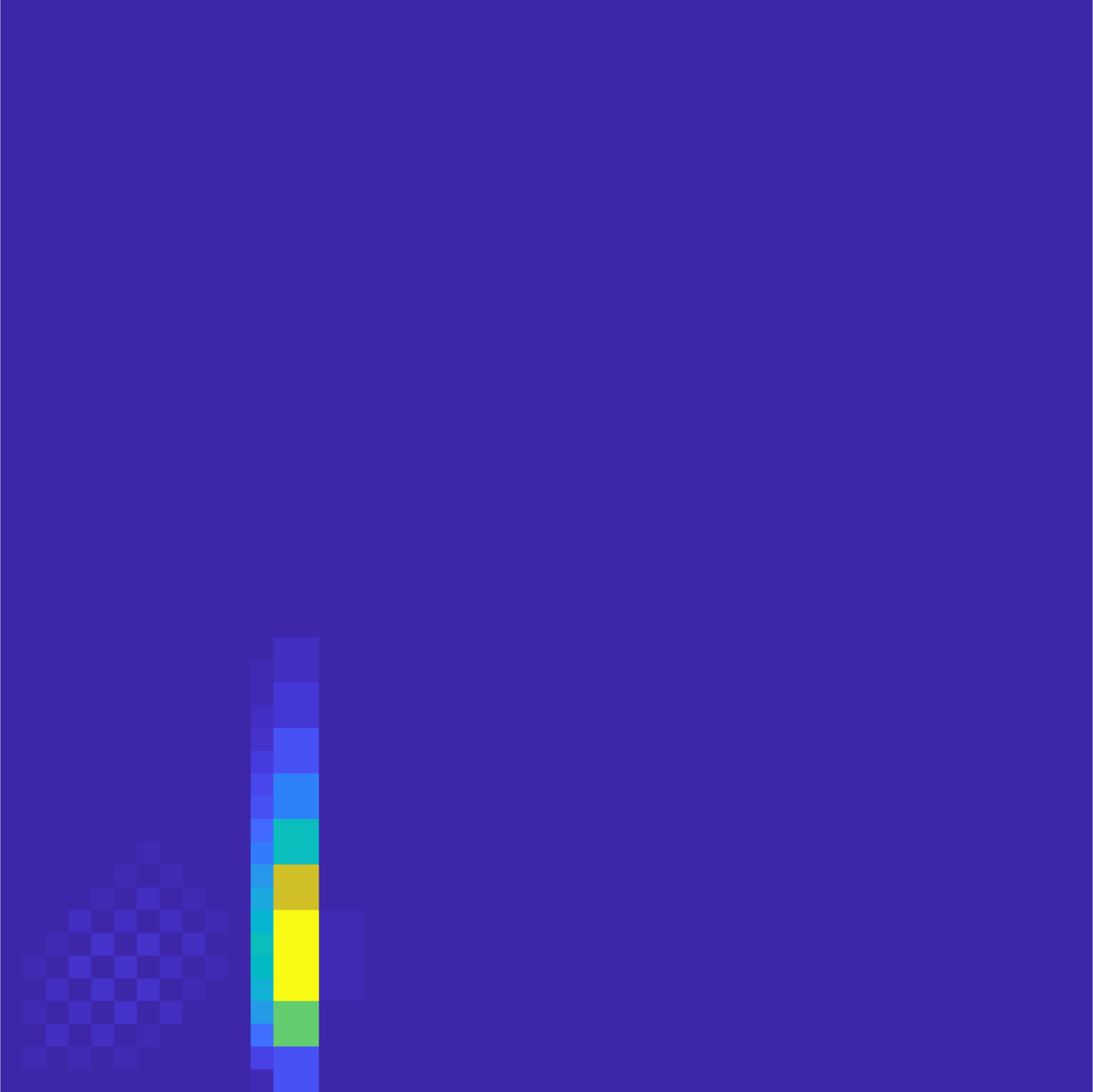};
		\end{groupplot}
	\end{tikzpicture}
	\hspace{4pt}
	\begin{tikzpicture}
		\begin{groupplot}[
			group style={
				group size=2 by 1,
				horizontal sep=0.1 cm,
			},
			width=0.32\textwidth,
			axis equal image,
			xlabel={$x_1$},
			ylabel={$x_2$},
			xtick = {0.0, 4.0, 8.0},
			ytick = {0.0, 4.0, 8.0},
			xmin=0, xmax=8,
			ymin=0, ymax=8
			]
			\nextgroupplot[ ylabel={}, ytick=\empty, xlabel={}, xtick=\empty]
			\addplot graphics [xmin=0, xmax=8.0, ymin=0, ymax=8.0] {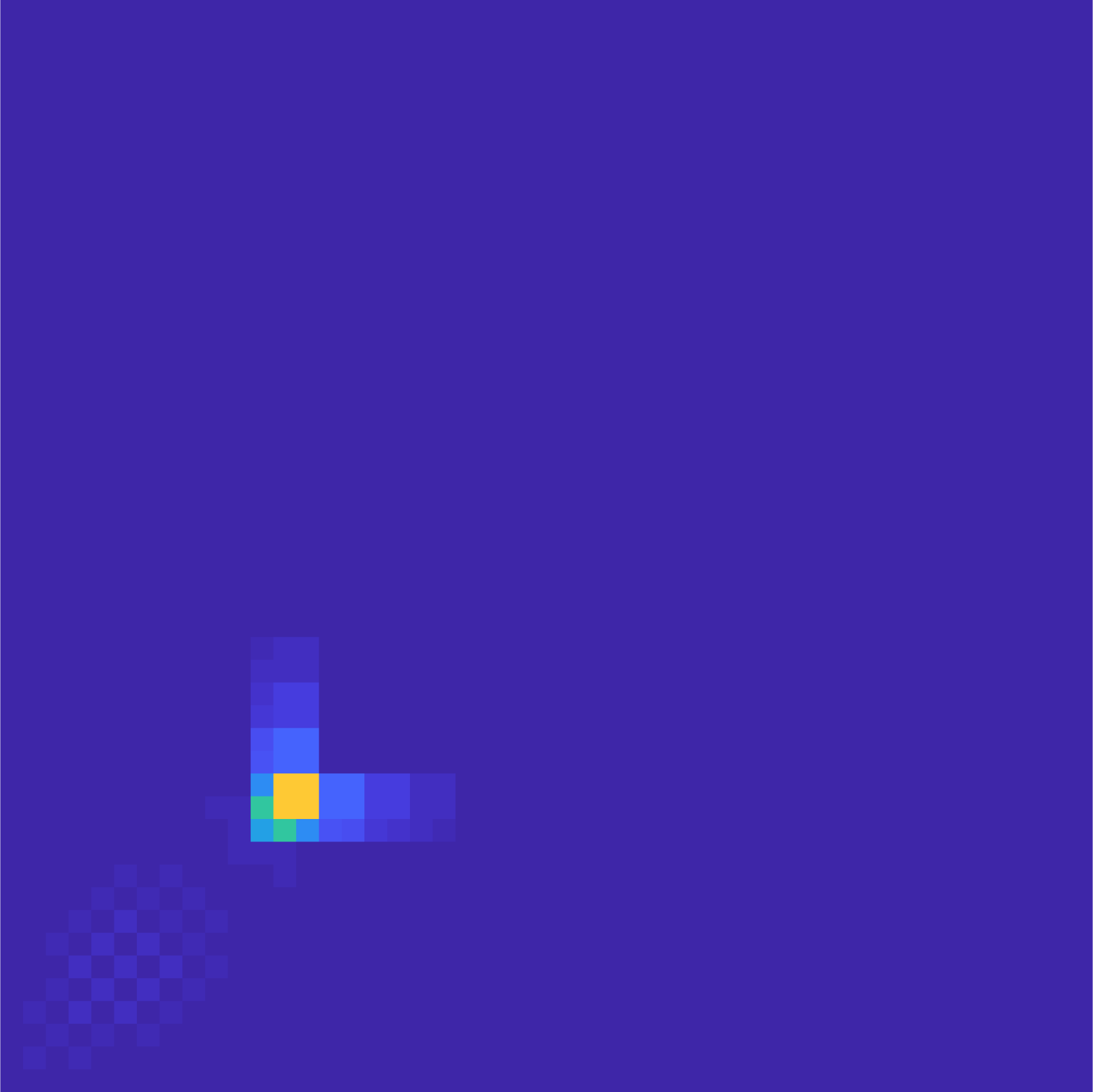};
			
			\nextgroupplot[ ylabel={}, ytick=\empty, xlabel={}, xtick=\empty]
			\addplot graphics [xmin=0, xmax=8.0, ymin=0, ymax=8.0] {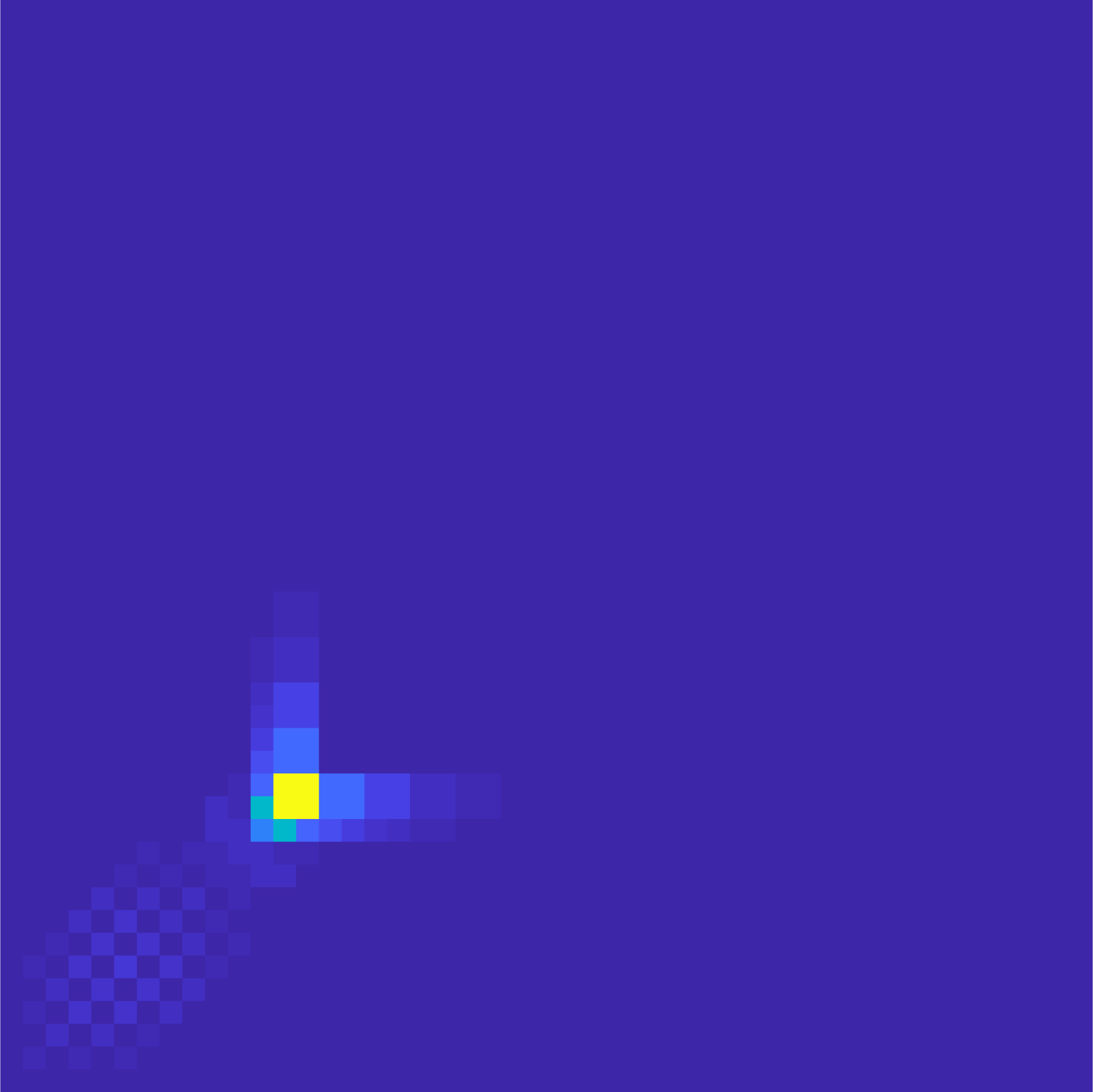};
		\end{groupplot}
	\end{tikzpicture}
	
	\vspace{8pt}
	\begin{tikzpicture}
		\begin{groupplot}[
			group style={
				group size=2 by 1,
				horizontal sep=0.1 cm,
			},
			width=0.32\textwidth,
			axis equal image,
			xlabel={$x_1$},
			ylabel={$x_2$},
			xtick = {0.0, 4.0, 8.0},
			ytick = {0.0, 4.0, 8.0},
			xmin=0, xmax=8,
			ymin=0, ymax=8
			]
			\nextgroupplot[ ylabel={}, ytick=\empty, xlabel={}, xtick=\empty]
			\addplot graphics [xmin=0, xmax=8.0, ymin=0, ymax=8.0] {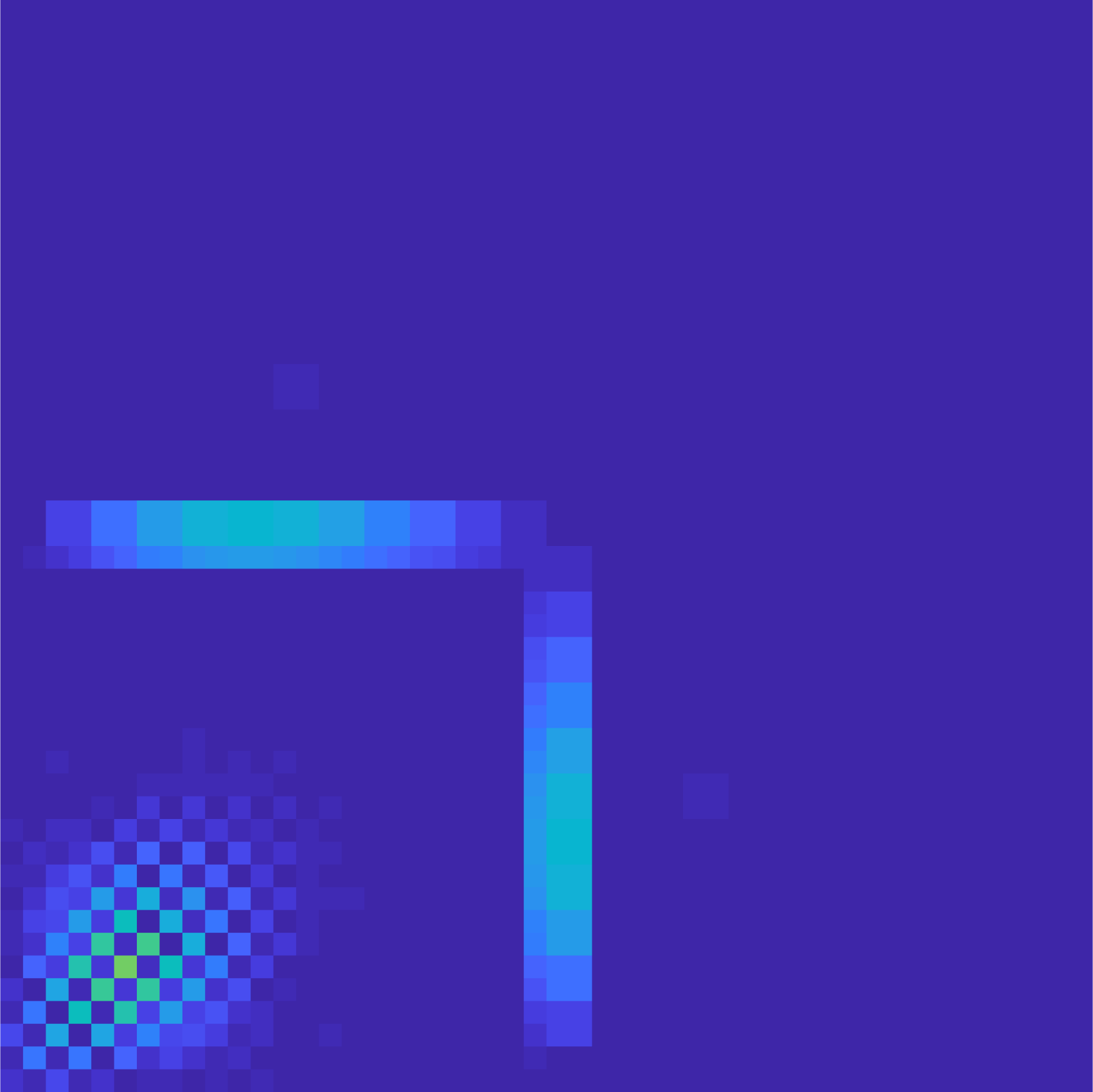};
			
			\nextgroupplot[ ylabel={}, ytick=\empty, xlabel={}, xtick=\empty]
			\addplot graphics [xmin=0, xmax=8.0, ymin=0, ymax=8.0] {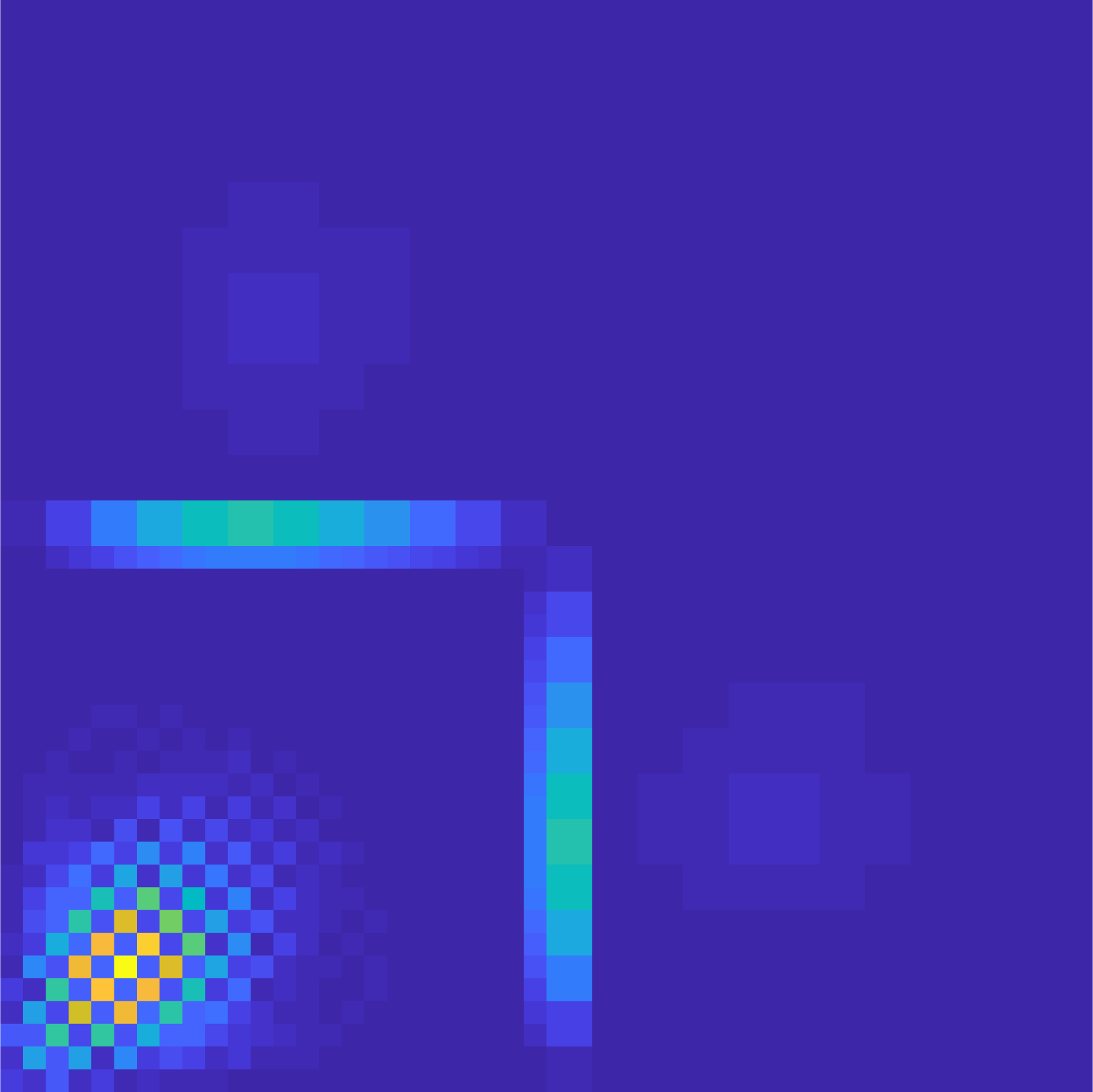};
		\end{groupplot}
	\end{tikzpicture}
	\hspace{4pt}
	\begin{tikzpicture}
		\begin{groupplot}[
			group style={
				group size=2 by 1,
				horizontal sep=0.1 cm,
			},
			width=0.32\textwidth,
			axis equal image,
			xlabel={$x_1$},
			ylabel={$x_2$},
			xtick = {0.0, 4.0, 8.0},
			ytick = {0.0, 4.0, 8.0},
			xmin=0, xmax=8,
			ymin=0, ymax=8
			]
			\nextgroupplot[ ylabel={}, ytick=\empty, xlabel={}, xtick=\empty]
			\addplot graphics [xmin=0, xmax=8.0, ymin=0, ymax=8.0] {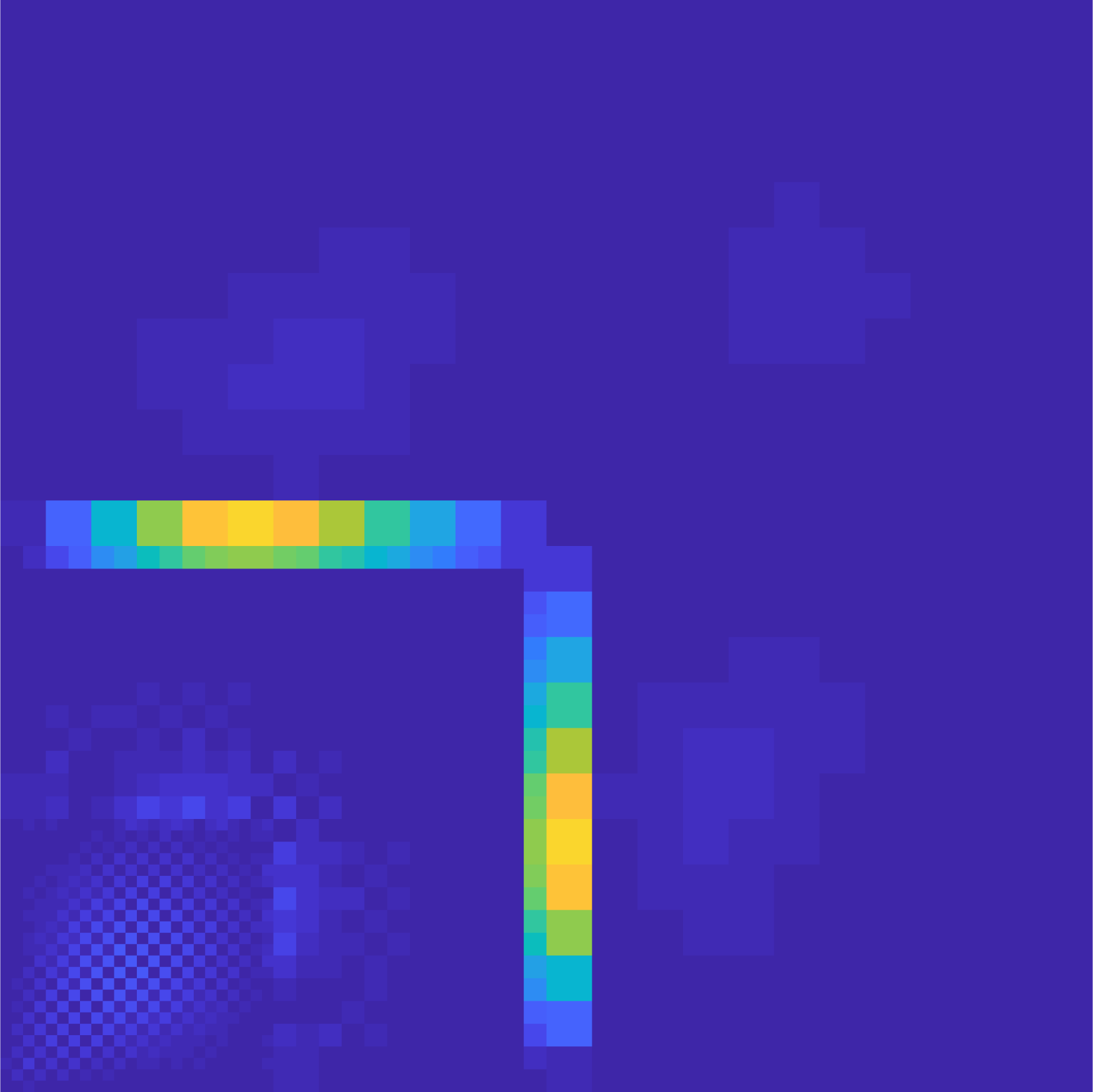};
			
			\nextgroupplot[ ylabel={}, ytick=\empty, xlabel={}, xtick=\empty]
			\addplot graphics [xmin=0, xmax=8.0, ymin=0, ymax=8.0] {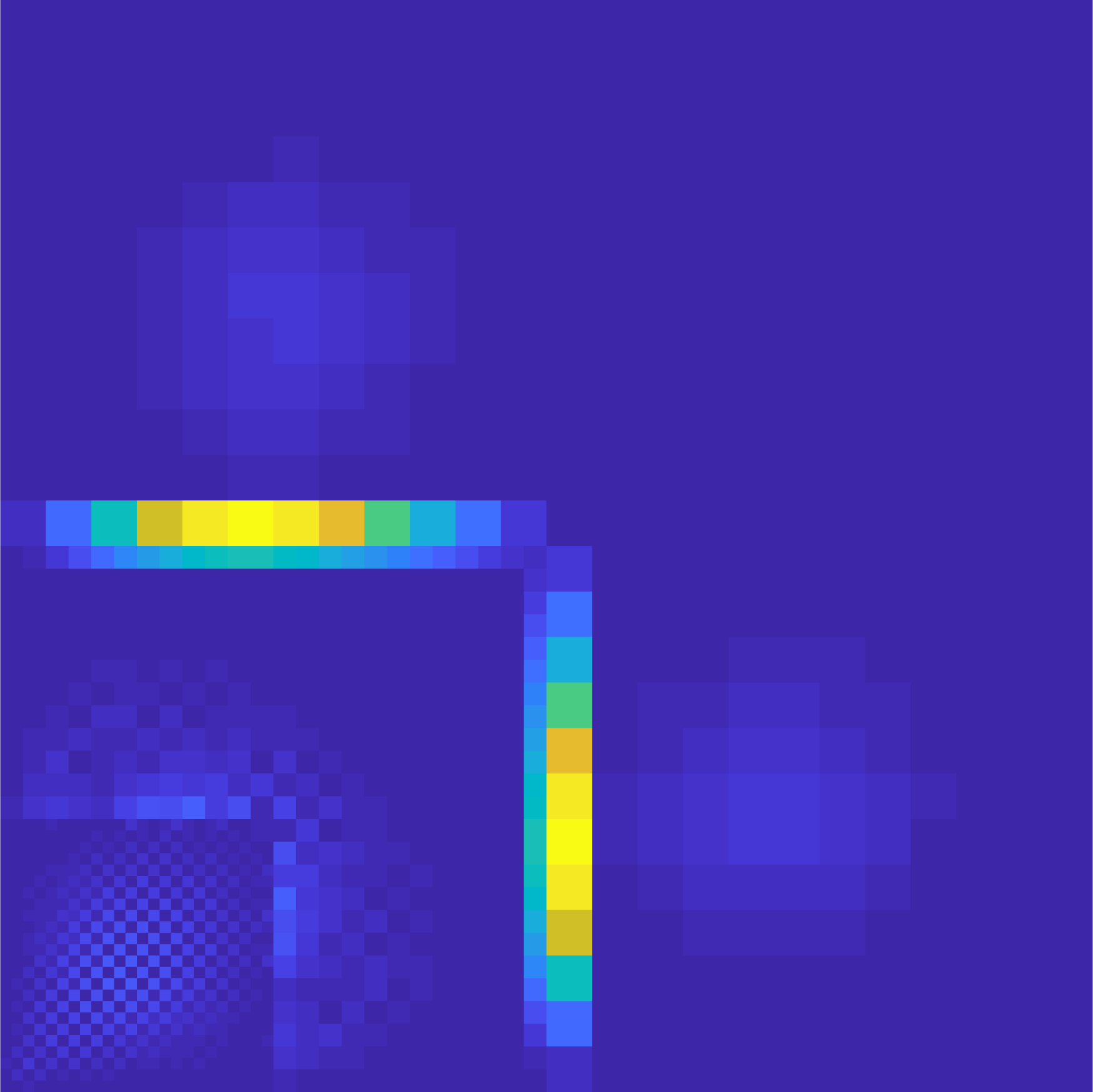};
		\end{groupplot}
	\end{tikzpicture}
	
	\vspace{8pt}
	\begin{tikzpicture}
		\begin{groupplot}[
			group style={
				group size=2 by 1,
				horizontal sep=0.1 cm,
			},
			width=0.32\textwidth,
			axis equal image,
			xlabel={$x_1$},
			ylabel={$x_2$},
			xtick = {0.0, 4.0, 8.0},
			ytick = {0.0, 4.0, 8.0},
			xmin=0, xmax=8,
			ymin=0, ymax=8
			]
			\nextgroupplot[ ylabel={}, ytick=\empty, xlabel={}, xtick=\empty]
			\addplot graphics [xmin=0, xmax=8.0, ymin=0, ymax=8.0] {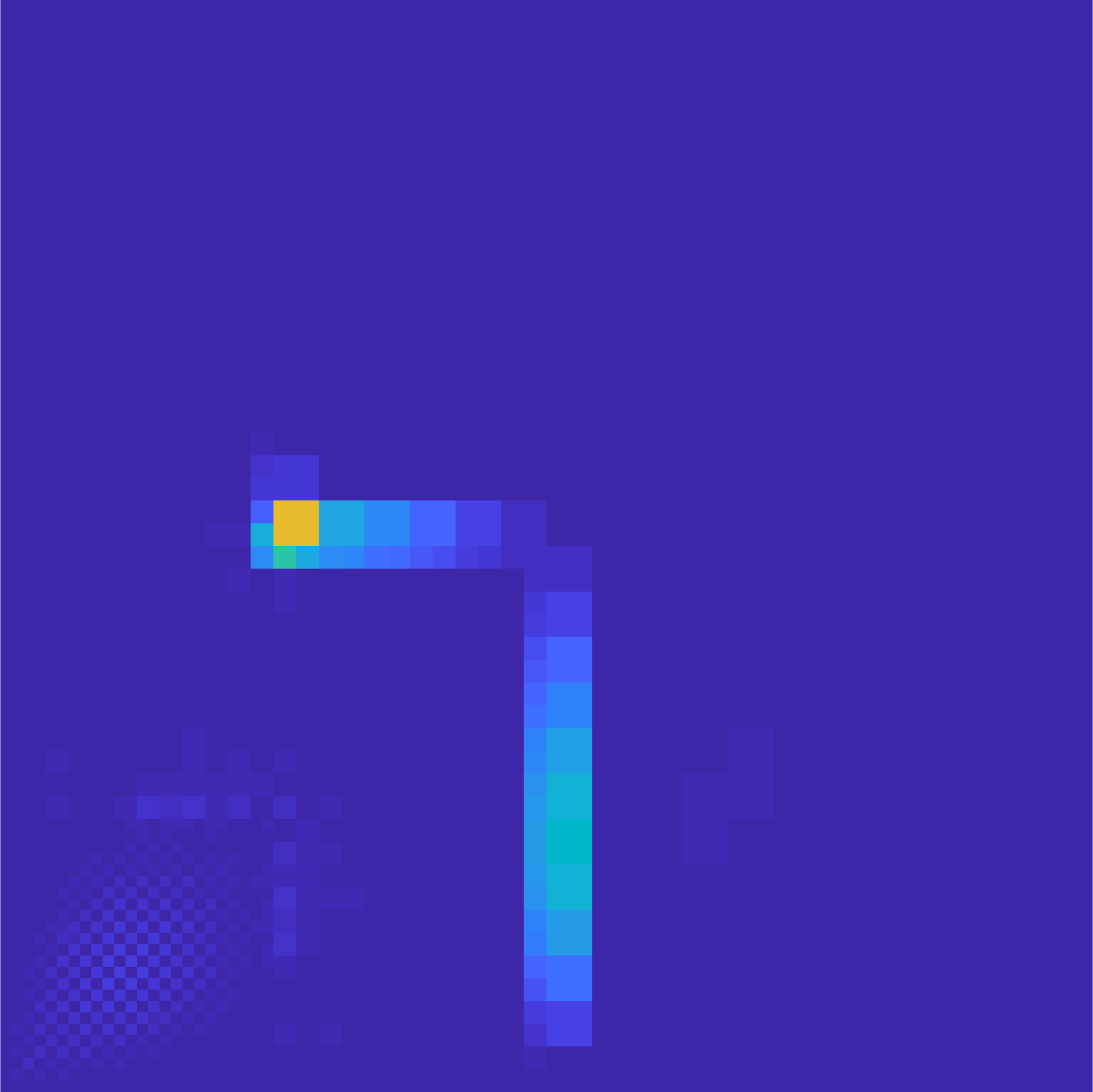};
			
			\nextgroupplot[ ylabel={}, ytick=\empty, xlabel={}, xtick=\empty]
			\addplot graphics [xmin=0, xmax=8.0, ymin=0, ymax=8.0] {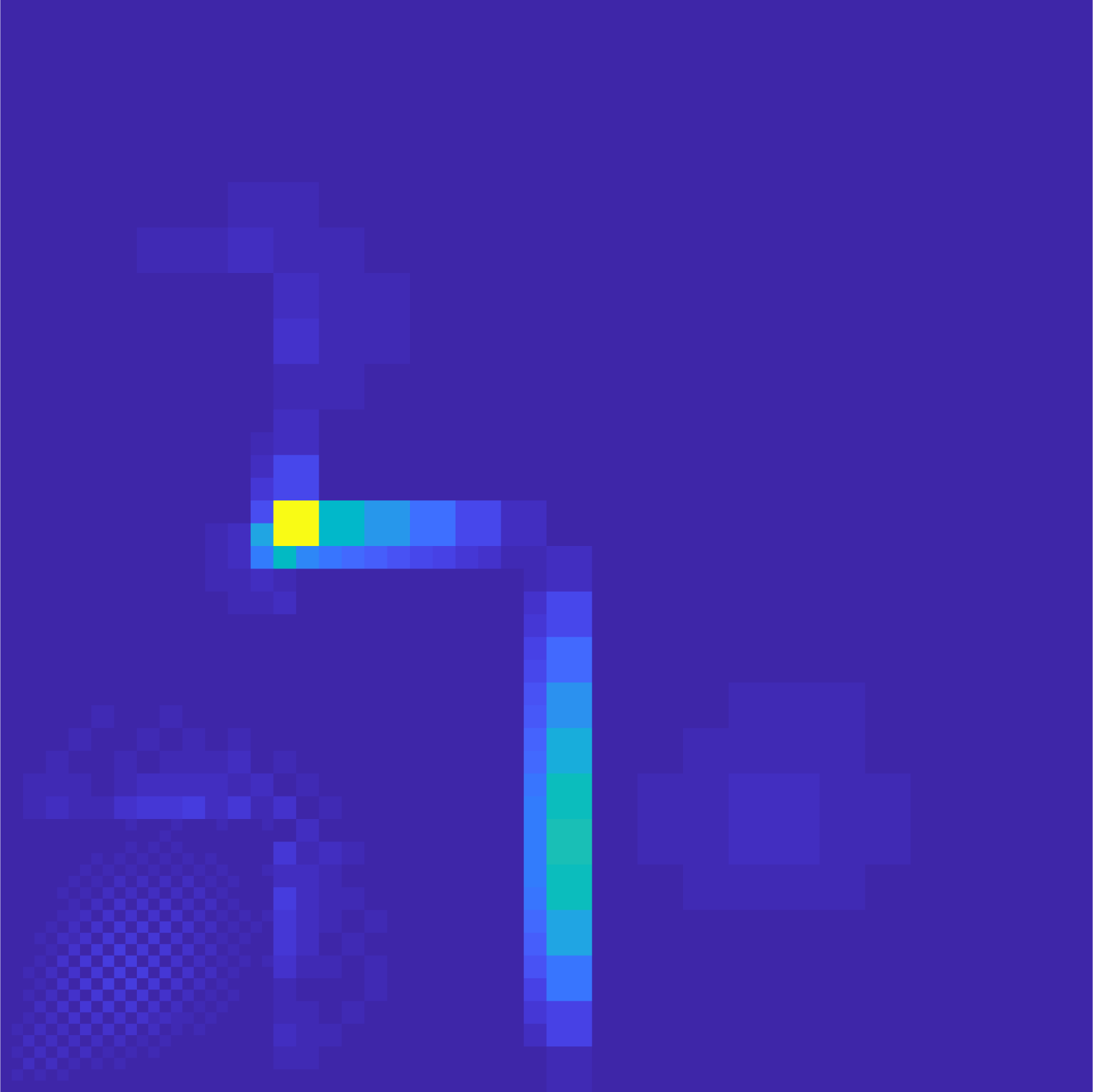};
		\end{groupplot}
	\end{tikzpicture}
	\hspace{4pt}
	\begin{tikzpicture}
		\begin{groupplot}[
			group style={
				group size=2 by 1,
				horizontal sep=0.1 cm,
			},
			width=0.32\textwidth,
			axis equal image,
			xlabel={$x_1$},
			ylabel={$x_2$},
			xtick = {0.0, 4.0, 8.0},
			ytick = {0.0, 4.0, 8.0},
			xmin=0, xmax=8,
			ymin=0, ymax=8
			]
			\nextgroupplot[ ylabel={}, ytick=\empty, xlabel={}, xtick=\empty]
			\addplot graphics [xmin=0, xmax=8.0, ymin=0, ymax=8.0] {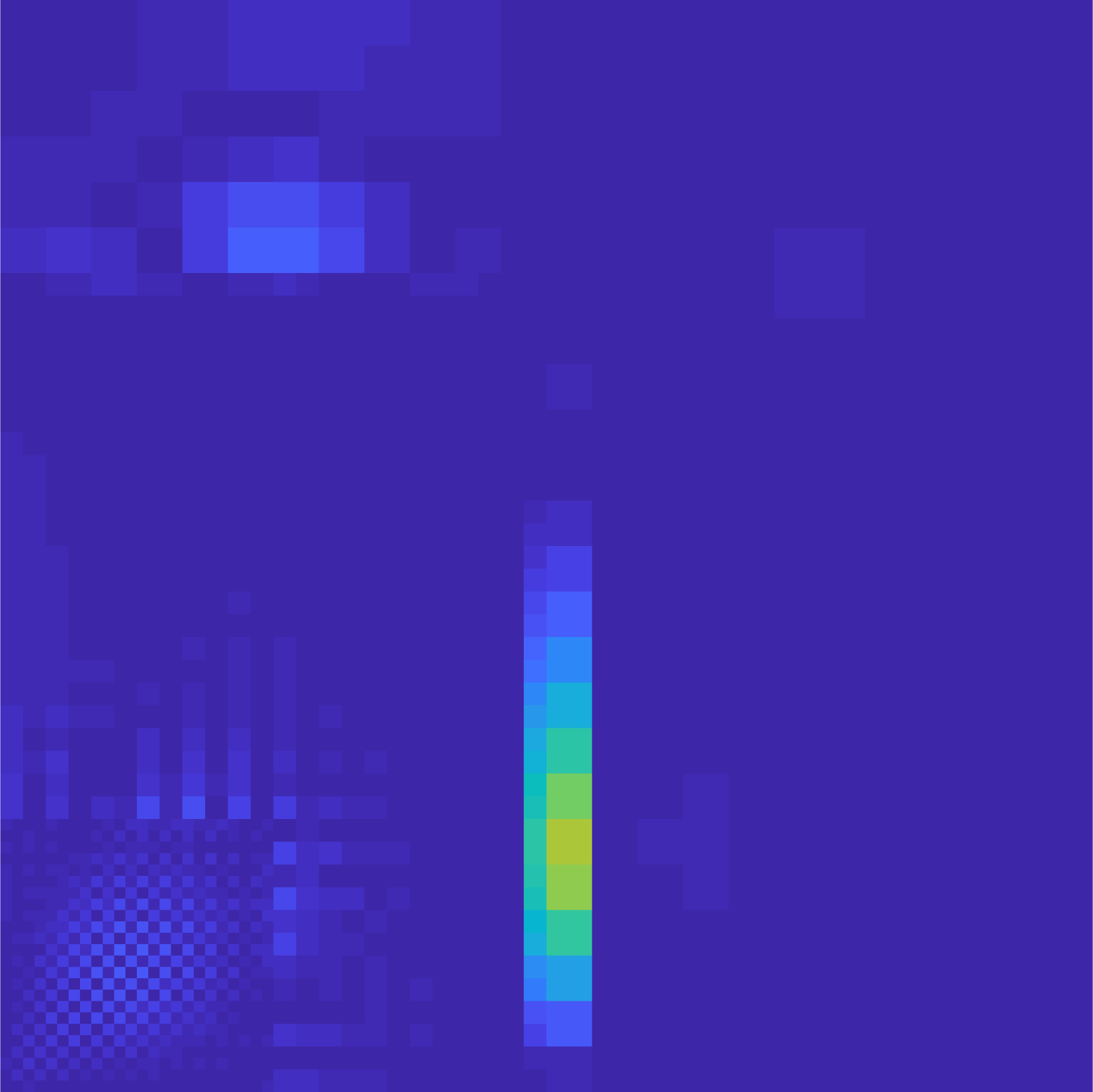};
			
			\nextgroupplot[ ylabel={}, ytick=\empty, xlabel={}, xtick=\empty]
			\addplot graphics [xmin=0, xmax=8.0, ymin=0, ymax=8.0] {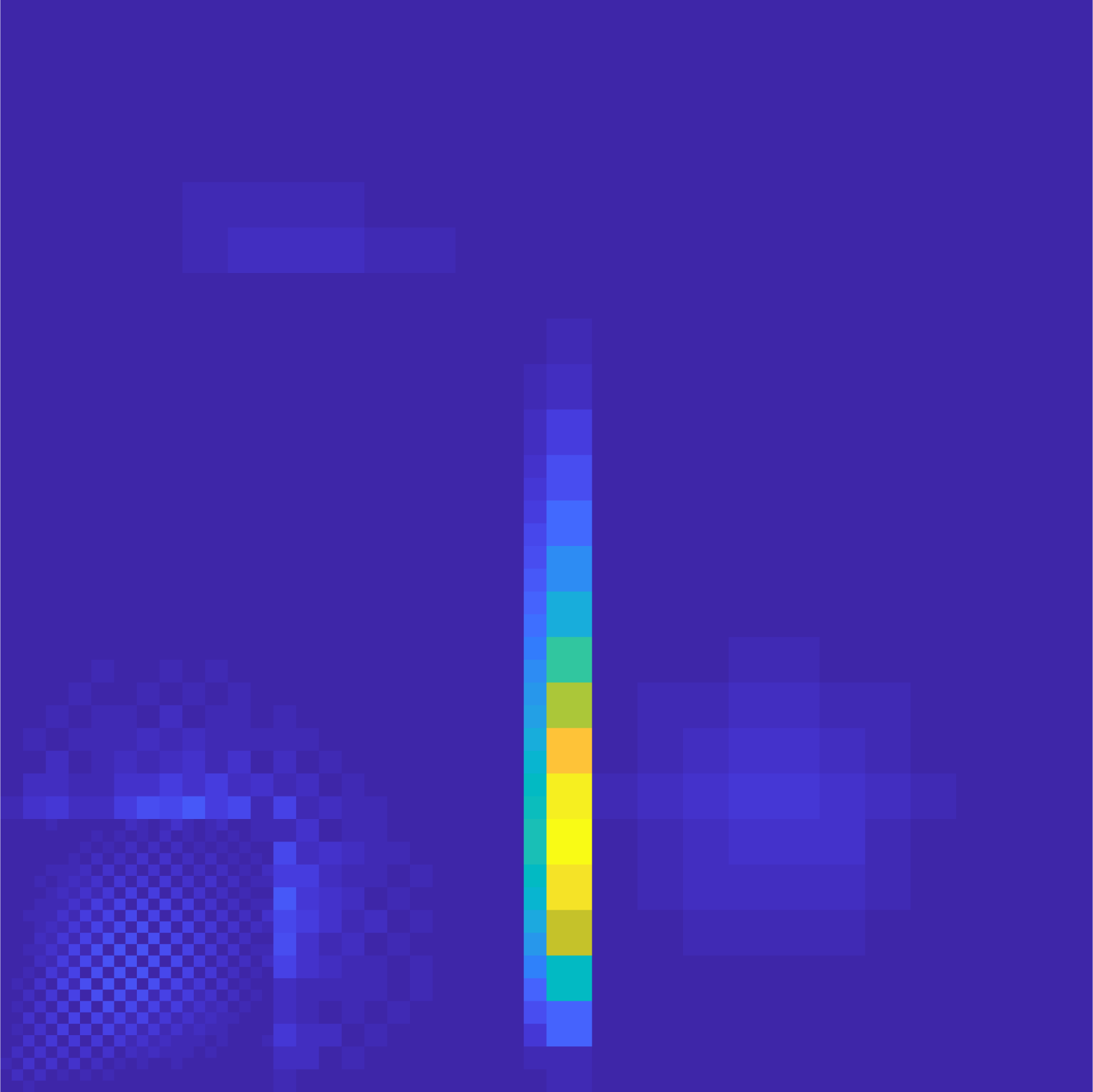};
		\end{groupplot}
	\end{tikzpicture}
	\caption{The localized DWR error estimate (\textit{left}) and true error (\textit{right}) at the first eight CG-GL adaptation
	iterations (\textit{left-to-right, top-to-bottom}) for the Poisson problem at $\Dcal_\mathrm{extrap}$ (Section~\ref{sec:rslt:poi:limited}).
	Each pair of figures corresponds to a different CG-GL iteration and use its own colorbar, scaled to its range to highlight similarity
	between the proposed error estimate and actual error.}
	\label{fig:error_contour_1}
\end{figure}

\begin{figure}
	\begin{tikzpicture}
\begin{groupplot} [
group style={group size = 3 by 1, horizontal sep = 1.5cm}]
\nextgroupplot[width=0.33\textwidth, xlabel=$r$, ymax=1.0485558278707041, xmax=11, ylabel=$e_\mathrm{sln}$, xmin=0, ymin=0.017113529346093007,ymode=log]
\addplot [color=blue, mark=o,thick]
coordinates {
( 0.00000000e+00,  4.39415721e-01)
( 1.00000000e+00,  4.51276221e-01)
( 2.00000000e+00,  4.11587530e-01)
( 3.00000000e+00,  4.04595910e-01)
( 4.00000000e+00,  3.89326978e-01)
( 5.00000000e+00,  2.14141709e-01)
( 6.00000000e+00,  1.92539942e-01)
( 7.00000000e+00,  1.07533271e-01)
( 8.00000000e+00,  9.90337694e-02)
( 9.00000000e+00,  1.88770240e-02)
( 1.00000000e+01,  1.87107072e-02)
( 1.10000000e+01,  1.86974577e-02)};\label{poi_limitedtraing_int}

\addplot [color=blue, mark=x,thick]
coordinates {
( 0.00000000e+00,  9.98624598e-01)
( 1.00000000e+00,  6.25744308e-01)
( 2.00000000e+00,  5.25279294e-01)
( 3.00000000e+00,  3.77154169e-01)
( 4.00000000e+00,  1.78719185e-01)
( 5.00000000e+00,  1.78720256e-01)
( 6.00000000e+00,  1.66722286e-01)
( 7.00000000e+00,  1.34024999e-01)
( 8.00000000e+00,  1.18294411e-01)
( 9.00000000e+00,  5.93256849e-02)
( 1.00000000e+01,  1.98511242e-02)
( 1.10000000e+01,  1.80142414e-02)};\label{poi_limitedtraing_ext}

\nextgroupplot[width=0.33\textwidth, xlabel=$r$, ymax=1.0498745773167717, xmax=11, ylabel=$e_\mathrm{qoi}$, xmin=0, ymin=0.0002952946380670536,ymode=log]
\addplot [color=blue, mark=o, thick, forget plot]
coordinates {
( 0.00000000e+00,  4.89800922e-01)
( 1.00000000e+00,  4.72104157e-01)
( 2.00000000e+00,  4.08037185e-01)
( 3.00000000e+00,  3.95243626e-01)
( 4.00000000e+00,  3.77888241e-01)
( 5.00000000e+00,  1.62306696e-01)
( 6.00000000e+00,  1.20137750e-01)
( 7.00000000e+00,  4.18047210e-02)
( 8.00000000e+00,  3.12040614e-02)
( 9.00000000e+00,  2.53582229e-03)
( 1.00000000e+01,  2.37764241e-03)
( 1.10000000e+01,  2.38729459e-03)};

\addplot [color=blue, mark=x, thick, forget plot]
coordinates {
( 0.00000000e+00,  9.99880550e-01)
( 1.00000000e+00,  2.54919060e-01)
( 2.00000000e+00,  1.76224696e-01)
( 3.00000000e+00,  6.82549782e-02)
( 4.00000000e+00,  1.31178611e-02)
( 5.00000000e+00,  1.30352423e-02)
( 6.00000000e+00,  1.13751909e-02)
( 7.00000000e+00,  8.21794615e-03)
( 8.00000000e+00,  6.39856433e-03)
( 9.00000000e+00,  1.77155328e-03)
( 1.00000000e+01,  3.59448409e-04)
( 1.10000000e+01,  3.10836461e-04)};

\nextgroupplot[width=0.33\textwidth, xlabel=$r$, ymax=9.607222471916558e-05, xmax=10, ylabel=$e_\mathrm{est}$, xmin=0, ymin=1.1723320474307565e-07,ymode=log]
\addplot [color=blue, mark=o, thick, forget plot]
coordinates {
( 0.00000000e+00,  9.40585272e-06)
( 1.00000000e+00,  4.03797714e-05)
( 2.00000000e+00,  5.46367031e-05)
( 3.00000000e+00,  3.95272234e-05)
( 4.00000000e+00,  1.66597314e-05)
( 5.00000000e+00,  7.14267305e-06)
( 6.00000000e+00,  1.65401930e-05)
( 7.00000000e+00,  2.05952051e-05)
( 8.00000000e+00,  1.78854074e-05)
( 9.00000000e+00,  1.60303401e-05)
( 1.00000000e+01,  5.37321301e-05)};

\addplot [color=blue, mark=x, thick, forget plot]
coordinates {
( 0.00000000e+00,  9.14973569e-05)
( 1.00000000e+00,  8.98985352e-06)
( 2.00000000e+00,  7.54249853e-06)
( 3.00000000e+00,  6.12929419e-06)
( 4.00000000e+00,  6.14707018e-06)
( 5.00000000e+00,  1.86711328e-07)
( 6.00000000e+00,  1.77453479e-07)
( 7.00000000e+00,  1.23403373e-07)
( 8.00000000e+00,  1.33866028e-07)
( 9.00000000e+00,  4.10359159e-07)
( 1.00000000e+01,  5.17410097e-07)};

\end{groupplot}\end{tikzpicture}
	\caption{The mean error of the CG-GL method applied to the Poisson problem over the test set $\Dcal_\mathrm{interp}$
	(\ref{poi_limitedtraing_int}) and $\Dcal_\mathrm{extrap}$ (\ref{poi_limitedtraing_ext}) as a function of refinement.
	Because the CG-GL method is initialized with $\Omega_l = \emptyset$ and $N_l=0$, $r= 0$ corresponds to a traditional ROM.}
\label{fig:limited_train_UJE}
\end{figure}
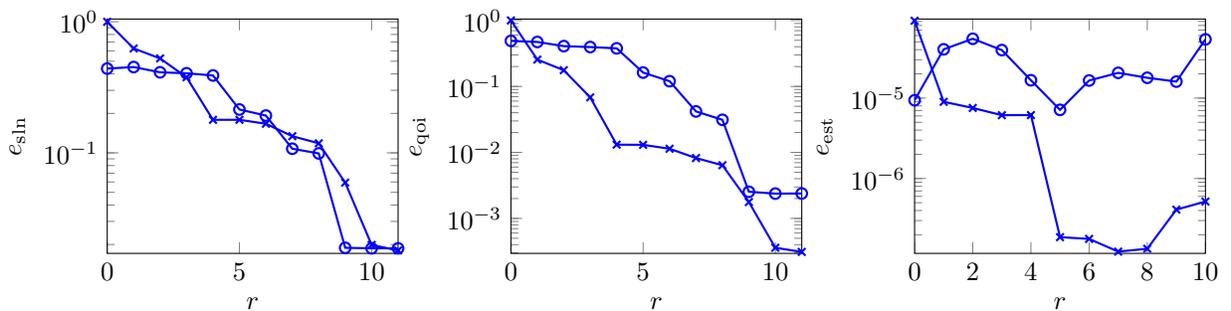

Next, we study the performance of the ROM and CG-GL methods with limited training. The training set is defined as $\Dcal_{\mathrm{train}}\coloneqq\{(1,1,2,6), (1,1,6,2), (1,1,6,6)\}$, which not only reduces the total training parameters compared with the previous case but also breaks the symmetry by removing snapshots with the peak in the lower left of the domain. This asymmetry manifests in the POD modes (Figure~\ref{fig:poi_pod_1}) and, due to the limited training, we do not truncate the POD basis ($k=3$). For this study, we choose the test set to be $\Dcal_\mathrm{test} = \Dcal_\mathrm{interp} \cup \Dcal_\mathrm{extrap}$, where
$\Dcal_\mathrm{interp} = \{(1,1,14/3,14/3)\}$ and $\Dcal_\mathrm{extrap} = \{(1,1,0.5,0.5)\}$. Similar to the previous section, we use the ROM to initialize the CG-GL method ($\Omega_l = \emptyset$, $N_l = 0$) with 16 patches of 36 quadratic quadrilateral elements.
At the interpolated parameter, the peak (location and magnitude) is incorrectly predicted by the ROM and gradually transformed to the true peak (Figure~\ref{fig:poi_interp_example_sol_sequential}) using the error estimation framework (Figure~\ref{fig:error_contour_0}) and CG-GL adaptation. At the extrapolated parameter, the ROM prediction is completely useless; however, an accurate approximation quickly emerges (Figure~\ref{fig:poi_extrap_example_sol_sequential}) from the adaptive CG-GL framework guided by the error estimates (Figure~\ref{fig:error_contour_1}). These results are quantified in Figure~\ref{fig:limited_train_UJE}, which shows the CG-GL error estimation and adaptation framework is able to effectively reduce the approximation error, even when the initial ROM prediction is poor.

\subsubsection{Effect of feature locality}
\label{sec:rslt:poi:local}
\begin{figure}
	\centering
	\begin{tikzpicture}
		\begin{groupplot}[
			group style={
				group size=3  by 1,
				horizontal sep=1cm
			},
			width=0.35\textwidth,
			axis equal image,
			xlabel={$x_1$},
			ylabel={$x_2$},
			xtick = {0.0, 4.0, 8.0},
			ytick = {0.0, 4.0, 8.0},
			xmin=0, xmax=8,
			ymin=0, ymax=8
			]
			\nextgroupplot[title={$\mu=(10,   0.1,    0.5,     0.5)$}]
			\addplot graphics [xmin=0, xmax=8.0, ymin=0, ymax=8.0] {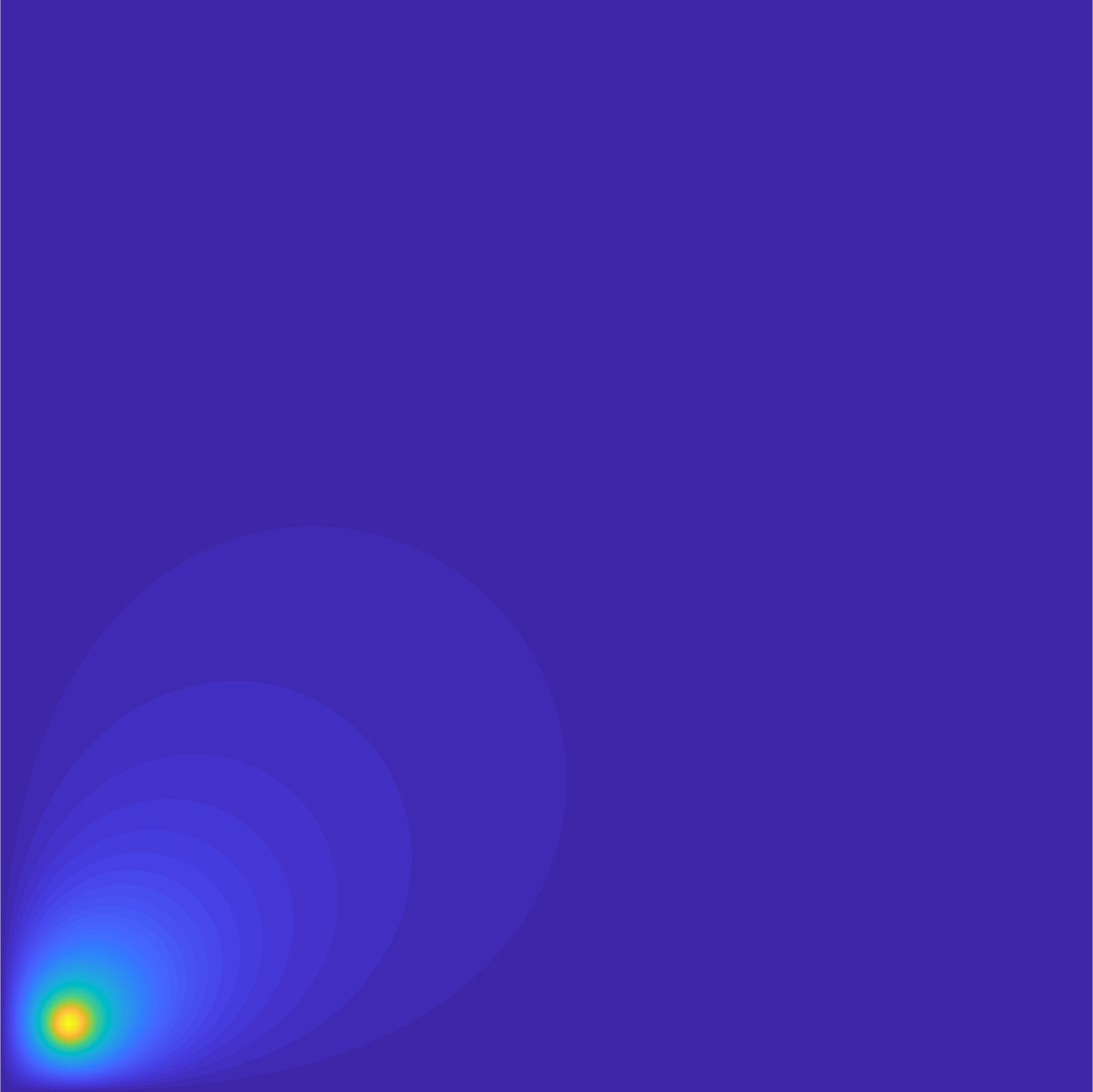};
			
			\nextgroupplot[title={$\mu=(10, 0.5,  0.5,  0.5)$}, ylabel={}, ytick=\empty]
			\addplot graphics [xmin=0, xmax=8.0, ymin=0, ymax=8.0] {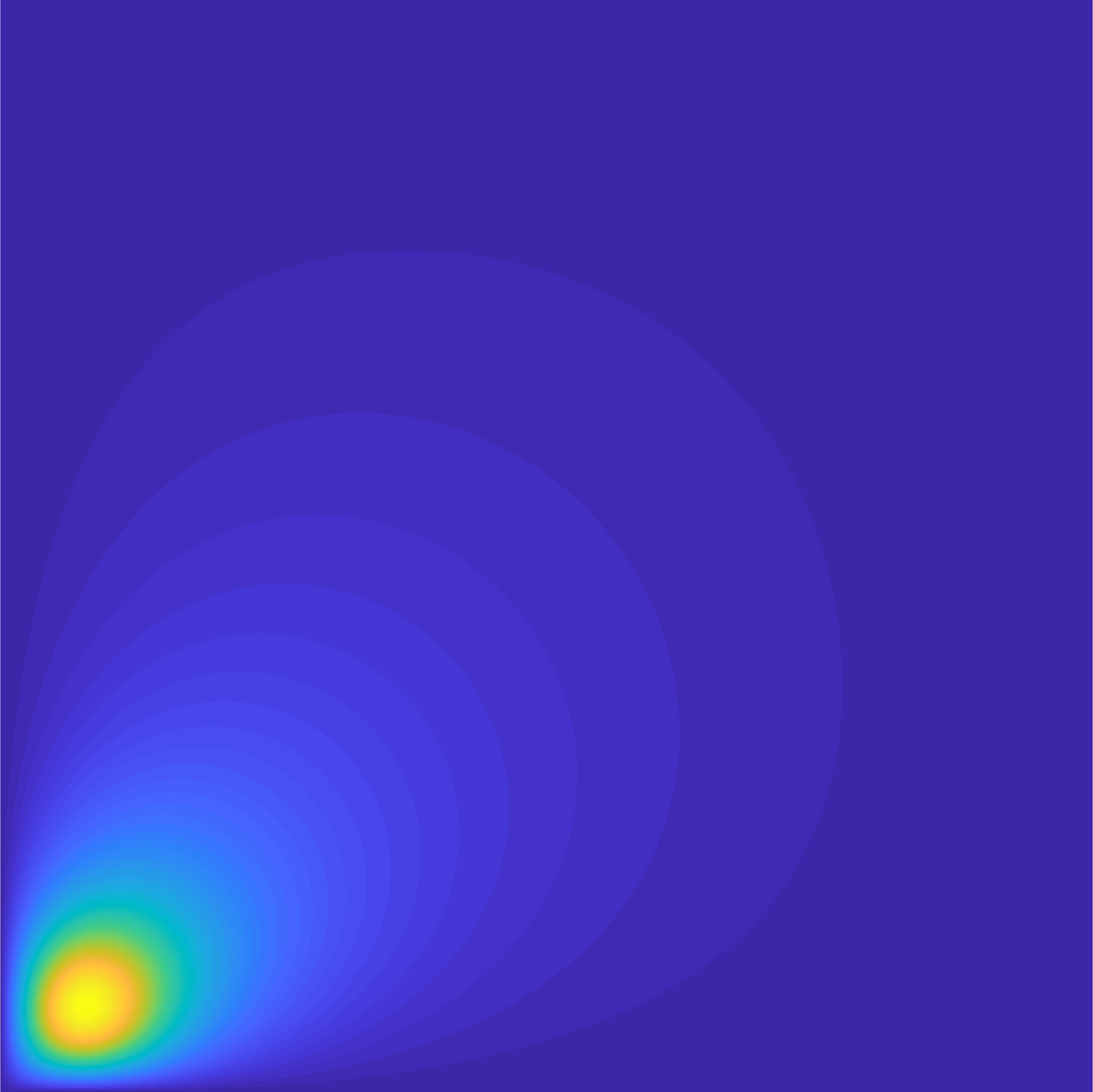};
			
			\nextgroupplot[title={$\mu=(10,  1,  0.5,  0.5)$}, ylabel={}, ytick=\empty]
			\addplot graphics [xmin=0, xmax=8.0, ymin=0, ymax=8.0] {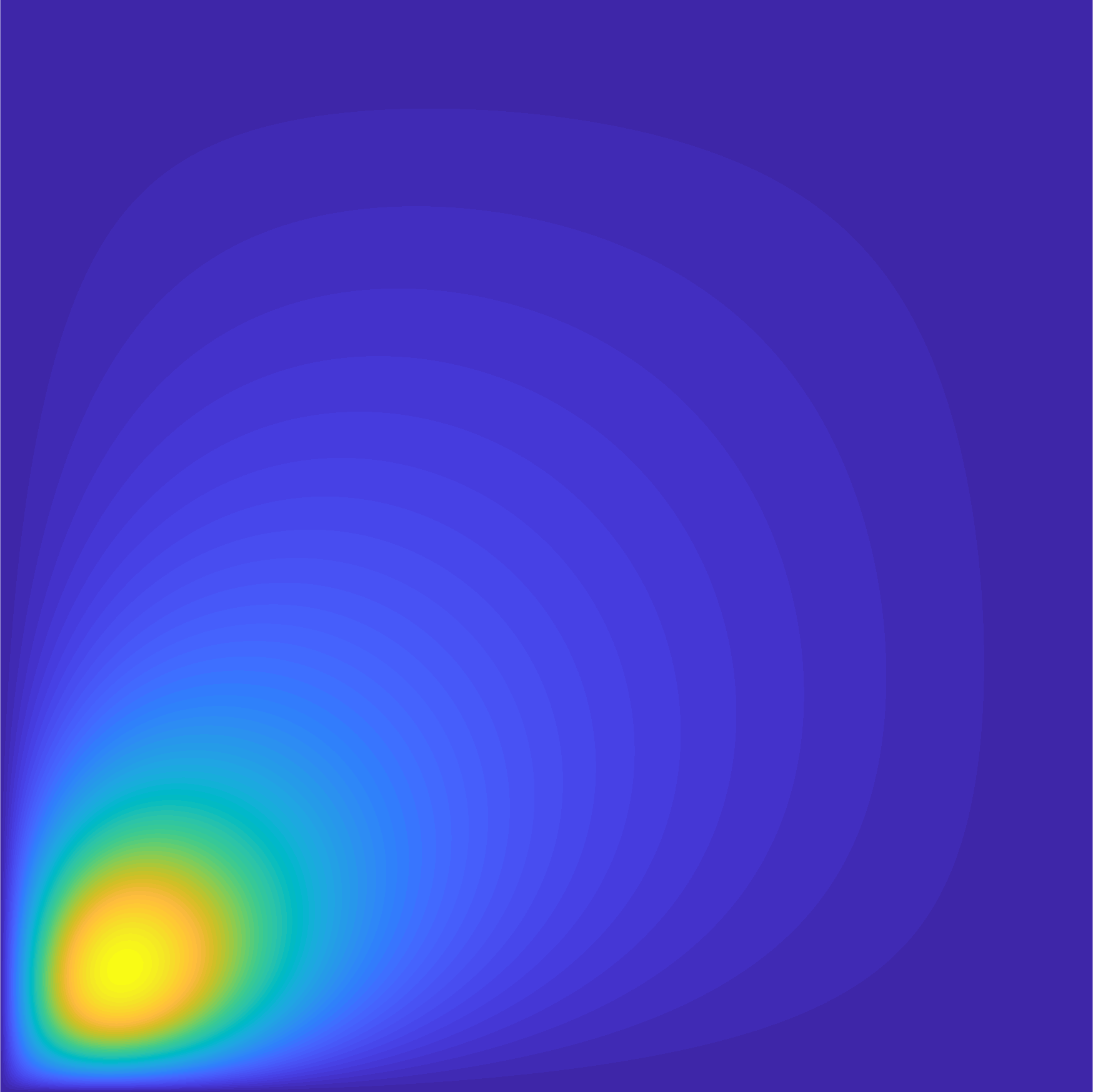};
		
		\end{groupplot}
	\end{tikzpicture}
	\caption{Solution of Poisson equation (\ref{eqn:poi}) for three different values of $\mu_2$.}
	\label{fig:poi_mu1}
\end{figure}

\begin{figure}
	\centering
	\begin{tikzpicture}
		\begin{groupplot}[
			group style={
				group size=4 by 1,
				horizontal sep=0.5cm
			},
			xmajorgrids=true,
			ymajorgrids=true,
			width=0.35\textwidth,
			axis equal image,
			xtick = {0, 2,4,6, 8},
			ytick = {0, 2,4,6, 8},
			xticklabel=\empty,
			yticklabel=\empty,
			xmin=0, xmax=8,
			ymin=0, ymax=8,
			axis on top,
			grid style={line width=1pt, draw=gray}
			]
			
			\nextgroupplot[title={}, xmajorgrids=false, ymajorgrids=false]
			\addplot graphics [xmin=0, xmax=8.0, ymin=0, ymax=8.0] {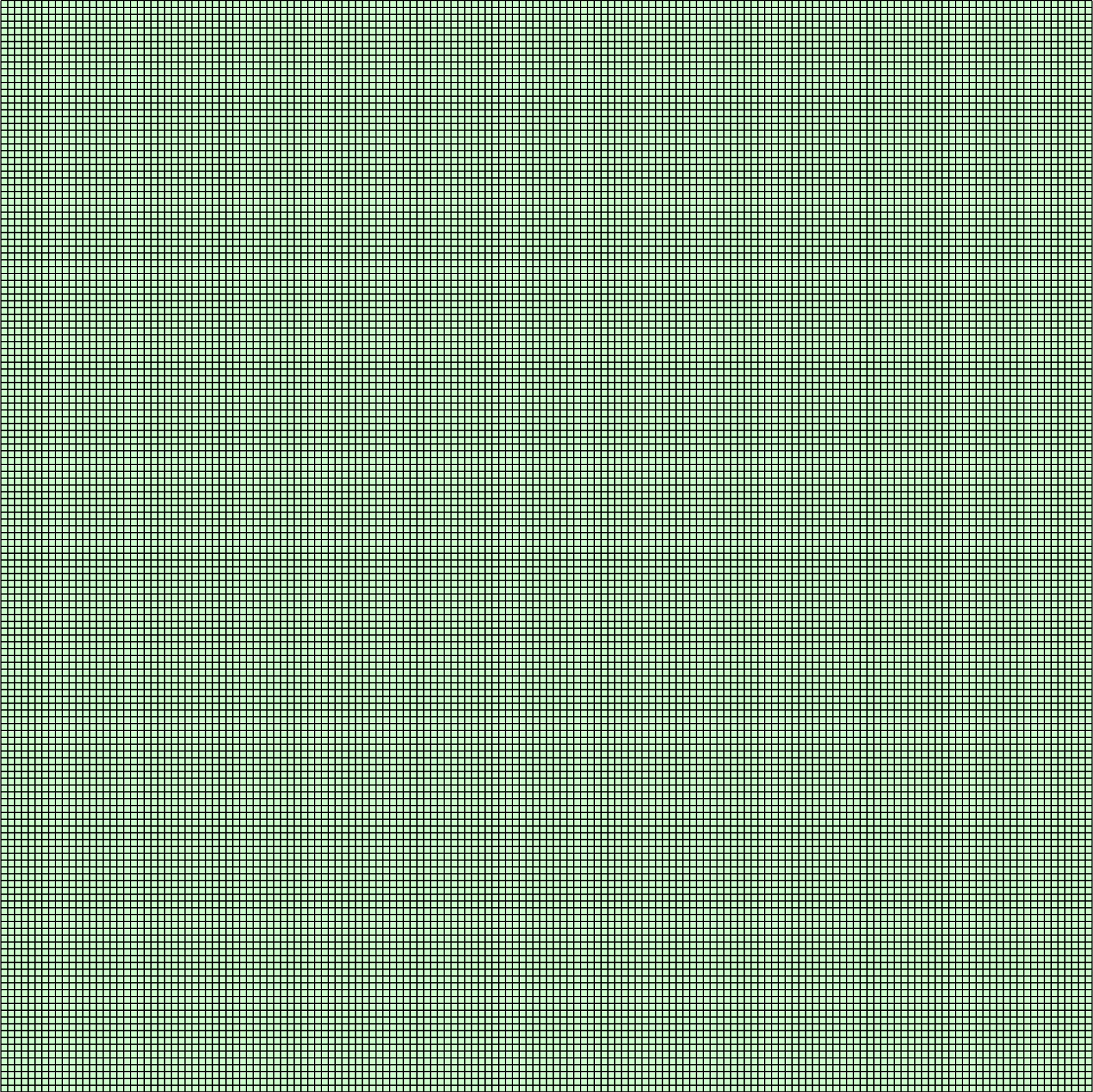};
		
			\nextgroupplot[title={},ylabel={}]
			\addplot graphics [xmin=0, xmax=8.0, ymin=0, ymax=8.0] {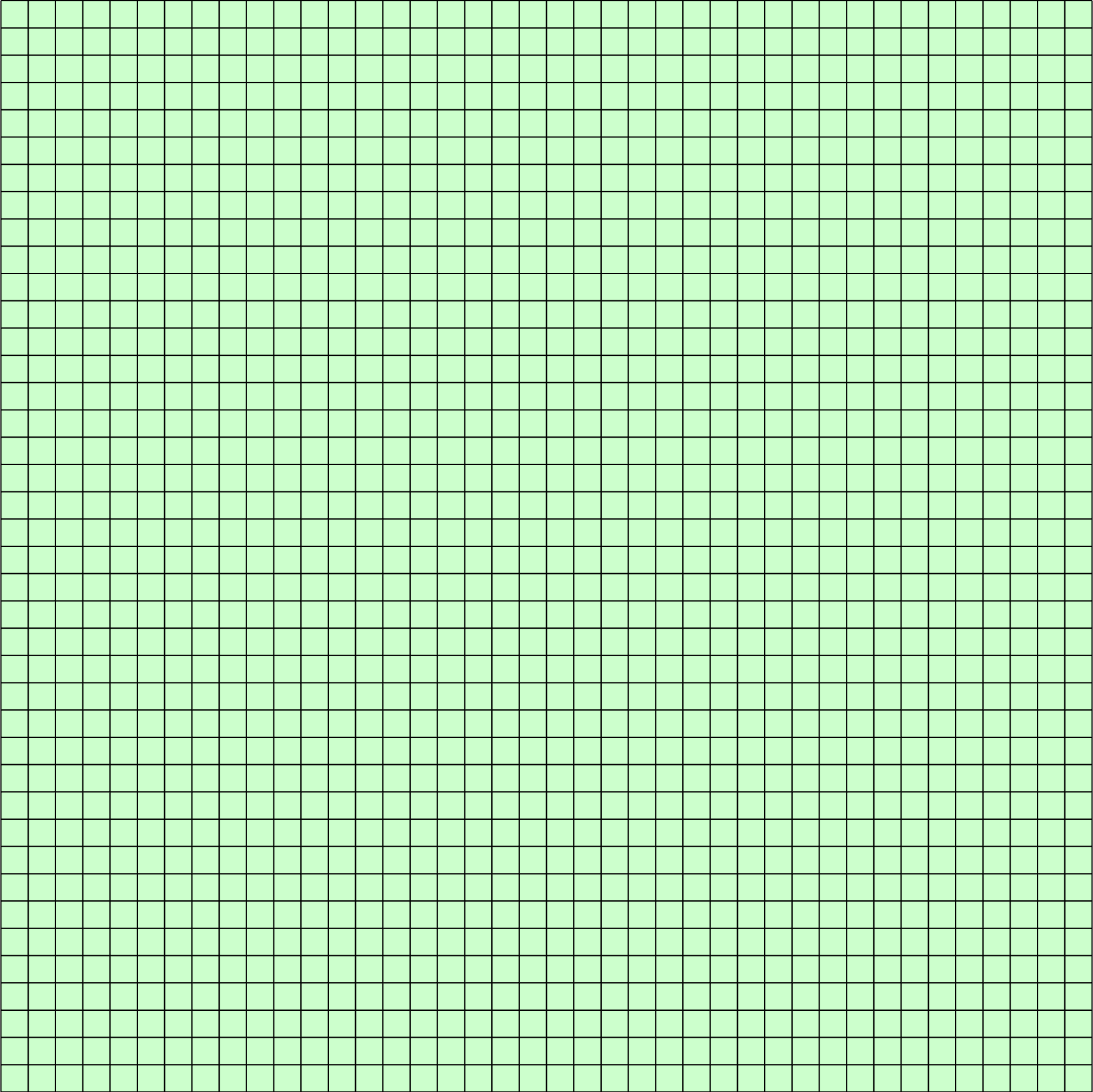};
			
			\nextgroupplot[title={}, ylabel={}, ytick={}]
			\addplot graphics [xmin=0, xmax=8.0, ymin=0, ymax=8.0] {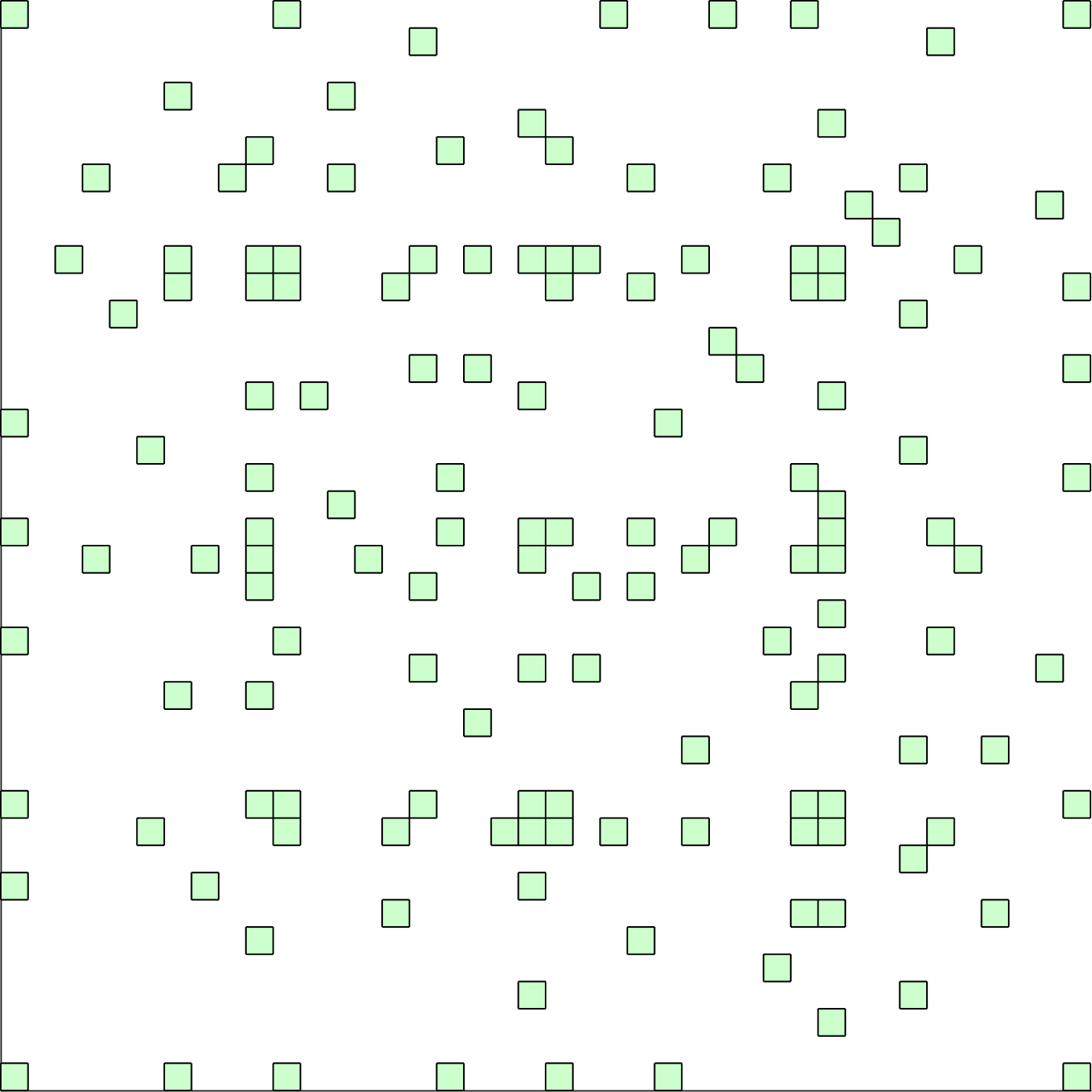};
			
			\nextgroupplot[title={}, ylabel={}, ytick={}, xmajorgrids=false, ymajorgrids=false]
			\addplot graphics [xmin=0, xmax=8.0, ymin=0, ymax=8.0] {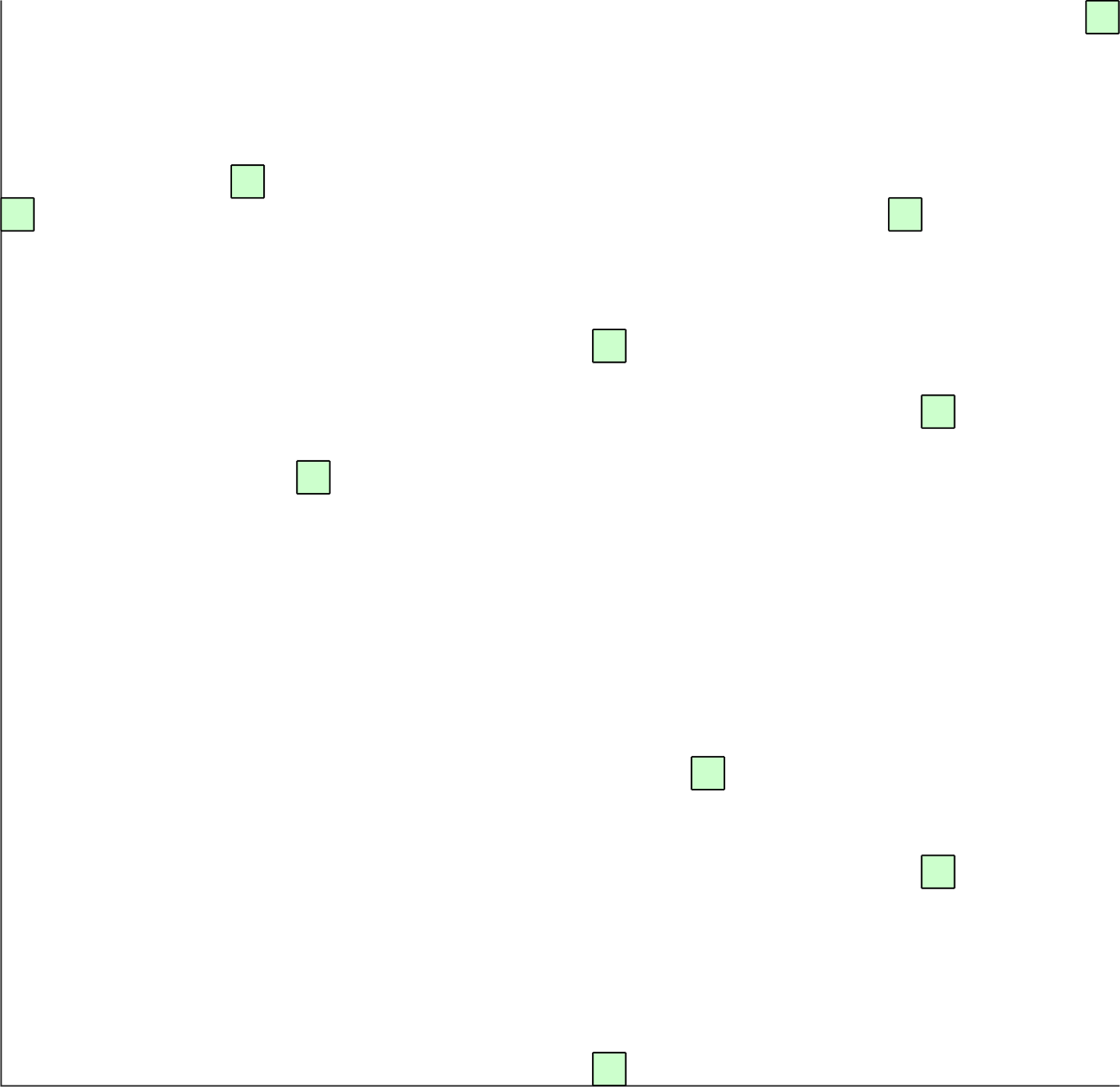};
		\end{groupplot}
	\end{tikzpicture}
	\caption{The mesh for the reference solution (\textit{left}), initial CG-GL patch meshes (\textit{middle left}), sample meshes for each CG-GL patch (\textit{middle right}), and sample mesh for the ROM (\textit{right}).}
	\label{fig:meshes_EQP}
\end{figure}

\begin{figure}
	\centering
	\begin{tikzpicture}
		\begin{groupplot}[
			group style={
				group size = 3 by 1,
				horizontal sep=0.5cm
			},
			xmajorgrids=true,
			ymajorgrids=true,
			width=0.35\textwidth,
			axis equal image,
			xtick = {0, 2,4,6, 8},
			ytick = {0, 2,4,6, 8},
			xticklabel=\empty,
			yticklabel=\empty,
			xmin=0, xmax=8,
			ymin=0, ymax=8,
			axis on top,
			grid style={line width=1pt, draw=gray}
			]
			
			\nextgroupplot[title={}]
			\addplot graphics [xmin=0, xmax=8.0, ymin=0, ymax=8.0] {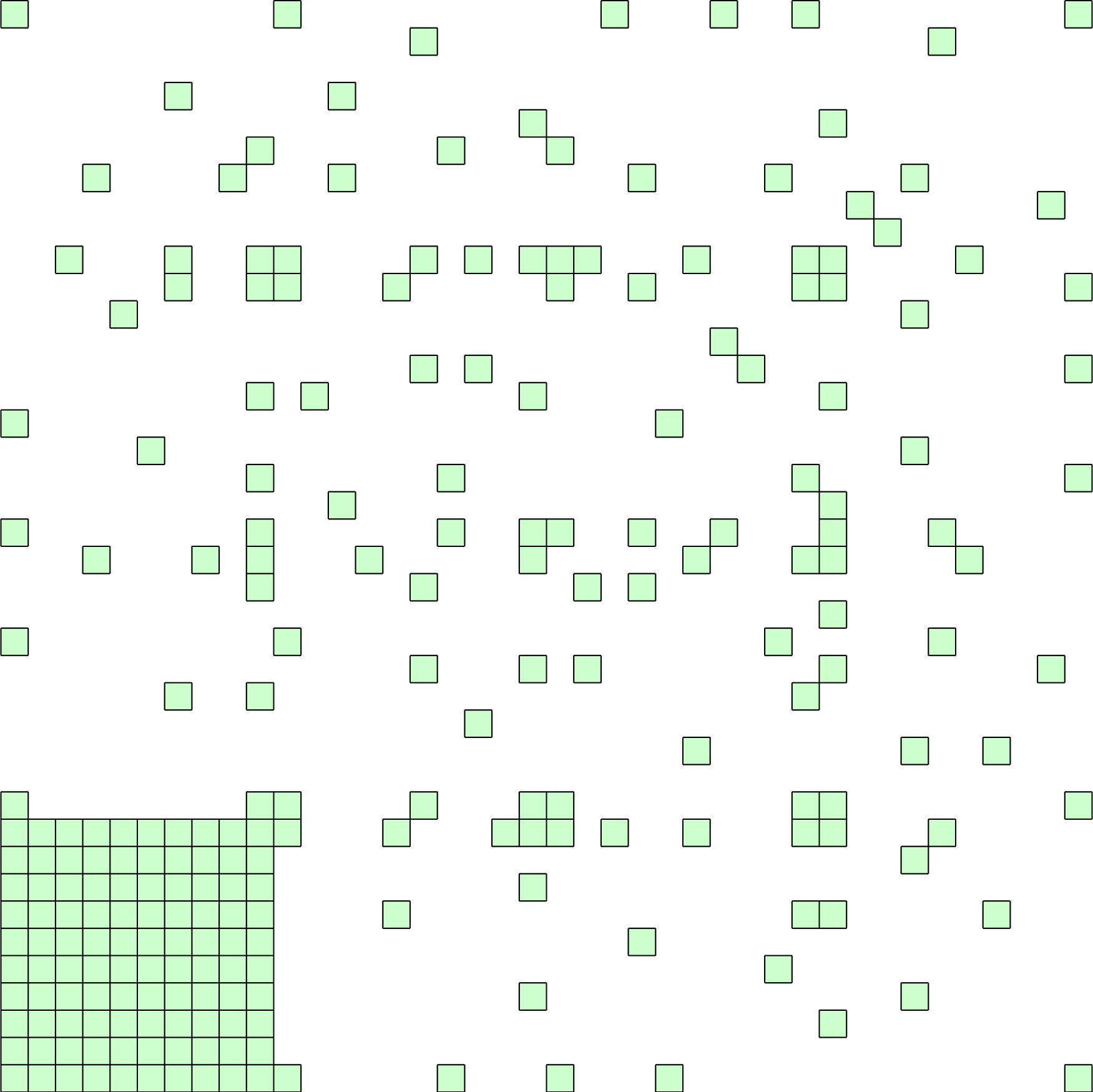};
			
			\nextgroupplot[title={}, ylabel={}, ytick={}]
			\addplot graphics [xmin=0, xmax=8.0, ymin=0, ymax=8.0] {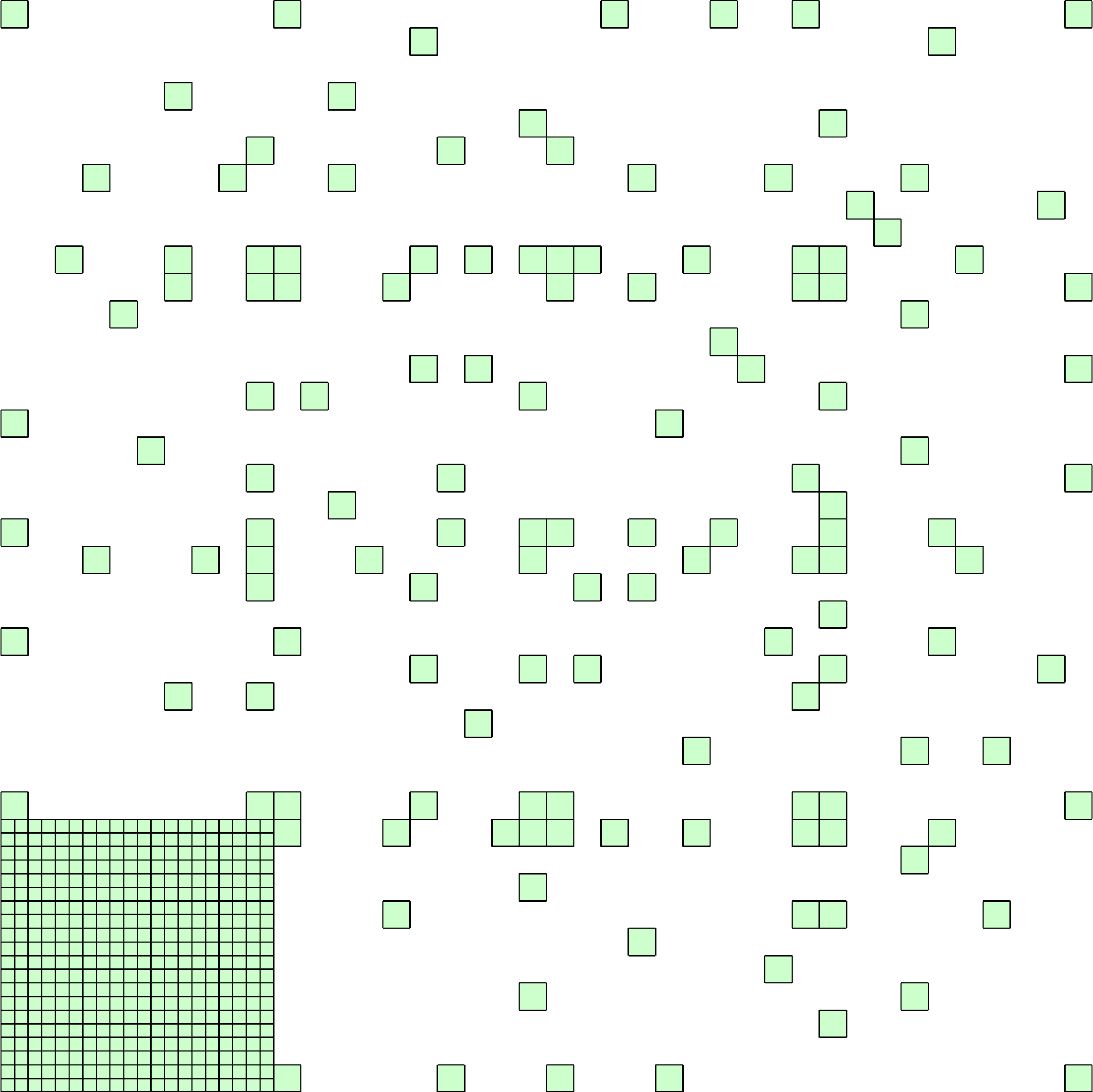};
			
			\nextgroupplot[title={}, ylabel={}, ytick={}]
			\addplot graphics [xmin=0, xmax=8.0, ymin=0, ymax=8.0] {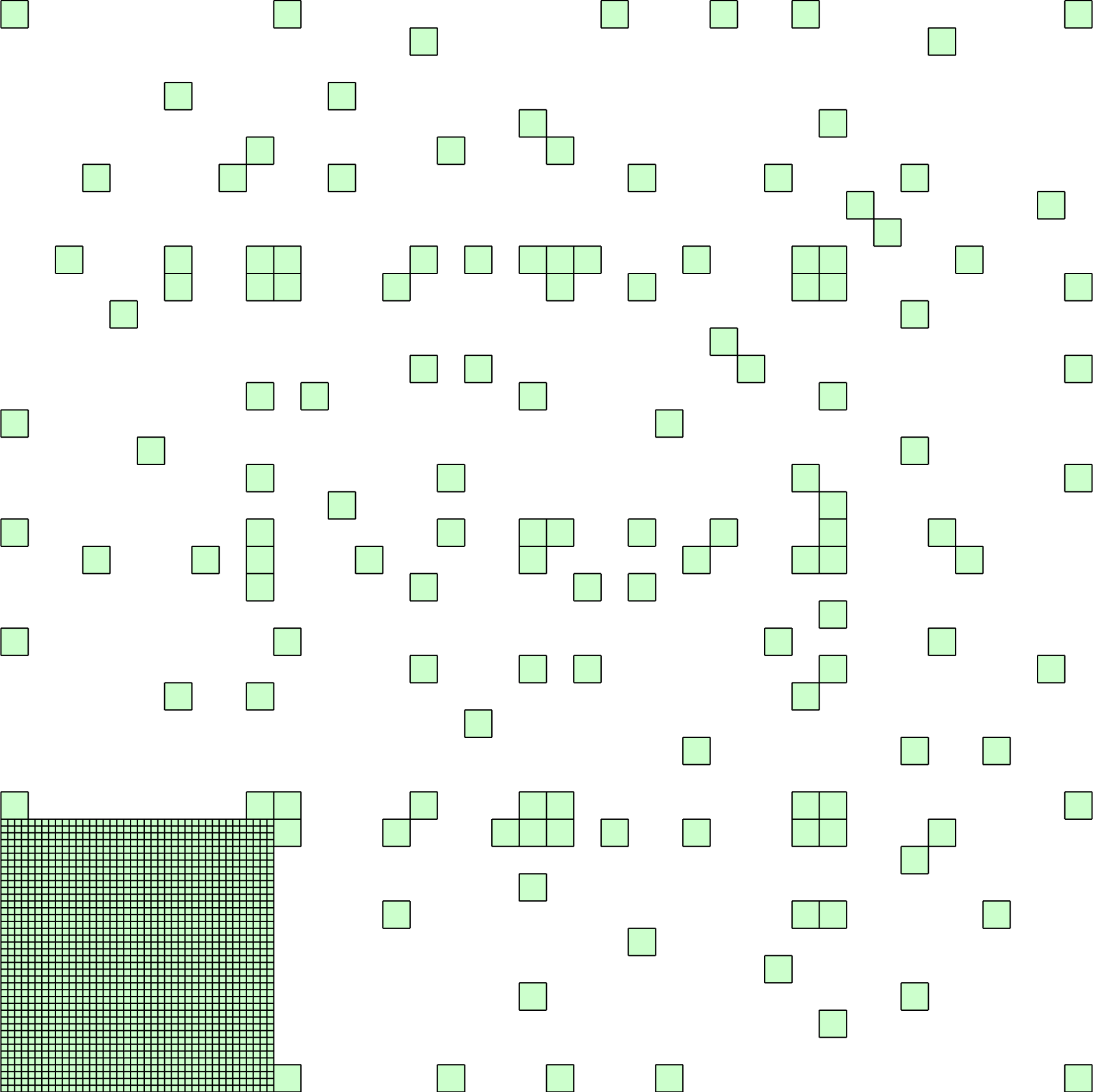};
		\end{groupplot}
	\end{tikzpicture}
	\caption{The sequence of discretizations used by the adaptive CG-GL method (\textit{left-to-right}) for the Poisson problem study in Section~\ref{sec:rslt:poi:local} at $\mu_\mathrm{test}   =  (10, 0.1, 0.5, 0.5)$.}
	\label{fig:poi_local0}
\end{figure}

\begin{figure}
	\centering
	\begin{tikzpicture}
		\begin{groupplot}[
			group style={
				group size = 3 by 2,
				horizontal sep=0.5cm,
				vertical sep=0.5cm
			},
			xmajorgrids=true,
			ymajorgrids=true,
			width=0.35\textwidth,
			axis equal image,
			xtick = {0, 2,4,6, 8},
			ytick = {0, 2,4,6, 8},
			xticklabel=\empty,
			yticklabel=\empty,
			xmin=0, xmax=8,
			ymin=0, ymax=8,
			axis on top,
			grid style={line width=1pt, draw=gray}
			]
			
			\nextgroupplot[title={}, xlabel = {}]
			\addplot graphics [xmin=0, xmax=8.0, ymin=0, ymax=8.0] {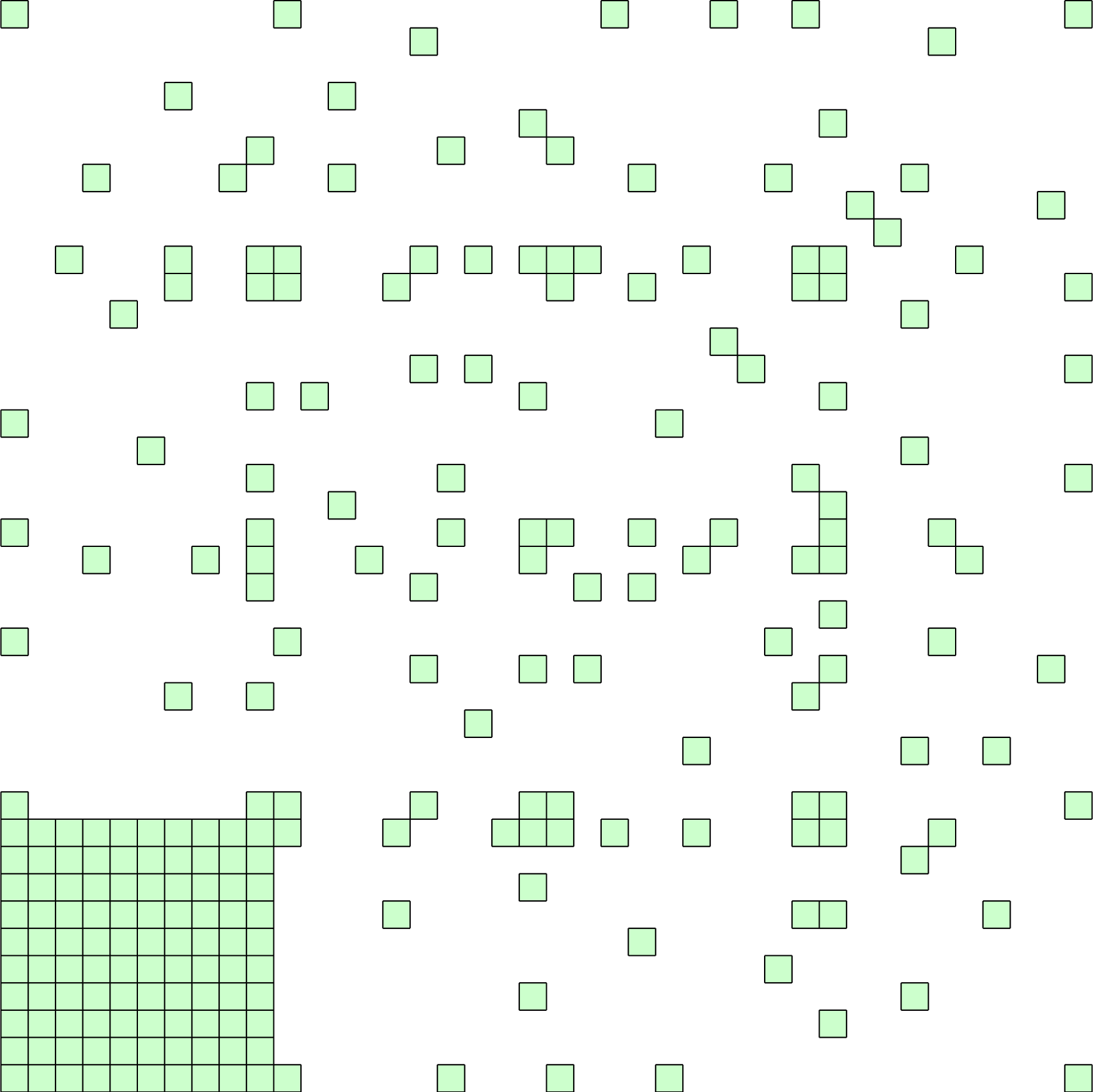};
			
			\nextgroupplot[title={}, xlabel = {}, ylabel={}, ytick={}]
			\addplot graphics [xmin=0, xmax=8.0, ymin=0, ymax=8.0] {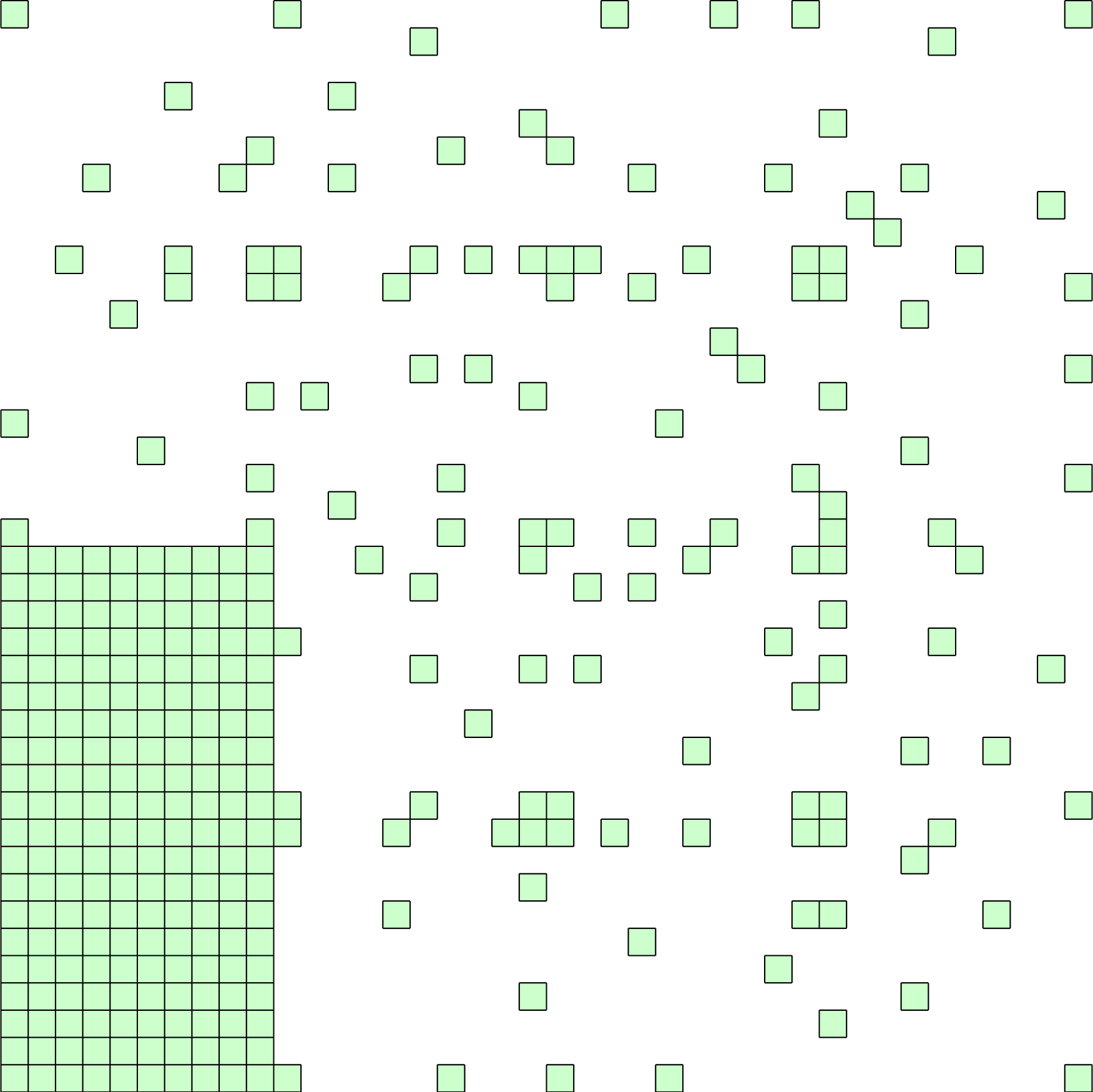};
			
			\nextgroupplot[title={}, ylabel={},xlabel = {},  ytick={}]
			\addplot graphics [xmin=0, xmax=8.0, ymin=0, ymax=8.0] {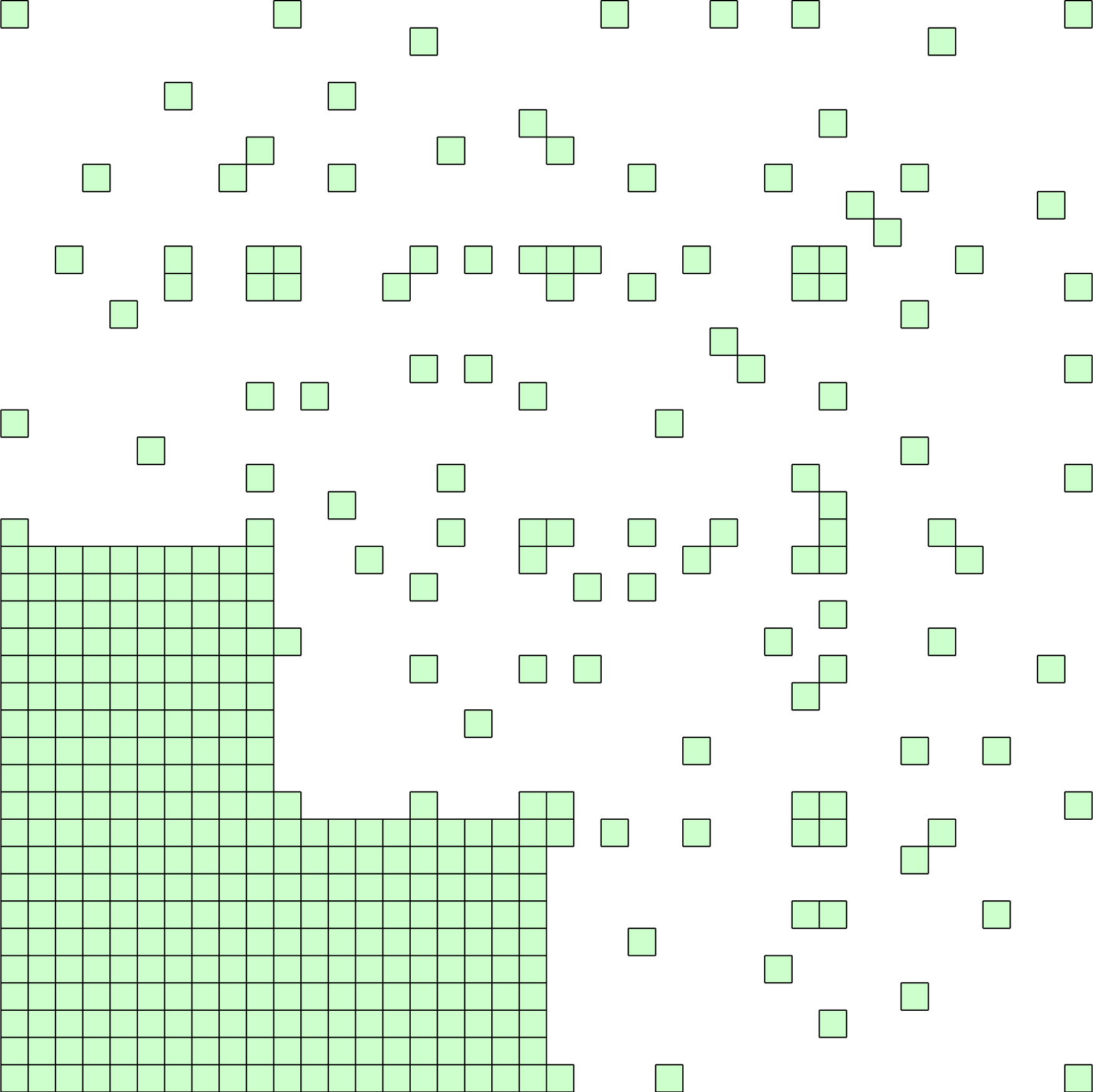};
			
			\nextgroupplot[title={}]
			\addplot graphics [xmin=0, xmax=8.0, ymin=0, ymax=8.0] {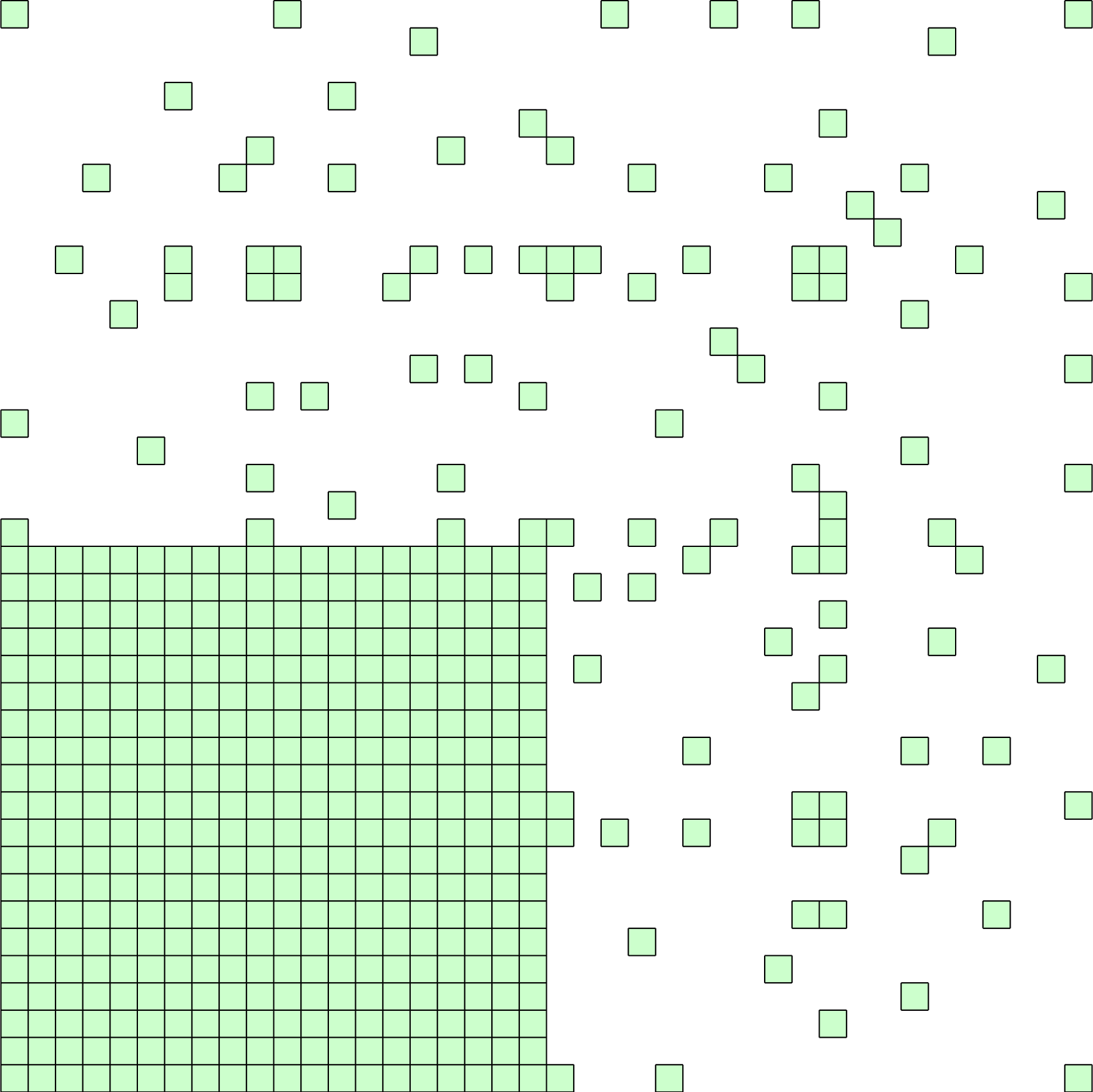};
			
			\nextgroupplot[title={}, ylabel={}, ytick={}]
			\addplot graphics [xmin=0, xmax=8.0, ymin=0, ymax=8.0] {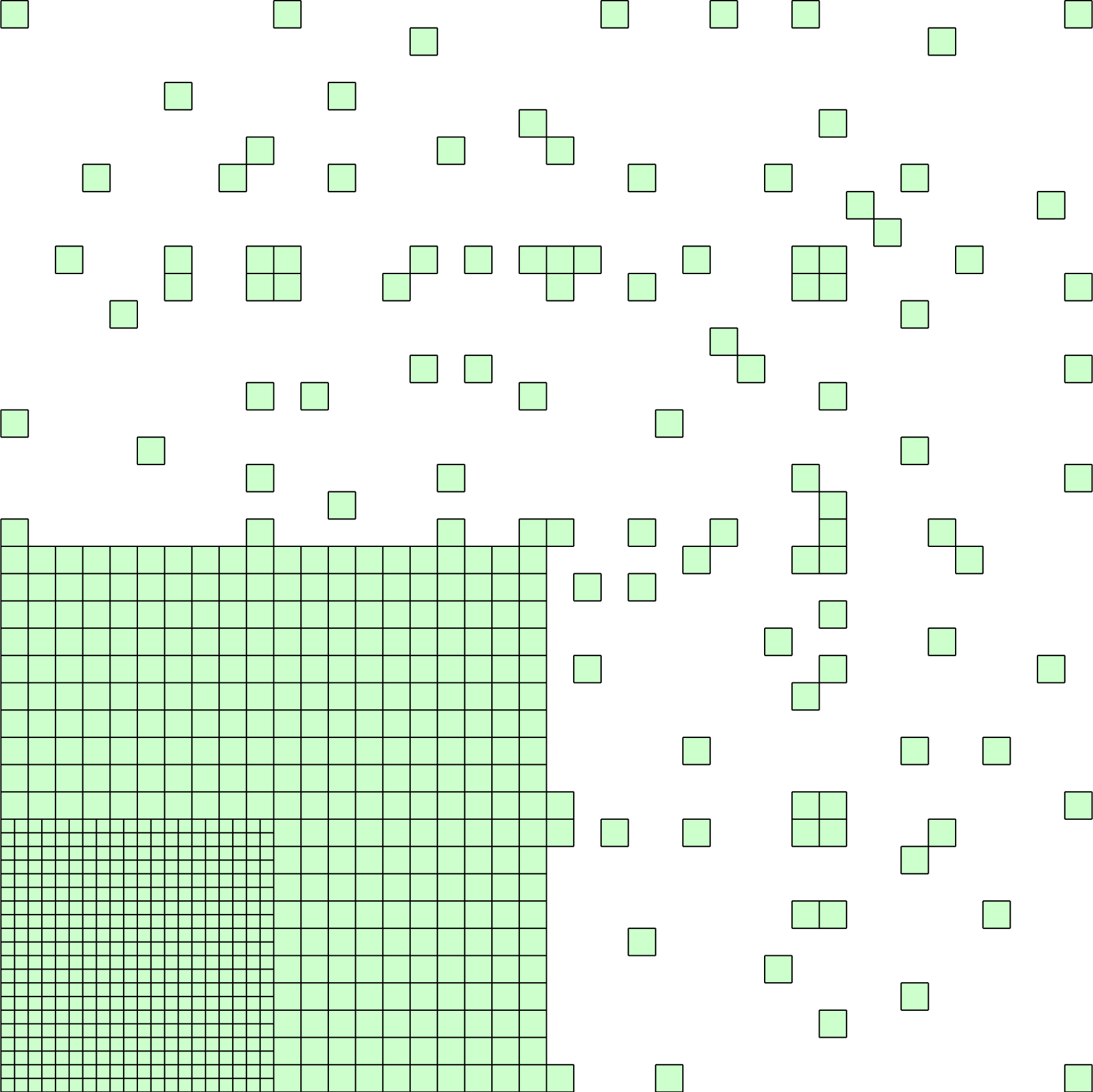};
			
			\nextgroupplot[title={}, ylabel={}, ytick={}]
			\addplot graphics [xmin=0, xmax=8.0, ymin=0, ymax=8.0] {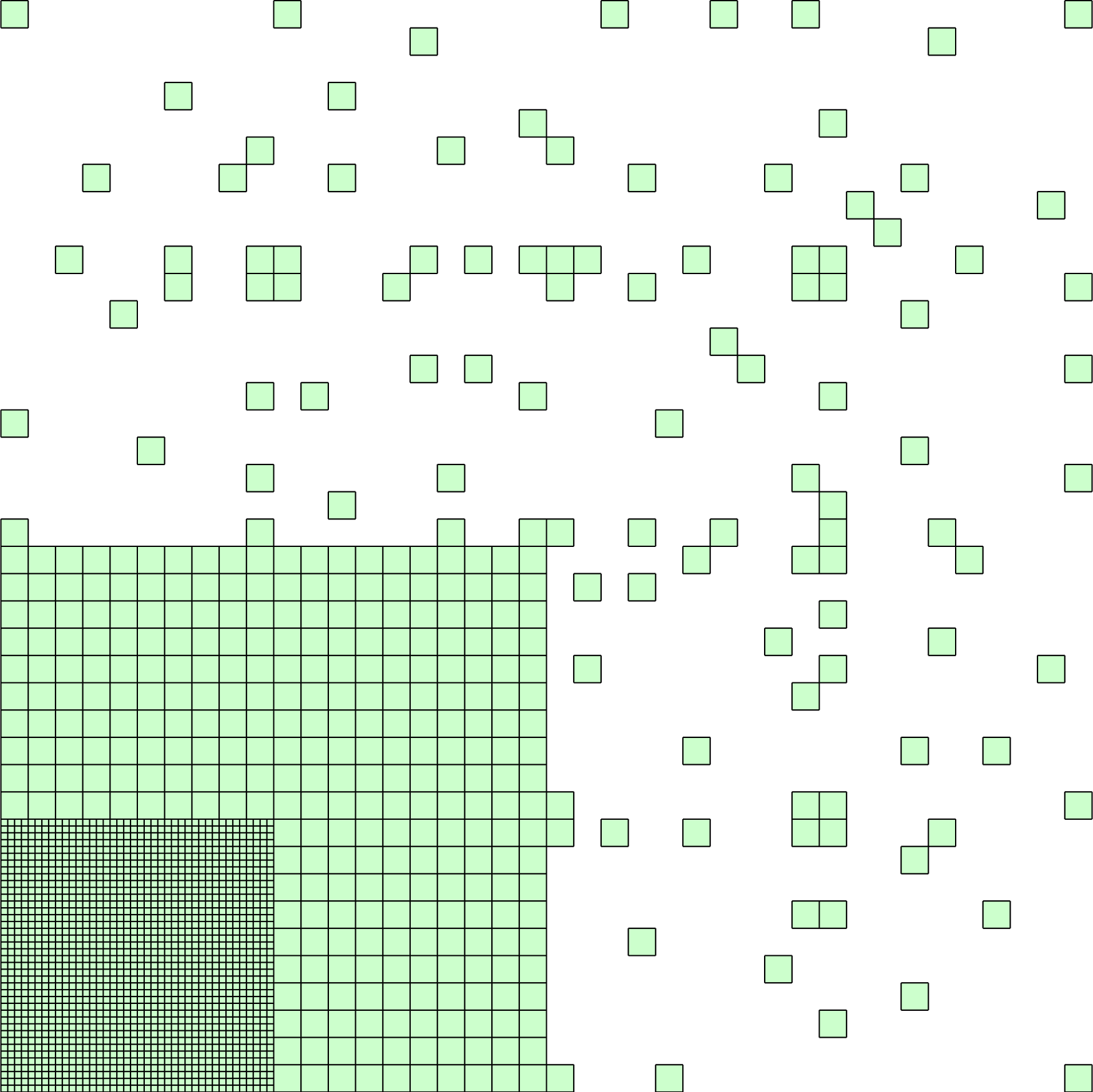};

		\end{groupplot}
	\end{tikzpicture}
	\caption{The sequence of discretizations used by the adaptive CG-GL method (\textit{left-to-right}, \textit{top-to-bottom}) for the Poisson problem study in Section~\ref{sec:rslt:poi:local} at $\mu_\mathrm{test}   =  (10, 0.5, 0.5, 0.5)$.}
	\label{fig:poi_local1}
\end{figure}

\begin{figure}
	\centering
	\begin{tikzpicture}
		\begin{groupplot}[
			group style={
				group size = 4 by 3,
				horizontal sep=0.5cm,
				vertical sep=0.5cm
			},
			xmajorgrids=true,
			ymajorgrids=true,
			width=0.35\textwidth,
			axis equal image,
			xtick = {0, 2,4,6, 8},
			ytick = {0, 2,4,6, 8},
			xticklabel=\empty,
			yticklabel=\empty,
			xmin=0, xmax=8,
			ymin=0, ymax=8,
			axis on top,
			grid style={line width=1pt, draw=gray}
			]
			
			\nextgroupplot[title={}, xlabel = {}]
			\addplot graphics [xmin=0, xmax=8.0, ymin=0, ymax=8.0] {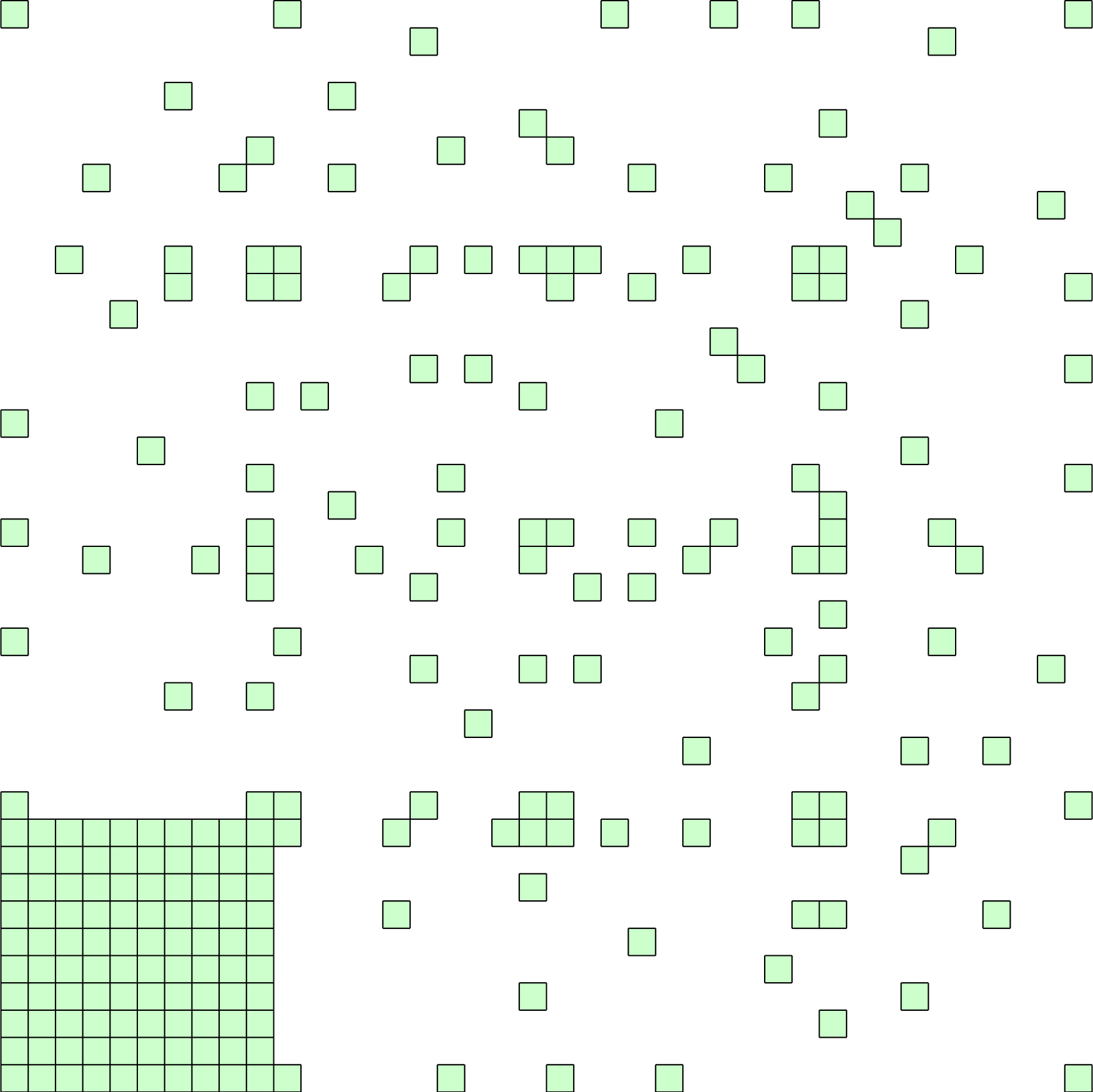};
			
			\nextgroupplot[title={}, xlabel = {}, ylabel={}, ytick={}]
			\addplot graphics [xmin=0, xmax=8.0, ymin=0, ymax=8.0] {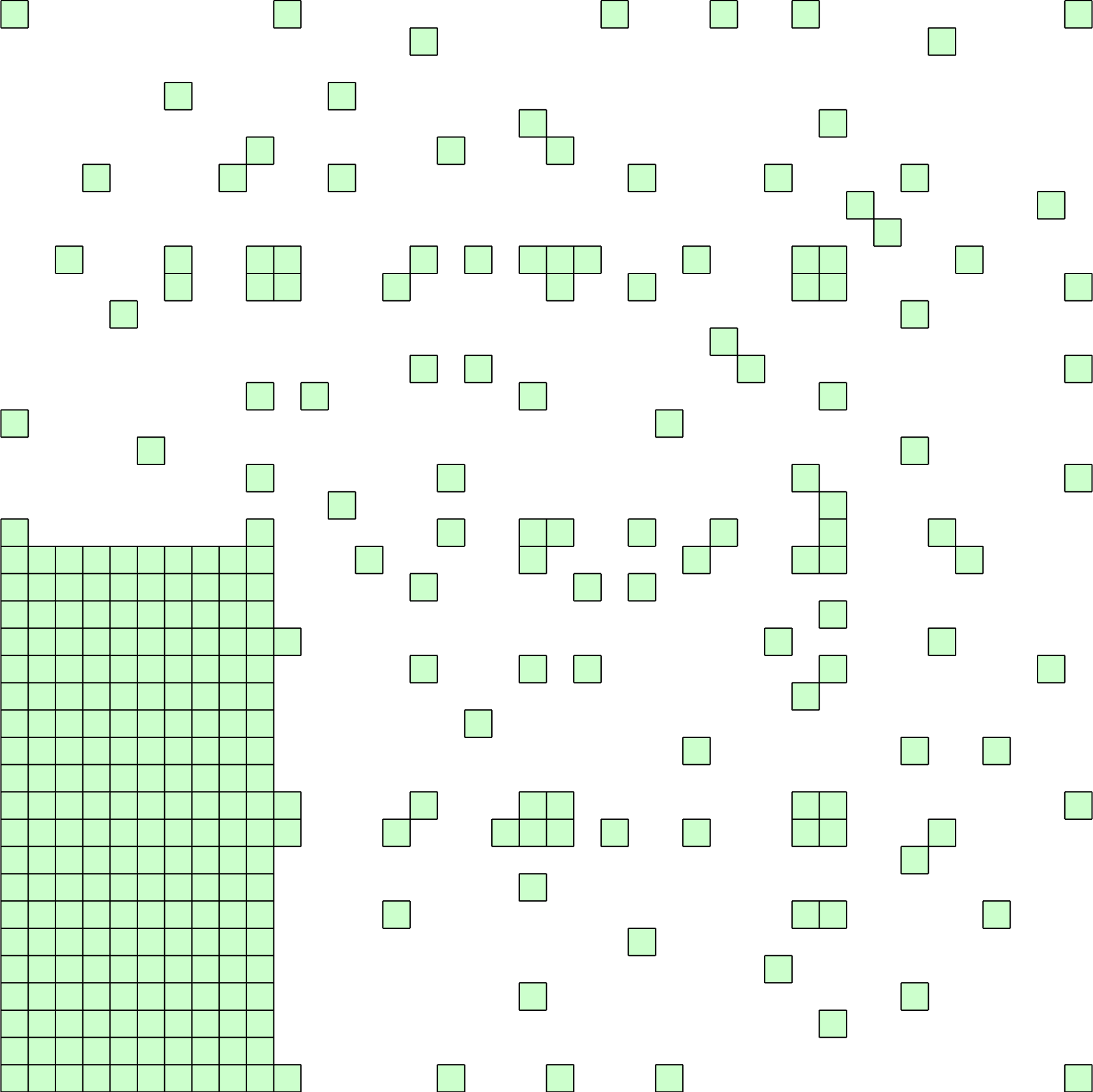};
		
			\nextgroupplot[title={}, xlabel = {}, ylabel={}, ytick={}]
			\addplot graphics [xmin=0, xmax=8.0, ymin=0, ymax=8.0] {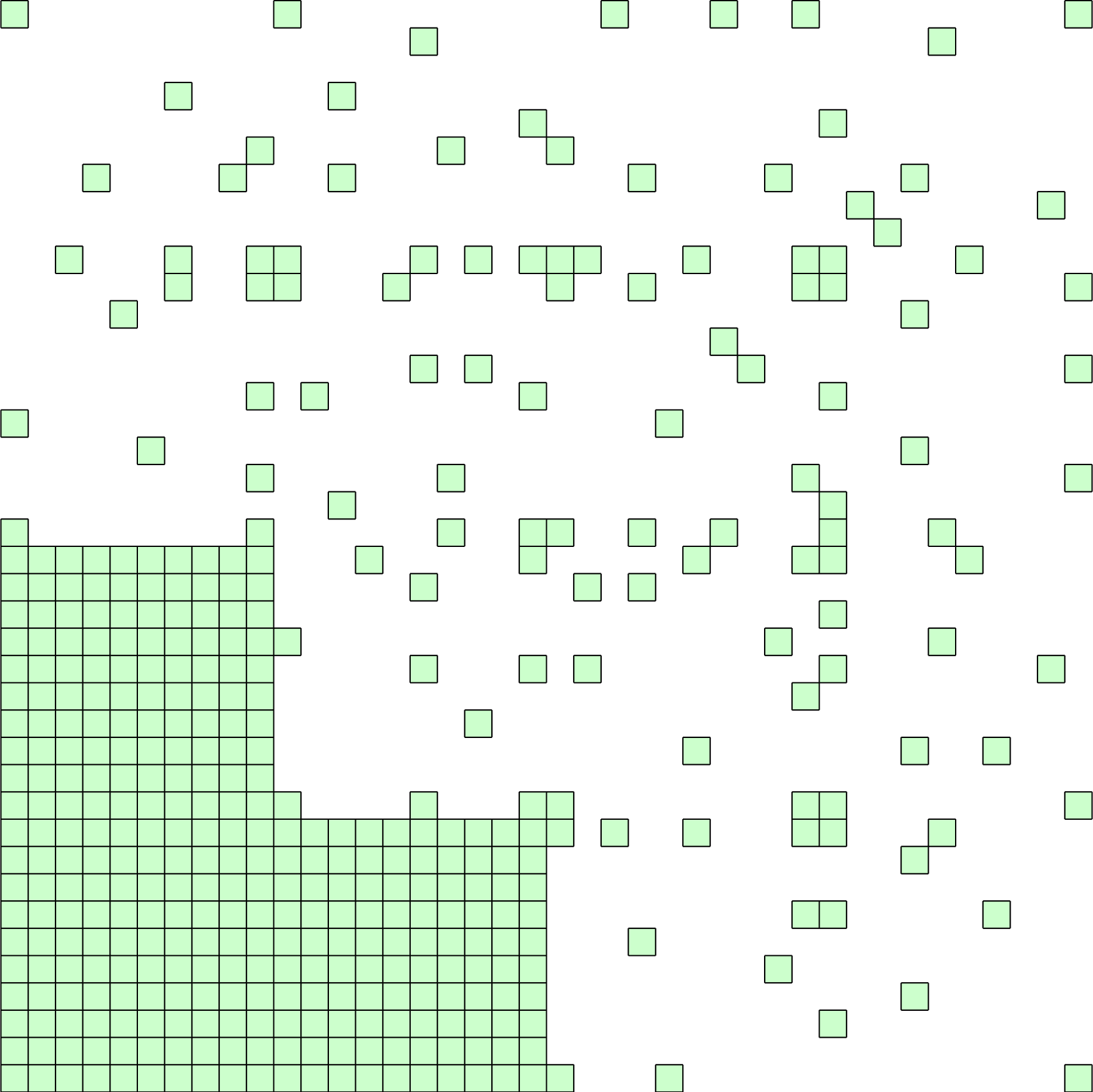};
			
			\nextgroupplot[title={}, ylabel={},xlabel = {},  ytick={}]
			\addplot graphics [xmin=0, xmax=8.0, ymin=0, ymax=8.0] {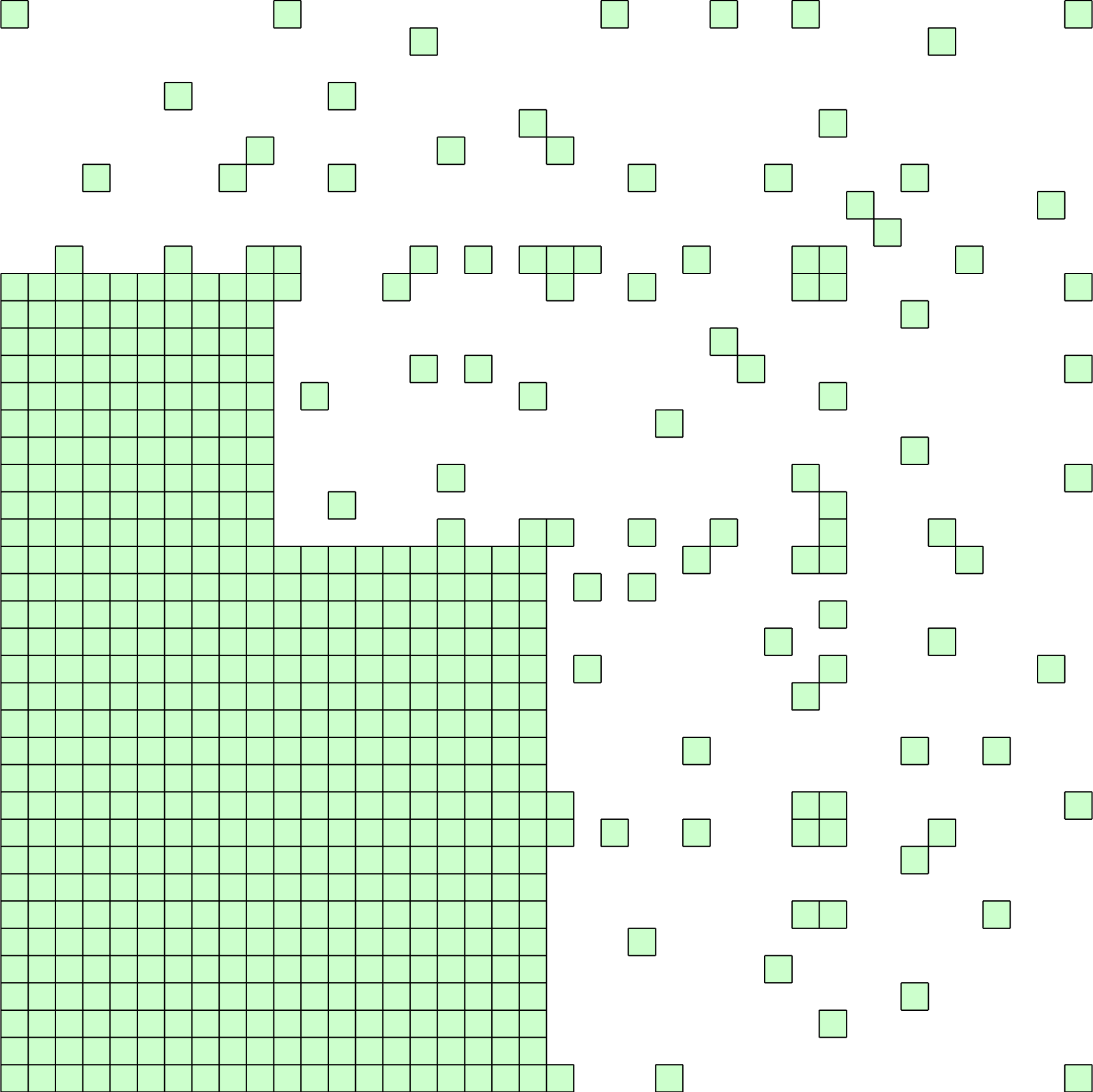};
			
			\nextgroupplot[title={}, xlabel = {}]
			\addplot graphics [xmin=0, xmax=8.0, ymin=0, ymax=8.0] {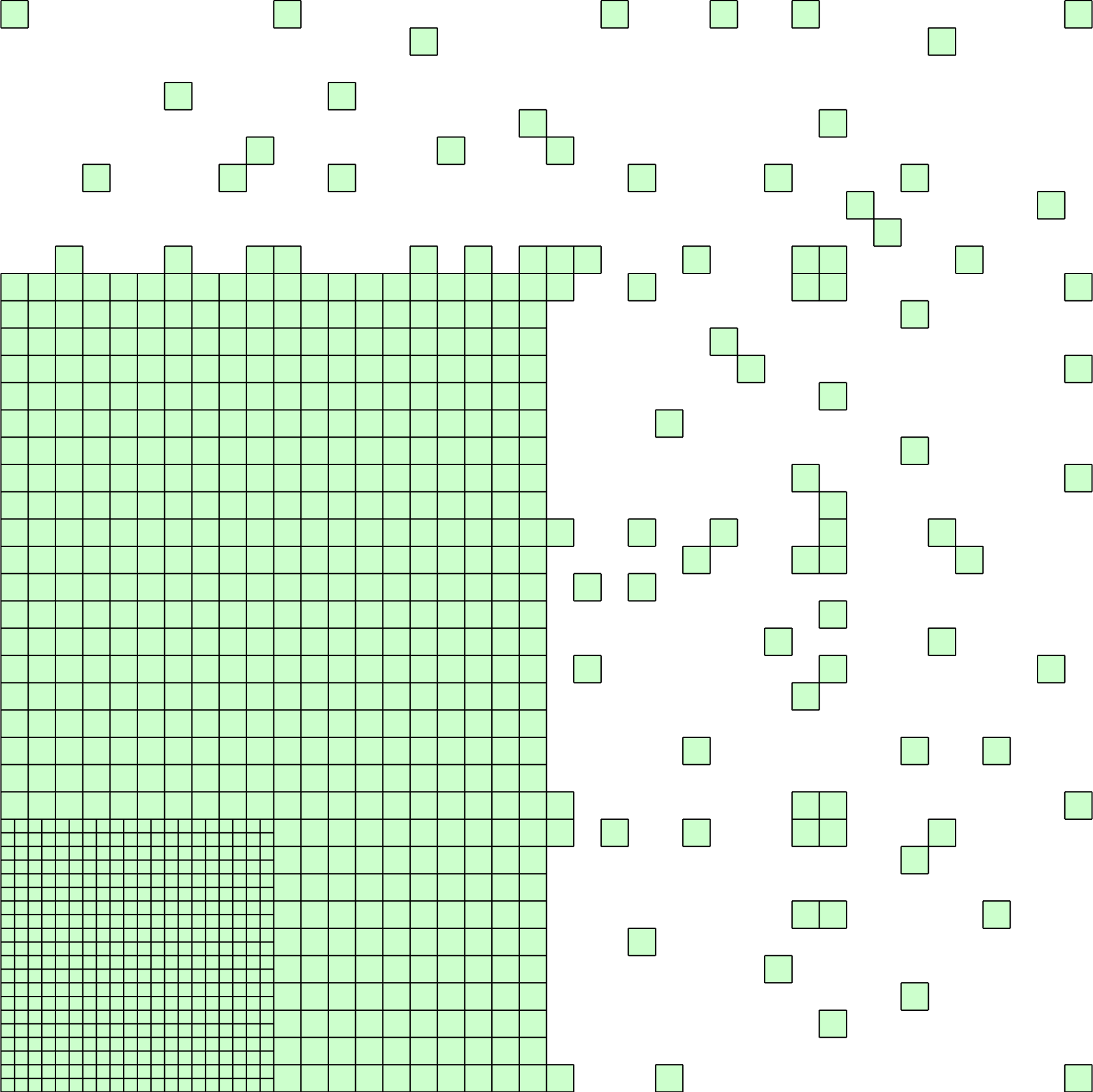};
			
			\nextgroupplot[title={}, xlabel = {}, ylabel={}, ytick={}]
			\addplot graphics [xmin=0, xmax=8.0, ymin=0, ymax=8.0] {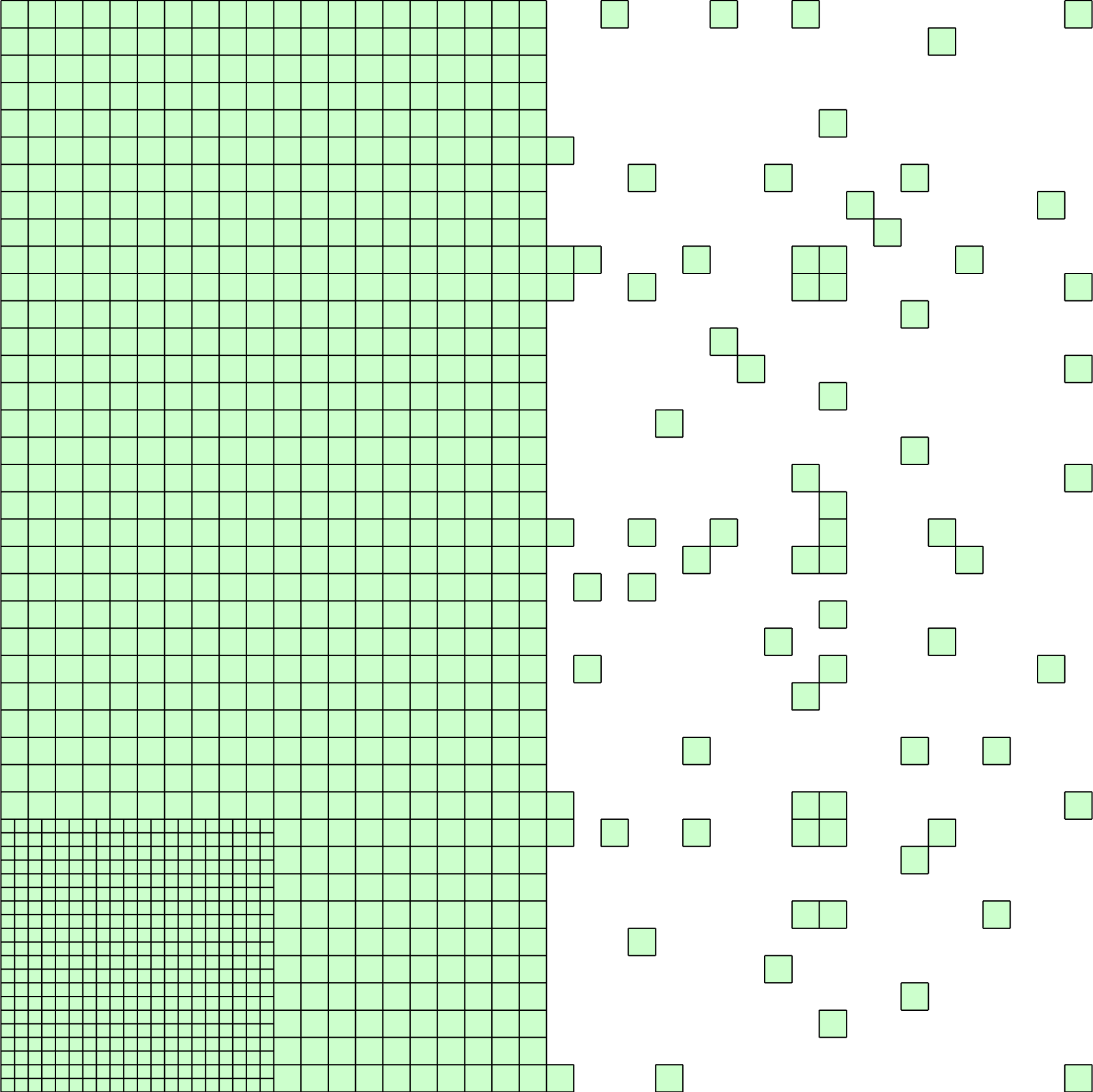};
			
			\nextgroupplot[title={}, xlabel = {}, ylabel={}, ytick={}]
			\addplot graphics [xmin=0, xmax=8.0, ymin=0, ymax=8.0] {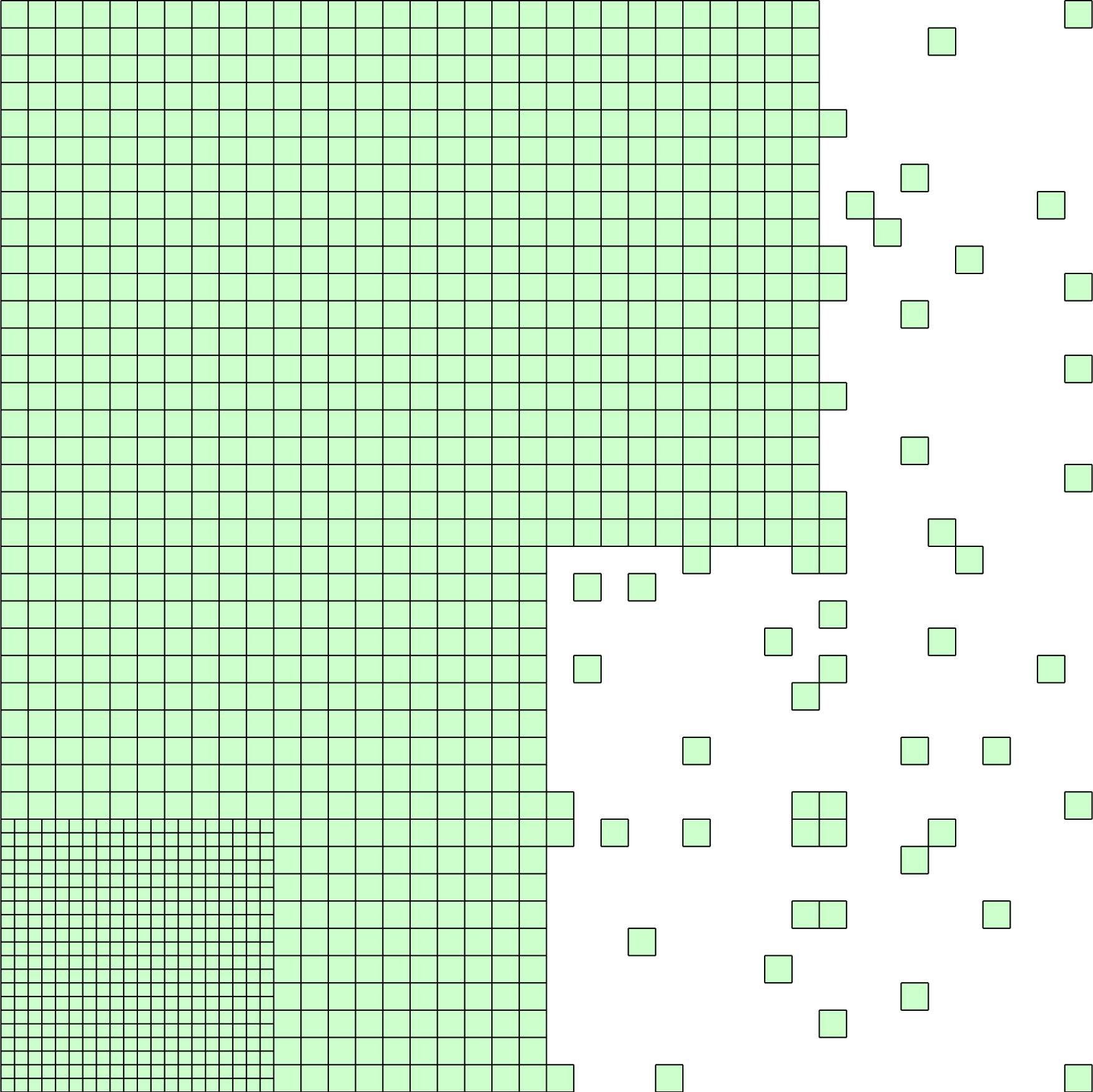};
			
			\nextgroupplot[title={}, ylabel={},xlabel = {},  ytick={}]
			\addplot graphics [xmin=0, xmax=8.0, ymin=0, ymax=8.0] {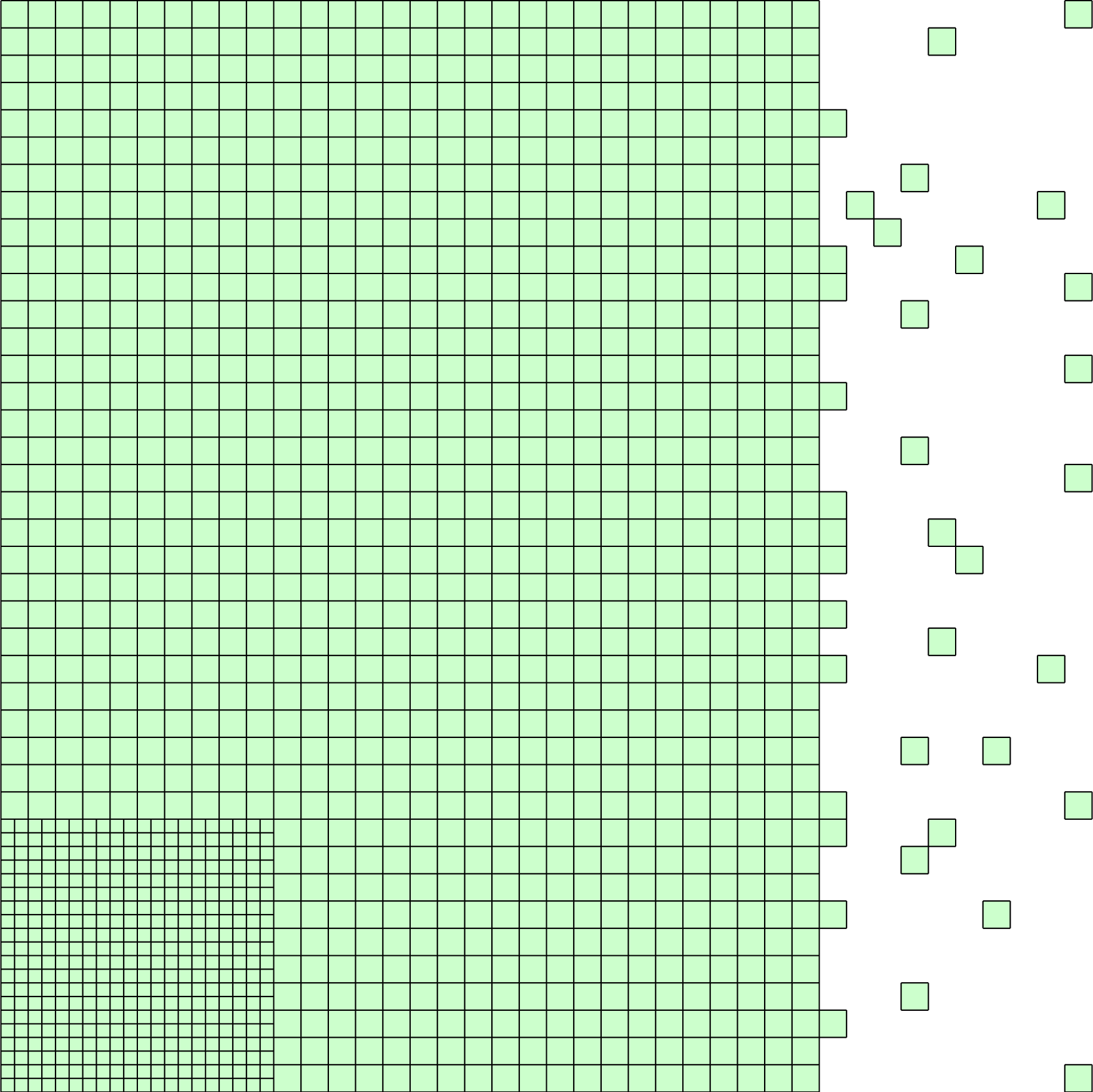};
			
			\nextgroupplot[title={}]
			\addplot graphics [xmin=0, xmax=8.0, ymin=0, ymax=8.0] {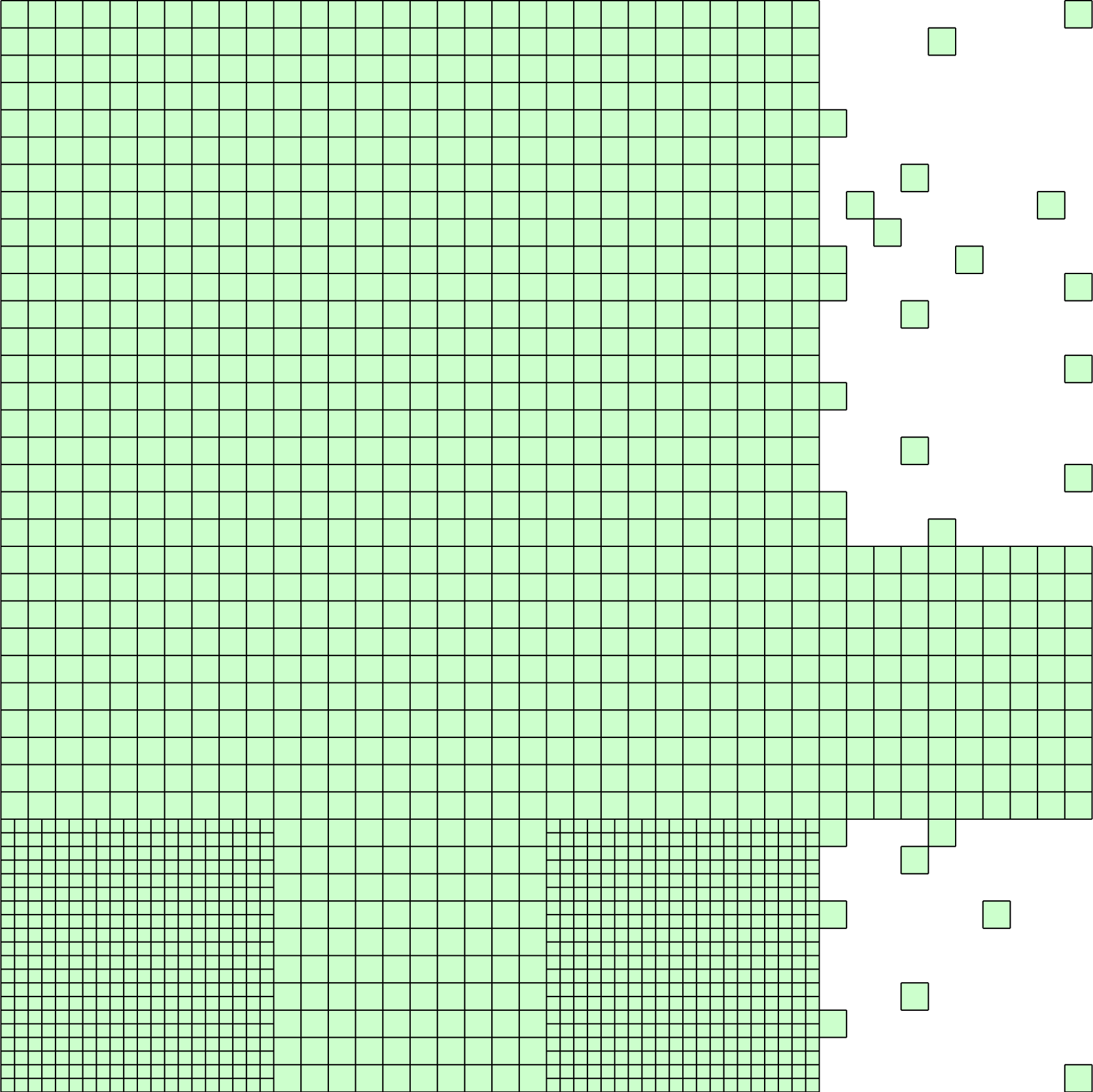};
			
			\nextgroupplot[title={}, ylabel={}, ytick={}]
			\addplot graphics [xmin=0, xmax=8.0, ymin=0, ymax=8.0] {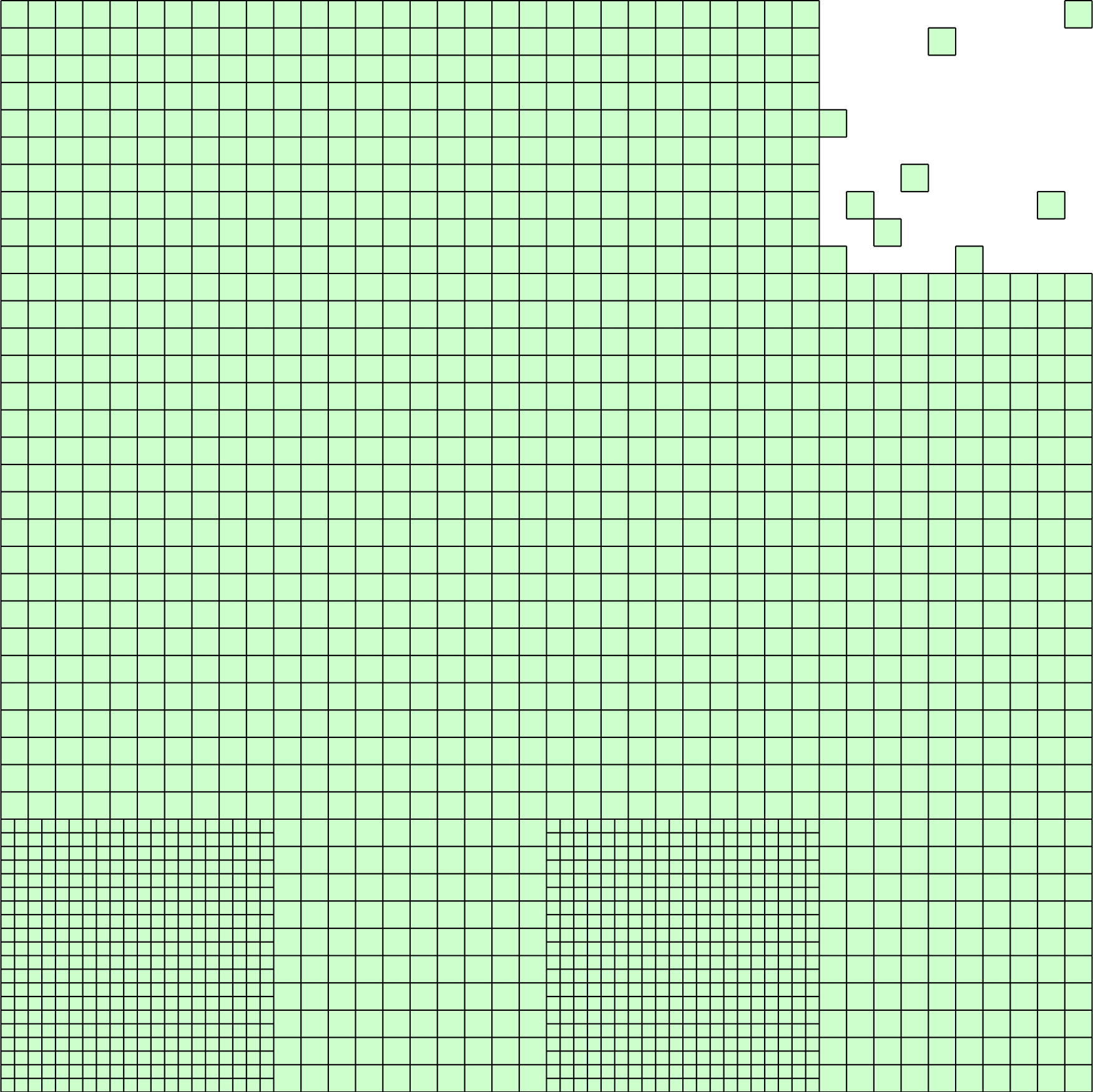};
			
			\nextgroupplot[title={}, ylabel={}, ytick={}]
			\addplot graphics [xmin=0, xmax=8.0, ymin=0, ymax=8.0] {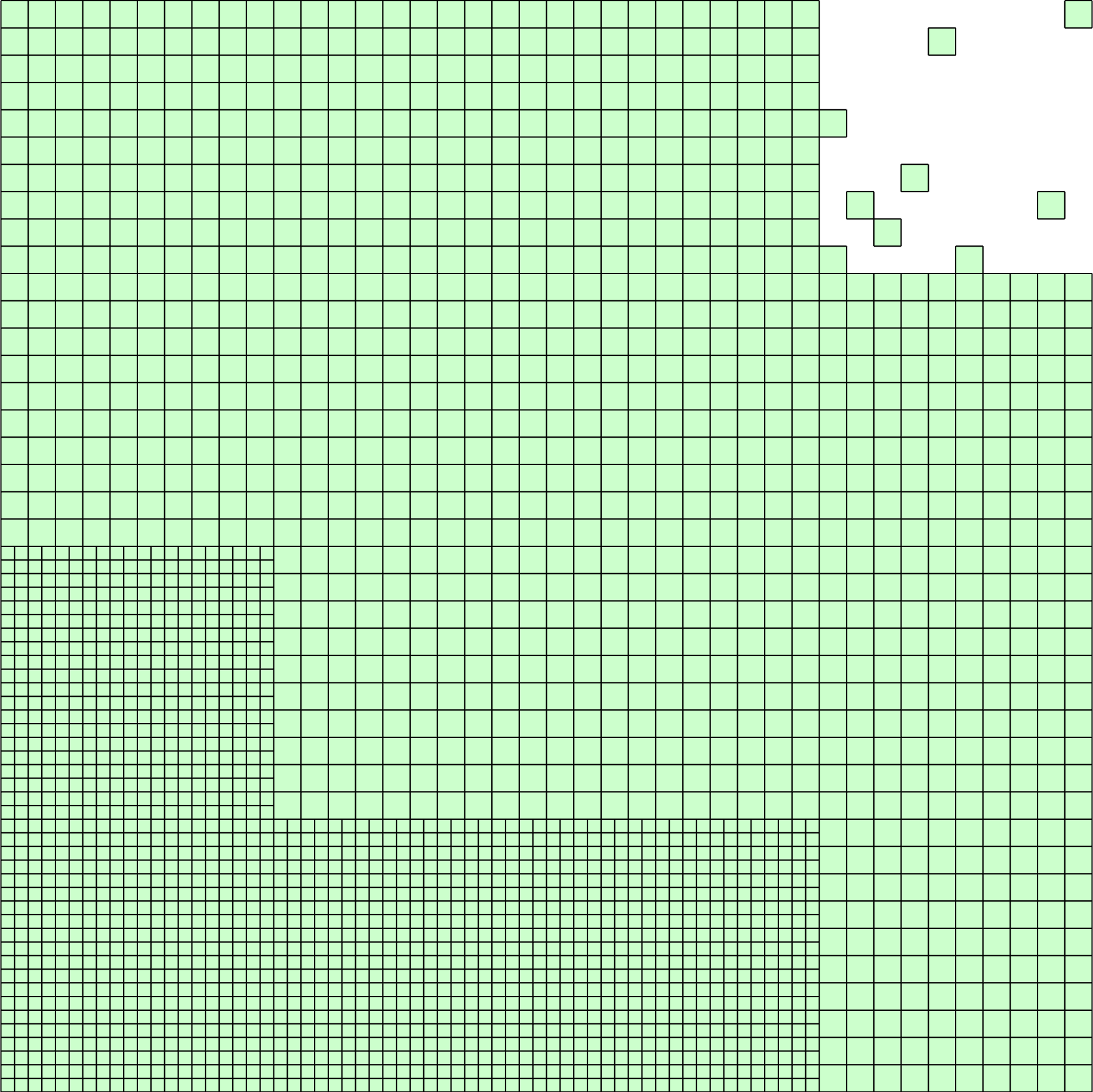};
			
			\nextgroupplot[title={}, ylabel={},  ytick={}]
			\addplot graphics [xmin=0, xmax=8.0, ymin=0, ymax=8.0] {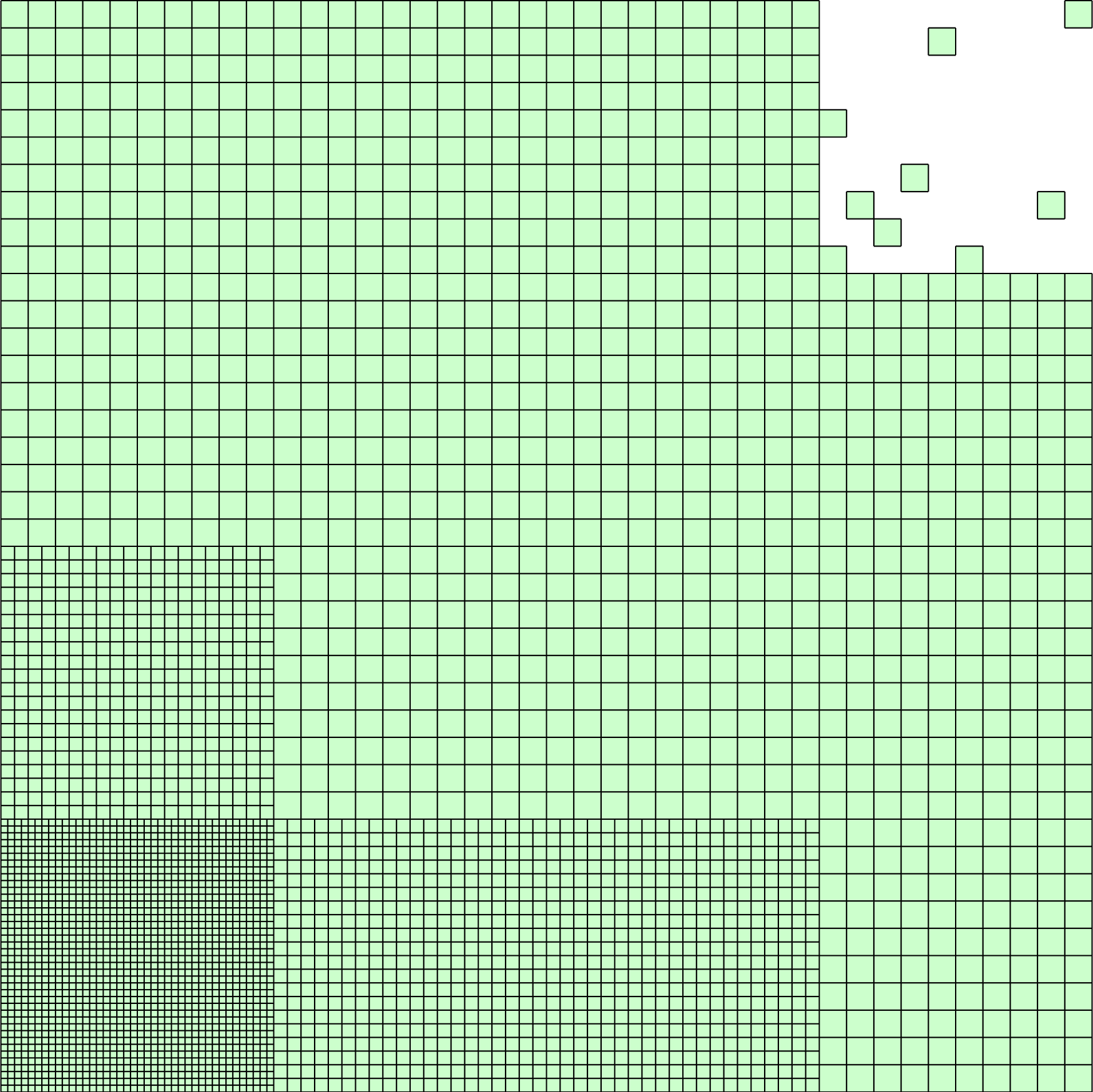};
			
		\end{groupplot}
	\end{tikzpicture}
		\caption{The sequence of discretizations used by the adaptive CG-GL method (\textit{left-to-right}, \textit{top-to-bottom}) for the Poisson problem study in Section~\ref{sec:rslt:poi:local} at $\mu_\mathrm{test}   =  (10, 1, 0.5, 0.5)$.}
		\label{fig:poi_local2}
\end{figure}

\begin{figure}
	\begin{tikzpicture}
\begin{groupplot} [
group style={group size = 3 by 2, horizontal sep = 1.5cm, vertical sep=0.5cm}]

\nextgroupplot[width=0.33\textwidth, grid=both, ymax=1.05, xmin=0, ymode=log, ymin=0.00903082679930243, title={$\mu_\mathrm{test} = (10, 0.1, 0.5, 0.5)$}, ylabel={$e_\mathrm{qoi}$}, xticklabel=\empty]
\addplot [color=magenta, mark=o]
coordinates {
( 0.00000000e+00,  1.00000000e+00)
( 1.00000000e+00,  1.15263017e-02)
( 2.00000000e+00,  1.03847477e-02)
( 3.00000000e+00,  9.50613347e-03)};\label{poi_localness_ej}

\addplot [color=red, mark=*, mark size=1]
coordinates {
( 0.00000000e+00,  9.97979751e-01)
( 1.00000000e+00,  1.06463730e-02)
( 2.00000000e+00,  1.03847322e-02)};\label{poi_localness_ejesti}

\nextgroupplot[width=0.33\textwidth, grid=both, ymax=1.05, xmin=0, ymode=log, ymin=0.0022365888213299197, title={$\mu_\mathrm{test} = (10, 0.5, 0.5, 0.5)$}, xticklabel=\empty]
\addplot [color=magenta, mark=o]
coordinates {
( 0.00000000e+00,  1.00000000e+00)
( 1.00000000e+00,  5.04627066e-02)
( 2.00000000e+00,  3.44890462e-02)
( 3.00000000e+00,  1.25959989e-02)
( 4.00000000e+00,  2.41025758e-03)
( 5.00000000e+00,  2.35766344e-03)
( 6.00000000e+00,  2.35430402e-03)};\label{poi_localness_ej}

\addplot [color=red, mark=*, mark size=1]
coordinates {
( 0.00000000e+00,  9.99944047e-01)
( 1.00000000e+00,  5.04592762e-02)
( 2.00000000e+00,  3.44856672e-02)
( 3.00000000e+00,  1.25926635e-02)
( 4.00000000e+00,  2.40689894e-03)
( 5.00000000e+00,  2.35766161e-03)};\label{poi_localness_ejesti}

\nextgroupplot[width=0.33\textwidth, ytick={1e-6,1e-4,1e-3,1e-2,1e-1,1}, grid=both, ymax=2.051772565530045, xmin=0, ymode=log, ymin=6.737239083460801e-08, title={$\mu_\mathrm{test} = (10, 1, 0.5, 0.5)$}, xticklabel=\empty]
\addplot [color=magenta, mark=o]
coordinates {
( 0.00000000e+00,  1.00000000e+00)
( 1.00000000e+00,  2.54996552e-01)
( 2.00000000e+00,  1.76228315e-01)
( 3.00000000e+00,  6.79959846e-02)
( 4.00000000e+00,  1.30552689e-02)
( 5.00000000e+00,  1.13831845e-02)
( 6.00000000e+00,  8.23079732e-03)
( 7.00000000e+00,  8.21886661e-03)
( 8.00000000e+00,  8.18675092e-03)
( 9.00000000e+00,  7.04319705e-03)
( 1.00000000e+01,  6.89624681e-03)
( 1.10000000e+01,  4.71917300e-03)
( 1.20000000e+01,  2.64033918e-03)
( 1.30000000e+01,  5.32750680e-05)
( 1.40000000e+01,  3.97424214e-05)
( 1.50000000e+01,  3.97394017e-05)
( 1.60000000e+01,  8.27405409e-06)
( 1.70000000e+01,  1.22659604e-06)
( 1.80000000e+01,  1.06916989e-06)
( 1.90000000e+01,  9.11949003e-07)
( 2.00000000e+01,  7.09183061e-08)};\label{poi_localness_ej}

\addplot [color=red, mark=*, mark size=1]
coordinates {
( 0.00000000e+00,  9.99987789e-01)
( 1.00000000e+00,  2.54996064e-01)
( 2.00000000e+00,  1.76227749e-01)
( 3.00000000e+00,  6.79952978e-02)
( 4.00000000e+00,  1.30544657e-02)
( 5.00000000e+00,  1.13823820e-02)
( 6.00000000e+00,  8.23032436e-03)
( 7.00000000e+00,  8.21920382e-03)
( 8.00000000e+00,  8.18701867e-03)
( 9.00000000e+00,  7.04337416e-03)
( 1.00000000e+01,  6.89651167e-03)
( 1.10000000e+01,  4.71939152e-03)
( 1.20000000e+01,  2.64075836e-03)
( 1.30000000e+01,  5.38308755e-05)
( 1.40000000e+01,  4.02330112e-05)
( 1.50000000e+01,  4.02572444e-05)
( 1.60000000e+01,  8.82758438e-06)
( 1.70000000e+01,  1.75233425e-06)
( 1.80000000e+01,  1.33171762e-06)
( 1.90000000e+01,  9.49343036e-07)};\label{poi_localness_ejesti}

\nextgroupplot[width=0.33\textwidth, grid=both, xmin=0, ymode=log, xlabel=$r$, ylabel={Cost}]

\addplot [color=blue, mark=x]
coordinates {
( 0.00000000e+00,  1.06332529e-02)
( 1.00000000e+00,  3.16916924e-02)
( 2.00000000e+00,  6.75626777e-02)
( 3.00000000e+00,  2.01500600e-01)};\label{poi_localness_c}

\nextgroupplot[width=0.33\textwidth, grid=both, ymin=1e-2, xmin=0, ymode=log, xlabel=$r$]

\addplot [color=blue, mark=x]
coordinates {
( 0.00000000e+00,  1.53496776e-02)
( 1.00000000e+00,  6.47689308e-02)
( 2.00000000e+00,  7.80977093e-02)
( 3.00000000e+00,  9.26441605e-02)
( 4.00000000e+00,  1.10401407e-01)
( 5.00000000e+00,  1.47039523e-01)
( 6.00000000e+00,  2.20080115e-01)};

\nextgroupplot[width=0.33\textwidth, grid=both, xmin=0, ymode=log, xlabel=$r$]

\addplot [color=blue, mark=x]
coordinates {
( 0.00000000e+00,  3.85179404e-03)
( 1.00000000e+00,  1.77142408e-02)
( 2.00000000e+00,  3.73428896e-02)
( 3.00000000e+00,  6.27908084e-02)
( 4.00000000e+00,  9.41169394e-02)
( 5.00000000e+00,  1.33282246e-01)
( 6.00000000e+00,  1.77407327e-01)
( 7.00000000e+00,  2.45783276e-01)
( 8.00000000e+00,  3.16585034e-01)
( 9.00000000e+00,  3.99696201e-01)
( 1.00000000e+01,  4.85544780e-01)
( 1.10000000e+01,  5.79808066e-01)
( 1.20000000e+01,  6.78378732e-01)
( 1.30000000e+01,  7.85519928e-01)
( 1.40000000e+01,  8.95093067e-01)
( 1.50000000e+01,  1.02870768e+00)
( 1.60000000e+01,  1.17420567e+00)
( 1.70000000e+01,  1.31308078e+00)
( 1.80000000e+01,  1.47377582e+00)
( 1.90000000e+01,  1.65855682e+00)
( 2.00000000e+01,  1.95406911e+00)};

\end{groupplot}
\end{tikzpicture}
	\caption{The true QoI error (\ref{poi_localness_ej}), QoI error estimate (\ref{poi_localness_ejesti}), and cumulative (non-dimensional) computational cost (\ref{poi_localness_c}) at each test parameter for the Poisson problem study in Section~\ref{sec:rslt:poi:local}. The non-dimensional computational cost is defined as the cumulative cost (wall time) of the adaptive CG-GL method normalized by the cost to compute the reference solution.}
	\label{fig:poi_error_estimation_and_cost}
\end{figure}
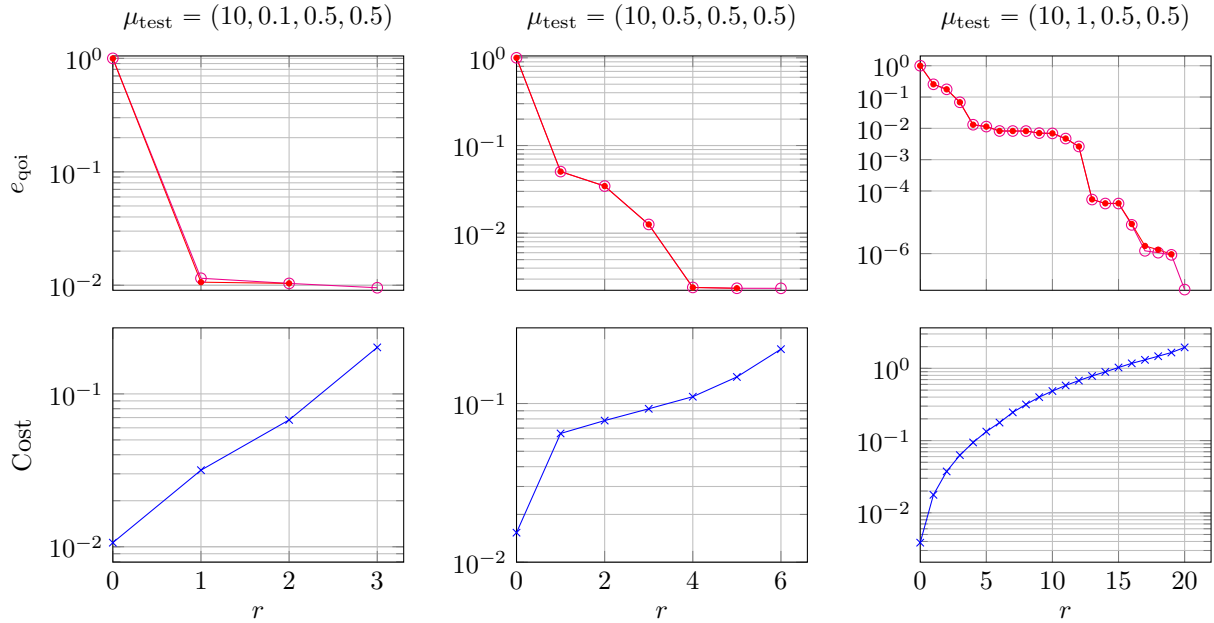

Next we study the impact of feature locality on the performance of the adaptive CG-GL framework by considering variations to $\mu_2$ (Figure~\ref{fig:poi_mu1}). The training set remains the same as the previous section, and the new testing set is
\begin{equation}
\Dcal_\mathrm{test}=\{(10, 0.1, 0.5, 0.5),  (10, 0.5, 0.5, 0.5), (10, 1, 0.5, 0.5) \}.
\end{equation}
Because $\sigma = 0.1$ leads to a steep feature in the PDE solution, we require increase resolution. As such,  the conservation law is discretized using $25600$ elements of polynomial degree $p=2$, and split the domain into $16$ patches, each discretized using $100$  elements of polynomial degree $p=2$ (Figure~\ref{fig:meshes_EQP}). Similar to the previous sections, the global basis and empirical quadrature weights for the ROM and CG-GL adaptions are trained on the same set ($\Dcal_{\mathrm{train}}$) with tolerances $\delta_\mathrm{dv} = \delta_\mathrm{rp} = 10^{-8}$. This leads to a sample mesh for the ROM (of the entire domain) with $10$ elements and sample meshes for the CG-GL discretization with $147$ elements ($9.2$ elements per patch, on average) (Figure~\ref{fig:meshes_EQP}). Figure~\ref{fig:poi_local0}-\ref{fig:poi_local2} show the sequence of grids used by the adaptive CG-GL method for $\mu_2 = 1,~0.5,~0.1$, respectively. As the locality of the feature increases, fewer patches are required to reach small errors. Finally, we see the adaptive CG-GL framework leads to a user-specified tradeoff between accuracy and computational cost that leverages information from global basis functions (Figure~\ref{fig:poi_error_estimation_and_cost}).

\subsection{Incompressible, viscous flow over backward-facing step}
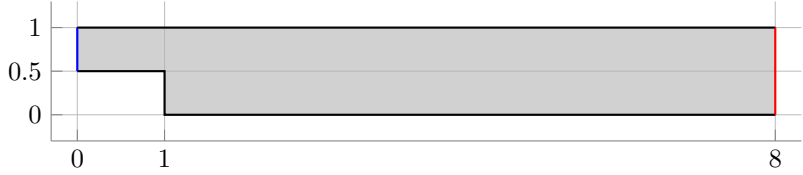
\begin{figure}
\centering
\begin{tikzpicture}
\begin{axis}[
axis equal image,
axis line style={gray},
axis x line*=bottom,
axis y line*=left,
width=0.7\textwidth,
xtick={0, 1, 8},
ytick={0, 0.5, 1},
grid=major,
ymax=1.3,
xmax=8.3,
xmin=-0.3,
ymin=-0.3]
\addplot [opacity=0.6, fill=black!30!white, opacity=0.6, forget plot]
coordinates {
( 0.00000000e+00,  5.00000000e-01)
( 1.00000000e+00,  5.00000000e-01)
( 1.00000000e+00,  0.00000000e+00)
( 8.00000000e+00,  0.00000000e+00)
( 8.00000000e+00,  1.00000000e+00)
( 0.00000000e+00,  1.00000000e+00)
( 0.00000000e+00,  5.00000000e-01)};

\addplot [thick, color=black]
coordinates {
( 0.00000000e+00,  5.00000000e-01)
( 1.00000000e+00,  5.00000000e-01)
( 1.00000000e+00,  0.00000000e+00)
( 8.00000000e+00,  0.00000000e+00)};\label{line:cyl0:wall}

\addplot [thick, color=black, forget plot]
coordinates {
( 0.00000000e+00,  1.00000000e+00)
( 8.00000000e+00,  1.00000000e+00)};

\addplot [thick, color=red]
coordinates {
( 8.00000000e+00,  0.00000000e+00)
( 8.00000000e+00,  1.00000000e+00)};\label{line:cyl0:out}

\addplot [thick, color=blue]
coordinates {
( 0.00000000e+00,  5.00000000e-01)
( 0.00000000e+00,  1.00000000e+00)};\label{line:cyl0:inlet}

\end{axis}
\end{tikzpicture}
\caption{Schematic of the backward-facing step with boundaries: $\Gamma_\mathrm{in}$ (\ref{line:cyl0:inlet}), $\Gamma_\mathrm{out}$ (\ref{line:cyl0:out}) and $\Gamma_\mathrm{wall}$ (\ref{line:cyl0:wall}).}
\label{fig:geons}
\end{figure}

\begin{figure}
	\centering
	\begin{tikzpicture}
		\begin{groupplot}[
			group style={
				group size=1 by 3,
				horizontal sep=0.5cm
			},
			width=0.8\textwidth,
			axis equal image,
			xlabel={$x_1$},
			ylabel={$x_2$},
			xtick = {0.0, 4.0, 8.0},
			ytick = {0.0, 0.5, 1.0},
			xmin=0, xmax=8,
			ymin=0, ymax=1.0
			]
			\nextgroupplot[title={$\mathrm{Re}=10$},xlabel={},xtick=\empty]
			\addplot graphics [xmin=0, xmax=8.0, ymin=0, ymax=1.0] {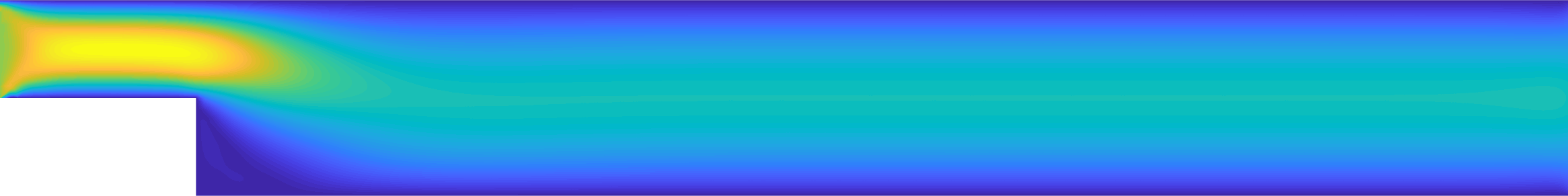};
			
			\nextgroupplot[title={$\mathrm{Re}=80$}, xlabel={},xtick=\empty]
			\addplot graphics [xmin=0, xmax=8.0, ymin=0, ymax=1.0] {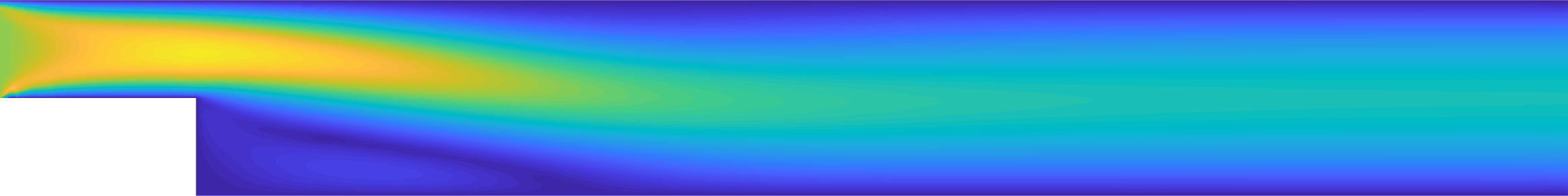};
			
			\nextgroupplot[title={$\mathrm{Re}=150$}]
			\addplot graphics [xmin=0, xmax=8.0, ymin=0, ymax=1.0] {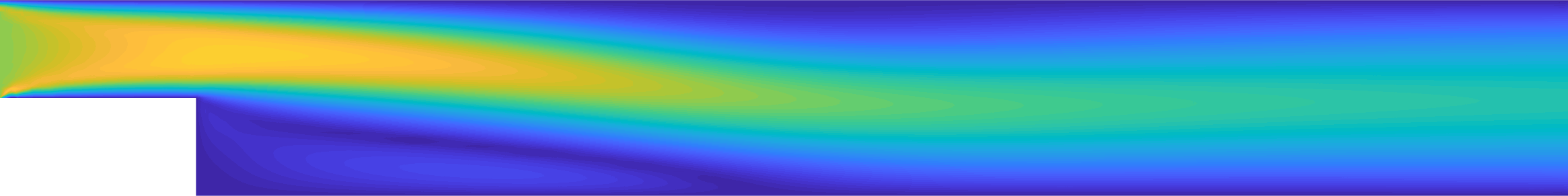};
		\end{groupplot}
	\end{tikzpicture}
	\colorbarMatlabParula{0}{0.3669}{0.7339}{1.1008}{1.4677}
	\caption{Velocity magnitude of viscous, incompressible flow over backward-facing step for
		three Reynolds numbers.}
	\label{fig:ns_vmag_3Re}
\end{figure}

Finally, we consider viscous, incompressible flow over a backward-facing step ($\Omega\subset\Rbb^2$) shown in Figure~\ref{fig:geons}. The governing equations are the steady, incompressible Navier-Stokes equations,
\begin{equation}
	(v\cdot\nabla)v - \nu\nabla^2v+\frac{1}{\rho}\nabla P=0,\quad\nabla\cdot v = 0\quad\text{in }\Omega,
	\label{eqn:NS}
\end{equation}
where $\nu\in\Rbb_{>0}$ is the kinematic viscosity of the fluid, $\rho\in\Rbb_{>0}$ is the density of the fluid, and  $v:\Omega\rightarrow\Rbb^2$ and $P:\Omega\rightarrow\Rbb$ are the velocity and pressure, respectively, of the fluid implicitly defined as the solution of \eqref{eqn:NS}. We set the boundary conditions for the boundaries shown in Figure~\ref{fig:geons} as,
\begin{equation}
	v = v_\mathrm{in}\quad\text{on }\Gamma_\mathrm{in},\quad \sigma\cdot n = 0\quad\text{on }\Gamma_\mathrm{out},\quad v = 0\quad\text{on }\Gamma_\mathrm{wall}
\end{equation}
where $v_\mathrm{in}=(1,0)$ is the inflow velocity, $n:\partial\Omega\rightarrow\Rbb^2$ is the outward unit normal to the boundary of
the domain. The rate of strain, $\epsilon:\Omega\rightarrow\Rbb^{2\times2}$ and stress, $\sigma:\Omega\rightarrow\Rbb^{2\times2}$, tensors are defined as 
\begin{equation}
	\epsilon = \frac{1}{2}(\nabla v+\nabla v ^T),\quad\sigma=2\rho\nu\epsilon + PI_2,
\end{equation}
where $I_2$ is the $2\times2$ identity matrix. In this case, we take $\rho=1$.  The Reynolds number, $\mathrm{Re}\in\Rbb_{>0}$, defined as $Re=\frac{1}{2\nu}$ (characteristic length based on the inlet, $L=0.5$), is used to parametrize the problem. The solution has a zone of re-circulation (Figure~\ref{fig:ns_vmag_3Re}) and reduction in the static pressure in that region, and the position of flow reattachment moves further from the wall  as the Reynolds number increases. The quantity of interest we choose for this problem is the integral of the pressure over the entire domain,
\begin{equation}
	q((v,P); \mu) \coloneqq \int_\Omega P \, dV.
\end{equation} 

\begin{figure}
	\centering
	\begin{tikzpicture}
		\begin{groupplot}[
			group style={
				group size=1 by 3,
				horizontal sep=0.7cm
			},
			xmajorgrids=true,
			ymajorgrids=true,
			width=0.8\textwidth,
			axis equal image,
			xlabel={$x_1$},
			ylabel={$x_2$},
			xtick = {0, 4.5, 8},
			ytick = {0, 0.5,1},
			xmin=0, xmax=8,
			ymin=0, ymax=1,
			axis on top,
			grid style={line width=0.5pt, draw=black}
			]
			
			\nextgroupplot[title={},xlabel={}]
			\addplot graphics [xmin=0, xmax=8.0, ymin=0, ymax=1] {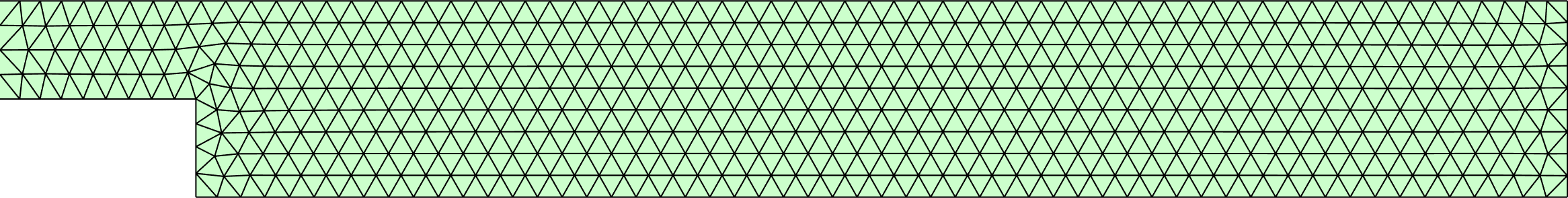};
			
			\nextgroupplot[title={}, xlabel={}]
			\addplot graphics [xmin=0, xmax=8.0, ymin=0, ymax=1] {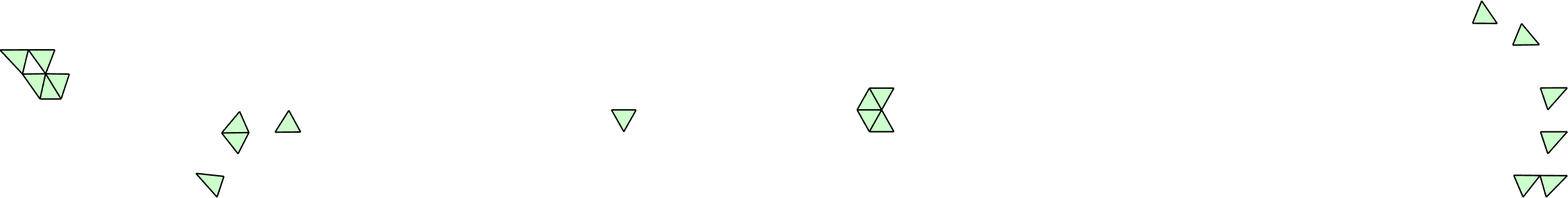};
			
			\nextgroupplot[title={}]
			\addplot graphics [xmin=0, xmax=0.35, ymin=0.5, ymax=0.75] {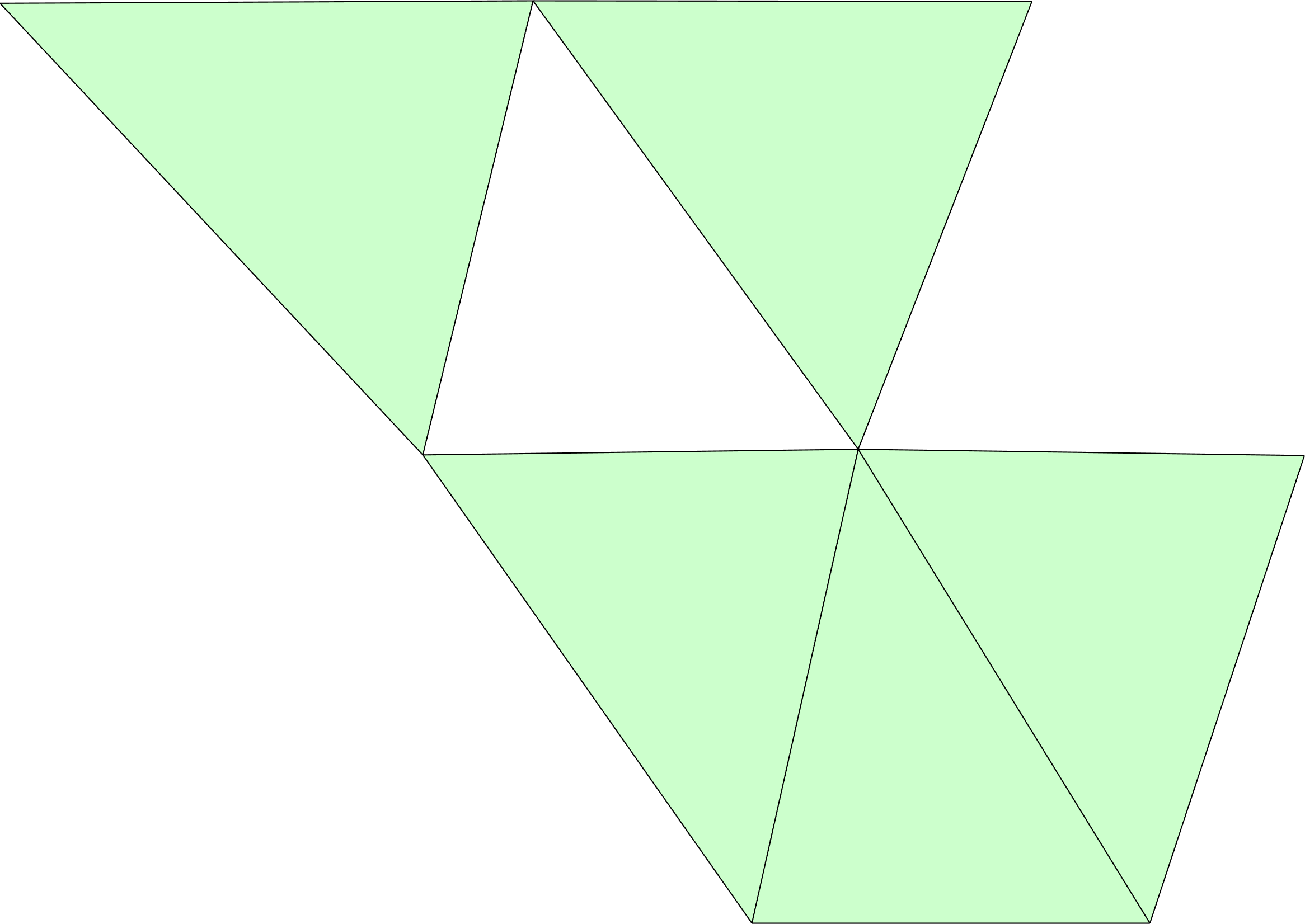};
		\end{groupplot}
	\end{tikzpicture}
	\caption{The mesh for the reference solution for the viscous flow problem (\textit{top}), the samples meshes for the CG-GL patches (\textit{middle}), and the sample mesh for the ROM (\textit{bottom}).}
	\label{fig:NSmesh}
\end{figure}

The system of conservation laws and quantity of interest are discretized with $\Pcal^2$-$\Pcal^1$ Taylor-Hood elements. Reference solutions are computed on a mesh consisting of 1067 triangular elements (Figure~\ref{fig:NSmesh}, \textit{top}) and the domain is partitioned into four patches for the CG-GL method.  We perform a targeted study to demonstrate the improvement brought by the CG-GL method over a traditional ROM solution with limited training. In particular, we consider the training set with only two parameters $\Dcal_\mathrm{train}=\{1, 130\}$ that is used train the ROM and CG-GL methods, which produces two global basis functions. Empirical quadrature weights are computed based on this training set over the entire domain for the ROM (Figure~\ref{fig:NSmesh}, \textit{middle}) and over each patch for the CG-GL method (Figure~\ref{fig:NSmesh}, \textit{bottom}) using $\Xi = \Dcal_\mathrm{test}$ and $\delta_\mathrm{dv} = \delta_\mathrm{rp} = 10^{-8}$. The predictive accuracy of both methods will be tested on the set $\Dcal_\mathrm{test}$ consisting of $50$ Reynolds numbers from the range $(1, 130)$.

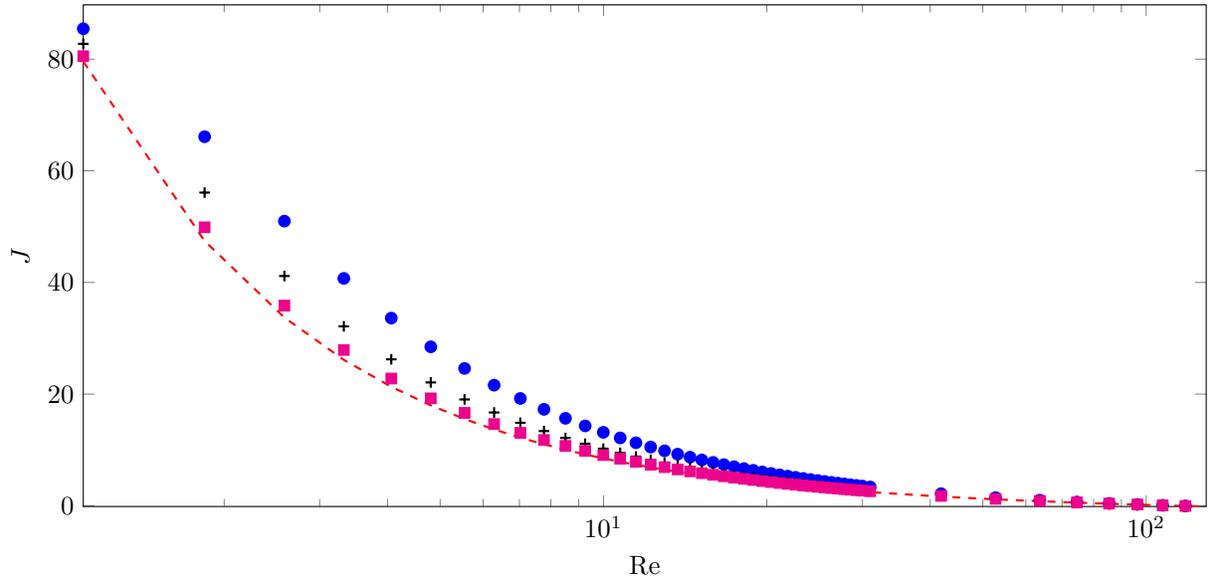
\begin{figure}
	\centering
	\begin{tikzpicture}
\begin{axis}[
width=1\textwidth,
height=0.5\textwidth,
xlabel=$\mathrm{Re}$,
ymax=89.71388853320312,
xmax=129.0,
ylabel=$J$,
xmin=1.1,
ymin=-0.14193481698889276,
xmode=log]
\addplot [mark options={solid, thick}, mark=*, mark size=2, only marks, color=blue]
coordinates {
( 1.10000000e+00,  8.54417986e+01)
( 1.84102564e+00,  6.60996071e+01)
( 2.58205128e+00,  5.09791716e+01)
( 3.32307692e+00,  4.07234392e+01)
( 4.06410256e+00,  3.36069253e+01)
( 4.80512821e+00,  2.84611038e+01)
( 5.54615385e+00,  2.45959247e+01)
( 6.28717949e+00,  2.15983420e+01)
( 7.02820513e+00,  1.92114033e+01)
( 7.76923077e+00,  1.72686760e+01)
( 8.51025641e+00,  1.56583110e+01)
( 9.25128205e+00,  1.43026636e+01)
( 9.99230769e+00,  1.31462901e+01)
( 1.07333333e+01,  1.21486144e+01)
( 1.14743590e+01,  1.12792963e+01)
( 1.22153846e+01,  1.05152167e+01)
( 1.29564103e+01,  9.83846297e+00)
( 1.36974359e+01,  9.23494858e+00)
( 1.44384615e+01,  8.69344850e+00)
( 1.51794872e+01,  8.20491058e+00)
( 1.59205128e+01,  7.76195611e+00)
( 1.66615385e+01,  7.35851150e+00)
( 1.74025641e+01,  6.98953291e+00)
( 1.81435897e+01,  6.65079748e+00)
( 1.88846154e+01,  6.33874330e+00)
( 1.96256410e+01,  6.05034528e+00)
( 2.03666667e+01,  5.78301793e+00)
( 2.11076923e+01,  5.53453847e+00)
( 2.18487179e+01,  5.30298552e+00)
( 2.25897436e+01,  5.08668977e+00)
( 2.33307692e+01,  4.88419407e+00)
( 2.40717949e+01,  4.69422075e+00)
( 2.48128205e+01,  4.51564491e+00)
( 2.55538462e+01,  4.34747228e+00)
( 2.62948718e+01,  4.18882079e+00)
( 2.70358974e+01,  4.03890523e+00)
( 2.77769231e+01,  3.89702429e+00)
( 2.85179487e+01,  3.76254971e+00)
( 2.92589744e+01,  3.63491696e+00)
( 3.00000000e+01,  3.51361744e+00)
( 3.10000000e+01,  3.35915955e+00)
( 4.18888889e+01,  2.15629951e+00)
( 5.27777778e+01,  1.45128976e+00)
( 6.36666667e+01,  9.88108496e-01)
( 7.45555556e+01,  6.60565171e-01)
( 8.54444444e+01,  4.16696383e-01)
( 9.63333333e+01,  2.28073533e-01)
( 1.07222222e+02,  7.78351214e-02)
( 1.18111111e+02, -4.46527454e-02)
( 1.29000000e+02, -1.46428280e-01)};\label{j_ns_0}

\addplot [mark options={solid, thick}, mark=+, mark size=2, only marks, color=black]
coordinates {
( 1.10000000e+00,  8.27389018e+01)
( 1.84102564e+00,  5.60959688e+01)
( 2.58205128e+00,  4.11449700e+01)
( 3.32307692e+00,  3.21339351e+01)
( 4.06410256e+00,  2.62257076e+01)
( 4.80512821e+00,  2.20866261e+01)
( 5.54615385e+00,  1.90374274e+01)
( 6.28717949e+00,  1.67025027e+01)
( 7.02820513e+00,  1.48592658e+01)
( 7.76923077e+00,  1.33681311e+01)
( 8.51025641e+00,  1.21374063e+01)
( 9.25128205e+00,  1.11044993e+01)
( 9.99230769e+00,  1.02252829e+01)
( 1.07333333e+01,  9.46778442e+00)
( 1.14743590e+01,  8.80828746e+00)
( 1.22153846e+01,  8.22883851e+00)
( 1.29564103e+01,  7.71560266e+00)
( 1.36974359e+01,  7.25775045e+00)
( 1.44384615e+01,  6.84668630e+00)
( 1.51794872e+01,  6.47550246e+00)
( 1.59205128e+01,  6.13858528e+00)
( 1.66615385e+01,  5.83132664e+00)
( 1.74025641e+01,  5.54990904e+00)
( 1.81435897e+01,  5.29114350e+00)
( 1.88846154e+01,  5.05234563e+00)
( 1.96256410e+01,  4.83123995e+00)
( 2.03666667e+01,  4.62588501e+00)
( 2.11076923e+01,  4.43461443e+00)
( 2.18487179e+01,  4.25598995e+00)
( 2.25897436e+01,  4.08876376e+00)
( 2.33307692e+01,  3.93184802e+00)
( 2.40717949e+01,  3.78429004e+00)
( 2.48128205e+01,  3.64525199e+00)
( 2.55538462e+01,  3.51399405e+00)
( 2.62948718e+01,  3.38986049e+00)
( 2.70358974e+01,  3.27226809e+00)
( 2.77769231e+01,  3.16069635e+00)
( 2.85179487e+01,  3.05467926e+00)
( 2.92589744e+01,  2.95379836e+00)
( 3.00000000e+01,  2.85767682e+00)
( 3.10000000e+01,  2.73491038e+00)
( 4.18888889e+01,  1.76117035e+00)
( 5.27777778e+01,  1.17163327e+00)
( 6.36666667e+01,  7.76277647e-01)
( 7.45555556e+01,  4.95123219e-01)
( 8.54444444e+01,  2.87869993e-01)
( 9.63333333e+01,  1.31615660e-01)
( 1.07222222e+02,  1.22605895e-02)
( 1.18111111e+02, -7.93536670e-02)
( 1.29000000e+02, -1.49405071e-01)};\label{j_ns_1}

\addplot [mark options={solid, thin}, mark=square*, mark size=2, only marks, color=magenta]
coordinates {
( 1.10000000e+00,  8.05265695e+01)
( 1.84102564e+00,  4.98777755e+01)
( 2.58205128e+00,  3.58508811e+01)
( 3.32307692e+00,  2.78939942e+01)
( 4.06410256e+00,  2.27869422e+01)
( 4.80512821e+00,  1.92367493e+01)
( 5.54615385e+00,  1.66275814e+01)
( 6.28717949e+00,  1.46297927e+01)
( 7.02820513e+00,  1.30512961e+01)
( 7.76923077e+00,  1.17727023e+01)
( 8.51025641e+00,  1.07159471e+01)
( 9.25128205e+00,  9.82785658e+00)
( 9.99230769e+00,  9.07096959e+00)
( 1.07333333e+01,  8.41813897e+00)
( 1.14743590e+01,  7.84921358e+00)
( 1.22153846e+01,  7.34892301e+00)
( 1.29564103e+01,  6.90548554e+00)
( 1.36974359e+01,  6.50966654e+00)
( 1.44384615e+01,  6.15412602e+00)
( 1.51794872e+01,  5.83295689e+00)
( 1.59205128e+01,  5.54135193e+00)
( 1.66615385e+01,  5.27535937e+00)
( 1.74025641e+01,  5.03170096e+00)
( 1.81435897e+01,  4.80763449e+00)
( 1.88846154e+01,  4.60084859e+00)
( 1.96256410e+01,  4.40938146e+00)
( 2.03666667e+01,  4.23155724e+00)
( 2.11076923e+01,  4.06593579e+00)
( 2.18487179e+01,  3.91127271e+00)
( 2.25897436e+01,  3.76648710e+00)
( 2.33307692e+01,  3.63063590e+00)
( 2.40717949e+01,  3.50289225e+00)
( 2.48128205e+01,  3.38252831e+00)
( 2.55538462e+01,  3.26890090e+00)
( 2.62948718e+01,  3.16143949e+00)
( 2.70358974e+01,  3.05963629e+00)
( 2.77769231e+01,  2.96303779e+00)
( 2.85179487e+01,  2.87123779e+00)
( 2.92589744e+01,  2.78387133e+00)
( 3.00000000e+01,  2.70060961e+00)
( 3.10000000e+01,  2.59423543e+00)
( 4.18888889e+01,  1.74693383e+00)
( 5.27777778e+01,  1.22438197e+00)
( 6.36666667e+01,  8.61962699e-01)
( 7.45555556e+01,  5.91833033e-01)
( 8.54444444e+01,  3.80863312e-01)
( 9.63333333e+01,  2.10787102e-01)
( 1.07222222e+02,  7.05721526e-02)
( 1.18111111e+02, -4.69164289e-02)
( 1.29000000e+02, -1.46536822e-01)};\label{j_ns_2}

\addplot [dashed, thick, color=red]
coordinates {
( 1.10000000e+00,  7.95912641e+01)
( 1.84102564e+00,  4.74211414e+01)
( 2.58205128e+00,  3.37189921e+01)
( 3.32307692e+00,  2.61295160e+01)
( 4.06410256e+00,  2.13086304e+01)
( 4.80512821e+00,  1.79751268e+01)
( 5.54615385e+00,  1.55325633e+01)
( 6.28717949e+00,  1.36657214e+01)
( 7.02820513e+00,  1.21923495e+01)
( 7.76923077e+00,  1.09997490e+01)
( 8.51025641e+00,  1.00144905e+01)
( 9.25128205e+00,  9.18668465e+00)
( 9.99230769e+00,  8.48125120e+00)
( 1.07333333e+01,  7.87280567e+00)
( 1.14743590e+01,  7.34252707e+00)
( 1.22153846e+01,  6.87616616e+00)
( 1.29564103e+01,  6.46273723e+00)
( 1.36974359e+01,  6.09363471e+00)
( 1.44384615e+01,  5.76202194e+00)
( 1.51794872e+01,  5.46239901e+00)
( 1.59205128e+01,  5.19029143e+00)
( 1.66615385e+01,  4.94202184e+00)
( 1.74025641e+01,  4.71454016e+00)
( 1.81435897e+01,  4.50529537e+00)
( 1.88846154e+01,  4.31213748e+00)
( 1.96256410e+01,  4.13324176e+00)
( 2.03666667e+01,  3.96704946e+00)
( 2.11076923e+01,  3.81222115e+00)
( 2.18487179e+01,  3.66759941e+00)
( 2.25897436e+01,  3.53217898e+00)
( 2.33307692e+01,  3.40508255e+00)
( 2.40717949e+01,  3.28554108e+00)
( 2.48128205e+01,  3.17287755e+00)
( 2.55538462e+01,  3.06649325e+00)
( 2.62948718e+01,  2.96585812e+00)
( 2.70358974e+01,  2.87049945e+00)
( 2.77769231e+01,  2.77999534e+00)
( 2.85179487e+01,  2.69396785e+00)
( 2.92589744e+01,  2.61207739e+00)
( 3.00000000e+01,  2.53401801e+00)
( 3.10000000e+01,  2.43426665e+00)
( 4.18888889e+01,  1.63875700e+00)
( 5.27777778e+01,  1.14732069e+00)
( 6.36666667e+01,  8.05972622e-01)
( 7.45555556e+01,  5.51025544e-01)
( 8.54444444e+01,  3.51478766e-01)
( 9.63333333e+01,  1.90398290e-01)
( 1.07222222e+02,  5.76440161e-02)
( 1.18111111e+02, -5.33377219e-02)
( 1.29000000e+02, -1.47056066e-01)};\label{ns_hdm}

\end{axis}
\end{tikzpicture}
	\caption{QoI as a function of Reynolds number for the viscous flow problem: reference solution (\ref{ns_hdm}), traditional ROM (\ref{j_ns_0}), and the CG-GL method with one (\ref{j_ns_1}) and two (\ref{j_ns_2}) adaptation iterations. }
	\label{fig:J_NS}
\end{figure}

\begin{figure}
	\centering
	\input{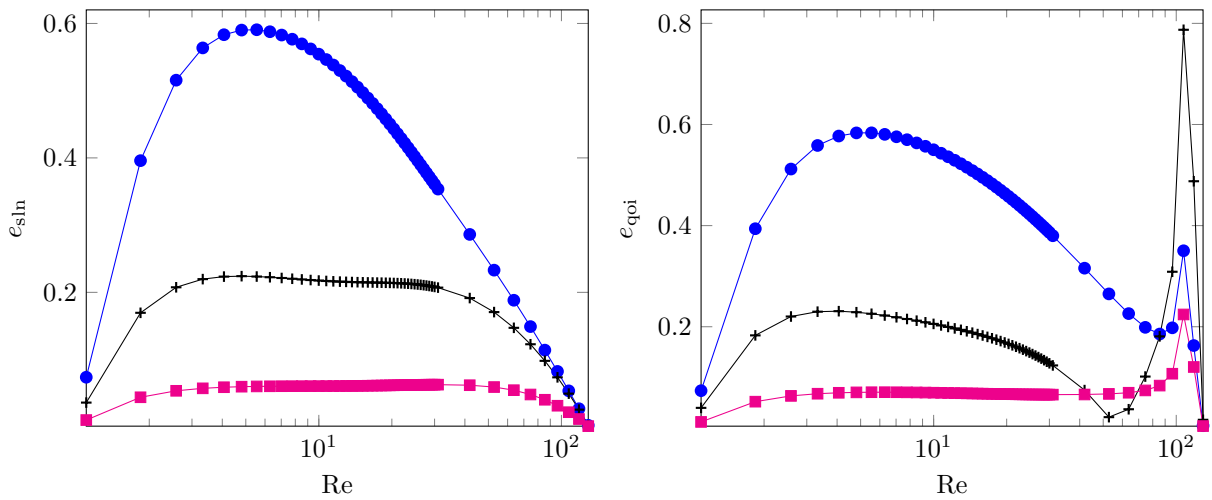}
	\caption{The relative solution (\textit{left}) and QoI (\textit{right}) error across test set $\Dcal_{\mathrm{test}}$ viscous flow problem: traditional ROM (\ref{ns_first_round}) and the CG-GL method with one (\ref{ns_second_round}) and two (\ref{ns_third_round}) adaptation iterations.}
	\label{fig:ns_errors}
\end{figure}

\begin{figure}
	\centering
	\begin{tikzpicture}
		\begin{groupplot}[
			group style={
				group size=1 by 3,
				horizontal sep=0.5cm
			},
			width=0.8\textwidth,
			axis equal image,
			xlabel={$x_1$},
			ylabel={$x_2$},
			xtick = {0.0, 4.0, 8.0},
			ytick = {0.0, 0.5, 1.0},
			xmin=0, xmax=8,
			ymin=0, ymax=1.0
			]
			\nextgroupplot[title={},xlabel={},xtick=\empty]
			\addplot graphics [xmin=0, xmax=8.0, ymin=0, ymax=1.0] {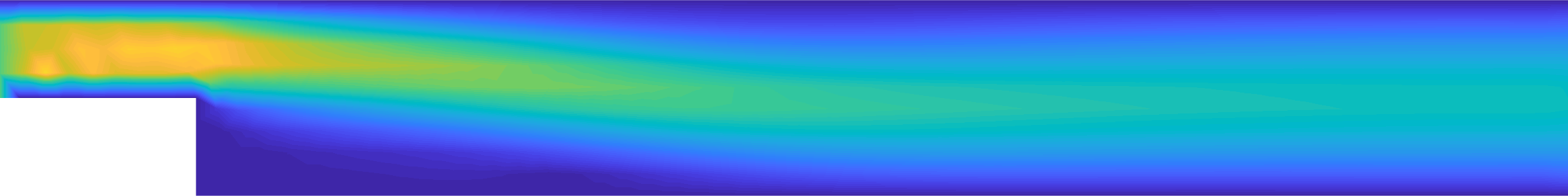};
			
			\nextgroupplot[title={}, xlabel={}, xtick=\empty]
			\addplot graphics [xmin=0, xmax=8.0, ymin=0, ymax=1.0] {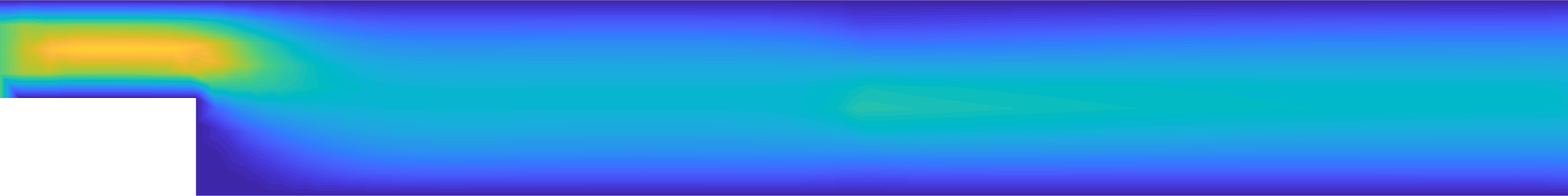};
			
			\nextgroupplot[title={}]
			\addplot graphics [xmin=0, xmax=8.0, ymin=0, ymax=1.0] {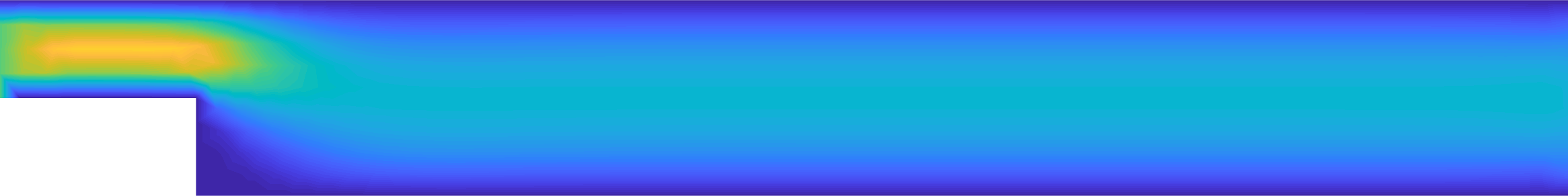};
		\end{groupplot}
	\end{tikzpicture}
	\colorbarMatlabParula{0}{0.4042}{0.8085}{1.2127}{1.6170}
	\caption{Velocity magnitude of the viscous flow problem at $\mathrm{Re}_\mathrm{worst}=4.8051$ using the ROM (\textit{top}), CG-GL method after two adaption iterations (\textit{middle}), and the reference solution (\textit{bottom}).}
	\label{fig:ns_test_vmag}
\end{figure}

The online performance is accessed by considering the solution and QoI error of the ROM and CG-GL method with respect to the reference solutions. Figure~\ref{fig:J_NS} shows the QoI as a function of the Reynolds number using the reference solution, ROM, and CG-GL after one and two adaption iterations. The value of the QoI predicted by the ROM is only accurate near the training points $\mathrm{Re}\in\{1,130\}$. On the other hand, enrichment with local basis functions makes the CG-GL predictions more accurate away from training configurations (Figure~\ref{fig:ns_errors}). All methods show a peak in the relative QoI error because the magnitude of the QoI is approaching zero; these peaks do not show up in the relative solution error. The worst-case solution error for the ROM occurs at $\mathrm{Re}_\mathrm{worst} = 4.0851$, where $e_\mathrm{sln}^\mathrm{ROM}(\mathrm{Re}_\mathrm{worst}) = 58.3\%$ and a single CG-GL adaption iteration has dropped the error to $e_\mathrm{sln}^\mathrm{CGGL}(\mathrm{Re}_\mathrm{worst}) = 7.02\%$. The ROM struggles to represent the flow separation and low-pressure regions in the left half of the domain, whereas these features are well-approximated by the CG-GL solution (Figure~\ref{fig:ns_test_vmag}).

%
%

%
%
\section{Conclusion}
\label{sec:conclu}
In this paper, we introduce a new numerical framework that is a hybrid between projection-based model
reduction and discretization methods for partial differential equations, where reliable \textit{a posteriori}
error estimation and adaptivity will be central in order to achieve robustness. The method uses basis functions
with global support to represent dominant features of the solution, similar to reduced-order models, to achieve a
high degree of accuracy with little computational cost, which are supplemented with locally supported basis functions,
such as those used in finite element methods, for additional accuracy when the pre-computed global modes are insufficient.
Two numerical experiments demonstrate the CG-GL method has better parametric robustness than classical model reduction,
which leads to more accurate predictions with limited training. The CG-GL method leverages global modes within an adaptive
framework with controllable tradeoff between accuracy and computational cost. Additional research is required to extend the approach to a discontinuous Galerkin (DG) setting so it can be used for convection-dominated problems.

\section*{Acknowledgments}
This material is based upon work supported by the Air Force Office of
Scientific Research (AFOSR) under award numbers FA9550-20-1-0236
and FA9550-22-1-0004. The content of this publication does not necessarily
reflect the position or policy of any of these supporters, and no official
endorsement should be inferred.


\bibliographystyle{plain}
\bibliography{biblio}

\end{document}